\documentclass[12pt,twoside,leqno]{article}
\usepackage{amsmath}
\usepackage{amssymb}
\usepackage{amsxtra}
\usepackage{amscd}
\usepackage{amsthm}
\usepackage[mathscr]{eucal}%\usepackage{eucal}

\theoremstyle{plain}

\newtheorem{sbthm}[subsubsection]{Theorem}
\newtheorem{sbprop}[subsubsection]{Proposition}
\newtheorem{sbcor}[subsubsection]{Corollary}
\newtheorem{sblem}[subsubsection]{Lemma}

\theoremstyle{definition}

\newtheorem{sbrem}[subsubsection]{Remark}

\newtheorem{sbpara}[subsubsection]{}

\newenvironment{pf}{\proof[\proofname]}{\endproof}

\begin{document}
\renewcommand{\thepropositionA}{}
\renewcommand{\thepropositionB}{}

\title{Classifying spaces of degenerating mixed Hodge structures, IV: 
The fundamental diagram}

\author
{Kazuya Kato, Chikara Nakayama, Sampei Usui}
%\subjclass{Primary 14M25; Secondary 14F20}

\maketitle
\renewcommand{\mathbb}{\bold}

\newcommand\Cal{\mathcal}
\newcommand\define{\newcommand}
\define\gp{\mathrm{gp}}%
\define\fs{\mathrm{fs}}%
\define\an{\mathrm{an}}%
\define\mult{\mathrm{mult}}%
\define\add{\mathrm{add}}%
\define\Ker{\mathrm{Ker}\,}%
\define\Coker{\mathrm{Coker}\,}%
\define\Hom{\mathrm{Hom}\,}%
\define\Ext{\mathrm{Ext}\,}%
\define\rank{\mathrm{rank}\,}%
\define\gr{\mathrm{gr}}%
\define\cHom{\Cal{Hom}}
\define\cExt{\Cal Ext\,}%

\define\cA{\Cal A}
\define\cC{\Cal C}
\define\cD{\Cal D}
\define\cO{\Cal O}
\define\cS{\Cal S}
\define\cM{\Cal M}
\define\cG{\Cal G}
\define\cH{\Cal H}
\define\cE{\Cal E}
\define\cF{\Cal F}
\define\cN{\Cal N}
\define\fF{\frak F}
\define\fg{\frak g}
\define\fh{\frak h}
\define\Dc{\check{D}}
\define\Ec{\check{E}}

\newcommand{\N}{{\mathbb{N}}}
\newcommand{\Q}{{\mathbb{Q}}}
\newcommand{\Z}{{\mathbb{Z}}}
\newcommand{\R}{{\mathbb{R}}}
\newcommand{\C}{{\mathbb{C}}}
\newcommand{\bN}{{\mathbb{N}}}
\newcommand{\bQ}{{\mathbb{Q}}}
\newcommand{\bF}{{\mathbb{F}}}
\newcommand{\bZ}{{\mathbb{Z}}}
\newcommand{\bP}{{\mathbb{P}}}
\newcommand{\bR}{{\mathbb{R}}}
\newcommand{\bC}{{\mathbb{C}}}
\newcommand{\bS}{{\bold{S}}}
\newcommand{\bbQ}{{\bar \mathbb{Q}}}
\newcommand{\ol}[1]{\overline{#1}}
\newcommand{\too}{\longrightarrow}
\newcommand{\respect}{\rightsquigarrow}
\newcommand{\compatible}{\leftrightsquigarrow}
\newcommand{\upc}[1]{\overset {\lower 0.3ex \hbox{${\;}_{\circ}$}}{#1}}
\newcommand{\Gmlog}{\bG_{m, \log}}%{\mathbb{G}_{m,\log}}
\newcommand{\Gm}{\bG_m}%{\mathbb{G}_m}
\newcommand{\ep}{\varepsilon}
\newcommand{\Spec}{\operatorname{Spec}}
\newcommand{\val}{{\mathrm{val}}} 
\newcommand{\n}{\operatorname{naive}}
\newcommand{\bs}{\operatorname{\backslash}}
\newcommand{\Gal}{\operatorname{{Gal}}}
\newcommand{\gal}{{\rm {Gal}}({\bar \Q}/{\Q})}
\newcommand{\galp}{{\rm {Gal}}({\bar \Q}_p/{\Q}_p)}
\newcommand{\gall}{{\rm{Gal}}({\bar \Q}_\ell/\Q_\ell)}
\newcommand{\wep}{W({\bar \Q}_p/\Q_p)}
\newcommand{\wel}{W({\bar \Q}_\ell/\Q_\ell)}
\newcommand{\Ad}{{\rm{Ad}}}
\newcommand{\BS}{{\rm {BS}}}
\newcommand{\even}{\operatorname{even}}
\newcommand{\End}{{\rm {End}}}
\newcommand{\odd}{\operatorname{odd}}
\newcommand{\GL}{\operatorname{GL}}
\newcommand{\np}{\text{non-$p$}}
\newcommand{\g}{{\gamma}}
\newcommand{\G}{{\Gamma}}
\newcommand{\Lam}{{\Lambda}}
\newcommand{\La}{{\Lambda}}
\newcommand{\lam}{{\lambda}}
\newcommand{\la}{{\lambda}}
\newcommand{\uL}{{{\hat {L}}^{\rm {ur}}}}
\newcommand{\uQp}{{{\hat \Q}_p}^{\text{ur}}}
\newcommand{\sel}{\operatorname{Sel}}
\newcommand{\dt}{{\rm{Det}}}
\newcommand{\Sig}{\Sigma}
\newcommand{\fil}{{\rm{fil}}}
\newcommand{\SL}{{\rm{SL}}}
\newcommand{\spl}{{\rm{spl}}}
\newcommand{\st}{{\rm{st}}}
\newcommand{\Isom}{{\rm {Isom}}}
\newcommand{\Mor}{{\rm {Mor}}}
\newcommand{\bg}{\bar{g}}
\newcommand{\id}{{\rm {id}}}
\newcommand{\cone}{{\rm {cone}}}
\newcommand{\al}{a}
\newcommand{\ChL}{{\cal{C}}(\La)}
\newcommand{\Image}{{\rm {Image}}}
\newcommand{\toric}{{\operatorname{toric}}}
\newcommand{\torus}{{\operatorname{torus}}}
\newcommand{\Aut}{{\rm {Aut}}}
\newcommand{\Qp}{{\mathbb{Q}}_p}
\newcommand{\barQp}{{\mathbb{Q}}_p}
\newcommand{\Qpur}{{\mathbb{Q}}_p^{\rm {ur}}}
\newcommand{\Zp}{{\mathbb{Z}}_p}
\newcommand{\Zl}{{\mathbb{Z}}_l}
\newcommand{\Ql}{{\mathbb{Q}}_l}
\newcommand{\Qlur}{{\mathbb{Q}}_l^{\rm {ur}}}
\newcommand{\F}{{\mathbb{F}}}
\newcommand{\eps}{{\epsilon}}
\newcommand{\epsLa}{{\epsilon}_{\La}}
\newcommand{\epsLaVxi}{{\epsilon}_{\La}(V, \xi)}
\newcommand{\epsOLaVxi}{{\epsilon}_{0,\La}(V, \xi)}
\newcommand{\Qplin}{{\mathbb{Q}}_p(\mu_{l^{\infty}})}
\newcommand{\otimesQplin}{\otimes_{\Qp}{\mathbb{Q}}_p(\mu_{l^{\infty}})}
\newcommand{\galFl}{{\rm{Gal}}({\bar {\Bbb F}}_\ell/{\Bbb F}_\ell)}
\newcommand{\gallur}{{\rm{Gal}}({\bar \Q}_\ell/\Q_\ell^{\rm {ur}})}
\newcommand{\galFF}{{\rm {Gal}}(F_{\infty}/F)}
\newcommand{\galFv}{{\rm {Gal}}(\bar{F}_v/F_v)}
\newcommand{\galF}{{\rm {Gal}}(\bar{F}/F)}
\newcommand{\epsVxi}{{\epsilon}(V, \xi)}
\newcommand{\epsOVxi}{{\epsilon}_0(V, \xi)}
\newcommand{\sig}{{\sigma}}
\newcommand{\ga}{{\gamma}}
\newcommand{\del}{{\delta}}
\newcommand{\Vss}{V^{\rm {ss}}}
\newcommand{\Bst}{B_{\rm {st}}}
\newcommand{\Dpst}{D_{\rm {pst}}}
\newcommand{\Dcrys}{D_{\rm {crys}}}
\newcommand{\DdR}{D_{\rm {dR}}}
\newcommand{\Fin}{F_{\infty}}
\newcommand{\Kla}{K_{\lambda}}
\newcommand{\Ola}{O_{\lambda}}
\newcommand{\Mla}{M_{\lambda}}
\newcommand{\Det}{{\rm{Det}}}
\newcommand{\Sym}{{\rm{Sym}}}
\newcommand{\LaSa}{{\La_{S^*}}}
\newcommand{\cX}{{\cal {X}}}
\newcommand{\MHG}{{\frak {M}}_H(G)}
\newcommand{\tauMla}{\tau(M_{\lambda})}
\newcommand{\Fvur}{{F_v^{\rm {ur}}}}
\newcommand{\Lie}{{\rm {Lie}}}
\newcommand{\cB}{{\cal {B}}}
\newcommand{\cL}{{\cal {L}}}
\newcommand{\cW}{{\cal {W}}}
\newcommand{\fq}{{\frak {q}}}
\newcommand{\cont}{{\rm {cont}}}
\newcommand{\SC}{{SC}}
\newcommand{\Om}{{\Omega}}
\newcommand{\dR}{{\rm {dR}}}
\newcommand{\crys}{{\rm {crys}}}
\newcommand{\hatSig}{{\hat{\Sigma}}}
\newcommand{\rdet}{{{\rm {det}}}}
\newcommand{\ord}{{{\rm {ord}}}}
\newcommand{\BdR}{{B_{\rm {dR}}}}
\newcommand{\BdRO}{{B^0_{\rm {dR}}}}
\newcommand{\Bcrys}{{B_{\rm {crys}}}}
\newcommand{\Qw}{{\mathbb{Q}}_w}
\newcommand{\barkappa}{{\bar{\kappa}}}
\newcommand{\cP}{{\Cal {P}}}
\newcommand{\cZ}{{\Cal {Z}}}
\newcommand{\oppLa}{{\Lambda^{\circ}}}
\newcommand{\bG}{{\mathbb{G}}}
\newcommand{\br}{{{\bold r}}}
\newcommand{\triv}{{\rm{triv}}}
\newcommand{\sub}{{\subset}}
\newcommand{\LD}{{D^{\star,\mild}_{\SL(2)}}}
\newcommand{\LbD}{{D^{\star}_{\SL(2)}}}
\newcommand{\dbDv}{{D^{\star}_{\SL(2),\val}}}
\newcommand{\nspl}{{{\rm nspl}}}
\newcommand{\nilp}{{{\rm nilp}}}
\newcommand{\lval}{{[\val]}}
\newcommand{\mild}{{{\rm{mild}}}}
\newcommand{\lan}{\langle}
\newcommand{\ran}{\rangle}
\newcommand{\sat}{{{\rm sat}}}
\newcommand{\Map}{{{\rm Map}}}
\newcommand{\lpr}{\langle \Phi \rangle}

\let\t=\tilde
\def\xb{\overline{x}}
\let\op=\oplus
\let\x=\times
\def\abtoric{\operatorname{|toric|}}
\def\be{\bold e}
\def\tp{\normalsize\prod}
\def\ts{\normalsize\sum}

\begin{abstract}
  We complete the construction of the fundamental diagram of various partial compactifications 
of the moduli spaces of mixed Hodge structures with polarized graded quotients. 
  The diagram includes 
the space of nilpotent orbits, the space of SL(2)-orbits, and the space of Borel--Serre 
orbits. We give amplifications of this fundamental diagram, and amplify the relations of these spaces. We describe how this work is useful to understand asymptotic behaviors of Beilinson regulators and of local height parings in degeneration. We discuss \lq\lq mild degenerations'' in which regulators converge. 
\end{abstract}

\section*{Contents}

\noindent \S\ref{s:intro}. Introduction

\noindent \S\ref{s:pre}. Preliminaries

\noindent \S\ref{s:new}. The new space $D^{\star}_{\SL(2)}$  of $\SL(2)$-orbits

\noindent \S\ref{s:val}. Valuative Borel--Serre orbits and valuative SL(2)-orbits 

\noindent \S\ref{s:newval}. New spaces $D^{\sharp}_{\Sigma,[:]}$ and $D^{\sharp}_{\Sigma,[\val]}$ of nilpotent orbits

\noindent \S\ref{s:dia}. Mild nilpotent orbits and the space $D^{\diamond}_{\SL(2)}$ of $\SL(2)$-orbits

\noindent \S\ref{s:NSB}. Complements

\noindent \S\ref{s:Ex}. Relations with asymptotic behaviors of regulators and local height pairings

\noindent \S\ref{s:co}. Corrections to \cite{KU}, supplements to Part III.

\setcounter{section}{-1}
\section{Introduction}\label{s:intro}

\renewcommand{\thefootnote}{\fnsymbol{footnote}}
\footnote[0]{Primary 14C30; 
Secondary 14D07, 32G20} 

\subsection{The fundamental diagram and its amplification}
\begin{sbpara} Let $D$ be the period domain which classifies mixed Hodge structures with polarized graded quotients 
with respect to 
the weight filtration (\cite{Gr}, \cite{U1}), with fixed Hodge numbers of graded quotients. In Part I--Part III (\cite{KNU2}) 
of this series of papers, we constructed extended period domains 
in the diagram  
$$\begin{matrix}
&&&& D_{\SL(2),\val}&\overset{\eta}{\underset{\subset}\to} &D_{\BS,\val}\\
&&&&\downarrow&&\downarrow\\
 \Gamma \bs D_{\Sig,\val}& \leftarrow & D_{\Sig,\val}^{\sharp}&\overset {\psi} \to &D_{\SL(2)}&&D_{\BS}\\
\downarrow &&\downarrow&&&&\\
 \Gamma \bs D_{\Sig}&\leftarrow &D_{\Sig}^{\sharp}&&&&
\end{matrix}$$
which we call the fundamental diagram, as the mixed Hodge versions of the extended period domains in  \cite{KU} for the pure case. We have constructed the maps in the diagram except the map $\eta$.   In this Part IV of our series of papers, we define the injective  map $\eta$. There is a big issue concerning this map $\eta$, which did not appear in the pure case,  as we explain below soon. In this Part IV, we amplify this fundamental diagram as in \ref{afd1} and \ref{afd2} below, and we remedy the issue as a result of the amplification. 

\end{sbpara}

\begin{sbpara}
The spaces in the fundamental diagram in 0.1.1 are topological spaces, the  right six spaces have $D$ as dense open sets, and the left two spaces have the quotient $\Gamma \bs D$ of $D$ by a discrete group $\Gamma$ as  dense open subsets. These left two spaces have sheaves of holomorphic functions extending that of $\Gamma \bs D$ (though these spaces need not be complex analytic spaces) and have log structures, and the right four spaces have sheaves of real analytic functions extending that of $D$ (though these spaces need not be real analytic spaces) and have log structures. The maps in the fundamental diagram except $\eta$ respect these structures. 

  Among these eight spaces, the main spaces are the three spaces $\Gamma\bs D_{\Sig}$ (the space of nilpotent orbits), $D_{\SL(2)}$ (the space of $\SL(2)$-orbits), and $D_{\BS}$ (the space of Borel--Serre orbits).
  We defined and studied  $D_{\BS}$  in Part I, 
$D_{\SL(2)}$ in Part II, and $\Gamma \bs D_{\Sig}$  in Part III. 
  The other five spaces appear to help the connection of these three spaces. 
  
  The map $\psi$ in the center of the fundamental diagram connects the four spaces in the world of nilpotent orbits on the left with the world of $\SL(2)$-orbits. We call $\psi$ the CKS map, for it is obtained in the pure case by using the work of Cattani-Kaplan-Schmid \cite{CKS} on the relation between nilpotent orbits and $\SL(2)$-orbits.
    
However, to connect the world of SL(2)-orbits and the world of Borel--Serre orbits on the right, 
the map $\eta$ has the following defect.
Though the map $\eta$ is a natural map and is continuous in the pure case (\cite{KU}), a big issue is that 
 in the mixed case, the map $\eta$ is not necessarily continuous (see Section 3.5). 

\end{sbpara}

\begin{sbpara}

To remedy this issue and to amplify the connections of the spaces in the fundamental diagram, we will introduce new spaces
$$D^{\star}_{\SL(2)}\;\;\text{and}\;\; D^{\diamond}_{\SL(2)}\quad \text{in the world of $\SL(2)$-orbits (see Section 2), and}$$
$$D^{\sharp}_{\Sig,[\val]} \;\;\text{and}\;\; D^{\sharp}_{\Sig,[:]}\quad \text{in the world of nilpotent orbits (see Section 4)}.$$

These spaces are topological spaces, and the first two have sheaves of real analytic functions and log structures. 
They have the following special properties. 

The space $D^{\star}_{\SL(2)}$ has better relations to Borel--Serre orbits than $D_{\SL(2)}$ (see Section 3.4), and this space remedies the above issue.
The spirit of the definition of $D^{\star}_{\SL(2)}$ (Section 2) is near  that of $D_{\BS}$.

As is shown in Section 5, the space $D^{\diamond}_{\SL(2)}$ has better relations to nilpotent orbits of
\lq\lq mild degeneration'' (see \ref{ss:mild} for the meaning of mildness) than $D_{\SL(2)}$, though among $D_{\SL(2)}$, $D^{\star}_{\SL(2)}$ and $D^{\diamond}_{\SL(2)}$, $D_{\SL(2)}$ is the best for the relation with general nilpotent orbits.

In the pure case, we have
$$D_{\SL(2)}=D^{\star}_{\SL(2)}=D^{\diamond}_{\SL(2)}.$$

The space $D^{\sharp}_{\Sig,[\val]}$ has a nice relation to $D_{\SL(2),\val}$ (see Section 4), which $D^{\sharp}_{\Sig,\val}$ does not have. 

The space $D^{\sharp}_{\Sig,[:]}$ is a quotient of $D^{\sharp}_{\Sig,[\val]}$ and also a quotient of $D^{\sharp}_{\Sig,\val}$, and has a nice relation to $D_{\SL(2)}$ (see Section 4).

The symbols $\star$ and $\diamond$ are used to express that the spaces are  shiny like stars and diamonds in the  relations to Borel--Serre orbits and  nilpotent orbits, respectively. The symbol $[:]$ is used because $D^{\sharp}_{\Sig, [:]}$ is regarded as a space of ratios. The symbol $[\val]$ similar to $[:]$ is used because $D^{\sharp}_{\Sig,[\val]}$ is the valuative space associated to $D^{\sharp}_{\Sig,[:]}$ for a certain log structure. 

Actually, as is explained in Part II, $D_{\SL(2)}$  has two structures $D_{\SL(2)}^I$ and $D_{\SL(2)}^{II}$ of a topological space with sheaves of real analytic functions and log structures. Everything in this Introduction is true for $D^{II}_{\SL(2)}$. 
\end{sbpara}

\begin{sbpara}\label{afd1} By using the above spaces, we have the following amplified fundamental diagram
and supplemental amplifications in \ref{afd3} and \ref{afd2}, which connect the \lq\lq three worlds'' better.
$$\begin{matrix}
&&&&&& D^{\star}_{\SL(2),\val}&\overset{\eta^{\star}}{\underset{\subset}\to}&D_{\BS,\val}\\
&&&&&&\downarrow && \downarrow\\
&&&&D^{\sharp}_{\Sigma,\lval}&
\overset{\psi}\to & D_{\SL(2),\val}& &D_{\BS}\\
&&&&\downarrow &&\downarrow&&\\
 \Gamma \bs D_{\Sig,\val}& \leftarrow & D_{\Sig,\val}^{\sharp}&\to &D^{\sharp}_{\Sig,[:]} & \overset {\psi} \to &D_{\SL(2)}&&\\
\downarrow &&&&\downarrow&&&&\\
 \Gamma \bs D_{\Sig}&&\leftarrow &&D_{\Sig}^{\sharp}&&&&
\end{matrix}$$
This diagram is commutative and  the maps respect the structures of the spaces. 
As indicated in this diagram, 
the valuative space $D^{\star}_{\SL(2),\val}$ associated to $D^{\star}_{\SL(2)}$  has  an injective morphism $\eta^{\star}: D^{\star}_{\SL(2),\val}\to D_{\BS,\val}$ (Theorem \ref{SL2BS}), which is an improved version of $\eta$, and a proper surjective morphism $D^{\star}_{\SL(2),\val}\to D_{\SL(2),\val}$ (Theorem \ref{0thm}). Here morphism means a morphism of topological spaces endowed with  sheaves of real analytic functions and with log structures. 
As also indicated in the diagram, the CKS map $\psi: D^{\sharp}_{\Sig,\val}\to D_{\SL(2)}$ factors  as $D^{\sharp}_{\Sig,\val}\to D^{\sharp}_{\Sig, [:]}\to D_{\SL(2)}$, and we have a continuous map $\psi: D^{\sharp}_{\Sig, [\val]}\to D_{\SL(2),\val}$ (Theorem \ref{valper}). 
\end{sbpara}

\begin{sbpara}\label{afd3} In the case $\Sig$ is the fan $\Xi$ of all rational nilpotent cones of rank $\leq 1$, we have $$D^{\sharp}_{\Xi,[\val]}= D^{\sharp}_{\Xi,[:]}= D^{\sharp}_{\Xi,\val}= D^{\sharp}_{\Xi}.$$ Furthermore in this case, we have  a CKS map  $\psi: D^{\sharp}_{\Xi}\to D^{\star}_{\SL(2),\val}$, and hence the three worlds are connected directly by 
$$ D^{\sharp}_{\Xi}
\overset{\psi}\to   D^{\star}_{\SL(2),\val}
\overset{\eta^{\star}}{\underset{\subset}\to} D_{\BS,\val}.$$ See Theorem \ref{rk1}. 
  \end{sbpara}
  
  \begin{sbpara}\label{star+-}
  
   As is described above, the spaces $D_{\SL(2)}$, $D^{\star}_{\SL(2)}$ and $D_{\BS}$ are related via their associated valuative spaces $D_{\SL(2),\val}$, $D^{\star}_{\SL(2),\val}$ and $D_{\BS,\val}$. The associated valuative space is a kind of a projective limit of blowing-ups. In Section 2, we will construct also spaces $D^{\star,+}_{\SL(2)}$, $D^{\star,-}_{\SL(2)}$ and $D^{\star,\BS}_{\SL(2)}$ which are related to $D^{\star}_{\SL(2)}$ via kinds of blowing-ups and blowing-downs and which work as bridges between  $D_{\SL(2)}$, $D^{\star}_{\SL(2)}$ and $D_{\BS}$ before going to the valuative spaces. See Section 2.

  \end{sbpara}
  
\begin{sbpara}

A nilpotent orbit appears as the limit of a variation of mixed 
Hodge structure in degeneration. SL(2)-orbits are simpler objects and Borel--Serre orbits are further simpler. 
The theory of SL(2)-orbits (\cite{Scw}, \cite{CKS} for the pure case and \cite{Pe}, \cite{KNU1} for the mixed case) tells that, roughly speaking, an SL(2)-orbit is associated to a nilpotent orbit, and we can
read real analytic behaviors of the degeneration better by looking at  the simpler object SL(2)-orbit. 
The map $\psi$  gives the  SL(2)-orbit associated to a nilpotent orbit. 

We hope that the above extended period domains and their relations are useful in the study of degeneration of mixed Hodge structures.

Actually, as illustrated in Section \ref{relBK} below and in Section 7, 
our theory has an application to the study \cite{BK} of asymptotic behaviors of degenerations of Beilinson regulators and local height pairings. In these subjects, the asymptotic behaviors are understood by degeneration of mixed Hodge structures.

\end{sbpara}

\subsection{Mild degenerations}\label{ss:mild}

\begin{sbpara}

We will define 
the subsets
$$D^{\mild}_{\Sig}\subset D_{\Sig}, \quad D^{\star,\mild}_{\SL(2)}\subset D^{\star}_{\SL(2)},
\quad D^{\mild}_{\BS}\subset D_{\BS}$$ of elements with mild degenerations.

Any element of $D^{\diamond}_{\SL(2)}$ is regarded as having mild degeneration. 

\end{sbpara}

\begin{sbpara}
Let $D^{\mild}_{\Sigma}$ be the subset of $D_{\Sig}$ consisting of all points $p$ satisfying the following condition. 

\medskip

For any element $N$ of the monodormy cone associated to $p$, there is a splitting of $W$ which is compatible with $N$.
(The splitting  can depend on $N$ and need not have any relation  with the Hodge filtration).

\medskip

Denote the subset $D_{\BS}^{(A)}$
of $D_{\BS}$ (Part I, 8.1) by $D^{\mild}_{\BS}$. 
Let $D^{\mild}_{\BS,\val}\subset D_{\BS,\val}$ be the inverse image of $D^{\mild}_{\BS}$. There is also a subset $D^{\star,\mild}_{\SL(2)}$ of $D^{\star}_{\SL(2)}$ consisting of $A$-orbits (see Section 2) whose inverse image $D^{\star,\mild}_{\SL(2),\val}$ in $D^{\star}_{\SL(2),\val}$ coincide with the inverse image of $D^{\mild}_{\BS,\val}$ under $\eta^{\star}$. We have also the mild parts of $D_{\Sigma,[:]}^{\sharp}$ and $D^{\sharp}_{\Sigma,[\val]}$, i.e., let $D_{\Sigma,[:]}^{\sharp,\mild}$ and $D^{\sharp,\mild}_{\Sigma,[\val]}$ be the inverse images of $\Gamma\bs D^{\mild}_{\Sigma}$ in $D_{\Sigma,[:]}^{\sharp}$ and in $D_{\Sigma,[\val]}^{\sharp}$, respectively.

All these mild parts $\Gamma \bs D^{\mild}_{\Sigma}$, 
$D_{\BS}^{\mild}$,  $\dots$, etc.\ are  open sets of $\Gamma \bs D_{\Sig}$, $D_{\BS}$,  $\dots$ etc., respectively.

\end{sbpara}

\begin{sbpara}\label{afd2}
For mild degenerations,
we can replace the upper right part of the amplified fundamental diagram by the following 
commutative diagram (maps respect structures of the spaces) which contain the space $D^{\diamond}_{\SL(2)}$ and its associated valuative space $D^{\diamond}_{\SL(2),\val}$ (Theorem \ref{diathm}).
$$\begin{matrix}
  D^{\sharp,\mild}_{\Sigma,\lval}&
\overset{\psi}\to  & D^{\diamond}_{\SL(2),\val} & \to &D^{\star,\mild}_{\SL(2),\val}&\overset{\eta^{\star}}{\underset{\subset}\to} & D^{\mild}_{\BS, \val}\\
\downarrow &&\downarrow&&\downarrow &&\downarrow \\
D^{\sharp,\mild}_{\Sig,[:]} &\overset{\psi}\to & D^{\diamond}_{\SL(2)} &\to & D^{\star,\mild}_{\SL(2)}&& D^{\mild}_{\BS}\\
&&&&\downarrow&&\\
&&&&D_{\SL(2)}&&
\end{matrix}$$
\end{sbpara}

\begin{sbpara}\label{app} In the applications of our work as in Section 7, the following part of the fundamental diagrams in Section 0.1  and Section \ref{ss:mild} becomes important.
$$\begin{matrix} D_{\Sig,[:]}^{\sharp,\mild}&\to& D^{\diamond}_{\SL(2)}\\ \cap
&&\downarrow\\
D_{\Sig,[:]}^{\sharp}&\to& D_{\SL(2)}
\end{matrix}$$
Via this diagram, we can understand degeneration of mixed Hodge structure in the space $D_{\SL(2)}$, and understand mild degeneration better in $D^{\diamond}_{\SL(2)}$. The right vertical arrow is usually not injective, and hence $D^{\diamond}_{\SL(2)}$ can tell informations about mild degeneration which is lost in $D_{\SL(2)}$. This is explained in Section \ref{relBK} below, and in Section \ref{s:Ex} more precisely. 
\end{sbpara}

\subsection{Relations with regulators and local height pairings}\label{relBK}
We illustrate the relations of this work with the work \cite{BK}.

\begin{sbpara}\label{K2intro} 
Let $S$ be a smooth curve over $\C$ and let 
 $f:X\to S$ be a proper surjective morphism from a smooth algebraic variety  $X$. 
  Let $0\in S$ be a point, 
and assume that $X\smallsetminus X_0\to S\smallsetminus \{0\}$ is smooth and $X$ is of semistable reduction at $0\in S$. 
 
 For $Z\in K_n(X\smallsetminus X_0)$ $(n\geq 1)$, the asymptotic behavior of the regulator of the restriction $Z(t)\in K_n(X_t)$ of $Z$ to $X_t$ $(t\in S
\smallsetminus\{0\}$,  $t\to 0)$ is studied in \cite{BK} by using our theory of degeneration of MHS. For each $r\geq 0$, $Z$ defines a variation of mixed Hodge structure $H_Z$ on $S\smallsetminus\{0\}$ with an exact sequence 
$0\to H^m(X/S)(r) \to H_Z \to \Z\to 0$, where $m=2r-n-1$, $H^m(X/S)$ is the $m$-th direct image $R^mf_*\Z$ on $S\smallsetminus \{0\}$ with Hodge filtration, and $(r)$ is the Tate twist. The ($r$-th) regulator of $Z(t)$ is determined by the fiber $H_Z(t)$ of $H_Z$ at $t$. 
\end{sbpara}

\begin{sbpara} We describe how our theory is related to this subject. Our description in the rest of Section 0.3 is rough and imprecise. More precise matters are described in Section \ref{ss:reg} and details are given in \cite{BK}. 

We have the period map 
$$(S\smallsetminus
\{0\}) \times K_n(X\smallsetminus X_0) \to \Gamma\bs D\quad (t, Z)\mapsto \text{class}(H_Z(t)).$$ By Part III, this extends to 
$$S\times K_n(X\smallsetminus X_0) \to \Gamma \bs D_{\Xi}, \quad S^{\log}\times K_n(X\smallsetminus
X_0) \to \Gamma \bs D^{\sharp}_{\Xi}$$
where $S^{\log}$ is the space associated to $S$ defined in \cite{KN}. If $Z$ comes from $K_n(X)$, then $H_Z$ has mild degeneration at $0\in S$ (\ref{KXmild}).  
The diagram in \ref{app} produces the following commutative diagram. 
$$\begin{matrix}  S^{\log}\times K_n(X) &\to &\Gamma\bs D^{\sharp,\mild}_{\Xi} &\to  & \Gamma\bs D^{\diamond}_{\SL(2)} \\
\downarrow && \cap &&\downarrow\\
S^{\log} \times K_n(X\smallsetminus
X_0)&\to &\Gamma\bs D^{\sharp}_{\Xi} &\to& \Gamma \bs D_{\SL(2)} 
\end{matrix}$$
\end{sbpara}

\begin{sbpara}  
We can prove that for  $Z\in K_n(X)$, the regulator of $Z(t)$ converges when $t\to 0$ (Theorem \ref{thm2}). In fact, this is a consequence of the fact that the period map $S\smallsetminus
\{0\}\to \Gamma \bs D\;;\;t\mapsto \text{class}(H_Z(t))$ induced by $Z$ extends to a continuous map $S^{\log}\to \Gamma\bs D^{\diamond}_{\SL(2)}$ as indicated by the upper row of the above diagram. We recover the limit of the regulator of $Z(t)$ for $t\to 0$ from the image of a point $b$ of $S^{\log}$ over $0$
in $\Gamma \bs D^{\diamond}_{\SL(2)}$. 
On the other hand, for $Z\in K_n(X\smallsetminus X_0)$ which need not come from $K_n(X)$, the regulator of $Z(t)$ need not converge when $t\to 0$, and the image of $b$ in $\Gamma \bs D_{\SL(2)}$ tells how rapidly it 
 diverges. When $Z$ comes from $K_n(X)$, 
 the image of $b$ in $\Gamma \bs D_{\SL(2)}$ has smaller information than the image of $b$ in $\Gamma \bs D^{\diamond}_{\SL(2)}$, and cannot tell the limit of the regulator of $Z(t)$. 
\end{sbpara}

\begin{sbpara}\label{htintro} We have a similar story for the asymptotic behavior of the local height pairing (at the Archimedean place). This is introduced in Section \ref{ss:reBK}. 
\end{sbpara}

\subsection{Organization of this paper, acknowledgements}

\begin{sbpara} The organization of this paper is as follows.

Section 1 is a preparation.

In Section \ref{s:new}, we consider the new space $D^{\star}_{\SL(2)}$ 
of $\SL(2)$-orbits.

In Section \ref{s:val}, we consider the spaces $D_{\SL(2),\val}$ and $D^{\star}_{\SL(2),\val}$ of valuative $\SL(2)$-orbits and the space $D_{\BS,\val}$ of valuative Borel--Serre orbits.

In Section \ref{s:newval}, we consider the new spaces $D^{\sharp}_{\Sig,[\val]}$ and $D^{\sharp}_{\Sig, [:]}$ in the world of nilpotent orbits, and improve CKS maps by using these spaces.

In Section \ref{s:dia}, we consider the new space $D^{\diamond}_{\SL(2)}$ 
of $\SL(2)$-orbits and construct mild CKS maps.

In Section \ref{s:NSB}, we give complementary results on properties of extended period domains, on relations of nilpotent orbits, SL(2)-orbits, and Borel-Serre orbits, and on extended period maps. 

 In Section \ref{s:Ex}, we illustrate the relations to the work \cite{BK} and give examples. 
 
 In the appendix Section \ref{s:co}, we give corrections to \cite{KU} and supplements to Part III.  Sections A.1 and  Section A.3 in this appendix are directly related to Section \ref{ss:Lthm} of this Part IV. 
\end{sbpara}

\begin{sbpara} We thank Spencer Bloch.  Theorem \ref{diathm}  and Theorem \ref{thm2} in this Part IV were obtained in  joint efforts with him related to the joint work \cite{BK}.

  The first author is partially supported by NFS grant DMS 1001729.
  The second author is partially supported by JSPS. KAKENHI (C) No. 18540017, (C) No. 22540011.
  The third author is partially supported by JSPS. KAKENHI  (B) No. 23340008.
\end{sbpara}

\noindent 

\section{Preliminaries}\label{s:pre}

\subsection{The setting}\label{setting}
  We recall the basic setting and the notation used throughout this series of papers.

\begin{sbpara}\label{hodge}
  We fix $\Lambda=(H_0, W, (\langle\;,\;\rangle_w)_w, (h^{p,q})_{p, q})$, where 

\medskip

$H_0$ is a finitely generated free $\bZ$-module,

$W$ is a finite increasing rational filtration on $H_{0, \bR} = \bR\otimes H_0$,

$\langle\;,\;\rangle_w$ for each $w\in \bZ$ is a rational nondegenerate $\bR$-bilinear form  $\gr^W_w\times \gr^W_w\to \bR$ which is symmetric if $w$ is even and is anti-symmetric if $w$ is odd, 

$h^{p,q}$ is a nonnegative integer given for each $(p, q)\in \bZ^2$,

\noindent
satisfying the following conditions (1)--(3).

\smallskip

(1) $\sum_{p, q} h^{p,q}= \text{rank}_\bZ(H_0)$,
\smallskip

(2) $\sum_{p+q=w} h^{p,q}= \dim_\bR(\gr^W_w)$ for any $w\in \bZ$,
\smallskip

(3) $h^{p,q}=h^{q,p}$ for any $(p, q)$. 
\end{sbpara}

\begin{sbpara}
Let $D$ be the classifying space of gradedly polarized mixed 
Hodge structures in \cite{U1} associated to the data fixed in \ref{hodge}.
As a set, $D$ consists of all increasing filtrations $F$ on $H_{0,\bC} = \bC \otimes H_0$ such that $(H_0, 
W, (\langle\;,\;\rangle_w)_w, F)$ is a gradedly polarized mixed Hodge structures with $\dim_\bC F^p(\gr^W_{p+q})/F^{p+1}(\gr^W_{p+q}) 
= h^{p,q}$ for all $p, q$. 

The space $D$ is an open subset of a simpler complex analytic manifold  $\Dc$  (Part I, 1.5) which is defined by dropping the condition of positivity for $\langle\;,\;\rangle_w$ in the definition of $D$.  

\end{sbpara}

\begin{sbpara}
%As in Part I, 1.6, 
For $A=\bZ, \bQ, \bR, \bC$, let 
$G_A$ be the group of the $A$-automorphisms of $(H_{0,A}, W)$ whose $\gr^W_w$ 
are compatible with $\langle\;,\;\rangle_w$ for all $w$.
  Here $H_{0,A} = A \otimes H_0$.  
  Then $G_\bC$ (resp.\ $G_\bR$)  acts on $\Dc$ (resp.\ $D$). 
%As in Part I, 1.7, 
For $A= \bQ, \bR, \bC$, let $\fg_A$ be the associated Lie algebra of $G_A$. 
%, which is 
%identified with the set of all $A$-endomorphisms 
%$N$ of $(H_{0,A}, W)$ whose $\gr^W_w$ satisfies 
%$\lan\gr^W_w(N)(x),y\ran_w+\lan x,\gr^W_w(N)(y)\ran_w=0$ for all $w, x, y$. 

%As in Part I, 1.6 and 1.7, 
Let $G_{A,u}=\{\g\in G_A\;|\;\gr^W(g)=1\}$,
$\fg_{A,u}=\{N\in \fg_A\,|\,\gr^W(N)=0\}$. Then $G_A/G_{A,u}$ 
is isomorphic to $G_A(\gr^W):=\prod_w G_A(\gr^W_w)$ and $\fg_A/\fg_{A,u}$ 
is isomorphic to $\fg_A(\gr^W):=\prod_w\fg_A(\gr^W_w)$, where $G_A(\gr^W_w)$ (resp.\ $\fg_A(\gr^W_w)$) is \lq\lq the  $G_A$ (resp.\ $\fg_A$) for $\gr^W_w$''. 
\end{sbpara}

\begin{sbpara}
  For each $w\in \bZ$, let 
$D(\gr^W_w)$ be the $D$ for the graded quotient $((H_0\cap W_w)/(H_0\cap W_{w-1}), \lan\;,\;\ran_w, (h^{p,q})_{p+q=w})$. 
  Let $D(\gr^W) = \prod_{w\in \bZ} \; D(\gr^W_w).$
  Then the canonical morphism 
$$D \to D(\gr^W); F \mapsto F(\gr^W):= (F(\gr^W_w))_{w\in \bZ}$$ is surjective.
\end{sbpara}

\begin{sbpara}\label{splW}
  Let $W'$ be a finite increasing filtration on $H_{0, \bR}$. 
  A {\it splitting} of $W'$ is an 
isomorphism 
$$s\colon \gr^{W'}:=\normalsize\bigoplus_w \gr^{W'}_w \overset \simeq \to H_{0,\bR}$$
of $\bR$-vector spaces such that for any $w\in \bZ$ and $v \in
\gr^{W'}_w$, $s(v) \in W'_w$ and $v = (s(v)\bmod W'_{w-1})$.

  Let $\spl(W')$ be the set of all splittings of $W'$. 

  Consider the case $W'=W$. 
  Then $\spl(W)$ is regarded as a $G_{\bR, u}$-torsor. 

  Let $D_{\spl}:= \{s(F) \;|\; s\in \spl(W),\, F\in D(\gr^W)\} \subset D$ be 
the subset of {\it $\bR$-split} elements.
Here $s(F)^p:= s(\bigoplus_w F_{(w)}^p)$ for $F = (F_{(w)})_w \in D(\gr^W)$.

Then, $D_{\spl}$ is a real analytic closed submanifold of $D$, and we have a real analytic isomorphism
$\spl(W) \times D(\gr^W) \overset \sim \to D_{\spl}$, $(s, F)\mapsto s(F)$.

Let $D_{\nspl}:=D\smallsetminus D_{\spl}$.
\end{sbpara}

\subsection{Canonical splitting of the weight filtration and the invariant $\delta$ of non-splitting}

 \begin{sbpara}\label{grsd} We review the fact that the weight filtration of an $\R$-mixed Hodge structure has a canonical splitting over $\R$ (which does not split the Hodge filtration except the case of an $\R$-split mixed Hodge structure) and the fact that there is an important map $\delta$  which tells how the $\R$-mixed Hodge structure is far from $\R$-split. 
  We review that we have an isomorphism of real analytic manifolds $$D\overset{\cong}\to \{(F, s, \delta)\in D(\gr^W)\times \spl(W) \times\cL  \;|\; \delta\in \cL(F)\}$$
$$ x\mapsto (x(\gr^W), \spl_W(x), \delta_W(x))$$ 
(Part II, Proposition 1.2.5) by using the canonical splitting $\spl_W(x)$ of $W$ associated to $x$ and the invariant $\delta_W(x)$ of non-splitting associated to $x$. 
$\cL$ and $\cL(F)$ are explained in \ref{cL(F)}, $\delta_W(x)$ is explained in \ref{II,1.2.2}, and $\spl_W(x)$ is explained in \ref{II,1.2.3} below. 
  
  \end{sbpara}
  
  \begin{sbpara}\label{cL(F)}

  Let $\cL= W_{-2}\mathrm{End}_{\bR}(\gr^W)$ be the set of all $\bR$-linear maps $\delta:\gr^W\to\gr^W$ such that 
  $\delta(\gr^W_w)\subset\bigoplus_{w'\le w-2}\gr^W_{w'}$ for all $w\in\bZ$ (Part II, 1.2.1). 
This is a finite dimensional weighted $\bR$-vector space.

For $F\in D(\gr^W)$, let 
$\cL(F)$ be the weighted subspace of $\cL$ consisting of all elements  whose $(p,q)$-Hodge components for $F$ are $0$ unless $p<0$ and $q<0$. That is, $\cL(F)$ is the set of all $\delta\in \cL$ such that 
$\delta(H^{p,q}_F)\subset %\tsize
\bigoplus_{p'<p, q'<q}\; 
H^{p',q'}_F\;\text{for all}\;p,q\in\bZ$. Here $H^{p,q}_F$ denotes the $(p,q)$-Hodge component of $F(\gr^W_{p+q})$
 (Part II, 1.2.1).

  \end{sbpara}

\begin{sbpara}\label{II,1.2.2}
  We explain $\delta_W(x)\in \cL(x(\gr^W))$. 

For $x \in D$, there is a unique pair of 
$s' \in \spl(W)$ and $\delta \in 
\cL(x(\gr^W))$ such that
$$
x = s'(\exp(i\delta)x(\gr^W))
$$
(\cite{CKS} (2.20)).
  We write $\delta_W(x)$ (or  $\delta(x)$) for this $\delta$. 
  
  \end{sbpara}
  
  \begin{sbpara}
Roughly speaking, $\delta_W(x)$ is the invariant of the mixed Hodge structure $x$ 
which measures how  $x$ is far from $D_{\spl}$ in $D$. We have  $\delta_W(x)=0$ if and only if $x\in D_{\spl}$ (\ref{splW}). 

This $\delta_W(x)$ plays important roles in our series of papers. It is related to the regulator in number theory and in arithmetic geometry 
as is discussed  in \cite{BK} and in Section 7 of this Part IV.  
 Hence we would like to propose to call $\delta_W(x)$ the regulator of the mixed Hodge structure $x$. 
\end{sbpara}

\begin{sbpara}\label{II,1.2.3} 
We explain $\spl_W(x)\in \spl(W)$.

Let $x \in D$, and let $s'\in \spl(W)$ and $\delta$ be as in \ref{II,1.2.2}.
Then the {\it canonical splitting} $s=\spl_W(x)$ of $W$ associated to $x$ is defined by
$$
s=s'\exp(\zeta),
$$
where $\zeta=\zeta(x(\gr^W), \delta)$ is a certain element of $\End_{\bR}(\gr^W)$ determined by $x(\gr^W)$ and $\delta=\delta_W(x)$ 
roughly as in the following way. 

Let $\delta_{p,q}$ $(p,q\in \bZ)$ be the $(p,q)$-Hodge component of $\delta$ with respect to $x(\gr^W)$.
Then the $(p,q)$-Hodge component $\zeta_{p, q}$ of $\zeta=\zeta(x(\gr^W), \delta)$ with respect to $x(\gr^W)$ is 
given as a certain universal Lie polynomial of $\delta_{p', q'}$ $(p', q'\in \bZ$,
$p'\leq -1$, $q'\leq -1)$. 
See \cite{CKS} (6.60), and Section 1 and Appendix of \cite{KNU1} for more explanations.

For $x\in D$, $x_{\spl}:=s(x(\gr^W))\in D_{\spl}$ with $s=\spl_W(x)$  is called the {\it associated $\R$-split mixed Hodge structure}. We have $x\in D_{\spl}$ if and only if $x=x_{\spl}$.

\end{sbpara}

\begin{sbpara}\label{liftac} We have the following action of the  group $\prod_{w\in \Z} \Aut_\R(\gr^W_w)$  on $D$, which we call the {\it lifted action}. 
For $a=(a_w)_w\in \prod_{w\in \Z} \Aut_\R(\gr^W_w)$, $a$ sends $x\in D$ to $x'\in D$ 
which is characterized by $x'(\gr^W_w)=a_wx(\gr^W_w)$, $\spl_W(x')=\spl_W(x)$, 
and $\delta_W(x')=\Ad(a)\delta_W(x)$. In other words, $a$ sends the Hodge filtration 
$F\in D$ to the Hodge filtration $s_F a(s_F^{-1}(F))$ where $s_F:=\spl_W(F)$ and 
$s_F^{-1}(F)$ denotes 
the filtration on $\gr^W_{\C}=\prod_w \gr^W_{w,\C}$ induced by $F$ via $s_F^{-1}: H_{0,\C}\overset{\cong}\to \gr^W_{\C}$. 

This lifted action will be used in Section 2. 

\end{sbpara}

\subsection{Spaces with real analytic structures and with fs log structures with sign}
\label{ss:B}
This is essentially a review of Section 3.1 of Part II.

\begin{sbpara}\label{2.5.1} Endow $\R^n$ $(n\geq 0)$ with the sheaf $\cO_{\R^n}$ of real analytic functions.

Let $\cB_\R'$ be the category of locally ringed spaces $S$ over $\R$ satisfying the following condition (i) locally on $S$. 

\smallskip

(i) There are $n\geq 0$ and a morphism $\iota:S\to \R^n$ of locally ringed spaces over $\R$ 
such that $\iota$ is injective, the topology of $S$ coincides with the topology induced from that of $\R^n$, and the map $\iota^{-1}(\cO_{\R^n}) \to \cO_S$ is surjective. 

\smallskip

For an object $S$ of $\cB'_\R$, we often call the structural sheaf $\cO_S$  the sheaf of real analytic functions on $S$ (though $S$ need not be a real analytic space). 

\smallskip

 Let $\cC_\R$ be the category of locally ringed spaces $S$ over $\R$ satisfying the following condition (ii).

\smallskip

(ii) For any open set $U$ of $S$ and for any $n\geq 0$, the canonical map $\text{Mor}(U, \R^n)\to \cO_S(U)^n$ is bijective. 
\end{sbpara}

\begin{sbpara}
We have $$\cB'_\R\subset \cC_\R.$$ For the proof, see Part II, Lemma 3.1.2. 

\end{sbpara}

\begin{sbpara}\label{value} For a topological field $K$ and for a locally ringed space $S$ over $K$, the following three conditions (i)--(iii) are equivalent.

\smallskip

(i) For any $s\in S$, the map $K\to \cO_{S,s}/m_s$ ($m_s$ denotes the maximal ideal of $\cO_{S,s})$ is an isomorphism. Furthermore for any open set $U$ of $S$ and for any $f\in \cO_S(U)$, the map $U\to K\;;\; s\mapsto f(s)$  is continuous. Here $f(s)$ denotes the image of $f$ in $\cO_{S,s}/m_s=K$. 
\smallskip

(ii) Let $\cO'_S$ be the sheaf on $S$ of all $K$-valued continuous functions. Then there is a homomorphism $\cO_S\to \cO'_S$ of sheaves of rings over $K$. 

\smallskip

(iii) Let $S'$ be the topological space $S$ endowed with the sheaf of all $K$-valued continuous functions. Then there is a morphism of locally ringed spaces $S'\to S$  over $K$ whose underlying map $S'\to S$ is the identity map.

\smallskip

If these equivalent conditions are satisfied, there is only one homomorphism $\cO_S\to \cO'_S$ of sheaves of rings over $K$, and  there is only one morphism $S'\to S$ of locally ringed spaces over $K$ lying over the identity map of $S$. 

These can be proved easily. 
\end{sbpara}

\begin{sbpara}
Note that objects of $\cC_\R$ satisfy the equivalent conditions in \ref{value} with $K=\R$.

\end{sbpara}

\begin{sbpara}\label{lsign}
Let $S$ be a locally ringed space over $\R$ satisfying the equivalent conditions in \ref{value} with $K=\R$. By a {\it log structure with sign} on $S$, we mean a log structure $M$ on $S$ endowed with a submonoid sheaf $M_{>0}$ of $M$ satisfying the following (i) and (ii).

\smallskip

(i)  $M_{>0}\supset \cO^\times_{S,>0}$. Here $\cO^\times_{S,>0}$ denotes the subgroup sheaf of $\cO_S^\times$ consisting of 
all local sections whose values are $>0$. 

(ii) The map $M_{>0} \times \{\pm 1\}\to M\;;\;(f, \varepsilon)\mapsto \varepsilon f$ is an isomorphism of sheaves. Here we regard $\{\pm 1\}\subset \cO_S^\times \subset M$. 

\medskip

Note that the map $\cO^\times_{S,>0}\times \{\pm 1\}\to \cO^\times_S\;;\;(f, \varepsilon)\mapsto \varepsilon f$ is an isomorphism. Indeed,  if $f \in \cO_S^\times$ has value $>0$ (resp.\ $<0$) at $s\in S$, then $f$ (resp.\ $-f$) belongs to $\cO_{S,>0}^{\times}$ on some open neighborhood of $s$. Hence this map is surjective. The injectivity is clear. 
\end{sbpara}

\begin{sbpara}\label{integr}

In Pat II, Section 3.1, we defined the notion log structure with sign in a more restrictive situation where $S$ is an object of $\cC_\R$ requiring $M$ is integral (that is, the canonical map $M\to M^{\gp}$ is injective), and the presentation of the definition there was  more complicated. So here we are improving the generality and the presentation of the definition. (But in this paper, we do not need this generalization.) 
If $M$ is integral, the present definition is equivalent to the definition in Part II,  3.1.5 
which uses  a subgroup sheaf $M^{\gp}_{>0}$. The relation with the present definition is that $M^{\gp}_{>0}$ in Part II, 3.1.5 is obtained from  $M_{>0}$ in the present definition as $M^{\gp}_{>0}= (M_{>0})^{\gp}$, and  $M_{>0}$ here is obtained from $M^{\gp}_{>0}$ there as $M_{>0}= M\cap M^{\gp}_{>0}$. To prove the equivalence, the non-trivial point is to show that

\smallskip

(1) $M_{>0}\cap \cO_S^\times = \cO^\times_{S,>0}$

\smallskip
\noindent 
for a log structure with sign in the present sense. We prove (1). If $f\in M_{>0}\cap \cO_S^\times$ has a value $<0$ at $s\in S$, $-f$ belongs to $\cO^\times_{S,>0}\subset M_{>0}$ on some open neighborhood of $s$, and this contradicts the condition (ii) in \ref{lsign}. Hence $f\in \cO^\times_{S,>0}$. 

Note that (1) implies the condition (3) in Part II, 3.1.5 on $M$, that is, the values of $f\in M_{>0}$ are $\geq 0$. (The values of $f$ mean the values of the image of $f$ in $\cO_S$.)  Indeed,  for $s\in S$, if  the image of $f$ in $M_s$ belongs to $\cO_{S,s}^\times$, then it belongs to $\cO^\times_{S, >0,s}$ by the above (1), and hence $f$ has value $>0$ at $s$. If the image of $f$ in $M_s$ does not belong to $\cO_{S,s}^\times$, then $f$ has value $0$ at $s$.

\end{sbpara}

\begin{sbpara}\label{cblog} 

Let $\cB'_\R(\log)$ be the category of objects of $\cB'_\R$  (\ref{2.5.1}) endowed with an fs log structure with sign. 

Let $\cC_\R(\sat)$ be the category of objects of $\cC_\R$ endowed with a saturated log structure with sign.

Here a log structure $M$ on a locally ringed space $S$ is said to be {\it saturated} if all stalks of $M$ are saturated in the following sense. We say 
a commutative monoid $\cS$ is saturated if it is integral (that is, the canonical map $\cS\to \cS^{\gp}$ is injective) and if for any $a\in \cS^{\gp}$ such that $a^n\in \cS\subset \cS^{\gp}$ for some integer $n\geq 1$, we have $a\in \cS$.

We have 
$$\cB'_\R(\log)\subset \cC_\R(\sat).$$
\end{sbpara}

\begin{sbpara}\label{2.3ex}
{\it Examples.}

\medskip

(1) {\it The object $\R^n_{\geq 0}$ of $\cB'_\R(\log)$. }
The sheaf $\cO$ of real analytic functions is the inverse image of the sheaf of real analytic functions on $\R^n$. The log structure $M$ with sign is  as follows. $M$ (resp.\ $M_{>0}$) is the multiplicative submonoid sheaf of $\cO$ generated by $\cO^\times$ (resp.\ $\cO^\times_{>0}$) and the coordinate functions $t_1, \dots, t_n$.  

\medskip

(2) {\it A real analytic manifold with corners}
(\cite{BS}, Appendix) is regarded as an object of $\cB'_\R(\log)$.
The log structure with sign is given as follows. Let $S$ be a real analytic manifold with corners and let $\cO$ be the sheaf of real analytic functions. If $S$ is an open set of $\R^n_{\geq 0}$ (endowed with the sheaf of real analytic functions), the log structure with sign $(M, M_{>0})$ is defined as the inverse image of that of $\R^n_{\geq 0}$. In this situation, the canonical map $M\to \cO$ is injective and hence $M$ and $M_{>0}$ are regarded as subsheaves of $\cO$.  In general, $S$ is locally isomorphic to an open set of $\R^n_{\geq 0}$, and the log structure with sign on $S$ induced from such isomorphism is independent of the choice of the isomorphism ($M$ and $M_{>0}$ are independent of the choice as subsheaves of $\cO$).

By this, we have (a real analytic manifold with corners) $=$ 
(an object of $\cB'_\R(\log)$ which is locally isomorphic to an open subobject of $\R^n_{\geq 0}$ $(n\geq 0)$). 

\medskip

(3) {\it The real toric variety $\Hom(\cS, \R^{\mult}_{\geq 0})$ for an fs monoid $\cS$.}
(Here $\R^{\mult}_{\geq 0}$ is the set $\R_{\geq 0}$ regarded as a 
multilpicative monoid.) This is also an object of $\cB'_\R(\log)$. (The above (1) is the case $\cS=\N^n$ of this (3).) 

The sheaf $\cO$ of real analytic functions is defined as follows. 
Take a surjective homomorphism $\N^n\to \cS$ of monoids for some $n\geq 0$. 
It gives an embedding $\Hom(\cS, \R^{\mult}_{\geq 0})\subset \R^n$. We say an {\it $\R$-valued function on an open set of $\Hom(\cS, \R^{\mult}_{\geq 0})$ is real analytic}
if it is locally a restriction of a real analytic function on an open set of $\R^n$. This defines $\cO$ and it is independent of the 
 the choice of the surjective homomorphism $\N^n \to \cS$. 
 
 The log structure $M$ is the one associated to the canonical embedding $\cS\to \cO$. $M_{>0}$ is the submonoid sheaf of $M$ generated by $\cO^\times_{>0}$ and the image of $\cS$.

\medskip

(4) {\it The compactified vector space.}
Let $V$ be a finite dimensional graded $\R$-vector space 
 $V=\bigoplus_{w\in \Z, w\leq -1} \; V_w$ of weight $\leq -1$. Then we have a real analytic manifold with corners $\bar V$ (Part I, Section 7). It is covered by two open sets $V$ and 
 $\bar V\smallsetminus \{0\}$. Here $V$ has the usual sheaf of real analytic functions and the trivial log structure, and $\bar V\smallsetminus \{0\}$ is described as follows. For $a\in \R_{>0}$ and $v\in V$, let $a\circ v = \sum_w a^wv_w\in V$ where $v_w$ denotes the component of $v$ of weight $w$.  By choosing a  real analytic closed submanifold $V^{(1)}$ of $V\smallsetminus \{0\}$ such that 
 $\R_{>0}\times V^{(1)}\to V\smallsetminus \{0\}\;;\;(a, v) \mapsto a\circ v$  is an isomorphism of real analytic manifolds, we have an isomorphism of real analytic manifolds with corners 
$$\R_{\geq 0}\times V^{(1)}\cong \bar V\smallsetminus \{0\}$$
extending the above isomorphism. We will denote this extended isomorphism as $(a,v)\mapsto a\circ v$. 

For example, in the cases $V=\cL$ and $V=\cL(F)$
(\ref{cL(F)}), we have the compactified vector spaces $\bar \cL$ and $\bar \cL(F)$, respectively. We can identify $\bar \cL(F)$ with the closure of $\cL(F)$ in $\bar \cL$.

\end{sbpara}

\begin{sbprop}\label{cv+0} Let $\cS$ be an fs monoid and consider the real toric variety $T:=\Hom(\cS, \R^{\mult}_{\geq 0})$. Then if $S$ is an object of $\cC_\R(\sat)$, we have a natural bijection between the set of all morphisms $S\to T$ in $\cC_\R(\sat)$ and the set of all homomorphisms $\cS\to \Gamma(S, M_{S,>0})$. 

\end{sbprop}

\begin{pf} Since $\cS\subset \Gamma(T, M_{T,>0})$, a morphism $S\to T$ induces $\cS\to \Gamma(S, M_{S,>0})$. 
It is easy to see that this correspondence is bijective. 
\end{pf}

\begin{sbpara}\label{pchar}

If $M$ is an fs log structure with sign, 
locally we have a chart $\cS\to M$ whose image is contained in $M_{>0}$. (Here $\cS$ is an fs monoid.)
In fact, if $\cS\to M$ is a chart, the composition $\cS\to M\cong M_{>0}\times \{\pm 1\}\to M_{>0}\subset M$ is also a chart. We will call such a chart $\cS\to M_{>0}$ a {\it positive chart}. 
\end{sbpara}

\begin{sbprop}\label{fiberpr} 

$(1)$  The category  $\cB'_\R(\log)$  has fiber products. 

\smallskip

$(2)$  A fiber product in $\cB'_\R(\log)$ is a fiber product in $\cC_\R(\sat)$.

  \end{sbprop}
   (1) is proved in Part II, Proposition 3.1.7. We give here a proof which proves both (1) and (2).

\begin{pf} 
For a 
diagram $S_1\to S_0\leftarrow S_2$ in $\cB'_\R(\log)$, locally on $S_0, S_1, S_2$, we can find fs monoids $\cS_0, \cS_1, \cS_2$ with homomorphisms $\cS_1\leftarrow \cS_0\to \cS_2$ and a morphism $\iota_j: S_j\to T_j:=  \Hom(\cS_j, \R^{\mult}_{\geq 0})$ of $\cB'_\R(\log)$ for each $j=0,1,2$, satisfying the following conditions (i) and (ii).

\smallskip

(i) The diagram
$$\begin{matrix}   S_1 &\to & S_0 &\leftarrow& S_2\\
\downarrow &&\downarrow && \downarrow \\
T_1&\to &T_0& \leftarrow & T_2\end{matrix}$$
is commutative. 

\smallskip

(ii) For each $j=0,1,2$, the  underlying map $S_j\to T_j$ of $\iota_j$ is injective, the topology and the log structure of $S_j$ with sign are induced from those of $T_j$, and the homomorphism $\iota_j^{-1}(\cO_{T_j}) \to \cO_{S_j}$ is surjective. 

\smallskip

This is proved by using positive charts (\ref{pchar}) on $S_0$, $S_1$, $S_2$ which are compatible. 

\smallskip

To prove \ref{fiberpr}, it is sufficient to prove that in this situation, we have the fiber product $S_3$ of $S_1\to S_0 \leftarrow S_2$ in $\cC_\R(\sat)$ which belongs to $\cB'_\R(\log)$.  Let $\cS_3$ be the pushout of the diagram $\cS_1\leftarrow \cS_0\to \cS_2$ in the category of fs monoids. This $\cS_3$ is obtained from the pushout $\cS_3'$ of 
$\cS_1\leftarrow \cS_0 \to \cS_2$ in the category of commutative monoids as follows. $\cS_3$ is the submonoid of $ (\cS_3')^{\gp}$ consisting of all elements $a$ such that for some integer $n\geq 1$, $a^n$ belongs to the submonoid of $(\cS_3')^{\gp}$ generated by the images of $\cS_1$ and $\cS_2$. Let $T_3$ be the real toric variety $\Hom(\cS_3, \R^{\mult}_{\geq 0})$, let $S_3'$ be the fiber product of $S_1\to S_0\leftarrow S_2$ 
in the category of topological spaces, and let $T'_3$ be the fiber product of $T_1\to T_0\leftarrow T_2$ which is identified with $\Hom(\cS'_3, \R^{\mult}_{\geq 0}) $ as a topological space.  As a topological space, we define $S_3$ as the fiber product of $S_3' \to T_3' \leftarrow T_3$. Let $\iota_3: S_3\to T_3$ be the canonical injection. We define the structure sheaf $\cO_{S_3}$ on $S_3$  as follows. For $j=0, 1,2$, let $I_j$ be the kernel of $\iota^{-1}_j(\cO_{T_j})\to \cO_{S_j}$. Let $I_3$ be the ideal of $\iota^{-1}_3(\cO_{T_3})$ generated by the images of $I_1$ and $I_2$. Define $\cO_{S_3}= \iota_3^{-1}(\cO_{T_3})/I_3$. Define the log structure with sign on $S_3$ as the inverse image of that of $T_3$. Then $S_3$ is clearly an object of $\cB'_\R(\log)$. 

We prove that $S_3$ is the fiber product of $S_1\to S_0\leftarrow S_2$ in $\cC_\R(\sat)$. By \ref{cv+0}, for an object $X$ of $\cC_\R(\sat)$ and $j=0,1,2,3$, a morphism $X\to S_j$ corresponds in one to one manner to a homomorphism $\cS_j\to \Gamma(X, M_{X,>0})$ such that the associated morphism $X\to T_j$ has the following two properties (1) and (2).
\medskip

(1) 
The image of the set $X$ in $T_j$ is contained in $S_j$.

\smallskip

(2) The image of $I_j$ in $\cO_X$ is zero. 
\smallskip

Since $\Gamma(X, M_{X,>0})$ is a saturated monoid, a homomorphism $\cS_3'\to \Gamma(X, M_{X,>0})$ and a homomorphism $\cS_3\to \Gamma(X, M_{X,>0})$ correspond in one to one manner. These prove that $S_3$ is the fiber product of $S_1\to S_0\leftarrow S_2$ in $\cC_\R(\sat)$. 
\end{pf}

  \begin{sbpara}\label{gest0} The proof of \ref{fiberpr} shows that the underlying topological space of a fiber product in $\cB'_\R(\log)$ need not be the fiber product of the underlying topological spaces. We consider this point.

   We call a homomorphism $\cS_0\to \cS_1$ of saturated commutative monoids (\ref{cblog}) 
 {\it universally saturated} if for any commutative monoid $\cS_2$ and any homomorphism $\cS_0\to \cS_2$, the pushout of $\cS_1\leftarrow \cS_0 \to \cS_2$ in the category of commutative monoids is saturated.

 For a morphism $S_1\to S_0$ of $\cB'_\R(\log)$, we say $f$ is universally saturated if for any $s_1\in S_1$ with image $s_0$ in $S_0$, the homomorphism $M_{S_0,s_0}\to M_{S_1,s_1}$ is universally saturated. (The last condition is equivalent to the condition that the homomorphism $(M_{S_0}/\cO^\times_{S_0})_{s_0} \to (M_{S_1}/\cO^\times_{S_1})_{s_1}$ is universally saturated.)

 The following can be proved easily: 
 Let  $f:S_1\to S_0$ be a morphism in $\cB'_\R(\log)$. Let the triple of homomorphisms $\cS_j\to M_{S_j}$ $(j=0,1)$ and $h:\cS_0\to \cS_1$ be a chart of $f$. Then, if $h$ is universally saturated, $f$ is universally saturated. Conversely, if $f$ is universally saturated, then locally on $S_0$ and $S_1$, there are positive charts (\ref{pchar}) and a homomorphism $h$ of charts as above such that $h$ is universally saturated. 
 
 \end{sbpara}
 
 \begin{sblem}
 Let $S_1\to S_0$ be a universally saturated morphism in
 $\cB'_\R(\log)$, let $S_2\to S_0$ be a morphism in $\cB'_\R(\log)$, and let $S_3$ be the fiber product of 
 $S_1\to S_0\leftarrow S_2$ in the category 
 $\cB'_\R(\log)$. Then the underlying topological space of 
 $S_3$ is the fiber product of the underlying topological spaces of $S_j$ $(j=0,1,2)$. 
 
 \end{sblem}
 
 This follows from the proof of \ref{fiberpr}.

 \begin{sbprop}\label{gest5}

 \smallskip
 
 (1) For $r\geq 1$, the homomorphism $\N\to \N^r\;;\;m\mapsto (m,m,\dots,m)$ is universally saturated.
 
 \smallskip
 
 (2) For any saturated commutative monoid $\cS$, the homomorphisms $\{1\}\to \cS$ and $\cS\to \{1\}$ are universally saturated.
 
 \smallskip
 
 (3) Let $\cS_j$ $(j=0,1,2)$ be saturated commutative monoids, let $\cS_0 \to \cS_1$ be a universally saturated homomorphism, let $\cS_0\to \cS_2$ be a homomorphism, and let $\cS_3$ be the pushout of $\cS_1\leftarrow \cS_0\to \cS_2$ in the category of commutative monoids. Then the homomorphism $\cS_2\to \cS_3$ is universally saturated.
 
 \smallskip
 
 (4) Let $\cS_j\to \cS_j'$ $(j=1, \dots, n)$ be universally saturated homomorphisms of saturated commutative monoids. Then the homomorphism $\prod_{j=1}^n \cS_j \to \prod_{j=1}^n \cS_j'$ is universally saturated.
 
  \smallskip
 
 (5) A homomorphism $\cS\to \cS'$ of saturated commutative monoids is universally saturated if and only if the induced homomorphism $\cS/\cS^\times \to \cS'/(\cS')^\times$ is unversally saturated.
 
 \smallskip
 
 (6) For a saturated commutative monoid $\cS$ and for $a\in \cS$, the canonical homomorphism $\cS\to \cS[1/a]$ is universally saturated. Here $\cS[1/a]$ denotes the submonoid $\{xa^{-n}\;|\;x\in \cS, n\geq 0\}$ of $\cS^{\gp}$. 
   \end{sbprop} 
   
   \begin{pf}
 The proofs of (1), (2), (5), (6) are easy. (3) is evident. We can prove (4) by induction on $n$  as follows. We may assume $n\geq 2$. Then the homomorphism between products in (4) is the composition $(\prod_{j=1}^{n-1} \cS_j) \times \cS_n \to (\prod_{j=1}^{n-1} \cS_j') \times \cS_n \to (\prod_{j=1}^{n-1} \cS_j')\times \cS_n'$ in which the first homomorphism is universally saturated by induction on $n$ and by (3), and the second homomorphism is universally saturated by (3). \end{pf}

   \begin{sbcor}\label{gest6} For a  diagram $S_1\to S_0\leftarrow S_2$ in $\cB'_\R(\log)$, the underlying topological space of  the fiber product is the fiber product of the underlying topological spaces in the following cases (i) and (ii).

  (i) The case where at least one of $S_1\to S_0$ and $S_2\to S_0$ is strict. 
  
  \medskip
  
  Here for a morphism $f: X\to Y$ of locally ringed spaces with log structures, we say $f$ is strict if the log structure of $X$ coincides with the inverse image of the log structure of $Y$ via $f$. 
  
  \medskip
  
  (ii) The case where the log structure of $S_0$ is trivial.

   \end{sbcor}

The following will be used for many times in this paper.

\begin{sbpara}\label{embstr} Let $X$ be an object of $\cB'_\R(\log)$ and let $Y$ be a subset of $X$. Assume that the following condition (C) is satisfied. 

(C) The  homomorphism from $\cO_X$ to the sheaf of $\R$-valued continuous functions on $X$ is injective.

Then we
have a structure on $Y$ as an object of $\cB'_\R(\log)$, which also satisfies (C), 
as follows. The topology of $Y$ is the one as a subspace of $X$. $\cO_Y$ is the sheaf of $\R$-valued functions on $Y$ which are locally restrictions of functions in $\cO_X$. The log structure with sign is the pullback 
of that of $X$. 

For an object $S$ of $\cB'_\R(\log)$ which satisfies (C), the map $\Mor(S, Y)\to \Mor(S,X)$ is injective and the image coincides with $\{f\in \Mor(S, X)\;|\; f(S)\subset Y\}$. 

\end{sbpara}

\subsection{Review on toric geometry}

We recall  toric varieties over a field and the real toric varieties associated to fans, by comparing them. 

\begin{sbpara}\label{rvtoric} Let $L$ be a finitely generated free abelian group and let $N:=\Hom(L, \Z)$. We will denote the group law of $L$ multiplicatively and that of $N$ additively.  

For a  rational finitely generated sharp cone $\sig$ in $N_\bR$, define an fs monoid $\cS(\sigma)$ by
$$\cS(\sig):=\{l \in L\,|\; l(\sigma)\geq 0\}.$$

For a rational fan  $\Sig$ in $N_\bR$, 
we have a toric variety $\toric_k(\Sig)$ over a field $k$ associated to $\Sig$ which is an fs log scheme over $k$, and a real toric variety $|\toric|(\Sig)$ which is an object of $\cB'_\R(\log)$. We review these. 

\end{sbpara}

\begin{sbpara}\label{torsig1}  The toric variety $\toric_k(\Sig)$ over $k$ is described as 
$$\toric_k(\Sig)= \bigcup_{\sig\in \Sig}  \Spec(k[\cS(\sigma)])\quad \text{(an open covering)}$$
where $k[\cS(\sig)]$ denotes the semigroup algebra of $\cS(\sigma)$ over $k$ and $\Spec(k[\cS(\sig)])$ is endowed with the standard log structure.  

It represents the contravariant functor from the category $(\fs/k)$ of fs log schemes over $k$ to the category of sets, which sends $S$ to the set of all homomorphisms $h:L\to M_S^{\gp}$ satisfying the following condition. 

\smallskip
(C) Let $s\in S$. Then there exists $\sig\in \Sig$ such that for any homomorphism $a:(M_S/\cO^\times_S)_{\bar s}\to \N$, the homomorphism $a\circ h: L\to \Q$ belongs to $\sig$. Here $\bar s$ is a geometric point over $s$.

\smallskip

Note that this condition is equivalent to the following condition.
%%% add a reference

\medskip

(C$^\prime)$ \'Etale locally on $S$, there is $\sig\in \Sig$ such that $h(\cS(\sigma)) \subset M_S$. 

\smallskip

The set $\toric_k(\Sig)(k)$ of all $k$-rational points of $\toric_k(\Sig)$ is identified with the set of pairs $(\sig, h)$ consisting of $\sig\in \Sig$ and a homomorphism $h:\cS(\sig)^\times \to k^\times$. 
The point corresponding to this pair is the element of $\Spec(k[\cS(\sig)])(k)=\Hom(\cS(\sig), k)$ 
 which sends $a\in \cS(\sig)^\times$ to $h(a)$ and sends $a\in \cS(\sig)\smallsetminus \cS(\sig)^\times$ to $0$. 

\end{sbpara}

\begin{sbpara}\label{torsig2} The real toric variety $|\toric|(\Sig)$ is described as 
$$|\toric|(\Sig)=\bigcup_{\sig \in \Sig} \Hom(\cS(\sig), \R^{\mult}_{\geq 0})\quad \text{(an open covering)}$$ where  $\Hom(\cS(\sig), \R^{\mult}_{\geq 0})$ is regarded as an object of $\cB'_\R(\log)$ as in \ref{2.3ex} (3). It represents the contravariant functor from $\cC_\R(\sat)$ to the category of sets, which sends $S$ to the set of all homomorphisms $h:L\to M_{S,>0}^{\gp}$ satisfying the following condition. 

\smallskip
(C) Let $s\in S$. Then there exists $\sig\in \Sig$ such that for any homomorphism $a:(M_S/\cO^\times_S)_{\bar s}\to \N$, the homomorphism $a\circ h: L\to \Q$ belongs to $\sig$. 

\smallskip

Note that this condition is equivalent to the following condition.
%%% add a reference

\smallskip

(C$^\prime)$
Locally on $S$, there is $\sig\in \Sig$ such that $h(\cS(\sigma)) \subset M_{S,>0}$.

\smallskip

The set $|\toric|(\Sig)$ is identified with the set of pairs $(\sig, h)$ consisting of $\sig\in \Sig$ and a homomorphism $h:\cS(\sig)^\times \to \R_{>0}$. 
The point corresponding to this pair is the element of $\Hom(\cS(\sig), \R^{\mult}_{\geq 0})$ which sends $a\in \cS(\sig)^\times$ to $h(a)$ and sends $a\in \cS(\sig)\smallsetminus \cS(\sig)^\times$ to $0$. 
By this understanding, we can regard $|\toric|(\Sig)$ as a 
closed subset of $\toric_\R(\Sig)(\R)$. 

\end{sbpara}

\begin{sbpara}\label{torsig3}
The set $|\toric|(\Sig)$ is also identified with the set of all pairs $(\sig, Z)$ consisting of $\sig\in \Sig$ and a subset $Z$ of $\Hom(L, \R^{\mult}_{>0})$ which is a $\Hom(L/\cS(\sig)^\times, \R^{\mult}_{>0})$-orbit.

  In fact, $(\sig, Z)$ corresponds to $(\sig, h)$ in \ref{torsig2} where 
   $h$ is the restriction of any element of $Z$ to $\cS(\sig)^\times$.

\end{sbpara}

\begin{sbpara}

 If $\Sig$ is finite and $\Sig'$ is a rational finite subdivision of $\Sig$, we have a proper surjective morphism 
 $\toric_k(\Sig')\to \toric_k(\Sig)$. In the case $k=\R$, this induces a morphism
  $|\toric|(\Sig')\to |\toric|(\Sig)$ which is proper and surjective.

\end{sbpara}

\begin{sbpara}\label{logmod}

A morphism $S'\to S$ in the category $(\fs/k)$ (resp.\  $\cB'_\R(\log)$) is called a {\it log modification} if locally on $S$, there are a homomorphism $\cS\to M_S$ (resp.\ $\cS\to M_{S,>0}$) with $\cS$ a sharp fs monoid and a rational finite subdivision $\Sig'$ of the fan $\Sig$ of all faces of the cone $\Hom(\cS, \R^{\add}_{\geq 0})\subset \Hom(\cS^{\gp}, \R^{\add})$ such that $S'$ is isomorphic over $S$ to $S\times_{\toric_k(\Sig)} \toric_k(\Sig')$ (resp.\ $S\times_{|\toric|(\Sig)} |\toric|(\Sig')$). 

The underlying map of topological spaces of a 
 log modification is proper and surjective.
 
 \end{sbpara}

\begin{sbpara}\label{logmod2} We introduce a functor $[\Sig]$ associated to a fan $\Sig$, and consider its relation to log modification. 

Let $L$ and $N$ be as in \ref{rvtoric}. For a rational fan $\Sig$  in $N_\R$, let $[\Sig]$ be the contravariant functor from $(\fs/k)$ 
 (resp.\  $\cB'_\R(\log)$) to the category of sets which sends $S$ to the set of all homomorphisms $h:L\to M^{\gp}_S/\cO^\times_S$ satisfying the condition (C) in \ref{torsig1} (resp.\ \ref{torsig2}). 
 In the present situation, (C) is equivalent to (C$^\prime$) with $M_S$ (resp.\ $M_{S,>0}$) replaced by $M_S/\cO_S^\times$.

Let $S$ be an object of $(\fs/k)$ (resp.\ $\cB'_\R(\log)$) and assume that we are given $h\in [\Sig](S)$. This induces a continuous map $S\to \Sig$ which sends $s\in S$ to the unique cone $\sig\in \Sig$ such that $\cS(\sig)\subset L$ coincides with the inverse image of $(M_S/\cO^\times_S)_s$ under $L\to (M^{\gp}_S/\cO^\times_S)_s$.

 Assume $\Sig$ is finite and let $\Sig'$ be a rational finite subdivision of $\Sig$. Then we have a morphism of functors $[\Sig']\to [\Sig]$. 
The contravariant functor 
$\Mor(-, S)\times_{[\Sig]} [\Sig']$ from $(\fs/k)$ (resp.\ $\cB'_\R(\log)$) to the category of sets is represented by a log modification $S'\to S$. In fact, locally  on $S$, $h:L\to M^{\gp}_S/\cO^\times_S$ lifts to a morphism $S\to \toric_k(\Sig)$ (1.4.2) (resp.\ $S\to |\toric|(\Sig)$ (1.4.3)) and this functor is represented by $S \times_{\toric_k(\Sig)} \toric_k(\Sig')$ (resp.\ $S\times_{|\toric|(\Sig)} |\toric|(\Sig')$).

\end{sbpara}

\begin{sbpara}\label{gest1} This 1.4.8 will be used in Section 2.4--Section 2.6. 
    
      Let $\cS_1$ be an fs monoid, let $T:=\Hom(\cS_1, \R^{\mult}_{>0})$, and let $Z$ be a $T$-torsor.

  The purpose of this \ref{gest1} is to introduce an object $\bar Z$ of $\cB'_\R(\log)$ and to give set-theoretical descriptions (1) below of log modifications of $\bar Z$.

      Let $\bar T:=\Hom(\cS_1, \R^{\mult}_{\geq 0}) \supset T$, and let $\bar Z:=Z\times^T \bar T$. 
  
  We regard $\bar Z$ as an object of $\cB'_\R(\log)$ as follows. Take $\br\in Z$. Then we have bijection $T\to Z\;;\;t\mapsto t\br$ and this induces a bijection $\bar T\to \bar Z$. Via the last bijection from the real toric variety $\bar T$, we obtain a structure of $\bar Z$ as an  
  object of $\cB'_\R(\log)$. This structure is independent of the choice of $\br$.

  We prepare notation. For $s\in \bar Z$, we define a subgroup $T(s)$ of $T$ and a $T(s)$-orbit $Z(s)$ inside $Z$ as follows.
  In the case $Z=T$ and hence $\bar Z=\bar T$, $s$ is a homomorphism $\cS_1\to \R^{\mult}_{\geq 0}$. In this case, let $T(s)$ be the subgroup of $T=\Hom(\cS_1, \R^{\mult}_{>0})$ consisting of all elements which kill $s^{-1}(\R_{>0})\subset \cS_1$, and let $Z(s)\subset T$ be the set of all elements of $\cS_1\to \R^{\mult}_{>0}$ whose restriction to $s^{-1}(\R_{>0})$ coincides with the homomorphism induced by $s$. Then $Z(s)$ is a $T(s)$-orbit.  In general, take $\br\in Z$,  
 consider the induced isomorphism $\bar Z \cong \bar T$,  let $t$ be the image of $s$ in $\bar T$, let $T(s):=T(t)$,  
 and let $Z(s)$ be the $T(s)$-orbit in $Z$ corresponding to the $T(t)$-orbit $Z(t)$ in $T$ via the isomorphism $Z\cong T$. Then $T(s)$ and $Z(s)$ are independent of the choice of $\br$.

  Consider $L$ and $N$  in \ref{rvtoric}, let $\sig$ be a rational finitely generated sharp cone in $N_\R$, and let $\Sig$ be the fan of  all faces of $\sig$.   Assume that we are given a homomorphism $\cS(\sig)\to \cS_1$. 
  
  Then we have a morphism of functors $\Mor(-, \bar Z)\to [\Sig]$ where $[\Sig]$ is as in \ref{logmod2}. This morphism is obtained as follows. The homomorphism $\cS(\sig) \to \cS_1$ induces $\Mor(-, \bar T)\to [\Sig]$.   
  Take $\br\in Z$. Then $\br$ gives an isomorphism $\bar Z\cong \bar T$ and hence the composite morphism $\Mor(-, \bar Z) \cong \Mor(-, \bar T) \to [\Sig]$. This composite morphism is independent of the choice of $\br$.

  Assume further that the homomorphism $\cS(\sig)\to \cS_1$ is universally saturated (\ref{gest0}). 
  
  Let $\Sig'$ be a rational finite subdivision of $\Sig$, and let $E$ be the log modification of $\bar Z$ which represents the fiber product $\Mor(-, \bar Z)\times_{[\Sig]} [\Sig']$  (\ref{logmod2}). We give a description of $E$ as a set.
  
  For $s\in \bar Z$ and for $\sig'\in \Sig'$ such that the image $\tau$ of $s$ in $\Sig$ coincides with the image of $\sig'$ in $\Sig$, let $T(s,\sig')$ 
  be the subgroup of $T(s)$  
 consisting of all elements whose image in $\Hom(L/\cS(\tau)^\times, \R^{\mult}_{>0})$ is contained in its subgroup $\Hom(L/\cS(\sig')^\times, \R^{\mult}_{>0})$. 
  Then we have:

  \smallskip
  
  (1) There is a canonical bijection between $E$ and the set of all triples $(s, \sig', Z')$ where $s\in \bar Z$, $\sig'$ is an element of $\Sig'$ whose image in $\Sig$ coincides with the image of $s$ in $\Sig$, and $Z'$ is a $T(s, \sig')$-orbit in $Z(s)$. 
 
  \smallskip

 In fact, if $Z=T$, then $E=  \bar T \times_{|\toric|(\Sig)} |\toric|(\Sig')$ and hence the bijection is given by \ref{torsig3}. 
   In general, for $\br\in Z$, if $t$ denotes the image of $s$ under the isomorphism $\bar Z\cong \bar T$, we have $T(s, \sig')=T(t, \sig')$ and the isomorphism $Z\cong T$ sends a $T(t, \sig')$-orbit in $T$ to a $T(s,\sig')$-orbit in $Z$, and the induced composite bijection from the set of triples $(s, \sig', Z')$ to $E$ is independent of the choice of $\br$.

  \end{sbpara}

\section{The new space $D^{\star}_{\SL(2)}$  of $\SL(2)$-orbits}\label{s:new}

  In Part II, we defined and studied the space $D_{\SL(2)}$ of SL(2)-orbits. 
  Here we introduce a variant $D^{\star}_{\SL(2)}$. It is an object of the category $\cB'_\R(\log)$ (\ref{cblog}).

  Recall that $D_{\BS}$ is  an object of  $\cB'_\R(\log)$, $D_{\SL(2)}$ has two structures $D^I_{\SL(2)}$ and $D^{II}
 _{\SL(2)}$ as objects of $\cB'_\R(\log)$, and  the identity map of $D_{\SL(2)}$ gives a morphism $D^I_{\SL(2)}\to D^{II}_{\SL(2)}$ of $\cB'_\R(\log)$. We will relate the three spaces $D^{\star}_{\SL(2)}$, $D^{II}_{\SL(2)}$ and $D_{\BS}$ in the following way. 
  These three spaces are not  connected directly, but as we will see in this Section 2, they are connected as in the diagram
  $$\begin{matrix}  D^{\star,+}_{\SL(2)}& \to &D^{\star}_{\SL(2)}&  \to &D^{\star,-}_{\SL(2)}&\leftarrow  & D^{\star,\BS}_{\SL(2)}\\ \downarrow &&&&&&\downarrow\\
  D^{II}_{\SL(2)}&&&&&& D_{\BS}\end{matrix}$$
  in $\cB'_\R(\log)$ in which the horizontal arrows are log modifications (\ref{logmod}) and the left vertical arrow is proper surjective.

  As will be seen in Section 3, this diagram will induce morphisms $$D^{\star}_{\SL(2),\val} \to D^{II}_{\SL(2),\val}, \quad D^{\star}_{\SL(2),\val}\to D_{\BS,\val}$$ of associated valuative spaces, which appeared 
  in Introduction, for log modifications induce isomorphisms of the associated valuative spaces $D^{\star,+}_{\SL(2),\val} \overset{\cong}\to D^{\star}_{\SL(2),\val}  \overset{\cong}\to D^{\star,-}_{\SL(2),\val}\overset{\cong}\leftarrow   D^{\star,\BS}_{\SL(2),\val}$.
  
  In the pure case, the arrows in $D_{\SL(2)}\leftarrow D^{*,+}_{\SL(2)}\to D^{\star}_{\SL(2)}\to D^{\star,-}_{\SL(2)}$ are isomorphisms. 
  
  In Section 2.1, we review SL(2)-orbits in the pure situation. In Section 2.2, we continue reviews on Part II. 
  In Section 2.3, we define the spaces $D^{\star}_{\SL(2)}$ and $D^{\star,-}_{\SL(2)}$. After preparations in Section 2.4, we connect $D^{\star}_{\SL(2)}$ and $D^{II}_{\SL(2)}$ in Section 2.5 by introducing the space $D^{\star,+}_{\SL(2)}$, and we connect $D^{\star,-}_{\SL(2)}$ and $D_{\BS}$ in Section 2.6 by introducing the space $D^{\star,\BS}_{\SL(2)}$. In Section 2.7, we show that our spaces of SL(2)-orbits belong to a full subcategory $\cB'_\R(\log)^+$ of $\cB'_\R(\log)$ consisting of nice objects.

\subsection{Review on SL(2)-orbits in the pure case}

 Let the setting be as in \ref{setting}, and assume that we are in the pure situation of weight $w$.

\begin{sbpara}\label{sim1}  

In this pure case, an {\it $\SL(2)$-orbit in $n$ variables} means a pair $(\rho, \varphi)$, where $\rho$ is a homomorphism $\SL(2,\C)^n\to G(\C)$ of algebraic groups defined over $\R$ and $\varphi$ is a holomorphic map ${\mathbb P}^1(\C)^n\to \Dc$, satisfying
\begin{align*} &\varphi(gz)=\rho(g)\varphi(z)\quad \text{for}\;\;g\in \SL(2,\C)^n\;\;\text{and}\;\;z\in {\mathbb P}^1(\C)^n,\\
&\varphi(\frak h^n) \subset D\quad (\text{$\frak h$ is the upper half plane $\{x+iy\;|\;x, y\in \R, y>0\}$}),\\
&\rho_*(\text{fil}^p_z({\frak s}{\frak l}(2,\C)^n)) \subset \text{fil}^p_{\varphi(z)}({\frak g}_{\C})\quad (z \in {\mathbb P}^1(\C)^n, p\in \Z).
\end{align*}
 Here $\rho_*$ denotes the homomorphism ${\frak s}{\frak l}(2, \C)^n\to {\frak g}_{\C}$ of Lie algebras induced by $\rho$, and $\text{fil}^{\bullet}_z$ and $\fil^{\bullet}_{\varphi(z)}$  are filtrations given by $z$ and $\varphi(z)$, respectively (see Part II,  2.1.2).

\end{sbpara}

\begin{sbpara}\label{pure2} Let $(\rho, \varphi)$ be an $\SL(2)$-orbit in $n$ variables. Define the {\it associated homomorphisms} $\tau, \tau^{\star}: {\mathbb G}^n_{m,\R}\to \Aut_\R(H_{0,\R})$ of algebraic groups as
$$\tau^{\star}(t)= \rho(g_1, \dots, g_n)\quad \text{where}\; t=(t_j)_{1\leq j\leq n}\; \text{and}$$
$$g_j =\begin{pmatrix} 1/\prod_{k=j}^n t_k & 0 \\ 0 & \prod_{k=j}^n t_k\end{pmatrix},$$
$$\tau(t)=(\prod_{j=1}^n t_j)^w\cdot \tau^{\star}(t).$$
The image of the homomorphism $\tau^{\star}$ is contained in $G_\R$.

For $1\leq j\leq n$, we define the increasing filtration $W^{(j)}$ on $H_{0,\R}$ as follows. 

We have $H_{0,\R}= \bigoplus_{1\leq j\leq n, k\in \Z} \; H_{0,\R}(j, k)$ where $H_{0,\R}(j, k)$ is the part of $H_{0,\R}$ on which the action $\tau$ of ${\mathbb G}_{m,\R}^n$ is given by $(t_{\ell})_{1\leq\ell \leq n} \mapsto t_j^k$. 

Define $W^{(j)}$ by $W^{(j)}_k= \bigoplus_{k'\leq k} \; H_{0, \R}(j,k')$. 

We call $W^{(j)}$ $(1\leq j\leq n)$ the {\it associated weight filtrations}. 

\end{sbpara}

\begin{sbpara}\label{pure1} Let $(\rho, \varphi)$ be an $\SL(2)$-orbit in $n$ variables. 

For $1\leq j\leq n$, the following conditions (i)--(iii) are equivalent.

\smallskip

(i) The $j$-th component $\SL(2, \C)\to G(\C)$ of $\rho$ is trivial.

\smallskip

(ii)  $\varphi$ factors through the projection ${\mathbb P}^1(\C)^n \to {\mathbb P}^1(\C)^{n-1}$ which removes the $j$-th component.

\smallskip

(iii) Either $j\geq 2$ and $W^{(j)}= W^{(j-1)}$, or $j=1$ and $W^{(1)}=W$ (that is, $W^{(1)}_w=H_{0,\R}$ and $W^{(1)}_{w-1}=0$).

\end{sbpara}

\begin{sbpara}\label{sl2eq1} We consider the following equivalence relation on $\SL(2)$-orbits. 

We say an $\SL(2)$-orbit in $n$ variables $(\rho, \varphi)$ is non-degenerate if there is no $j$ $(1\leq j\leq n)$ which satisfies the equivalent conditions in \ref{pure1}. 

For a non-degenerate  $\SL(2)$-orbit $(\rho,\varphi)$ in $n$ variables and for a non-degenerate $\SL(2)$-orbit $(\rho', \varphi')$ in $n'$ variables, $(\rho, \varphi)$ and $(\rho', \varphi')$ are equivalent if and only if $n=n'$ and there is $t\in \R^n_{>0}$ such that $$\rho'(g)=\tau^{\star}(t) \rho(g) \tau^{\star}(t)^{-1}, \quad \varphi'(z)=\tau^{\star}(t)\varphi(z)$$
for any $g\in \SL_2(\C)^n$ and $z\in {\bf P}^1(\C)^n$. 
Here $\tau^{\star}$ is the homomorphism associated to $(\rho,\varphi)$ in \ref{pure2}. We have the same equivalence relation when we replace $\tau^{\star}(t)$ in the above by $\tau(t)$  in \ref{pure2}   associated to $(\rho,\varphi)$.

Any $\SL(2)$-orbit uniquely factors through a non-degenerate $\SL(2)$-orbit, called the associated non-degenerate $\SL(2)$-orbit, which is described as below. 
Two $\SL(2)$-orbits are equivalent if and only if their associated non-degenerate $\SL(2)$-orbits are equivalent in the above sense. 

For an $\SL(2)$-orbit $(\rho, \varphi)$ in $n$ variables, the associated non-degenerate $\SL(2)$-orbit $(\rho', \varphi')$ is as follows. 
Let $J=\{a(1), \dots, a(r)\}$ $(a(1)<\dots <a(r))$ be the set of $j$ $(1\leq j\leq n)$ such that the $j$-th component of $\rho$ is non-trivial. Then $(\rho',\varphi')$ is the $\SL(2)$-orbit in $r$ variables  defined by
$$
\rho(g_1, \dots, g_n)=\rho'(g_{a(1)},\dots, g_{a(r)}) \quad \varphi(z_1,\dots, z_n)=\varphi'(z_{a(1)},\dots, z_{a(r)}).$$

This number $r$ is called the rank of the (equivalence class of the) SL(2)-orbit $(\rho,\varphi)$.

\end{sbpara}

\begin{sbpara}\label{psl2}  The set $D_{\SL(2)}$ is defined as the set of all equivalence classes of $\SL(2)$-orbits $(\rho,\varphi)$ such that all members of the set of weight filtrations associated to $(\rho,\varphi)$ (\ref{pure2}) are rational (that is, defined already on $H_{0,\Q}$).

$D$ is embedded in $D_{\SL(2)}$ as the set of classes of SL(2)-orbits of rank $0$.

\end{sbpara}

\begin{sbpara}\label{sl2eq2}  Let $p\in D_{\SL(2)}$. We define objects
$$\tau^{\star}_p, \;\tau_p, \;Z(p), \;\cW(p)$$
associated to $p$.

Let $n$ be the rank of $p$. 
Let $(\rho,\varphi)$ be a non-degenerate $\SL(2)$-orbit which represents $p$.
The homomorphism 
 $\tau^{\star}$ (resp.\ $\tau$) (\ref{pure2}) associated to $(\rho,\varphi)$ depends only on the class $p$ (it does not depend on the choice of $(\rho,\varphi)$). We denote it as $\tau^{\star}_p$ (resp.\ $\tau_p$).

  The subset  $$\{\varphi((iy_j)_{1\leq j\leq n}) \;|\; y_j\in \R_{>0} \;(1\leq j\leq n)\}=
  \tau^{\star}(\R^n_{>0})\varphi({\bold i})=
  \tau(\R^n_{>0})\varphi({\bold i})\subset D$$
   $({\bold i}:=(i,\dots,i)\in \frak h^n)$ 
   depends only on the class $p$. We denote it as $Z(p)$ and call it the {\it torus orbit associated to $p$}. 
     
 The family $\{W^{(j)}\; |\; 1\leq j\leq n\}$ of weight filtrations associated to $(\rho,\varphi)$  (\ref{pure2}) depends only on the class $p$. Let $\cW(p)=\{W^{(j)}\;|\; 1\leq j\leq n\}$ and call it the set of weight filtrations associated to $p$.
 It consists of $n$ elements (Part II, Proposition 2.1.13).

\end{sbpara}

\begin{sbpara}\label{pure3} $D_{\SL(2)}$ has a structure as an object of $\cB'_\R(\log)$. For this, see Part II, Section 3.2.

A basic property of the topology of $D_{\SL(2)}$ is that, if $p\in D_{\SL(2)}$ is the class of an SL(2)-orbit $(\rho,\varphi)$, $p$ is the limit of $\varphi(iy_1, \dots, iy_n)\in  D$ where $y_j\in \R_{>0}$ and $y_j/y_{j+1}\to \infty$ ($1\leq j\leq n$, $y_n$ denotes $1$).

\end{sbpara}

\subsection{Reviews on $D_{\SL(2)}(\gr^W)$ and $D_{\SL(2)}(\gr^W)^{\sim}$}

We now consider the mixed Hodge situation. We review the spaces $D_{\SL(2)}(\gr^W)$ and $D_{\SL(2)}(\gr^W)^{\sim}$ considered in Part II and prepare notation which we will use later. 

  Actually there was an error concerning the definition of $D_{\SL(2)}(\gr^W)^{\sim}$ in Part II. We correct it in \ref{correct}.

\begin{sbpara}
Let $$D_{\SL(2)}(\gr^W):=\prod_{w\in \Z} D_{\SL(2)}(\gr^W_w)$$
where $D_{\SL(2)}(\gr^W_w)$ denotes the space $D_{\SL(2)}$ (Section 2.1) for the graded quotient $\gr^W_w$. 
\end{sbpara}

\begin{sbpara}\label{sim2}

The set $D_{\SL(2)}(\gr^W)^{\sim}$ is defined as follows (cf.\ Part II, 3.5.1). 

By an {\it $\SL(2)$-orbit on $\gr^W$ of rank $n$}, we mean a family $(\rho_w, \varphi_w)_{w\in \Z}$ of $\SL(2)$-orbits $(\rho_w, \varphi_w)$ on $\gr^W_w$
 in $n$ variables  in the sense of \ref{sim1} satisfying the following condition (1).

\smallskip

(1) For each $1\leq j\leq n$, there is a $w\in \Z$ such that 
 the $j$-th component of $\rho_w$ is non-trivial. 

\smallskip

 The equivalence relation is defined as follows.

 For an $\SL(2)$-orbit $(\rho_w, \varphi_w)_w$ on $\gr^W$ of rank $n$, the homomorphisms $\tau, \tau^{\star}: {\mathbb G}_{m,\R}^n \to \Aut_\R(\gr^W_w)$ associated to the $\SL(2)$-orbit $(\rho_w, \varphi_w)$ in $n$ variables of weight $w$   
 for $w\in \Z$ (\ref{pure2}) define homomorphisms 
$${\tau}, {\tau}^{\star}: {\mathbb G}_{m,\R}^n \to \prod_{w\in \Z} \Aut_\R(\gr^W_w)$$
of algebraic groups, respectively. 

An $\SL(2)$-orbit $(\rho_w, \varphi_w)_w$ on $\gr^W$ of rank $n$ 
and an $\SL(2)$-orbit $(\rho'_w, \varphi'_w)_w$ on $\gr^W$ of rank $n'$ 
are {\it equivalent} if and only if $n'=n$ and $(\rho'_w(g))_w= \tau^{\star}(t)(\rho_w(g))_w\tau^{\star}(t)^{-1}$, $(\varphi'_w(z))_w=\tau^{\star}(t)(\varphi_w(z))_w$ for some $t\in \R^n_{>0}$. 
 (We have the same equivalence relation when we replace  $\tau^{\star}$ here by $\tau$.)

The set $D_{\SL(2)}(\gr^W)^{\sim}$ is defined as the set of all equivalence classes of $\SL(2)$-orbits $(\rho_w, \varphi_w)_w$ on $\gr^W$ such that the weight filtrations on $\gr^W_w$ associated to $(\rho_w,\varphi_w)$ are rational (i.e., defined over $\Q$) for any $w\in \Z$.

\end{sbpara}

\begin{sbrem}\label{correct} In the definition of $D_{\SL(2)}(\gr^W)^{\sim}$ in Part II, 3.5.1, we forgot to put the condition of the rationality of the associated weight filtrations. This error does not affect the rest of Part II. 

\end{sbrem}

\begin{sbpara}\label{sim10}
We have the embedding
$$D(\gr^W)\overset{\subset}\to D_{\SL(2)}(\gr^W)^{\sim}$$
by identifying $D(\gr^W)$ with the set of $\SL(2)$-orbits on $\gr^W$ of rank $0$.

We have a map 
$$D_{\SL(2)}(\gr^W)^{\sim}\to D_{\SL(2)}(\gr^W)\;;\;p\mapsto (p(\gr^W_w))_w$$
which sends the class $p$ of $(\rho_w,\varphi_w)_w$ to $(\text{the class $p(\gr^W_w)$ of $(\rho_w, \varphi_w)$})_w$.

\end{sbpara}

\begin{sbpara}\label{simst2} 

For $p\in D_{\SL(2)}(\gr^W)^{\sim}$, we define a finite set $\overline{\cW}(p)$ of increasing filtrations on $\gr^W=\prod_w \gr^W_w$ as follows. Let $(\rho_w,\varphi_w)_w$ be an $\SL(2)$-orbit on $\gr^W$ in $n$ variables which represents $p$, let $W^{(w, j)}$ $(w\in \Z, 1\leq j\leq n)$ be the $j$-th weight filtration on $\gr^W_w$ associated to the $\SL(2)$-orbit $(\rho_w,\varphi_w)_w$ on $\gr^W_w$ in $n$ variables, and let $W^{(j)}=\bigoplus_w W^{(w,j)}$. Let $\overline{\cW}(p):=
\{W^{(j)} \; |\; 1\leq j \leq n\}$. Then $\overline{\cW}(p)$ is independent of the choice of the representative $(\rho_w,\varphi_w)_w$ of $p$. 

By an {\it admissible set of weight filtrations on $\gr^W$} (Part II, 3.2.2),
 we mean a set of increasing filtrations on $\gr^W$ which coincides with the set $\overline{\cW}(p)$  of weight filtrations associated to some point $p$ of $D_{\SL(2)}(\gr^W)^{\sim}$.     

 An admissible set $\Phi$ of weight filtrations on $\gr^W$ has a natural structure of a totally ordered set (given by the {\it variance} of $W'(\gr^W)$ for $W'\in \Phi$\;; see Part II, 2.1.11 and 2.1.13). For any $p\in D_{\SL(2)}(\gr^W)^{\sim}$ of rank $n$ such that $\Phi=\overline{\cW}(p)$, if $(W^{(j)})_{1\leq j\leq n}$ denotes the family of weight filtrations associated to $p$, 
$W^{(j)}\leq W^{(k)}$ for this order if and only if $j\leq k$. By using this ordering, we will identify $\Phi$ with the totally ordered set $\{1,\dots,n\}$. By this, we will identify $\bG_m^{\Phi}$, $\Z^{\Phi}$, etc. with $\bG_m^n$, $\Z^n$, etc.

 Let $\overline{\cW}$ be the set of all admissible sets of weight filtrations on $\gr^W$.

 Let $\cW(\gr^W_w)$ be the set of all admissible sets of weight filtrations on $\gr^W_w$, that is,  $\cW(\gr^W_w)=\{\cW(p)\;|\; p\in D_{\SL(2)}(\gr^W_w)\}$ (\ref{sl2eq2}). 
 
  We have a map 
 $$\overline{\cW}\to \prod_w \cW(\gr^W_w)\;;\; \Phi\mapsto (\Phi(w))_w,$$
 where $ \Phi(w):=\{W'(\gr^W_w)\; |\; W'\in \Phi, W'(\gr^W_w)\neq W(\gr^W_w)\}$. This map sends $\overline{\cW}(p)$ for $p\in D_{\SL(2)}(\gr^W)^{\sim}$ to $(\cW(p(\gr^W_w)))_w$. 

\end{sbpara}

\begin{sbpara}\label{simst2.5} For $\Phi\in \overline{\cW}$ and $Q=(Q(w))_w\in \prod_w \cW(\gr^W_w)$ such that $\Phi(w)\subset Q(w)$ for any $w\in \Z$, let
$$\bG_m^{\Phi} \to \prod_{w\in \Z} \bG_m^{Q(w)}$$
be the homomorphism which sends $(t_{W'})_{W'\in \Phi}$ to $(t'_{w,j})_{w\in \Z, j\in Q(w)}$, where $t'_{w,j}$ is the product of $t_{W'}$ for all elements $W'$ of $\Phi$ such that $W'(\gr^W_w)=j$. 

If $p\in D_{\SL(2)}(\gr^W)^{\sim}$ and $p'=(p(\gr^W_w))_w\in D_{\SL(2)}(\gr^W)$, for $\Phi=\overline{\cW}(p)$ and $Q(w)= \cW(p(\gr^W_w))$ $(w\in \Z)$,  $\tau^{\star}_p$ coincides with the composition $\bG_m^{\Phi} \to \prod_w \bG_m^{Q(w)}\to G_\R(\gr^W)$ where the first arrow is as above and the second arrow is $\tau^{\star}_{p'}$. 

\end{sbpara}

\begin{sbpara}\label{sim5} Let $p\in D_{\SL(2)}(\gr^W)^{\sim}$ (resp.\ $p\in D_{\SL(2)}(\gr^W)$). We define objects $$S_p,\; X(S_p)^+, \; A_p,\; B_p,\; \bar A_p, \;\bar B_p,\; \tau^{\star}_p,\; \tau_p,\; \tilde \tau^{\star}_p,\; \tilde \tau_p,\; Z(p)$$
associated to $p$. 

Let $$S_p= \bG_m^{\overline{\cW}(p)}\quad (\text{resp.}\;\;  \prod_w \bG_m^{\cW(p_w)}).$$ 

Then the character group $X(S_p)$ of $S_p$ is identified with $\prod_w \Z^{\overline{\cW}(p)}$ (resp.\ $\prod_w \Z^{\cW(p_w)}$). We define the submonoid $X(S_p)^+$ of $X(S_p)$ as the part corresponding to $\N^{\overline{\cW}(p)}$ (resp.\ $\prod_w \N^{\cW(p_w)}$). 

Let $A_p$ be the connected component in $S_p(\R)$ which contains the unit element. We identify $$A_p= \Hom(X(S_p), \R^{\mult}_{>0})=  \R_{>0}^{\overline{\cW}(p)} \; (\text{resp.}\ \prod_w \R_{>0}^{\cW(p_w)}).$$  
Let $$\bar A_p=\Hom(X(S_p)^+, \R^{\mult}_{\geq 0})=\R_{\geq 0}^{\overline{\cW}(p)}\;  (\text{resp.}\ \prod_w \R_{\geq 0}^{\cW(p_w)})
\supset A_p,$$ $$\bar B_p=\R_{\geq 0}\times \bar A_p\supset B_p: = \R_{>0}\times A_p.$$
We regard $\bar A_p$ and $\bar B_p$ as real toric varieties (\ref{2.3ex} (3)).

We define homomorphisms
$$\tau_p, \tau^{\star}_p: S_p \to \prod_w \Aut_{\R}(\gr^W_w)$$
of algebraic groups over $\R$ and a subset $Z(p)$ of $D(\gr^W)$.

Assume first $p\in D_{\SL(2)}(\gr^W)^{\sim}$. 
For an $\SL(2)$-orbit $(\rho_w, \varphi_w)_w$ on $\gr^W$ of rank $n$
which represents $p$, the associated homomorphisms 
$\tau, {\tau}^{\star}: S_p={\mathbb G}_{m,\R}^{\overline{\cW}(p)}=\bG_m^n \to \prod_w \Aut_\R(\gr^W_w)$ 
depend only on $p$. We denote $\tau$ as $\tau_p$ and $\tau^{\star}$ as ${\tau}^{\star}_p$. 
The set $$Z(p):=\{(\varphi_w(iy_1, \dots, iy_n))_w \;|\; y_j\in \R_{>0}\; (1\leq j\leq n)\}$$ $$=\{\tau_p(t)(\varphi_w({\bf i}))_w\; |\; t\in A_p\}=\{\tau^{\star}_p(t)(\varphi_w({\bf i}))_w\; |\; t\in A_p\}$$ $$\subset D(\gr^W)=\prod_{w\in \Z} D(\gr^W_w)\quad \quad  (\text{here}\;\;  {\bf i}=(i, \dots, i)\in {\frak h}^n)$$ depends only on $p$.

Next for $p\in D_{\SL(2)}(\gr^W)$, define $\tau_p$ and $\tau^{\star}_p$ as $(\tau_{p_w})_w$ and $(\tau^{\star}_{p_w})_w$, respectively, and let $Z(p)= \prod_w Z(p_w)$ (\ref{sl2eq2}).

Both for  $p\in D_{\SL(2)}(\gr^W)^{\sim}$ and for $p\in D_{\SL(2)}(\gr^W)$, 
we call $Z(p)$  the {\it torus orbit} of $p$. It is an $A_p$-torsor.

We define extended homomorphisms
$${\tilde \tau}_p, {\tilde \tau}_p^{\star}: {\mathbb G}_m\times S_p\to \prod_w \Aut_\R(\gr^W_w),$$ 
for $t_0\in\bG_m$ and $t\in S_p$, by 
$${\tilde \tau}_p(t_0,t)=(t_0^w)_w\tau_p(t)=\tau_p(t)(t_0^w)_w,$$ $$ {\tilde \tau}_p^{\star}(t_0,t)=(t_0^w)_w\tau^{\star}_p(t)=\tau^{\star}_p(t)(t_0^w)_w.$$ 
Here $(t_0^w)_w$ acts on $\gr^W_w$ as the multiplication by $t_0^w$.

\end{sbpara}

 \begin{sbpara}\label{simst3}
 Let $\Phi\in \overline{\cW}$.
 
 By a {\it splitting of $\Phi$} (Part II, 3.2.3), we mean a homomorphism 
 $\alpha=(\alpha_w)_w: \bG_m^{\Phi} \to \prod_w\Aut_\bR(\gr^W_w)$ of algebraic groups
 over $\R$ such that, for any $W'\in \Phi$ and $k\in \Z$, $W'_k$ coincides with the sum of the parts of $\gr^W$ of $\alpha$-weight $m$ for all $m\in \Z^{\Phi}$ such that $m(W')\leq k$.

 For a splitting $\alpha$ of $\Phi$, let $\alpha^{\star}: \bG_m^{\Phi}\to G_\R(\gr^W)$ be the homomorphism whose $G_{\bR}(\gr^W_w)$-component $\alpha^{\star}_w$ is 
$t=(t_j)_{j\in \Phi}\mapsto (\prod_{j\in \Phi} t_j)^{-w}\cdot \alpha_w(t)$. 

Note that the actions of $\alpha(t)$ and $\alpha^{\star}(t)$ $(t\in \R^{\Phi}_{>0})$ on $D(\gr^W)$ are the same.

 A splitting of $\Phi$ exists:
 If $p\in D_{\SL(2)}(\gr^W)^{\sim}$ and $\Phi=\overline{\cW}(p)$, ${\tau}_p$ is a splitting of $\Phi$. 
  In this case, for $ \alpha= {\tau}_p$, $\alpha^{\star}$ in the above coincides with ${\tau}_p^{\star}$ in \ref{sim5}.

Let $Q=(Q(w))_w\in \prod_w \cW(\gr^W_w)$. 

By a {\it splitting of $Q$}, we mean a family $\alpha=(\alpha_w)_w$ where $\alpha_w$ is a splitting of $Q(w)$.  Let $\alpha^{\star}=(\alpha^{\star}_w)_w$.

\end{sbpara}

\begin{sbpara}\label{simst3.5} Let $\Phi\in \overline{\cW}$. 

By a {\it distance to $\Phi$-boundary} (Part II, 3.2.4), we mean a real analytic map $\beta\colon D(\gr^W)\to \R_{>0}^{\Phi}$ 
such that $\beta(\alpha(t)x) = t\beta(x)$ $(t \in \bR^{\Phi}_{>0}, x \in D(\gr^W))$ 
for any splitting $\alpha$ of $\Phi$. (The last condition is equivalent to $\beta(\alpha^{\star}(t)x)= t\beta(x)$ $(t\in \bR^{\Phi}_{>0}, x\in D(\gr^W)$.)

A distance to $\Phi$-boundary exists (Part II, 3.2.5). 

Let $Q=(Q(w))_w\in \prod_w \cW(\gr^W_w)$. 

By a {\it distance} to $Q$-boundary, we mean a family $(\beta_w)_{w\in \Z}$ where $\beta_w$ is a distance to $Q(w)$-boundary for the pure situation $\gr^W_w$. 
\end{sbpara}

\begin{sbpara}\label{simst4} In Part II, we endowed $D_{\SL(2)}(\gr^W)$ and $D_{\SL(2)}(\gr^W)^{\sim}$ with structures as objects of $\cB'_\R(\log)$. These spaces satisfy the condition (C) in \ref{embstr}, that is, the sheaf of real analytic functions  is a sub-sheaf of the sheaf of all $\R$-valued continuous functions. $D_{\SL(2)}(\gr^W)$ is just the product of $D_{\SL(2)}(\gr^W_w)$ (\ref{pure3}) in $\cB'_\R(\log)$.  The canonical map $D_{\SL(2)}(\gr^W)^{\sim}\to D_{\SL(2)}(\gr^W)$ (\ref{sim10})  is a morphism in $\cB'_\R(\log)$ and it is a log modification (\ref{logmod}) as is explained in  Part II,  3.5.9 and 3.5.10.

We review some properties of these spaces. 
\end{sbpara}

\begin{sbpara}
Let $p\in D_{\SL(2)}(\gr^W)^{\sim}$ (resp.\ $p\in D_{\SL(2)}(\gr^W)$) and let $\br\in Z(p)$ (\ref{sim5}). Then $p$ is the limit of $\tau_p(t)\br=\tau^{\star}_p(t)\br$ where $t\in A_p$ tends to $0\in \bar A_p$. Here $0\in \bar A_p$ denotes $(0,\dots,0)\in \R_{\geq 0}^{\Phi}$ where $\Phi=\overline{\cW}(p)$ (resp.\ $\prod_w \R_{\geq 0}^{Q(w)}$ where $Q(w)=\cW(p(\gr^W_w)))$ (\ref{simst2}) in the identifications $\bar A_p=
 \R_{\geq 0}^{\Phi}$ (resp.\ $\prod_w \R_{\geq 0}^{Q(w)}$) $\supset A_p= 
 \R_{>0}^{\Phi}$ (resp.\ $\prod_w \R_{>0}^{Q(w)}$).

\end{sbpara}

\begin{sbpara}\label{phipart}

 For $\Phi\in \overline{\cW}$, let  $$D_{\SL(2)}(\gr^W)^{\sim}(\Phi)= 
  \{p\in D_{\SL(2)}(\gr^W)^{\sim} \;|\; \overline{\cW}(p)\subset \Phi\}.$$

    For $Q=(Q(w))_{w\in \Z} \in \prod_{w\in \Z} \cW(\gr^W_w)$, let $$D_{\SL(2)}(\gr^W)(Q)= \{p\in D_{\SL(2)}(\gr^W)\;|\; \cW(p_w)\subset Q(w)\;\text{for all}\;w\in \Z\}.$$ 
    
    Then $D_{\SL(2)}(\gr^W)^{\sim}(\Phi)$ (resp.\ $D_{\SL(2)}(\gr^W)(Q))$ is open in $D_{\SL(2)}(\gr^W)^{\sim}$ (resp.\ $D_{\SL(2)}(\gr^W)$). 
When $\Phi$ (resp.\ $Q$) moves, these open sets cover $D_{\SL(2)}(\gr^W)^{\sim}$ (resp.\ $D_{\SL(2)}(\gr^W)$).

If $\Phi\in \overline{\cW}$, $Q=(Q(w))_w\in \prod_w \cW(\gr^W_w)$ and $\Phi(w) \subset Q(w)$ (\ref{simst2}) for any $w\in \Z$, then the map $D_{\SL(2)}(\gr^W)^{\sim}\to D_{\SL(2)}(\gr^W)$ induces a map $D_{\SL(2)}(\gr^W)^{\sim}(\Phi)\to D_{\SL(2)}(\gr^W)(Q)$.

          \end{sbpara}

\begin{sbpara}\label{logst1}
Let $\Phi\in \overline{\cW}$ (resp.\ $Q=(Q(w))_w\in \prod_w\cW(\gr^W_w)$) and let $\beta$ be a distance to $\Phi$-boundary (resp.\ $Q$-boundary). Then the map $\beta$ extends uniquely to a morphism 
$$\beta: D_{\SL(2)}(\gr^W)^{\sim}(\Phi) \to \R_{\geq 0}^{\Phi}\quad (\text{resp.}\; D_{\SL(2)}(\gr^W)(Q) \to \prod_w \R_{\geq 0}^{Q(w)})$$ of $\cB'_\R(\log)$. 
The log structure with sign of $D_{\SL(2)}(\gr^W)^{\sim}(\Phi)$ (resp.\ $D_{\SL(2)}(\gr^W)(Q)$) coincides with the inverse image of the canonical log structure with sign of $\R^{\Phi}_{\geq 0}$ (resp.\ $\prod_w \R_{\geq 0}^{Q(w)}$) (\ref{2.3ex} (1)). 

For a distance $\beta$ to $\Phi$-boundary (resp.\ $Q$-boundary), each component $\beta_j$ $(j\in \Phi)$ (resp.\ $\beta_{w,j}$ $(w\in \Z$, $j\in Q(w))$) of $\beta$ is a section of the log structure $M_S$ where $S= D_{\SL(2)}(\gr^W)^{\sim}(\Phi)$ (resp.\ $D_{\SL(2)}(\gr^W)(Q)$). We have a chart $\N^{\Phi}\to M_S$ (resp.\ $\prod_w \N^{Q(w)} \to M_S$) defined as $m\mapsto \prod_j \beta_j^{m(j)}$ (resp.\ $m\mapsto \prod_w \beta_{w,j}^{m(w,j)}$). The induced homomorphism from $\N^{\Phi}$ (resp.\ $\prod_w \N^{Q(w)}$) to $ M_S/\cO^\times_S$ is independent of the choice of $\beta$. If $\Phi=\overline{\cW}(p)$ (resp.\ $Q(w)= \cW(p_w)$) for $p\in S$, this induces an isomorphism from $\N^{\Phi}$ (resp.\ $\prod_w \N^{Q(w)}$) to $(M_S/\cO^\times_S)_p$. 

If $\Phi(w) \subset Q(w)$ (cf. \ref{simst2}) for any $w\in \Z$, we have a commutative diagram 
$$\begin{matrix} \prod_w \N^{Q(w)} &\to & M_S/\cO^\times_S & (S:=D_{\SL(2)}(\gr^W))\\
\downarrow &&\downarrow & \\
\N^{\Phi} &\to & M_{S'}/\cO_{S'}^\times & (S':=D_{\SL(2)}(\gr^W)^{\sim})
\end{matrix}$$
where the left vertical arrow is the homomorphism induced from the homomorphism  $\bG_m^{\Phi}\to \prod_w \bG_m^{Q(w)}$ (\ref{simst2.5}) on the character groups. 
\end{sbpara}

\begin{sbpara}\label{bab} 
Let $\Phi\in \overline{\cW}$ (resp.\ $Q\in \prod_w \cW(\gr^W_w)$), and
let $\alpha$ be a splitting of $\Phi$ (resp.\ $Q$) and let $\beta$ be a distance to $\Phi$-boundary (resp.\ $Q$-boundary). Then the map
$D(\gr^W)\to D(\gr^W) \;;\; x\mapsto \alpha(\beta(x))^{-1}x = \alpha^{\star}(\beta(x))^{-1}x$
extends uniquely to a morphism 
$$b_{\alpha, \beta}: D_{\SL(2)}(\gr^W)^{\sim} \;(\text{resp.\ $D_{\SL(2)}(\gr^W)$}) \to D(\gr^W)$$
(Part II, Proposition 3.2.6).
\end{sbpara}

\subsection{The space $D^{\star}_{\SL(2)}$} 

We define the space $D^{\star}_{\SL(2)}$ comparing it with the space $D^{II}_{\SL(2)}$ which we defined in Part II. We also define a related space $D^{\star,-}_{\SL(2)}$.

\begin{sbpara}\label{symb1} Let $D^{\star}_{\SL(2)}$  (resp.\ $D^{\star,-}_{\SL(2)}$, resp.\  $D_{\SL(2)}$) be the set of all pairs $(p, Z)$ where $p\in D_{\SL(2)}(\gr^W)^{\sim}$ (resp.\ $D_{\SL(2)}(\gr^W)$, resp.\ $D_{\SL(2)}(\gr^W)^{\sim}$) and $Z$ is a subset of $D$ satisfying  the following two conditions (i) and (ii).  

Denote $\tau^{\star}_p$ (resp.\ $\tau^{\star}_p$, resp.\ $\tau_p$) by $\alpha_p$ and denote $\tilde \tau^{\star}_p$ (resp.\ $\tilde \tau^{\star}_p$, resp.\ $\tilde \tau_p$) by $\tilde \alpha_p$. 

\smallskip

(i)  $Z$ is either 

\smallskip
(i.A) an $\alpha_p(A_p)$-orbit in $D$, or 

\smallskip

(i.B) an $\tilde \alpha_p(B_p)$-orbit in $D_{\nspl}$ (\ref{splW})

\smallskip
\noindent
for the lifted action (\ref{liftac}).

\smallskip

(ii) The image of $Z$ in $D(\gr^W)$ coincides with the torus  orbit $Z(p)$ (\ref{sim5}) of $p$.
 
 \smallskip
 
We call an element $(p, Z)$ an {\it $A$-orbit} if it satisfies (i.A) and a {\it $B$-orbit} if it satisfies (i.B). 
This is similar to the case of $D_{\BS}$, which also consists of $A_P$-orbits and $B_P$-orbits 
for $\Q$-parabolic subgroups $P$ of $G_{\R}(\gr^W)$ (Part I, 5.1 and Definition 5.3). 

\end{sbpara}

\begin{sbpara}\label{slmap}

We  embed $D$ in $D^{\star}_{\SL(2)}$ (resp.\ $D^{\star,-}_{\SL(2)}$, resp.\ $D_{\SL(2)}$) by $F\mapsto (F(\gr^W), \{F\})$.

We have canonical maps $$D^{\star}_{\SL(2)}\to D_{\SL(2)}(\gr^W)^{\sim}, \;\;D^{\star,-}_{\SL(2)}\to D_{\SL(2)}(\gr^W),\;\;D_{\SL(2)}\to D_{\SL(2)}(\gr^W)^{\sim}$$
defined by $(p, Z) \mapsto p$.

We have a canonical map
$$D^{\star}_{\SL(2)}\to D^{\star,-}_{\SL(2)}\;;\;
(p, Z)\mapsto (p', Z'), \quad p':= (p(\gr^W_w))_w, \quad Z':= \tau^{\star}_{p'}(A_{p'})Z.$$
\end{sbpara}

\begin{sbpara}
The style of the definition of the set $D_{\SL(2)}$ in \ref{symb1} is slightly different from the one in Part II, Section 2.5. 
We explain the relation of the two styles. Let $(p, Z)\in D_{\SL(2)}$  in the present style, and let 
$(\rho_w, \varphi_w)_w$ be an $\SL(2)$-orbit on $\gr^W$ which represents $p$. 
If $(p, Z)$ is an $A$-orbit (\ref{symb1}), it
 is the class of $((\rho_w,\varphi_w)_w,\br)\in \cD'_{\SL(2), n}$ in Part II, 2.3.1 with $\br\in Z$. If $(p, Z)$ is a $B$-orbit (\ref{symb1}), it is 
the class of $((\rho'_w, \varphi'_w)_w, \br)\in \cD'_{\SL(2),n+1}$ in Part II, 2.3.1 where $\br \in Z$ and $\rho'_w$ (resp.\ $\varphi'_w$) is 
the composition $\SL(2,\R)^{n+1}\to \SL(2, \R)^n\to G_\R(\gr^W_w)$ (resp.\ 
${\mathbb P}^1(\C)^{n+1}\to {\mathbb P}^1(\C)^n \to D(\gr^W_w)$) of the projection to the last $n$ factors and $\rho_w$ (resp.\ $\varphi_w$).

\end{sbpara}

\begin{sbpara}\label{sm}

Let $D^{\star,\mild}_{\SL(2)}$ (resp.\ $D^{\star,-,\mild}_{\SL(2)}$) be the subset of $D^{\star}_{\SL(2)}$ (resp.\ $D^{\star,-}_{\SL(2)}$) consisting of all $A$-orbits. 

(We do not define the mild part of $D_{\SL(2)}$. The part of $A$-orbits in $D_{\SL(2)}$ does not fit our formulation of the mild part.)
\end{sbpara}

\begin{sbpara}\label{sit} Consider the following three situations (a)--(c).

(a) $\frak D:= D^{\star}_{\SL(2)}$, $\frak E= D_{\SL(2)}(\gr^W)^{\sim}$.

(b) $\frak D=D^{\star,-}_{\SL(2)}$,  $\frak E= D_{\SL(2)}(\gr^W)$.

(c) $\frak D= D_{\SL(2)}$,  $\frak E= D_{\SL(2)}(\gr^W)^{\sim}$.

We endow $\frak D$ with a structure of an object of $\cB'_\R(\log)$ as follows in \ref{2.3.6}--\ref{2.3.11}. In the situation (c), this coincides with the structure $D^{II}_{\SL(2)}$ treated in Part II. 

\end{sbpara}

\begin{sbpara}\label{2.3.6}
In the situations (a) and (c) (situation (b)) in \ref{sit}, 
for $\Phi\in \overline{\cW}$ (resp.\ $Q\in \prod_w \cW(\gr^W_w)$), let $\frak D(\Phi)$ (resp.\ $\frak D(Q)$) be the inverse image of $\frak E(\Phi)$ (resp.\ $\frak E(Q)$) (\ref{phipart}) in $\frak D$. 
\end{sbpara}

\begin{sbpara}\label{slspl} In the situations (a)--(c) in \ref{sit}, for $x=(p, Z)\in \frak D$, $\spl_W(\br)$ for $\br\in Z$ is independent of the choice of $\br$. We denote this $\spl_W(\br)$ $(\br\in Z)$ by $\spl_W(x)$.

\end{sbpara}

\begin{sbpara}
\label{2.3.9} 
 In the situations (a) and (c) (resp.\ situation (b)) in \ref{sit}, let $\Phi\in \overline{\cW}$ 
 (resp.\ $Q=(Q(w))_w\in \prod_w \cW(\gr^W_w))$,  let $\alpha$ be a splitting of $\Phi$ (resp.\ $Q$) (\ref{simst3}) and let $\beta$ be a distance to $\Phi$-boundary (resp, $Q$-boundary) (\ref{simst3.5}).

 In the situations (a) and (b) (resp.\ situation (c)), for $x\in D$, let $\delta_{\alpha,\beta}(x)\in \cL$ (\ref{cL(F)}) be $\Ad(\alpha^{\star}(\beta(p)))^{-1}\delta_W(x)$
  (resp.\ $\Ad(\alpha(\beta(p)))^{-1}\delta_W(x)$), where $p$ denotes the image of $x$ in $D(\gr^W)$ (for $\alpha^{\star}$, see \ref{simst3}). 
  Let $\frak D'=\frak D(\Phi)$ (resp.\ $\frak D(Q)$). Then, for $x=(p, Z)\in \frak D'$ and $\br\in Z$, $\delta_{\alpha,\beta}(\tau_p^{\star}(t)\br)$ (resp.\ 
  $\delta_{\alpha,\beta}(\tau_p(t)\br)$) converges in $\bar \cL$ (\ref{2.3ex} (4))
  when $t\in A_p$ tends to $0$ in $\bar A_p$, and the limit depends only on $x$ and is  independent of the choice of $\br$. We denote this limit by $\delta_{\alpha,\beta}(x)$. We have $\delta_{\alpha,\beta}(x) \in \bar \cL(b_{\alpha,\beta}(p))$, where $b_{\alpha,\beta}(p)$ is as in \ref{bab}.

 These $\delta_{\alpha,\beta}(x)$ and $b_{\alpha,\beta}(p)$ $(x=(p, Z))$ are  described as follows. In the situations (a) and (c) (resp.\ situation (b)), let $\alpha'$ and $(\alpha^{\star})'$ be the restrictions of $\alpha$ and $\alpha^{\star}$ (\ref{simst3}) to the subgroup $\bG_m^{\overline{\cW}(p)}$ (resp.\ $\prod_w \bG_m^{\cW(p_w)}$) 
 of $\bG_m^{\Phi}$ (resp.\ $\prod_w \bG_m^{Q(w)}$), respectively. 
 Since both $\alpha'$ and $\tau_p$ splits $\overline{\cW}(p)$ (resp.\ $(\cW(p_w))_w$), there is $u\in \prod_w \Aut_\R(\gr^W_w)$ such that for all $W'\in \overline{\cW}(p)$ (resp.\ for all $w\in \Z$ and all $W'\in \cW(p_w)$), $u$ preserves $W'$ and  induces the identity maps on $\gr^{W'}$, and such that 
 $$\tau_p(t) = u\alpha'(t)u^{-1}, \quad \tau^{\star}_p(t)=u(\alpha^{\star})'(t)u^{-1}$$
 for any $t\in \bG_m^{\overline{\cW}(p)}$ (resp.\ $\prod_w \bG_m^{\cW(p_w)}$). 
 Take $\br \in Z$ and let $\bar \br$ be the image of $\br$ in $Z(p)$ (\ref{sim5}). Then we have
 $$b_{\alpha,\beta}(p)=  b_{\alpha,\beta}(u^{-1}\bar \br),$$ $$\delta_{\alpha,\beta}(x)= \Ad(u\alpha^{\star}(\beta(u^{-1}\bar \br)))^{-1}\delta_W(\br) \;\;(\text{resp.}\; \Ad(u\alpha(\beta(u^{-1}\br)))^{-1}\delta_W(\br)). $$
 
 These are  shown in Part II, 3.3.9 in the situation (c). The proofs for the situations (a) and (b) are similar. 

\end{sbpara}

\begin{sbprop}\label{staran1} Consider the three situations in \ref{sit}.

In the situations (a) and (c), let $\Phi\in \overline{\cW}$ and $\frak D'=\frak D(\Phi)$, $\frak E'=\frak E(\Phi)$. In the situation (b), let $Q\in \prod_w \cW(\gr^W_w)$ and $\frak D'=\frak D(Q)$, $\frak E'=\frak E(Q)$.
In the situations (a) and (c) (resp.\ situation (b)), fix a
  splitting $\alpha$ of $\Phi$ (resp.\ $Q$) and a distance $\beta$ to $\Phi$-boundary (resp.\ $Q$-boundary).

Then we have a bijection 
$$\nu: \frak D'\to \{(p, s, \delta)\in \frak E'\times \spl(W) \times \bar \cL\;|\;
\delta\in \bar \cL(b_{\alpha,\beta}(p))\}
 $$
$(\bar \cL(\;)$ is as in \ref{2.3ex} (4))  defined as $x\mapsto (p, s, \delta)$, where $p$ is the image of $x$ in $\frak E$, $s=\spl_W(x)$ (\ref{slspl}), $\delta=\delta_{\alpha, \beta}(x)$ (\ref{2.3.9}).

  \end{sbprop}

\begin{pf} 
 The inverse map of $\nu$ is defined as $(p, s, \delta)\mapsto (p, Z)$ where $Z$ is as follows. 
Consider the situations (a) and (c) (resp.\ situation (b)).  Take $u\in G_\R(\gr^W)$ for $p$ as in \ref{2.3.9}. 
 
  In the case  $\delta\in \cL\subset \bar \cL$, 
  $Z$ is the  subset of  $D$ whose image in $D(\gr^W) \times \spl(W) \times \cL$ under the map in \ref{grsd} is 
 $$
 \{(\br, s, \Ad(u \alpha^{\star}(\beta(u^{-1}\br)))\delta) \;|\; \br \in Z(p)\}
 \quad 
 (\text{resp.}\; \{(\br, s, \Ad(u \alpha(\beta(u^{-1}\br)))\delta) \;|\; \br \in Z(p)\}).
 $$
 
  In the case $\delta=0\circ \delta'\in \bar \cL \smallsetminus \cL$ with $\delta'\in \cL\smallsetminus \{0\}$ (\ref{2.3ex} (4)), $Z$ is $\R_{>0}\circ Z'$ where $Z'$ is the above set $Z$  for $(p, s, \delta')$. 
 \end{pf}

\begin{sbprop}\label{staran2} Let the three situations be as in \ref{sit}.

$(1)$ In the situations (a) and (c) (resp.\ the situation (b)), endow  $\frak D(\Phi)$ (resp.\ $\frak D(Q)$) (\ref{2.3.6}) with a structure of an object of $\cB'_\R(\log)$ by using the bijection $\nu$ in Proposition \ref{staran1} (the target of $\nu$ is regarded as an object of $\cB'_\R(\log)$ by regarding it as $Y$ in $X= \frak E' \times \spl(W) \times \bar \cL$ in \ref{embstr}). Then this structure 
 is independent of the choice of $(\alpha,\beta)$. 

$(2)$ There is a unique structure on $\frak D$ as an object of $\cB'_\R(\log)$ such that for any $\Phi\in \overline \cW$ (resp.\ $Q\in \prod_w \cW(\gr^W_w)$), $\frak D':=\frak D(\Phi)$ (resp.\ $\frak D':=\frak D(Q)$) is an open subset and the restriction of this structure to $\frak D'$ coincides with the structure given in $(1)$. 
\end{sbprop}

\begin{pf} For the situation (c),  this follows from 
  Part II, Proposition 3.2.9 and Theorem 3.2.10. The proofs for the situations (a) and (b) are  similar. 
  \end{pf}
 
\begin{sbpara}\label{2.3.11} The structures of $D^{\star}_{\SL(2)}$,
 $D^{\star,-}_{\SL(2)}$, $D^{II}_{\SL(2)}$ as objects of $\cB'_\R(\log)$ are given by the situations (a), (b), (c) in Proposition \ref{staran2}, respectively.
 
In the situations (a)--(c) in \ref{sit},  the canonical map $\frak D\to \frak E$ is evidently a morphism of $\cB'_\R(\log)$. 
  \end{sbpara}

\begin{sbpara}\label{Lpart}  Via the bijection $\nu$ of Proposition \ref{staran1}, 
 $A$-orbits in $\frak D'$
 correspond to elements $(p,s, \delta)$ of the target of $\nu$  such that $\delta\in \cL$. 
 
 Hence the subset of $\frak D$ consisting all $A$-orbits is open in $\frak D$.

Elements $(p, Z)$ of $\frak D'$ such that $Z\subset D_{\spl}$ (1.1.5) correspond to elements $(p,s,\delta)$ of the target of $\nu$ 
 such that $\delta=0$. 
\end{sbpara}

\begin{sbpara}\label{thm+2}  Consider the situations (a)--(c) in \ref{sit}. In the situation (c), we consider the structure $D^{II}_{\SL(2)}$ of $D_{\SL(2)}$.

In Theorem \ref{ls1} below, we extend the result Part II, Theorem 3.4.4 on the local structure of $D^{II}_{\SL(2)}$ to all situations in \ref{sit}. This \ref{thm+2} is a preparation for it. 

Let $p\in \frak E$. We consider the local structure of $\frak D$ around the inverse image of $p$ in $\frak D$.

Consider the situations (a) and (c) (resp.\ situation (b)). 
Let $\Phi:=\overline{\cW}(p)$ (resp.\ $Q=(Q(w))_w$ with $Q(w):= \cW(p_w)$). Fix $\br\in Z(p)$.

Let $K_{\br}$ be the maximal compact subgroup of $G_\bR(\gr^W)$ associated to $\br$ (Part II, 3.4.1), and $K'_{\br}\subset K_{\br}$ be the isotropy subgroup of $G_\R(\gr^W)$ at $\br$. 

We use the notation in \ref{sim5}.
Let $R$ be an $\R$-subspace of ${\frak g}_\R(\gr^W)$ satisfying the following conditions (C1) and (C2).

\smallskip

(C1) 
${\frak g}_\R(\gr^W)= \Lie(\tau^{\star}(A_p)) \oplus R\oplus \Lie(K_{\br})$.

\smallskip

(C2)  $R=\sum_{m\in X(S_p)} R\cap (({\frak g}_\R)_m + ({\frak g}_\R)_{-m})$. Here $(-)_m$ denotes the part of weight $m$ for the adjoint action of $S_p$ via $\tau^{\star}_p$. (The definition of the part $(-)_m$ does not change if we replace $\tau^{\star}_p$ by $\tau_p$.)

\smallskip

 Let $S$ be an $\R$-subspace of $\Lie(K_{\br})$ such that $\Lie(K_{\br})=\Lie(K'_{\br})\oplus S$.

 For a subset $J$ of $\Phi$  (resp.\ for $J=(J(w))_{w\in \Z}$, $J(w)\subset Q(w)$), let $S_J$ be the subset of $S$ consisting of all elements $k$ such that 
  $\exp(k)\br \in (K_{\br}\cap G_{\bR,J}(\gr^W))\cdot \br$, where $G_{\bR,J}(\gr^W)$ is the subgroup of $G_{\bR}(\gr^W)$ consisting of all $g\in G_\bR(\gr^W)$ such that $g W' = W'$ for any $W'\in J$ (resp.\ for any $w\in \Z$ and any $W'\in J(w)$). 

\medskip

We define an object $Y$ of $\cB'_\R(\log)$ as follows. 
Let 
$$X= \bar A_p \times {\frak g}_\R(\gr^W)\times  {\frak g}_\R(\gr^W) \times {\frak g}_\R(\gr^W)\times S.$$ 
Note that $\bar A_p$ is $\R^{\Phi}_{\geq 0}$ (resp.\ $\prod_w \R_{\geq 0}^{Q(w)}$) (\ref{sim5}).

 Let $Y$ be the subset of $X$ consisting of all elements $(t,f,g,h,k)$ satisfying the following conditions (i)--(iv). 
In (ii) and (iv) below, let $J=\{j\in \Phi\;|\; t_j=0\}$ (resp.\ $J=(J(w))_{w\in \Z}$ with $J(w)=\{j\in Q(w)\;|\; t_{w,j}=0\}$). 
 
 For $\chi\in X(S_p)$, write $\chi=\chi_+(\chi_-)^{-1}$ with $\chi_+, \chi_-\in X(S_p)^+$ which are defined as follows. In the identification $X(S_p)= \prod_w \Z^{Q(w)}$, if we denote by $m(w,j)\in \Z$ the $(w,j)$-component 
  of $\chi$ for $w\in \Z$ and $j\in Q(w)$, then the $(w,j)$-component of $\chi_+$ is $\max(m(w,j),0)$ and the $(w,j)$-component of $\chi_-$ is $\max(-m(w,j),0)$. 

\smallskip
 
 (i) For any $\chi\in X(S_p)$, $t(\chi_+)g_{\chi}= t(\chi_-)f_{\chi}$ and $t(\chi_+)h_{\chi}= t(\chi_-)g_{\chi}$.

\smallskip

 Here $g_{\chi}$ etc. denotes the $\chi$-component for the adjoint action of $S_p$ via $\tau^{\star}_p$. $t(\chi_+), t(\chi_-) \in \R_{\geq 0}$ are defined by the understanding $\bar A_p= \Hom(X(S_p)^+, \R^{\mult}_{\geq 0})$. 
 
 \smallskip
 
(ii)  Let $\chi\in X(S_p)$. If $t(\chi_+)=0$, then $g_{\chi}=f_{\chi}=0$. If $t(\chi_-)=0$, then $g_{\chi}=h_{\chi}=0$.
 In other words, if $m(j)\in \Z$ for $j\in \Phi$ (resp.\ $m(w,j)\in \Z$ for $w\in \Z$ and $j\in Q(w)$) denotes the $j$ (resp.\ $(w,j)$)-component of $\chi$ in the identification $X(S_p)= \Z^{\Phi}$ (resp.\ $\prod_w \Z^{Q(w)}$), then $f_{\chi}=0$ unless $m(j)\leq 0$ for any $j\in J$ (resp.\ unless $m(w, j)\leq 0$ for any $w\in \Z$ and $j\in J(w)$), $g_m=0$ unless $m(j)=0$ for any $j\in J$ (resp.\ unless $m(w, j)= 0$ for any $w\in \Z$ and $j\in J(w)$), $h_m=0$ unless $m(j)\geq 0$ for any $j\in J$ (resp.\ unless $m(w, j)\geq 0$ for any $w\in \Z$ and $j\in J(w)$). 

\smallskip

(iii) $g_{\chi}\in R$ and $f_{\chi}+h_{\chi^{-1}}\in R$ for any $\chi\in X(S_p)$. 

\smallskip

(iv) $k\in S_J$. 

\smallskip

Regard $X$ as an object of $\cB'_\R(\log)$ in the natural way, and regard $Y\subset X$ as an object of $\cB'_\R(\log)$ by \ref{embstr}.

Let $$Y_0=\{(t,f,g,h,k)\in  Y\;|\; t\in A_p\}\subset Y.$$ 

\end{sbpara}

  \begin{sbthm}\label{ls1} Consider the three situations in \ref{sit}. Let the notation be as in \ref{thm+2}. 
  
  \smallskip
  (1) For a sufficiently small open neighborhood $U$ of $(0,0,0,0,0)$ in $Y$, there exists a unique open immersion $U\to \frak E$ in $\cB'_\R(\log)$ which sends $(t, f, g, h, k)\in U\cap Y_0$  to $$\exp(f){\tau}^{\star}_p(t) \exp(k) \br= \exp(f){\tau}_p(t) \exp(k) \br
 $$ of $D(\gr^W)\subset \frak E$.
  This morphism sends $(0,\dots, 0)\in Y$ to $p$.
  
  \smallskip
  (2) Let $\bar L=\bar \cL(\br)$, $L=\cL(\br)$. Then for a sufficiently small open neighborhood $U$ of $(0,0,0,0,0)$ in $Y$, there exists a unique open immersion $U\times \spl(W) \times \bar L\to \frak D$ in $\cB'_\R(\log)$ having the following property. In the situations 
 (a) and (b) (resp.\ situation (c)), 
 it sends
   $(t, f, g, h, k,s,\delta)\in Y\times \spl(W) \times L$, where $(t,f,g,h,k) \in U\cap Y_0$, $s\in \spl(W)$, and $\delta\in L$, to the element of $D$ whose image in 
$D(\gr^W) \times \spl(W) \times \cL$ under the isomorphism \ref{grsd} is 
$$(\exp(f){\tau}^{\star}_p(t) \exp(k) \br, s, \Ad(\exp(f)\tau^{\star}_p(t)\exp(k))\delta)$$ $$ (\text{resp}. \;
(\exp(f){\tau}_p(t) \exp(k) \br, s, \Ad(\exp(f)\tau_p(t)\exp(k))\delta)).$$ 

\smallskip

(3) For a sufficiently small open neighborhood $U$ of $(0,0,0,0,0)$ in $Y$, the diagram 
 $$\begin{matrix}  U\times \spl(W)\times \bar L &\to& \frak D\\
 \downarrow &&\downarrow \\
 U&\to& \frak E
 \end{matrix}$$
 is cartesian in $\cB'_\R(\log)$ and in the category of topological spaces.

\smallskip
(4) In the situations (a) and (c) (resp.\ situation (b)), the image of the map in (1) is contained in $\frak E(\Phi)$ (resp.\ $\frak E(Q)$) and the image of the map in (2) is contained in $\frak D(\Phi)$ (resp.\ $\frak D(Q)$), where $\Phi=\overline{\cW}(p)$ (resp.\ $Q=(\cW(p_w))_w$). 

\smallskip

(5) The underlying maps of the morphisms in (1) and (2) are described as in \ref{thm+4} below. 

 \end{sbthm}
  
  \begin{pf}
  In the situation (c), this is given in Part II, Theorem 3.4.4 and  3.4.12.  The proofs for the situations (a) and (b) are similar. 
\end{pf}

\begin{sbpara}\label{thm+4} The maps in (1) and (2) in Theorem \ref{ls1}  are induced from the maps 
$$Y \to \frak E,\quad Y\times \spl(W) \times \bar L \to \frak D,$$
respectively, defined as follows.

The first map sends $(t, f, g, h, k) \in Y$ to the following element $p'\in \frak E$:
  
  Assume we are in the situations (a) and (c) (resp.\ situation (b)). 
  Let $J=\{j\in \Phi\;|\; t_j=0\}$ (resp.\ $J=(J(w))_{w\in \Z}$ where $J(w)= \{j\in Q(w)\;|\; t_{w,j}=0\}$). 
   Define $p_J\in \frak E$ as follows. 
 Let $n=\sharp(\Phi)$ (resp.\ $n(w)= \sharp(Q(w))$ for $w\in \Z$).  Let $(\rho, \varphi)$ be the SL(2)-orbit on $\gr^W$ which represents $p$ (resp.\  $(\rho_w,\varphi_w)$ for $w\in \Z$  be the  $\SL(2)$-orbit on $\gr^W_w$ in $n(w)$ variables which represents $p_w$) such that $\br=\varphi(i,\dots,i)$ (resp.\ $\br_w=\varphi_w(i,\dots, i)$).
  Write $J=\{j_1,\dots, j_m\}$, $j_1<\dots< j_m$ (resp.\ $J(w)=\{j_{w,1}, \dots, j_{w,m(w)}\}$, $j_1<\dots<j_{m(w)}$). Then $p_J$ is the class of the following $\SL(2)$-orbit $(\rho',\varphi')$ on $\gr^W$ of rank $m$ (resp.\ the family $(\rho'_w,\varphi'_w)_w$ of $\SL(2)$-orbits in $m(w)$ variables). 
  $$\rho'(g_1, \dots, g_m)= \rho(g_1',\dots, g_n'), \quad \varphi'(z_1,\dots, z_m)= \varphi(z'_1,\dots,z'_{n(w)})$$
  $$(\text{resp.}\; \rho_w'(g_1, \dots, g_{m(w)})= \rho_w(g_1',\dots, g_{n(w)}), \quad \varphi'_w(z_1,\dots, z_{m(w)})= \varphi_w(z'_1,\dots,z'_{n(w)}).)$$
  Here $g'_j= g_k$ and $z'_j=z_k$  where $k$ is the smallest among integers $a$  such that $1\leq a\leq m$ (resp.\ $1\leq a\leq m(w)$) and $j\leq j_a$ if such $a$ exists, and $g'_j=1$ and $z'_j=i$ if such $a$ does not exist.

  Let $A'$ be the set of all elements $t'$ of $A_p$ 
  such that $t'_j=t_j$ for any $j\in \Phi\smallsetminus J$ (resp.\ $t'_{w,j}=t_{w,j}$ for any $w\in \Z$ and any $j\in Q(w)\smallsetminus J(w)$). 
  
  Then 
 $$p'= \exp(f)\tau_p(t')\exp(k)p_J$$ with $t'\in A'$. This $p'$ is independent of the choice of $t'\in A'$.
 
 Next the second map $Y\times \spl(W) \times \bar L\to \frak D$ sends $(t,f,g,h,k,s,\delta)$ to $(p', Z)\in \frak D$
 where $p'$ is as above and $Z\subset D$ is as follows. Consider the situations (a) and (b) (resp.\ situation (c)).

 If $\delta\in L\subset \bar L$,   $Z$ is the subset of $D$ whose image under the embedding $D\to D(\gr^W) \times \spl(W) \times \bar \cL$ in \ref{grsd} is the set of elements  $$(\exp(f) \tau^{\star}_p(t') \exp(k)\br, s, \Ad(\exp(f)\tau^{\star}_p(t') \exp(k)) \delta)$$ $$ (\text{resp.}\; 
  (\exp(f) \tau_p(t') \exp(k)\br, s, \Ad(\exp(f)\tau_p(t') \exp(k)) \delta))$$ where $t'$ ranges over all elements of $A'$. 
   If $\delta\in \bar L\smallsetminus L$ and $\delta=0\circ \delta^{(1)}$ for $\delta^{(1)}\in L\smallsetminus \{0\}$ (\ref{2.3ex} (4)), $Z$ is the subset of $D$ whose image under the embedding $D\to D(\gr^W) \times \spl(W) \times \bar \cL$ in \ref{grsd} is the set of elements 
  $$(\exp(f) \tau^{\star}_p(t') \exp(k)\br, s, \Ad(\exp(f)\tau^{\star}_p(t') \exp(k)) (c\circ \delta^{(1)})) $$ $$ (\text{resp.}\; (\exp(f) \tau_p(t') \exp(k)\br, s, \Ad(\exp(f)\tau_p(t') \exp(k)) (c\circ \delta^{(1)})))$$ where $t'$ ranges over all elements of $A'$ and $c$ ranges over all elements of $\R_{>0}$.

\end{sbpara}

  \medskip
  
The part for $D_{\SL(2)}$ of the following Proposition is Part II, Theorem 3.5.15.

  \begin{sbprop}\label{Lbund} Consider the situations in \ref{sit}. Fix any $F\in D(\gr^W)$ and let $\bar L=\bar \cL(F)$ (\ref{2.3ex}, (4)).  Then $\frak D$ is an $\bar L$-bundle over $\frak E\times \spl(W)$
   as an object of $\cB'_\R(\log)$.

Consequently, the map $\frak D\to \frak E \times \spl(W)$ is proper.
 \end{sbprop}

     \begin{pf} This follows from \ref{ls1}. 
    \end{pf}    
    
        Note that $\bar \cL(F)$ for all $F\in D(\gr^W)$ are isomorphic to each other as objects of $\cB'_\R(\log)$.

  \begin{sbprop}\label{2stars} The map $D_{\SL(2)}^{\star}\to D_{\SL(2)}^{\star,-}$ (\ref{slmap}) is a morphism of $\cB'_\R(\log)$. The following diagram is cartesian in $\cB'_\R(\log)$ and  also cartesian in the category of topological spaces. $$\begin{matrix}  
  D^{\star}_{\SL(2)} & \to & D^{\star,-}_{\SL(2)}\\\downarrow &&\downarrow\\
  D_{\SL(2)}(\gr^W)^{\sim} &\to & \;D_{\SL(2)}(\gr^W).
  \end{matrix}
  $$
  \end{sbprop}
  
     \begin{pf} We deduce this from Theorem \ref{ls1}. Let $p\in D_{\SL(2)}(\gr^W)^{\sim}$ and let $p'$ be the image of $p$ in $D_{\SL(2)}(\gr^W)$. Take $R$ and $S$ for the situation (b) in \ref{sit} as in \ref{thm+2} by using $p'$ as $p$ in \ref{thm+2}, and write this $R$ by $R'$. Let $C$ be an $\R$-subspace of $\fg_\R(\gr^W)$ such that $\Lie(\tau^{\star}_{p'}(A_{p'}))$ is the direct sum of $\Lie(\tau^{\star}_p(A_p))$ and $C$. Let $R=C\oplus R'$. Then $R$ and $S$ satisfy the conditions on $R$ and $S$ in \ref{thm+2} 
     for the situation (a) in \ref{sit} and for $p$. The homomorphism $S_p \to S_{p'}$ (\ref{simst2.5}) induces a homomorpjhism $X(S_{p'})^+\to X(S_p)^+$ and hence a morphism $\bar A_p=\Hom(X(S_p)^+, \R^{\mult}_{\geq 0}) \to \bar A_{p'}=\Hom(X(S_{p'})^+, \R^{\mult}_{\geq 0})$. Let $Y$ be the $Y$ in \ref{thm+2} defined by $(p, R, S)$ for the situation (a) in \ref{sit}, and let $Y'$ be the $Y$ in \ref{thm+2} defined by $(p', R', S)$ for the situation (b) in \ref{sit}. 
     
     For $(t,f,g,h,k)\in Y$, since $g\in R=C\oplus R'$, we can write $g=c+g'$ with $c\in C$ and $g'\in R'$ in a unique way, and we have $(t', f', g', h', k)\in Y'$ where $t'$ is the image of $t$ in $\bar A_{p'}$ and $f'=f-c$, $h'=h-c$. 
    We have a morphism $Y\to Y'$ which sends $(t,f,g,h,k)\in Y$ to $(t't'', f',g',h',k)\in Y'$ where $t''$ is the unique element of $A_{p'}$ such that $\tau^{\star}_{p'}(t'')=\exp(c)$. For a sufficiently small open neighborhood $U$ of $(0,0,0,0,0)$ in $Y$ and for a sufficiently small open neighborhood $U'$ of $(0,0,0,0,0)$ in $Y'$ such that the image of $U$ in $Y'$ is contained in $U'$, we have commutative diagrams
     $$\begin{matrix} U&\to & \frak E \\ \downarrow & &\downarrow \\ U' & \to & \frak E'\end{matrix} \quad\quad \begin{matrix}U\times \spl(W) \times \bar L &\to& \frak D\\ \downarrow &&\downarrow \\U'\times \spl(W) \times \bar L&\to & \frak D'\end{matrix}$$
     where $\frak E=D_{\SL(2)}(\gr^W)^{\sim}$, $\frak E'=D_{\SL(2)}(\gr^W)$, $\frak D=D^{\star}_{\SL(2)}$, $\frak D'=D^{\star,-}_{\SL(2)}$. This reduces Proposition \ref{2stars} to Theorem \ref{ls1}. 
     \end{pf}

\subsection{Basic facts on $\SL(2)$-orbits and Borel-Serre orbits}\label{ss:basic}

This Section \ref{ss:basic} is a preparation for the rest of Section 2. In \ref{BS1}--\ref{BS3}, we review the space $D_{\BS}$ defined and studied in Part I, and then in \ref{objass}--\ref{gest8} we give some basic facts about the spaces $D^{\star}_{\SL(2)}$, $D^{\star,-}_{\SL(2)}$,  $D^{II}_{\SL(2)}$, and $D_{\BS}$,.

\begin{sbpara}\label{BS1}

We review the definition of the set $D_{\BS}$ shortly (see Part I for details). 

Parabolic subgroups play central roles in the theory of Borel-Serre spaces. Following \cite{BS}, for a linear algebraic group $Z$ over a field, we call an algebraic subgroup $P$ of $Z$ a {\it parabolic subgroup} if it is geometrically connected and $Z/P$ is a projective variety. 

In our setting, there are  bijections 
$$\{\text{$\Q$-parabolic subgroup of $G$}\} \leftrightarrow \{\text{$\Q$-parabolic subgroup of $G(\gr^W)$}\}$$ 
$$
\leftrightarrow \{\text{family $(P_w)_{w\in \Z}$ of $\Q$-parabolic subgroups $P_w$ of $G(\gr^W_w)$}\}.$$
The bijection from the last set to the second set is given by $(P_w)_w\mapsto \prod_w P_w$, and the bijection from the second set to the first set is given by taking the inverse image under $G_\R\to G_\R(\gr^W)$.

 Let $P$ be a $\Q$-parabolic subgroup of $G_\R(\gr^W)$. Let $P_u$ be the unipotent radical of $P$, 
 let $S_P$ be the largest $\Q$-split torus in the center of $P/P_u$, and let $A_P$ (resp.\ $B_P$) be the connected component including $1$ of the topological group $S_P(\R)$ (resp.\ $(\bG_m\times S_P)(\R)$). 
 
 For each $p\in D(\gr^W)$, we have a canonical homomorphism $S_P\to P$ 
 of algebraic groups over $\R$ such that the composition $S_P\to P\to P/P_u$ 
 is the identify map, which we call the {\it Borel-Serre lifting at $p$} and denote
  by $t\mapsto t_p$. This $t_p$ is characterized by the following two properties.
  
  \smallskip
  
  (i) The image of $t_p$ in $P/P_u$ coincides with $t$.

  \smallskip

  (ii) $\theta_{K_p}(t_p)= t_p^{-1}$ where $\theta_{K_p}: G_\R(\gr^W)\to G_\R(\gr^W)$ denotes the Cartan involution associated to the maximal compact subgroup $K_p$ (cf.\ Part I, 2.1) of $G_\R(gr^W)$ associated to $p$. 
  
  \smallskip

  We have the following action of $B_P$ on $D$, which we call the Borel-Serre action and denote as $(b, F)\mapsto b\circ F$ $(b\in B_P$, $F\in D)$. For $b=(c,a)\in B_P$ with $c\in \R_{>0}$ and $a\in A_P$, we define $b\circ F:=  (c^w)_wa_{F(\gr^W)} F$, where $a_{F(\gr^W)}$ is the Borel-Serre lifting of $a$ at $F(\gr^W)$, 
   $(c^w)_w$ is the element of $\prod_w \Aut_\R(\gr^W_w)$ which acts on $\gr^W_w$ as the multiplication by $c^w$, and $(c^w)_wa_{F(\gr^W)}$ acts on $D$ by the lifted action \ref{liftac}. 
   
    The action of $A_P$ on $D$ and the action of $B_P$ on $D_{\nspl}$ are fixed point free.

 $D_{\BS}$ is defined as  the set of pairs $(P, Z)$ where $P$ is a $\Q$-parabolic subgroup of $G_\R(\gr^W)$ and $Z$ is either
 
 \smallskip
 
 (i) an $A_P$-orbit in $D$ or 
 
 \smallskip
 
 (ii) a $B_P$-orbit in $D_{\nspl}$
 
 \smallskip
 \noindent
 for the Borel-Serre action.
 
In the case (i), we call  $(P, Z)$ an {\it $A_P$-orbit}. In the case (ii), we call $(P, Z)$ a {\it $B_P$-orbit}.

We denote by $D^{\mild}_{\BS} $ the subset  of $D_{\BS}$ consisting of $A_P$-orbits. This subset was written as $D^{(A)}_{\BS}$ in Part I. 

\end{sbpara}

\begin{sbpara}\label{BS2} We review the structure of $D_{\BS}$ as an object of $\cB'_\R(\log)$ (actually it is a real analytic manifold with corners). 

For a $\Q$-parabolic subgroup $P$ of $G_\R(\gr^W)$, let
$$D_{\BS}(P)= \{(Q, Z)\in D_{\BS}\;|\; Q\supset P\}.$$
Then $D_{\BS}(P)$ forms an open covering of $D_{\BS}$ when $P$ varies. $D_{\BS}$ is also covered by the open sets $D^{\mild}_{\BS}$ (\ref{BS1}) and $D_{\BS,\nspl}$ where $D_{\BS,\nspl}$ denotes the subset of $D_{\BS}$ consisting of all elements $(P, Z)$ such that $Z\subset D_{\nspl}$.

The structures of  $D^{\mild}_{\BS}(P): =D_{\BS}(P)\cap D^{\mild}_{\BS}$ and $D_{\BS,\nspl}(P):= D_{\BS}(P)\cap D_{\BS,\nspl}$ as objects of $\cB'_\R(\log)$ are described as follows.

Let $X(S_P)$ be the character group of $S_P$,  and let 
$\Delta(P) \subset X(S_P)$ be the set of simple roots (\cite{BS}).  This set $\Delta(P)$  is characterized by the following two properties (i) and (ii).

\smallskip

(i) Let $n$ be the rank of $S_P$. Then $\Delta(P)$ is of order $n$ and generates $ \Q\otimes X(S_P)$ over $\Q$.

\smallskip

(ii) Let $X(S_P)^+$ be the submonoid of $X(S_P)$ generated by $\Delta(P)$. 
Lift $S_P$ to a subtorus of $P$. Then $X(S_P)^+$ coincides with the submonoid of $X(S_P)$ generated by $\chi^{-1}$ where $\chi$ ranges over all elements of $X(S_P)$ which appear in the adjoint action of $S_P$ on $\Lie(P)$.

\smallskip
Define a real toric variety  (\ref{2.3ex} (3)) $\bar A_P$ and $\bar B_P$ as
$$\bar A_P:=\Hom(X(S_P)^+, \R^{\mult}_{\geq 0})= \R_{\geq 0}^{\Delta(P)}\supset A_P=\Hom(X(S_P), \R^{\mult}_{>0})= \R_{>0}^{\Delta(P)},$$  
$$B_P:=\R_{\geq 0} \times \bar A_P\supset B_P= \R_{>0}\times A_P.$$

For a $\Q$-parabolic subgroup $Q$ of $G_\R(\gr^W)$ with $Q\supset P$, there is a canonical injection $\Delta(Q) \to \Delta(P)$, and $Q\mapsto \Delta(Q)\subset \Delta(P)$ is a bijection from the set of all $\Q$-parabolic subgroups of $G_\R(\gr^W)$ such that $Q\supset P$ to the set of all subsets of $\Delta(P)$. This is explained as follows. 

For such $Q$, we have $Q_u\subset P_u$, the composition $S_Q\to Q/Q_u\to Q/P_u$ is injective, and the image of this composite map is contained in $S_P\subset P/P_u\subset Q/P_u$. Hence $A_Q$ is regarded as a subgroup of $A_P$. There is a unique injection $\Delta(Q)\to \Delta(P)$ such that the composition $\R_{>0}^{\Delta(Q)}\cong A_Q\subset A_P=\R_{>0}^{\Delta(P)}$ coincides with the map $f\mapsto g$ where $g(j)=f(j)$ for $j\in \Delta(Q)$ and $g(j)=1$ for $j\in \Delta(P)\smallsetminus \Delta(Q)$.

We have bijections
$$D^{\mild}_{\BS}(P) \cong D\times^{A_P} \bar A_P, \quad
D_{\BS,\nspl}(P)\cong D_{\nspl} \times^{B_P} \bar B_P$$
which sends the element  $(Q, Z)$ of $D^{\mild}_{\BS}(P)$ (resp.\ $D_{\BS,\nspl}(P)$) to the class of $(z, h)$ (resp.\ $(z, \tilde h)$) where $z\in Z$ and $h\in \bar A_P=\R_{\geq 0}^{\Delta(P)}$ (resp.\ $\tilde h=(0, h)\in \bar B_P=\R_{\geq 0}\times \R_{\geq 0}^{\Delta(P)}$) is defined by $$h(j)=0\;\text{for}\;  j\in \Delta(Q)\subset \Delta(P), \quad h(j)=1\;\text{for}\; j\in \Delta(P)\smallsetminus \Delta(Q).$$
The right-hand sides of these bijections are regarded as objects of $\cB_\R'(\log)$ (Part I, Section 8) as is explained below, and the left-hand sides have the structures as objects of $\cB'_\R(\log)$ for which these bijections are isomorphisms of $\cB'_\R(\log)$.  

There is a closed real analytic sub-manifold $D^{(1,A)}$ (resp.\ $D^{(1,B)}$) of $D$ (resp.\ $D_{\nspl}$) such that 
we have an isomorphism $A_P\times D^{(1,A)}\overset{\cong}\to D$ (resp.\ $B_P\times D^{(1,B)}\overset{\cong}\to D_{\nspl}$), $(a, F)\mapsto a\circ F$, of real analytic manifolds. This induces a bijection $\bar A_P \times D^{(1,A)}\to D\times^{A_P} \bar A_P$ (resp.\ $\bar B_P\times D^{(1,B)} \to D_{\nspl}\times^{B_P} \bar B_P$) and by this, $D\times^{A_P} \bar A_P$ (resp, $D_{\nspl}\times^{B_P} \bar B_P$) has a structure of an object of $\cB'_\R(\log)$. This structure is independent of the choice of $D^{(1,A)}$ (resp.\ $D^{(1,B)}_{\nspl}$).

\end{sbpara}

\begin{sbpara}\label{BS3} The definition of the set $D_{\BS}$ can be rewritten in the style which is similar to the  definitions of the spaces of $\SL(2)$-orbits in Section 2.3.

Let $D_{\BS}(\gr^W)=\prod_{w\in \Z} D_{\BS}(\gr^W_w)$ where $D_{\BS}(\gr^W_w)$ is the space $D_{\BS}$ for the graded quotient $\gr^W_w$. For $p=(P_w, Z_w)_{w\in \Z}\in D_{\BS}(\gr^W)$, we denote $\prod_{w\in \Z} Z_w\subset D(\gr^W)$ as $Z(p)$. We call $Z(p)$ the {\it torus orbit} of $p$ and we call $\prod_{w\in \Z} P_w\subset G_\R(\gr^W)$ the $\Q$-parabolic subgroup of $G_\R(\gr^W)$ associated to $p$. Then,  $D_{\BS}$ is understood as the set of pairs $(p, Z)$ where $p\in D_{\BS}(\gr^W)$ and $Z$ is a subset of $D$ satisfying the following conditions (i) and (ii). 

\smallskip

(i)  $Z$ is either 

\smallskip

(i.A) an $A_P$-orbit in $D$ for the Borel--Serre action, or 

\smallskip

(i.B) a $B_P$-orbit in $D_{\nspl}$ for the Borel--Serre action. 

\smallskip

Here $P$ is the $\Q$-parabolic subgroup of $G_\R$ associated to $p$.
\smallskip

(ii) The image of $Z$ in $D(\gr^W)$ coincides with the torus  orbit $Z(p)$ of $p$.

\end{sbpara}

\begin{sbpara}\label{sit2} In the rest of this Section 2.4, we consider the situations (a)--(c) in \ref{sit} and also the situation 

\smallskip

(d)  $\frak D= D_{\BS}$, $\frak E=D_{\BS}(\gr^W)$. 

\end{sbpara}

\begin{sbpara}\label{objass}

For $x\in \frak D$, we define objects
$$S_x,\; X(S_x)^+, \; T(x),\;  \bar T(x), \; Z(x),\;  \bar Z(x)$$
associated to $x$.

In the situations (a)--(c), write
  $x=(p, Z)$  $(p\in \frak E$, $Z\subset D)$. In the situation (d), write $x=(P, Z)$. 
  
 In the situations (a)--(c),  let $S_x=S_p$ if $x$ is an $A$-orbit, and let $S_x=\bG_m\times S_p$ if $x$ is a $B$-orbit (\ref{sim5}, \ref{symb1}). In the situation (d), let $S_x=S_P$ if $x$ is an $A_P$-orbit, and let $S_x= \bG_m\times S_P$ if $x$ is a $B_P$-orbit (\ref{BS1}, \ref{BS2}). 
 
 We define a submonoid $X(S_x)^+$ of the character group $X(S_x)$ of $S_x$, as follows. In the situations (a)--(c), let $X(S_x)^+:=X(S_p)^+$ if $x$ is an $A$-orbit (\ref{sim5}), and let $X(S_x)^+:=\N\times X(S_p)^+\subset \Z\times X(S_p)=X(S_x)$ if $x$ is a $B$-orbit. In the situation (d), let $X(S_x)^+:=X(S_P)^+$ if $x$ is an $A_P$-orbit, and let $X(S_x)^+:=\N\times X(S_P)^+\subset \Z\times X(S_P)=X(S_x)$ if $x$ is a $B_P$-orbit, where $X(S_P)^+$ is as in \ref{BS2}.
 
 Let $T(x)$ be the connected component of $S_x(\R)$ containing the unit element. Let 
 $$\bar T(x) := \Hom(X(S_x)^+, \R^{\mult}_{\geq 0})\supset T(x) = \Hom(X(S_x)^+, \R^{\mult}_{>0}).$$
 We regard $\bar T(x)$ as a real toric variety.
 
 Define $Z(x):=Z$.  We call $Z(x)$ the torus orbit associated to $x$. 
 
 $T(x)$ acts on $Z(x)$ and $Z(x)$ is a $T(x)$-torsor. Let $\bar Z(x):= Z(x) \times^{T(x)} \bar T(x)$. 
 Then $\bar Z(x)$ has a unique structure of an object of $\cB'_\R(\log)$ such that for any $\br\in Z(x)$, the bijection $ \bar T(x)\to \bar Z(x)$ induced from the bijection $T(x)\to Z(x)\;;\;t\mapsto t\br$ becomes an isomorphism in $\cB'_\R(\log)$. We call $\bar Z(x)$ the extended torus obit associated to $x$. In \ref{gest2} below, we will embed $\bar Z(x)$ in $\frak D$ satisfying $x\in \bar Z(x)$.

\end{sbpara}

\begin{sbpara}\label{logstalk0}

This \ref{logstalk0} is a preparation for the next \ref{logstalk}. Consider the three situations in \ref{sit}.

In the situations (a) and (b) (resp.\ situation (c)), we have a global section $\beta_0^{\star}$ (resp.\ $\beta_0$) of $M_{\frak D}/\cO^\times_{\frak D}$ defined as follows. 

In the situations (a) and (c) (resp.\ situation (b)), let $\Phi \in \overline{\cW}$
(resp.\  $Q=(Q(w))_w \in \prod_w \cW(\gr^W_w)$), let $\frak D'= \frak D(\Phi)$ (resp.\ $\frak D'=\frak D(Q)$), 
let $\alpha$  be a splitting of $\Phi$ (resp.\ $Q$), and let $\beta$ be a distance to $\Phi$-boundary 
(resp.\ $Q$-boundary).
 Fix a real analytic closed sub-manifold $\cL^{(1)}$ of $\cL\smallsetminus \{0\}$ 
 such that $\R_{>0}\times \cL^{(1)}\to \cL\smallsetminus \{0\}\;;\; (a,\delta) \mapsto a\circ \delta$ 
 is an isomorphism of real analytic manifolds, and let 
 $\R_{\geq 0}\times \cL^{(1)}\overset{\cong}\to \bar \cL\smallsetminus \{0\}$ be the induced isomorphism in $\cB'_\R(\log)$.

Let $\frak D'_{\nspl}$ be the open subset of $\frak D'$ defined by $\delta\ne0$ via the bijection $\nu$ in Proposition \ref{staran1} associated to $(\alpha,\beta)$.
Then in  the situations (a) and (b) (resp.\ situation (c)), we have the composite morphism  $\frak D'_{\nspl}\to \bar \cL\smallsetminus \{0\} \cong \R_{\geq 0} \times \cL^{(1)} \to \R_{\geq 0}$ where the first arrow is $\nu$.
We denote this composite morphism $\cD_{\nspl}(\Phi) \to \R_{\geq 0}$ by $\beta_0^{\star}$ (resp.\ $\beta_0$). Then as is easily seen, this $\beta_0^{\star}$ (resp.\ $\beta_0$) belongs to $M_{\frak D_{\nspl}'}$, the class of $\beta_0^{\star}$ (resp.\ $\beta_0$) in $M_{\frak D_{\nspl}'}/\cO^\times_{\cD_{\nspl}'}$ is independent of the choices of $\alpha$, $\beta$, and $\cL^{(1)}$, this class extends uniquely to a section of 
$M_{\frak D'}/\cO^\times_{\frak D'}$ which is trivial on the part of $A$-orbits of $\frak D$, and this local section of $M_{\frak D}/\cO^\times_{\frak D}$ on $\frak D'=\frak D(\Phi)$ (resp.\ $\frak D'=\frak D(Q)$) extends, when $\Phi$ (resp.\ $Q$) moves, to a global section $\beta_0^{\star}$ (resp.\ $\beta_0)$ of $M_{\frak D}/\cO^\times_{\frak D}$ on $\frak D$ uniquely. 

\end{sbpara}

\begin{sbprop}\label{logstalk} Consider the four situations in \ref{sit2}. 
For  $x\in \frak D$, we  have a canonical isomorphism $$(M_{\frak D}/\cO^\times_{\frak D})_x\cong X(S_x)^+.$$ 

\end{sbprop}

\begin{pf} We first consider the situations (a)--(c).
 Write  $x=(p, Z)$. As in \ref{logst1}, we have a canonical isomorphism $(M_{\frak E}/\cO^\times_{\frak E})_p\cong X(S_p)^+$. In the case when $x$ is an $A$-orbit, we have  $(M_{\frak E}/\cO^\times_{\frak E})_p \overset{\cong}\to (M_{\frak D}/\cO^\times_{\frak D})_x$. If $x$ is a $B$-orbit, we have $\N\times (M_{\frak E}/\cO^\times_{\frak E})_p \overset{\cong}\to (M_{\frak D}/\cO^\times_{\frak D})_x$ where $1\in \N$ is sent to $\beta^{\star}_0$ in the situations (a) and (b) and to $\beta_0$ in the situation (c).

 We next consider the situation (d). Write $x=(P, Z)$. Assume first $x$ is an $A_P$-orbit. Consider the composite morphism $S:=D_{\BS}^{\mild}(P)\cong \bar A_P\times D^{(1,A)} \to \bar A_P=\R_{\geq 0}^{\Delta(P)}$ where the first isomorphism is as in \ref{BS2}. 
For $j\in \Delta(P)$, let $\beta_j: S\to \R_{\geq 0}$ be the $j$-component of this composite morphism. Then $\beta_j$ is a  section of $M_S$ and the class of $\beta_j$ in $M_S/\cO^\times_S$  is independent of the choice of $D^{(1,A)}$ in \ref{BS2}. We have a canonical isomorphism $X(S_x)^+ =\N^{\Delta(P)} \overset{\cong}\to (M_S/\cO^\times_S)_x$ which sends $m\in  \N^{\Delta(P)}$ to the class of $\prod_{j\in \Delta(P)} \beta_j^{m(j)}$. Assume next that $x$ is a $B_P$-orbit. 
Consider the composite morphism $S:=D_{\BS,\nspl}(P)\cong \bar B_P\times D^{(1,B)} \to \bar B_P=\R_{\geq 0}\times \R_{\geq 0}^{\Delta(P)}$ where the first isomorphism is as in \ref{BS2}. Let $\beta_0^{\BS}: S\to \R_{\geq 0}$ be the first component of this composite morphism, and 
for $j\in \Delta(P)$, let $\beta_j: S\to \R_{\geq 0}$ be the $j$-component of this composite morphism.
Then $\beta_0^{\BS}$ and $\beta_j$ $(j\in \Delta(P))$ are sections  of $M_S$ and their classes in $M_S/\cO^\times_S$ are independent of the choice of $D^{(1,B)}$ in \ref{BS2}. We have an isomorphism
$X(S_x)^+ \cong \N \times  \N^{\Delta(P)} \to (M_S/\cO^\times_S)_x$ which sends $(m_0, (m_j)_{j\in \Delta(P)})\in \N\times  \N^{\Delta(P)}$ to the class of $(\beta_0^{\BS})^{m_0}\cdot \prod_{j\in \Delta(P)} \beta_j^{m(j)}$.
\end{pf}

  \begin{sbpara}\label{gest2} Let the situations (a)--(d) be as in \ref{sit2}.

 Let $x\in \frak D$. 
     The inclusion map $Z(x) \to D$ extends uniquely to a morphism $$\bar Z(x) \to \frak D$$ of $\cB'_\R(\log)$. This morphism is described as follows. 
     
     Assume first we are in one of the situations (a)--(c). Write $x=(p, Z)$ and fix $\br\in Z(p)$. Consider the morphism
    $Y\times\spl(W)\times\bar L \to \frak D$ in \ref{thm+4} defined for $(p, \br, R, S)$ by fixing $R$ and $S$ \ref{thm+2}. Then the morphism $\bar Z(x) \to \frak D$ is the composite morphism $\bar Z(x) \to Y\times\spl(W)\times\bar L \to \frak D$ where the first morphism is as follows. 
   Let $F$ be an element of $Z(x)$ whose image under the embedding $D\to D(\gr^W)\times\spl(W) \times \cL$ is $(\br, s, \delta)$. Let $t\in \bar A_p$. Then the first morphism sends $(F, t)\in \bar Z(x)=Z(x) \times^{T(x)} \bar T(x)$ to $(t, 0,0,0,0, s, \delta)\in Y\times\spl(W)\times\bar L$ and if $x$ is a $B$-orbit, for $(c,t)\in \bar B_p$ $(c\in \R_{\geq 0})$, the first morphism sends $(F,(c,t))\in \bar Z(x)$ to  $(t, 0,0,0,0,s, c\circ \delta)\in Y\times\spl(W)\times\bar L$. 
   
   Next assume we are in the situation (d). Write $x=(P, Z)$. If $x$ is an $A_P$-orbit, this morphism $\bar Z(x)\to \frak D$ is the composition $\bar Z(x) = Z \times^{A_P} \bar A_P\subset  D\times^{A_P} \bar A_P \cong D_{\BS}^{\mild}(P)$. If $x$ is a $B_P$-orbit, this morphism is the composition $\bar Z(x) = Z \times^{B_P} \bar B_P \subset  D_{\nspl} \times^{B_P} \bar B_P \cong D_{\BS,\nspl}(P)$.

    This morphism $\bar Z(x) \to \frak D$ is injective and strict (\ref{gest6}), and sends  $0\in \bar Z(x)$ to $x$. Here $0$ denotes the class of $(\br, 0)$ where $\br\in Z(x)$ and $0\in \bar T(x)$ is the homomorphism $(M_{\frak D}/\cO^\times_{\frak D})_x \to \R^{\mult}_{\geq 0}$ which sends any non-trivial element of $(M_{\frak D}/\cO^\times_{\frak D})_x$ to $0$. (Then $0\in \bar Z(x)$ is independent of the choice of $\br$.) We will identify $\bar Z(x)$ with its image in $\frak D$, which coincides with 
  the closure of $Z(x)$ in $\frak D$. 
 
\end{sbpara}

    \begin{sbpara}\label{gest80} Consider the situations (a)--(d) as in \ref{sit2}. 
    
   In \ref{gest8} below, we give descriptions of log modifications of $\frak D$ as sets by using the extended torus orbit $\bar Z(x) \subset \frak D$ associated to $x\in \frak D$ (\ref{gest2}), which we will use in Section 2.5 and Section 2.6. 

    Let $U$ be an open set of $\frak D$. 
    
    Let $L$ and $N$ be as in \ref{rvtoric}, let $\Sig$ be a finite rational fan in $N_\R$, and let $\Sig'$ be a rational finite subdivision of $\Sig$. 
    
    Let $\Mor( -, U)\to [\Sig]$ be a morphism of functors (\ref{logmod2}) such that for any $x\in U$, if $\sig$ denotes the image of $x$ in $\Sig$ (\ref{logmod2}), the homomorphism $\cS(\sig) \to (M_U/\cO_U^\times)_x$ is universally saturated. For $x\in U$ and $\sig'\in \Sig'$ whose images in $\Sig$ coincide, we define a subgroup $T(x, \sig')$ of $T(x)= \Hom((M_U^{\gp}/\cO^\times_U)_x, \R^{\mult}_{>0})$ as follows. Let $\sig$ be the image in $\Sig$. Then the homomorphism $L\to (M_U^{\gp}/\cO^\times_U)_x$ factors through $L/\cS(\sig)^\times$. 
     $T(x,\sig')$ is the inverse image of $\Hom(L/\cS(\sig')^\times, \R^{\mult}_{>0})\subset \Hom(L/\cS(\sig)^\times, \R^{\mult}_{>0})$ in $T(x)$.

    Let $U'\to U$ be the log modification which represents the functor $\Mor(-, U) \times_{[\Sig]} [\Sig']$ (\ref{logmod2}).
    
    \end{sbpara}

    \begin{sblem}\label{gest8} 
    Let the notation and the assumptions be as in \ref{gest80}. 
    
     There exists a canonical bijection
    between $U'$ and the set of all triples $(x, \sig', Z')$ where $x\in U$, $\sig'$ is an element of $\Sig'$ whose image in $\Sig$ coincides with the image of $x$ in $\Sig$, and $Z'$ is an $T(x,\sig')$-orbit in $Z(x)$. 
    
    \end{sblem}

    \begin{pf} Let $x\in U$ and let $U''$ be the fiber product of $\bar Z(x) \to U\leftarrow U'$. Then the fiber on $x$ of $U'\to U$ coincides with the fiber on $x$ of $U''\to \bar Z(x)$. Since $U''$ represents the functor $\Mor(-,U)\times_{[\Sig]} [\Sig']$, this lemma follows from 
    \ref{gest1}. \end{pf}

\subsection{Relations with $D_{\SL(2)}$} 

We connect the spaces $D^{\star}_{\SL(2)}$ and $D^{II}_{\SL(2)}$ by introducing a new space $D^{\star,+}_{\SL(2)}$ of $\SL(2)$-orbits.

\begin{sbpara}\label{sl+def}  We define a log modification (\ref{logmod}) $$D^{\star,+}_{\SL(2)}\to D^{\star}_{\SL(2)}.$$

On $\frak D:=D^{\star}_{\SL(2)}$, there is a unique section $\beta_{\text{tot}}$ of $M_{\frak D}/\cO^\times_{\frak D}$ such that 
for any $\Phi\in \overline{\cW}$, the restriction of $\beta_{\text{tot}}$ to $\frak D(\Phi)$ coincides with the image of the product $\prod_{j\in \Phi} \beta_j$ in $M_{\frak D}/\cO_{\frak D}^\times$ where $\beta=(\beta_j)_{j\in \Phi}$ is a distance to $\Phi$-boundary. Let $\beta_0^{\star}$ be the section of $M_{\frak D}/\cO_{\frak D}^\times$ defined in \ref{logstalk0}. Consider the homomorphism
$\N^2\to M_S/\cO_S^\times\;;\; (a,b) \mapsto \beta_{\text{tot}}^a(\beta_0^{\star})^b$. 

Take $L=\Z^2$ in \ref{rvtoric}, let $\Sig$ be the fan of all faces of the cone $\R^2_{\geq 0}\subset N_\R^2=\R^2$, so we have a morphism $\Mor(-, \frak D) \to [\Sig]$. Let $\Sig'$ be the rational finite subdivision of $\Sig$ consisting of the cones $$\sigma_1:= \{(x,y)\in \R^2_{\geq 0}\; |\; x\geq y\}, \quad \sigma_2:=\{(x,y)\in \R^2_{\geq 0}\; |\; x\leq y\}$$
and their faces. Let $D^{\star,+}_{\SL(2)}$ be the log modification of $\frak D$ which represents the fiber product $\Mor(-,\frak D)\times_{[\Sig]} [\Sig']$ (\ref{logmod2}).

$D^{\star,+}_{\SL(2)}$ is covered by the open sets $D^{\star,+}_{\SL(2)}(\sig_j)$ for $j=1,2$ corresponding to the cone $\sig_j$, which represents $\Mor(-, \frak D)\times_{[\Sig]} [\text{face}(\sig_j)]$ where $\text{face}(\sig_j)$ denotes the fan of all faces of  $\sig_j$.  On the open set $U=D^{\star,+}_{\SL(2)}(\sig_1)$ (resp.\ $U=D^{\star,+}_{\SL(2)}(\sig_2)$), 
 the pull back of $\beta_{\text{tot}}/\beta_0^{\star}$ (resp.\ $\beta_0^{\star}/\beta_{\text{tot}}$) in $M^{\gp}_U/\cO^\times_U$ belongs to $M_U/\cO^\times_U$. 
\end{sbpara}

\begin{sbpara}
Since the restriction of $\beta_0^{\star}$ to $D^{\star,\mild}_{\SL(2)}$ is trivial, 
the canonical morphism $D^{\star,+}_{\SL(2)}\to D^{\star}_{\SL(2)}$ is an isomorphism over $D^{\star,\mild}_{\SL(2)}$, and hence $D^{\star,\mild}_{\SL(2)}$ is embedded in $D^{\star,+}_{\SL(2)}$ as an open set. Via this, $D\subset D^{\star,\mild}_{\SL(2)}$ is embedded in $D^{\star,+}_{\SL(2)}$ as an open set.

\end{sbpara}
 
  \begin{sbpara}\label{sl+def2}
 We describe $D^{\star,+}_{\SL(2)}$ as a set.  
 
  We have 
  $$\Sig= \{\tau_{1,2}, \tau_1, \tau_2, \tau_0\}, \quad \Sig'=\{\sig_1, \sig_2, \sig_0, \tau_1, \tau_2, \tau_0\}$$
  where 
 $$\tau_{1,2}:= \R^2_{\geq 0}, \;\; \tau_1:=\R_{\geq 0}\times \{0\}, \;\; \tau_2:= \{0\}\times \R_{\geq 0}, \;\;\tau_0:= \{(0,0)\}, \;\; \sig_0:= \{(x,x)\;|\; x\in \R_{\geq 0}\}.$$ So, $\Sig$  is the set of all faces of $\tau_{1,2}$, and 
 $\text{face}(\sig_j)=\{\sig_j, \tau_j, \sig_0, \tau_0\}$ for $j=1,2$.

The image of $x=(p, Z)\in D^{\star}_{\SL(2)}$ in $\Sig$  is $\tau_0$ if and only if $x\in D$, $\tau_1$ if 
and only if $x\in D^{\star,\mild}_{\SL(2)}\smallsetminus D$, $\tau_2$ if and only if $x$ is a $B$-orbit and $p\in D(\gr^W)$, and $\tau_{1,2}$ if and only if $x$ is a $B$-orbit and  $p\notin D(\gr^W)$.  

We apply \ref{gest8} to describe the log modification $D^{\star,+}_{\SL(2)}$ of $D^{\star}_{\SL(2)}$ as a set. 
For this, we show that the homomorphism 
$\N^2\to (M_{\frak D}/\cO_{\frak D}^\times)_x$ $(\frak D:=D^{\star}_{\SL(2)})$, given by $(\beta_{\text{tot}}, \beta_0^{\star})$ in \ref{sl+def}, is universally saturated for any $x\in \frak D$. If the image of $x$ in $\Sig$ is $\tau_0$ or $\tau_1$ (resp.\ $\tau_2$ or $\tau_{1,2}$), this homomorphism has the shape $\N^2\to \N^r\;;\;(a, b) \mapsto (b,\dots, b)$ (resp.\ $\N^2\to \N\times \N^r\;;\; (a,b)\mapsto (a, b,\dots, b)$) for some integer $r\geq 0$, and hence is universally saturated by Proposition \ref{gest5}.

  By Lemma \ref{gest8}, we have the following list of points of 
 $D^{\star,+}_{\SL(2)}$.
 
 \smallskip

(1) $(x, \tau_j, Z(x))$ ($x\in D^{\star}_{\SL(2)}$ and the image of $x$ in $\Sig$ is  $\tau_j$). (Here $j=0,1,2$.)

(2) $(x, \sig_j, Z(x))$ ($x\in D^{\star}_{\SL(2)}$ and the image of $x$ in $\Sig$ is  $\tau_{1,2}$).
 (Here $j=1,2$.)

(3) $(x, \sig_0, Z')$ ($x=(p, Z)\in D^{\star}_{\SL(2)}$, the image of $x$ in $\Sig$ is $\tau_{1,2}$, and $Z'$ is $\tau_p(A_p)$-orbit in $Z(x)$).

\smallskip

Actually, in (3), what Lemma \ref{gest8} directly tells is that a $\tau^{\star}_p(T(x,\sig_0))$-orbit $Z'$ in the ${\tilde \tau}_p^{\star}(B_p)$-orbit $Z(x)$ appears instead of a $\tau_p(A_p)$-orbit in $Z(x)$. But $\tau^{\star}_p(T(x,\sig_0))=\tau_p(A_p)$ inside ${\tilde \tau^{\star}_p}(B_p)$.

 \end{sbpara}

\begin{sbpara}\label{maps5}

We have a map $D^{\star,+}_{\SL(2)}\to D_{\SL(2)}$ defined as follows.

\smallskip

(1) $(x, \tau_j, Z)$ ($x=(p, Z)$ with image $\tau_j$ in $\Sig$ for $j=0,2$) and ($x, \sig_2, Z)$ $(x=(p, Z)$ with image $\tau_{1,2}$ in $\Sig$) 
are sent to $(p, Z)\in D_{\SL(2)}$. 

\smallskip

(2) $(x, \tau_1, Z)$ ($x=(p, Z)$ with image $\tau_1$ in $\Sig$) and $(x, \sig_1, Z)$ ($x=(p, Z)$ with image $\tau_{1,2}$ in $\Sig$) are sent to $(p, Z_{\spl})\in D_{\SL(2)}$.

Here $Z_{\spl}=\{F_{\spl}\;|\; F\in Z\}$ where $F_{\spl}$ is as in \ref{II,1.2.3}. 

\smallskip

(3) $(x, \sig_0, Z')$ ($x=(p, Z)$ with image $\tau_{1,2}$ in $\Sig$ and  $Z'$ is a $\tau_p(A_p)$-orbit inside $Z$) 
is sent to $(p, Z')\in D_{\SL(2)}$.

\end{sbpara}

\begin{sbthm}\label{0thm} 

 (1) The identity map of $D$ extends uniquely to a morphism $D^{\star,+}_{\SL(2)}\to D^{II}_{\SL(2)}$ in $\cB'_\R(\log)$. Its underlying map of sets is the map in \ref{maps5}. 
 This map is
 proper and surjective.

\smallskip
(2) Let $U$ be the open set $D^{II}_{\SL(2),\nspl}\cup D$ of $D^{II}_{\SL(2)}$. Then the inverse image of $U$ in $D^{\star,+}_{\SL(2)}$ coincides with the open set $D^{\star,+}_{\SL(2)}(\sig_2)$, and the induced morphism $D^{\star,+}_{\SL(2)}(\sig_2)\to U$ of $\cB'_\R(\log)$ is an isomorphism. 
\end{sbthm}

\begin{pf}  
We prove (1). It is sufficient to  prove that the map in \ref{maps5} is a 
 morphism $D^{\star,+}_{\SL(2)}\to D^{II}_{\SL(2)}$ of $\cB'_\R(\log)$. For an  admissible set of weight filtrations $\Phi$ on $\gr^W$, let $D^{\star,+}_{\SL(2)}(\Phi)\subset D^{\star,+}_{\SL(2)}$ be the inverse image of $D^{\star}_ {\SL(2)}(\Phi)\subset D^{\star}_{\SL(2)}$.  
 It is sufficient to prove that the induced map $D^{\star,+}_{\SL(2)}(\Phi)\to D^{II}_{\SL(2)}(\Phi)$ is a morphism in
 $\cB'_\R(\log)$.

Let $D^{\star,+}_{\SL(2),\nspl}\subset D^{\star,+}_{\SL(2)}$ be the inverse image of the open set $D^{\star}_{\SL(2),\nspl}$ of $D^{\star}_{\SL(2)}$. Then $D^{\star,+}_{\SL(2)}(\sig_1)$ is the union of the two open sets $D^{\star,\mild}_{\SL(2)}$ (which is embedded in $D^{\star,+}_{\SL(2)}$) and $D^{\star,+}_{\SL(2),\nspl}\cap D^{\star,+}_{\SL(2)}(\sig_1)$, and $D^{\star,+}_{\SL(2)}(\sig_2)$ is contained in $D^{\star,+}_{\SL(2),\nspl}$.

Take a splitting $\alpha$ of $\Phi$ and a distance $\beta$ to $\Phi$-boundary. 

First, the induced map  $D^{\star,\mild}_{\SL(2)}(\Phi)\to D^{II}_{\SL(2)}(\Phi)$ is a morphism in $\cB'_\R(\log)$ because this map is embedded in  a  commutative diagram
$$\begin{matrix} D^{\star,\mild}_{\SL(2)}(\Phi)& \overset{\subset}\to & D_{\SL(2)}(\gr^W)^{\sim}(\Phi) \times \spl(W) \times \cL\\
\downarrow && \downarrow \\
D^{II}_{\SL(2)}(\Phi) & \overset{\subset}\to & D_{\SL(2)}(\gr^W)^{\sim}(\Phi) \times \spl(W) \times \bar \cL
\end{matrix}$$
where the horizontal arrows are the maps $\nu$ in Proposition \ref{staran1} associated to $(\alpha,\beta)$ and 
the right vertical arrow is the morphism $(p, s, \delta)\mapsto (p, s, \sum_{w\leq -2} (\prod_{j\in \Phi} \beta_j(p))^{-w}\delta_w)$, and because the structure of $D^{II}_{\SL(2)}(\Phi)$ as an object of $\cB'_\R(\log)$ is induced from that of $D_{\SL(2)}(\gr^W)^{\sim}\times \spl(W)\times \bar \cL$ in the sense of \ref{embstr}. 

Next we consider the induced map $D^{\star,+}_{\SL(2),\nspl}(\Phi)\to D^{II}_{\SL(2)}(\Phi)$. 
Take a closed real analytic subset $\cL^{(1)}$ of $\cL\smallsetminus \{0\}$ such that 
$\R_{>0}\times \cL^{(1)}\to \cL\smallsetminus \{0\}\;;\;(a,\delta)\mapsto a \circ \delta$ is an isomorphism, and consider the induced isomorphism 
$\R_{\geq 0}\times \cL^{(1)}\overset{\cong}\to \bar \cL \smallsetminus \{0\}$.
 Let $\beta_0^{\star}: D^{\star}_{\SL(2),\nspl}(\Phi)\to \R_{\geq 0}$ be 
 the composition $D^{\star}_{\SL(2),\nspl}(\Phi) \to \bar \cL\smallsetminus \{0\} \cong \R_{\geq 0}\times \cL^{(1)} \to \R_{\geq 0}$ where the first arrow is induced by the map $\nu$ in Proposition \ref{staran1} associated to $(\alpha,\beta)$. 
For $j=1,2$, let 
 $U_j:= D^{\star,+}_{\SL(2),\nspl}(\Phi)\cap D^{\star,+}_{\SL(2)}(\sig_j)$. 
 Then when we regard $\beta_j$ $(j\in \Phi)$ and $\beta_0^{\star}$ as sections of $M_{U_j}$, then in $M^{\gp}_{U_j}$, $(\prod_{j\in \Phi}\beta_j)/\beta_0^{\star}$ belongs to $M_{U_1}$ and 
 $\beta_0^{\star}/\prod_{j\in \Phi} \beta_j$ belongs to $M_{U_2}$. Furthermore, $\beta_0^{\star}/\prod_{j\in \Phi} \beta_j$ on $U_2$ is the pull back of the section $\beta_0$ of the log structure of $D_{\SL(2),\nspl}(\Phi)$ which is defined as the composition $D_{\SL(2),\nspl}(\Phi) \to \bar \cL\smallsetminus \{0\} \cong \R_{\geq 0}\times \cL^{(1)}\to \R_{\geq 0}$ where the first arrow is induced by $\nu$ of Proposition \ref{staran1} associated to $(\alpha,\beta)$. 
 
 The induced maps $U_j\to D^{II}_{\SL(2)}(\Phi)$ for $j=1,2$ are morphisms because they are embedded in the commutative  diagrams 
 $$\begin{matrix} U_1 & \overset{\subset}\to & D_{\SL(2)}(\gr^W)^{\sim}(\Phi) \times \spl(W) \times (\R_{\geq 0}\times \cL^{(1)})\times \R_{\geq 0}\\
\downarrow && \downarrow \\
D^{II}_{\SL(2)}(\Phi) & \overset{\subset}\to & D_{\SL(2)}(\gr^W)^{\sim}(\Phi) \times \spl(W) \times \bar \cL,
\end{matrix}$$   
$$\begin{matrix} U_2 & \overset{\subset}\to & D_{\SL(2)}(\gr^W)^{\sim}(\Phi) \times \spl(W) \times  \cL^{(1)}\times \R_{\geq 0}\\
\downarrow && \Vert \\
D^{II}_{\SL(2),\nspl}(\Phi) & \overset{\subset}\to & D_{\SL(2)}(\gr^W)^{\sim}(\Phi) \times \spl(W) \times \cL^{(1)}\times \R_{\geq 0}.
\end{matrix}$$

Here in both diagrams, the lower horizontal arrows are induced by $\nu$ in Proposition \ref{staran1} associated to $(\alpha,\beta)$ and the isomorphism $\bar \cL \smallsetminus \{0\} \cong \cL^{(1)}\times \R_{\geq 0}$. 
In the first diagram,  the part $U_1\to \R_{\geq 0}\times \cL^{(1)}$ in the upper row is the composition $U_1\to D^{\star}_{\SL(2),\nspl}\to \bar \cL\smallsetminus \{0\} \cong \R_{\geq 0}\times \cL^{(1)}$,  the map from $U_1$ to the last $\R_{\geq 0}$ in the upper row is $(\prod_{j\in \Phi} \beta_j)/\beta_0^{\star}$, and the right vertical arrow is 
  $(p, s, t, \delta, t')\mapsto (p,s, \sum_{w\leq -2} (tt')^{-w}\delta_w)$.
 In the second diagram,  the part $U_2\to  \cL^{(1)}$ in the upper row is the composition
  $U_2\to D^{\star}_{\SL(2),\nspl}\to \bar \cL\smallsetminus \{0\} \cong \R_{\geq 0}\times \cL^{(1)}\to \cL^{(1)}$, 
   the map $U_2\to \R_{\geq 0}$ in the upper row is $\beta_0^{\star}/\prod_{j\in \Phi} \beta_j$, 
    and the right vertical arrow is the identity map.

  The surjectivity of $D^{\star,+}_{\SL(2)}\to D^{II}_{\SL(2)}$ is easily seen. The map is proper because  $D^{\star,+}_{\SL(2)}$ and $D^{II}_{\SL(2)}$ are proper over $D_{\SL(2)}(\gr^W)^{\sim} \times \spl(W)$.
  This completes the proof of (1). 

We prove (2). It is easy to check that the  inverse image of $U$ in $D^{\star,+}_{\SL(2)}$ is $D^{\star,+}_{\SL(2)}(\sig_2)$, and that the map $D^{\star,+}_{\SL(2)}(\sig_2)\to U$ is bijective. Hence for the proof of (2), it is sufficient to prove that the converse map $D^{II}_{\SL(2),\nspl}\to D^{\star,+}_{\SL(2)}(\sig_2)$ is a morphism in $\cB'_\R(\log)$. 
This is a morphism as is seen from the above last commutative diagram. (In the upper row of this diagram, the structure of the space of $U_2$ as an object of $\cB'_\R(\log)$ is induced from that of 
$D_{\SL(2)}(\gr^W)^{\sim}\times \spl(W) \times \cL^{(1)}\times \R_{\geq 0}$ in the sense of \ref{embstr}.) 
\end{pf}

\begin{sbpara}\label{lam} In the next Proposition \ref{twoSL2}, we consider when the identity map of $D$ extends to an isomorphism  $D^{\star}_{\SL(2)}\cong D^{II}_{\SL(2)}$ in $\cB'_\R(\log)$.

Let $\lambda: D_{\SL(2)}\to D^{\star}_{\SL(2)}$ be the map which coincides on $D_{\SL(2),\nspl}\cup D$ with the composition of morphisms $D^{II}_{\SL(2),\nspl}\cup D\cong D^{\star,+}_{\SL(2)}(\sigma_2) \to D^{\star}_{\SL(2)}$ in $\cB'_\R(\log)$ and 
which coincides on $D_{\SL(2),\spl}:=\{(p,Z)\in D_{\SL(2)}\; |\; Z\subset D_{\spl}\}$ with 
the composition of two isomorphisms $D_{\SL(2),\spl}\cong D_{\SL(2)}(\gr^W)^{\sim} \times \spl(W)\cong D^{\star}_{\SL(2),\spl}:= \{(p, Z)\in D^{\star}_{\SL(2)}\;|\;Z\subset D_{\spl}\}$ in $\cB'_\R(\log)$. 
\end{sbpara}

\begin{sbprop}\label{twoSL2} The following conditions {{\rm (i)}}--{{\rm (vii)}} are equivalent.

{{\rm (i)}} Either $D=D_{\spl}$ or $D_{\SL(2)}(\gr^W)=D(\gr^W)$.

{{\rm (ii)}} The identify map of $D$ extends to an isomorphism $D^{II}_{\SL(2)}\cong D^{\star}_{\SL(2)}$ in  $\cB'_\R(\log)$. 

{{\rm (iii)}} The identity map of $D$ extends to a homeomorphism $D^{II}_{\SL(2)}\cong D^{\star}_{\SL(2)}$. 

{{\rm (iv)}} The map $\lambda : D^I_{\SL(2)}\to D^{\star}_{\SL(2)}$ (\ref{lam}) is continuous. 

{{\rm (v)}} The identity map of $D$ extends to a continuous map $D^{\star}_{\SL(2)}\to D^{II}_{\SL(2)}$. 

{{\rm (vi)}} The map $D^{\star,\mild}_{\SL(2)}\to D_{\SL(2)}$ is injective. 

{{\rm (vii)}} The map $D_{\SL(2),\nspl} \to D^{\star}_{\SL(2)}$ is injective.
\end{sbprop}

\begin{pf} (i) $\Rightarrow$ (ii). If $\cL(F)=0$ (1.2.2) for any $F\in D(\gr^W)$, the isomorphism in
\ref{grsd} extends to  isomorphisms from 
$D^{II}_{\SL(2)}$ and $D^{\star}_{\SL(2)}$ onto $D_{\SL(2)}(\gr^W)^{\sim} \times \spl(W)$ in $\cB'_\R(\log)$. If $D_{\SL(2)}(\gr^W)= D(\gr^W)$, the isomorphism in \ref{grsd} extends to isomorphisms from $D^{II}_{\SL(2)}$ and $D^{\star}_{\SL(2)}$ onto $\{(F, s, \delta)\in D(\gr^W) \times \spl(W) \times \bar \cL\;|\; \delta\in {\bar \cL}(F)\}$ in $\cB'_\R(\log)$.

(ii) $\Rightarrow$ (iii). Clear.

(iii) $\Rightarrow $ (iv), (v), and (vi). Clear.

(v) $\Rightarrow$ (vii). 
 If (v) is satisfied, the composition $D^{II}_{\SL(2), \nspl}\to D^{\star}_{\SL(2)}\to D^{II}_{\SL(2)}$ will be the inclusion map. 

We prove (iv) $\Rightarrow$ (i), 
 (vi) $\Rightarrow$ (i), and (vii) $\Rightarrow$ (i). 
 
 In the rest of this proof, assume  $D\neq D_{\spl}$ and $D_{\SL(2)}(\gr^W) \neq D(\gr^W)$. That is, assume (i) does not hold. Then
there is $x= (p, Z)\in D^{\star,\mild}_{\SL(2)}$ with $p$ of rank $1$ such that $Z\subset  D_{\nspl}$. 
Let $x_{\spl}:= (p, Z_{\spl})\in D^{\star,\mild}_{\SL(2)}$.  
We have $x\neq x_{\spl}$.

We prove (iv) $\Rightarrow$ (i). Take $\br\in Z$. Then when $t\in \R_{>0}$ tends to $0$, $\tau^{\star}_p(t)\br$ converges to $x$ and $\tau^{\star}_p(t)\br_{\spl}$ converges to $x_{\spl}$ in $D^{\star}_{\SL(2)}$. 

\medskip

{\bf Claim.} Let $y= (p, Z_{\spl})\in D_{\SL(2)}$. Then when $t\in \R_{>0}$ converges to $0$, $\tau^{\star}_p(t)\br$ and $\tau^{\star}_p(t)\br_{\spl}$ converge to $y$ in $D^I_{\SL(2)}$.

\medskip

We prove Claim. Let $s:=\spl_W(p)=\spl_W(\br)$. By Part II, Proposition 3.2.12, it is sufficient to prove that when $t\in \R_{>0}$ tends to $0$, $(s \tau_p(t)s^{-1})^{-1}(s\tau^{\star}_p(t)s^{-1})\br$ and $(s\tau_p(t)s^{-1})^{-1}(s\tau^{\star}_p(t)s^{-1})\br_{\spl}$ converge to $\br_{\spl}$. 
The former is equal to $s (t^{-w})_w s^{-1}\br$ and hence  converges to $\br_{\spl}$. Here $(t^{-w})_w$ denotes the linear automorphism of $\gr^W=\prod_w \gr^W_w$ which acts on $\gr^W_w$ as the multiplication by $t^{-w}$. The latter is equal to $\br_{\spl}$. 
This proves Claim.

By Claim,  if the continuous map  $D^I_{\SL(2)}\to D^{\star}_{\SL(2)}$ exists, it should send $y$ to $x$ and also to $x_{\spl}\neq x$. A contradiction.

We prove (vi) $\Rightarrow $ (i).
The elements $x$ and $x_{\spl}$ of $D^{\star,\mild}_{\SL(2)}$ have the same image $(p, Z_{\spl}) \in D^{II}_{\SL(2)}$. Hence the map $D^{\star,\mild}_{\SL(2)}\to D_{\SL(2)}$ is not injective.

We prove (vii) $\Rightarrow$ (i). Take $\br \in Z$. 
Take $a\in \R_{>0}\smallsetminus \{1\}$ and let $\br'=a\circ \br$. 
 Then the elements  of $D_{\SL(2),\nspl}$ of the forms 
 $(p, \tau_p(\R_{>0})\br)$ and $(p, \tau_p(\R_{>0})\br')$ for the lifted action (\ref{liftac})
 are different but they have the same image
 $(p, \R_{>0}\circ Z)$  in $D^{\star}_{\SL(2)}$.  Hence the map $D_{\SL(2),\nspl}\to D^{\star}_{\SL(2)}$ is not injective. 
\end{pf}

\subsection{Relations with  $D_{\BS}$}

We connect the spaces $D^{\star,-}_{\SL(2)}$ and $D_{\BS}$ by introducing a new space $D^{\star,\BS}_{\SL(2)}$ of $\SL(2)$-orbits. 

\begin{sbpara}
For $Q=(Q(w))_w \in \prod_{w\in \Z} \cW(\gr^W_w)$, let $$G_\R(\gr^W)_{Q}:= \prod_w G_\R(\gr^W_w)_{Q(w)}\quad\text{where}$$ $$G_\R(\gr^W_w)_{Q(w)}:=\{g\in G_\R(\gr^W_w)\;|\; gW'=W'\;\text{for all}\; W'\in Q(w)\}.$$

Let $G_\R(\gr^W)_{Q,u}$ be the unipotent radical of $G_\R(\gr^W)_Q$. 
\end{sbpara}

\begin{sbpara}\label{SB1} 
Let $p\in D_{\SL(2)}(\gr^W)$. We define a set $\cP(p)$ of $\Q$-parabolic subgroups of $G_\R(\gr^W)$. 

Let $X(S_p)$ be the character group of the torus $S_p$ (\ref{sim5}) associated to $p$. 
For $\chi\in X(S_p)$, let $$\fg_\R(\gr^W)_{\chi}=\{v\in \fg_\R(\gr^W)\; |\; \Ad(\tau_p^{\star}(t))v= \chi(t)v \;\text{for all}\; t\in S_p\}.$$

Let $\cP(p)$ be the set of all $\Q$-parabolic subgroups $P$ of $G_\R(\gr^W)$ satisfying the following conditions (i) and (ii). 

\smallskip

(i) $P\supset G_\R(\gr^W)_Q$ and $P_u\supset G_\R(\gr^W)_{Q,u}$, where $Q= (\cW(p_w))_w$.

\smallskip

(ii) There is a subset $I$ of $X(S_p)$ such that $\Lie(P)= \sum_{\chi\in I} \fg_\R(\gr^W)_{\chi}$.

\end{sbpara}

\begin{sbpara}\label{2.6.3}  We define $D^{\star,\BS}_{\SL(2)}$ as a set.  

$D_{\SL(2)}^{\star,\BS}$ is the set of all triples  $(p, P, Z)$, where $p\in D^{\star,-}_{\SL(2)}(\gr^W)$, $P\in \cP(p)$, and $Z\subset D$  satisfying the following conditions (i) and (ii). Let $A_{p,P}\subset A_p$ be the inverse image of $A_P\subset P/P_u$ under the composite map $A_p \to G_\R(\gr^W)_Q/G_\R(\gr^W)_{Q,u}\to P/P_u$. Let $B_{p,P}= \R_{>0}\times A_{p,P}\subset B_p$.

\smallskip

(i) $Z$ is either an $\tau^{\star}_p(A_{p.P})$-orbit in $D$ or a $\tilde \tau^{\star}_p(B_{p,P})$-orbit in $D_{\nspl}$.

\smallskip

(ii) The image of $Z$ in $D(\gr^W)$ is contained in the torus orbit $Z(p)$. 

\smallskip
For $w\in \Z$, we denote by $D_{\SL(2)}(\gr^W_w)^{\BS}$ the set $D^{\star, \BS}_{\SL(2)}$ for $\gr^W_w$. Let $D_{\SL(2)}(\gr^W)^{\BS}:=\prod_w D_{\SL(2)}(\gr^W_w)^{\BS}$. 

We have an evident map $D^{\star,\BS}_{\SL(2)}\to D^{\star}_{\SL(2)}(\gr^W)^{\BS}$. 

\end{sbpara}

\begin{sbprop}\label{SB2} (1) We have a canonical map $$D_{\SL(2)}^{\star,\BS} \to D^{\star,-}_{\SL(2)}\;;\; (p,P,Z)\mapsto (p, \tau_p^{\star}(A_p)Z).$$

\smallskip

(2) We have a map
$$D_{\SL(2)}^{\star,\BS}\to D_{\BS}\;;\; (p, P, Z) \mapsto (P, A_P\circ Z).$$

Here $\circ$ denotes the Borel-Serre action with respect to $P$. 

\end{sbprop}

\begin{pf} (1) is clear.

We prove (2). It is sufficient to prove that, for $\br\in Z$ and $t\in A_{p,P}$, we have $\tau_p^{\star}(t)\br = (\tau^{\star}_p(t) \bmod P_u)\circ \br$. This follows from $\theta_{K_{\br}}(\tau^{\star}_p(t))= \tau^{\star}_p(t)^{-1}$ (\cite{KU0}, Lemma 3.8) where $\theta_{K_{\br}}$ denotes the Cartan involution $G_\R(\gr^W)\to G_\R(\gr^W)$ associated to the maximal compact subgroup $K_{\br}$ of $G_\R(\gr^W)$. 
\end{pf}

We give $D^{\star,\BS}_{\SL(2)}$ a structure of an object of $\cB'_\R(\log)$.

The following \ref{Pcone}--\ref{Pcone6} are preparations.

 \begin{sblem}\label{Pcone} Let $L$ and $N$ be as in \ref{rvtoric}.

 Let $R$ be a finite subset of $L$ such that $R^{-1}=R$  and such that the $\Q$-vector space $\Q\otimes L$ is generated by $R$.

 \smallskip
 
 $(1)$ Let $\sig$ be a rational finitely generated sharp cone in $N_\R$ and let $\cS(\sig)=\{l\in L\;|\;h(l) \geq 0\;\text{for all}\;h\in \sig\}$ be the corresponding fs submonoid of $L$ such that $\cS(\sig)^{\gp}=L$. Then $\sig$ satisfies the following condition $(i)$ if and only if $\cS(\sig)$ satisfies the following conditions $(ii.1)$ and $(ii.2)$. 
 
 \smallskip
 
 $(i)$ There exists  a subset $R'$ of $R$ such that $R=R'\cup (R')^{-1}$ and such that $$\sigma=\{h\in N_\R\;|\; h(l) \geq 0 \;\text{for all}\;  l\in R'\}.$$

  \smallskip
  
 $(ii.1)$ $R\subset \cS(\sig)\cup \cS(\sig)^{-1}$.  
 
 \medskip
 
 $(ii.2)$ For any $l\in \cS(\sig)$, there is an integer $n\geq 1$ such that $l^n$ belongs to the submonoid of $L$ generated by 
 $\cS(\sig)\cap R$.

\smallskip
 
 $(2)$ The set of all $\sig$ satisfying the condition $(i)$ in (1) is a rational fan whose support is the whole  $N_\R$. 
   
   \medskip
   
  $(3)$ Assume that we are given a subset $R^+$ of $R$ which generates $\Q\otimes L$ over $\Q$. Let  $\nu:= \{h\in N_{\R}\; |\; h(R^+)\subset \R_{\geq 0}\}$. Then $\sig$ as above such that $\sig\subset \nu$ form a rational finite subdivision of $\nu$.

 \end{sblem}
  
  \begin{pf} The proof of (1) is straightforwards. 
  
  We prove (2). Let $I$ be the set of all cones $\sig$ satisfying the condition (i) in (1). We first prove that $I$ is a fan.

  We prove that if $\sig_j\in I$ $(j=1,2)$, $\sig_1\cap \sig_2$ is a face of $\sig_1$. Let $R'_j\subset R$ and assume $R=R'_j\cup (R'_j)^{-1}$, $\sig_j= \{h\in N_\R\;|\; h(l)\geq 0 \;\text{for all}\; l\in R'_j\}$. Let $R'=R'_1\cup R'_2$. Then $\sig_1\cap \sig_2=\{h\in N_\R\;|\; h(l)\geq 0\; \text{for all}\; l\in R'\}$. Since $R'\smallsetminus R'_1\subset (R'_1)^{-1}$, $\sig_1\cap \sig_2$ is a face of $\sig_1$.

  We prove that if $\sig\in I$, any face $\tau$ of $\sig$ belongs to $I$.  
  Since $\tau$ is a face of $\sig$, we have $\cS(\tau)=\cS(\sig)[b^{-1}]
  = \{ab^{-n}\;|\; a\in \cS(\sig), n\geq 0\}$ for some $b\in \cS(\sig)$. 
  By the condition (ii.2) in (1) for  $\cS(\sig)$, there exists $n\geq 1$, $a_1, \dots, a_r\in\cS(\sig) \cap R$ 
  and $m(j)\geq 1$ $(1\leq j\leq r)$ such that $b^n = \prod_{j=1}^r  a_j^{m(j)}$. 
  We have $\cS(\tau)=\cS(\sig)[1/\prod_{j=1}^r a_j]$. For the set $R'\subset R$ such that $R=R'\cup (R')^{-1}$ and $\sig= \{h\in N_\R\;|\; h(l)\geq 0\;\text{for all}\; l\in R'\}$, 
  we have 
  $\tau=\{h\in N_{\R}\;|\; h(l)\geq 0\;\text{for all}\;l\in R'\cup\{a_1^{-1},\dots, a_r^{-1}\}\}$. Hence $\tau \in I$.
  
  These prove that $I$ is a fan. We show that $\bigcup_{\sig \in I} \sig=N_\R$. Let $h\in N_\R$. Let $R'=\{l\in R\;|\; h(l)\geq 0\}$. Then $R=R'\cup (R')^{-1}$. For  $\sig:=\{h'\;|\; h'(l) \geq 0\;\text{for all}\; l\in R'\}\in I$, we have $h\in \sig$. 
  
 These completes the proof of (2).
 
 We prove (3). By (2), we have $\nu=\bigcup_{\sig\in I} (\sig \cap \nu)$. It is sufficient to prove that $\sig \cap \nu \in I$ for any $\sig \in I$. For $R'\subset R$ such that $R=R'\cup (R')^{-1}$ and $\sig=\{h\in N_\R\;|\; h(l)\geq 0\;\text{for all}\; l\in R'\}$, we have $\sig\cap \nu= \{h\in N_\R\;|\; h(l)\geq 0 \;\text{for all}\; l\in R' \cup R^+\}\in I$. 
  \end{pf}

  \begin{sbpara}\label{Pcone2} Let $Q=(Q(w))_w\in \prod_w \cW(\gr^W_w)$. 
  Let 
  $L$  be the character group of $\prod_{w\in \Z} \bG_m^{Q(w)}$ and let $N=\Hom(L, \Z)$. We have the situation of \ref{rvtoric}. As in \ref{rvtoric}, we denote the group law of $L$ multiplicatively, though  $L$ is identified with $\prod_w \Z^{Q(w)}$.

  Let $\cP(Q)$ be the set of  
  all $\Q$-parabolic subgroups $P$ of $G_\R(\gr^W)$ satisfying the following  conditions (i) and (ii).

  \smallskip

  (i) $P\supset G_\R(\gr^W)_Q$.

  \smallskip
  
  (ii) 
  Take a splitting $\alpha=(\alpha_w)_w$ of $Q$. For $\chi\in L$, let $\fg_\R(\gr^W)_{\chi}$ be the part of $\fg_\R(\gr^W)$ on which the adjoint action of $\prod_w \bG_m^{Q(w)}$ via $\alpha^{\star}$ is given by $\chi$.  
Then there is a subset $I$ of $L$ such that $\Lie(P)= \sum_{\chi\in I} \fg_\R(\gr^W)_{\chi}$. 
  
    \smallskip
  Under the condition (i), the condition (ii) is independent of the choice of $\alpha$. This is because if $\alpha'$ is another splitting of $Q$, $\alpha'(t)=g\alpha(t)g^{-1}$ for some $g\in G_\R(\gr^W)_Q\subset P$.
  
  \end{sbpara}

  \begin{sbpara}\label{Pcone2.2} Let the notation be as in \ref{Pcone2}. 
  
  Taking a splitting $\alpha$ of $Q$, 
 define a subset $$R(Q)=\{\chi\in L\;|\;\fg_\R(\gr^W)_{\chi}\neq 0\}$$ 
 where $\fg_\R(\gr^W)_{\chi}$ is defined with respect to $\alpha$. This set is independent of the choice of $\alpha$ because all splittings of $Q$ are conjugate by elements of $G_\R(\gr^W)_Q$.
   
   Let $L^+=\prod_w \N^{Q(w)}\subset \prod_w \Z^{Q(w)}=L$. We will apply \ref{Pcone} by taking $R(Q)$ and $R(Q)\cap L^+$ as $R$ and $R^+$, respectively. We show that 
    $R^+$ generates the $\Q$-vector space $\Q\otimes L$, as is assumed in 
     \ref{Pcone} (3). Let $w\in \Z$ and take $p\in D_{\SL(2)}(\gr^W_w)$ such that $Q(w)= \cW(p_w)$. Let $n$ be the rank of $p$, take a representative of $p$ and let $N_1,\dots, N_n\in \fg_\R(\gr^W_w)$ be the monodromy logarithms of the representaive, and identify $Q(w)$ with $\{1,\dots, n\}$ (\ref{simst2}). Then $\Ad(\tau^{\star}_p(t))N_j= t_j^{-2}N_j$. Hence $R(Q(w))^+$ generates the $\Q$-vector space $\Q^{Q(w)}$. Hence $R(Q)^+$ generates the $\Q$-vector space $\Q\otimes L=\prod_w \Q^{Q(w)}$.

  Let $\cP'(Q)$ be the set 
  of all rational finitely generated sharp cones $\sig$ in $N_\R$ satisfying the following conditions $(i)$ and $(ii)$.

  \smallskip
  
  $(i)$ There is a subset $R'$ of $R(Q)$ such that $R(Q)=R'\cup (R')^{-1}$ and such that $\sig=\{h\in N_\R\;|\; h(\chi)\geq 0\;\text{for all}\;\chi\in R'\}$. 
  
  \smallskip
  
  $(ii)$ $\sig \subset \prod_w \R_{\geq 0}^{Q(w)}$ in $N_\R=\prod_w \R^{Q(w)}$. 
  \smallskip
  
  That is, $\cP'(Q)$ is the set of $\sig$ considered in \ref{Pcone} (3). Hence $\cP'(Q)$ is a rational fan in $N_\R$ and is a rational finite subdivision of the cone $\prod_w \R_{\geq 0}^{Q(w)}\subset N_\R$.
  
    \end{sbpara}

 \begin{sbpara}\label{Pcone2.5} Let the notation be as in \ref{Pcone2} and \ref{Pcone2.2}. 
   
 We have $\cP(Q) = \prod_w \cP(Q(w))$ where the element $(P_w)_w$ of the left  hand side corresponds to the element $\prod_w P_w$ of the right-hand side.

 We have $R(Q)=\prod_w R(Q(w))$ in $X(\prod_w \bG_m^{Q(w)})= \prod_w X(\bG_m^{Q(w)})$.

 We have  $\cP'(Q)= \prod_w \cP'(Q(w))$ where the element $(\sig_w)_w$ of the left  hand side corresponds to the element $\prod_w \sig_w$ of the right-hand side. 
 \end{sbpara}
  
  \begin{sbprop}\label{Pcone3} Let the notation be as in \ref{Pcone2} and \ref{Pcone2.2}.  
   For $P\in \cP(Q)$, let $$\sig_P=\{h\in N_\R\;|\;h(\chi)\geq 0\; \text{for all $\chi\in R(Q)$ such that $\fg_\R(\gr^W)_{\chi^{-1}}\subset \Lie(P)$}\}.$$ Then $\sig_P\in \cP'(Q)$ and we have a bijection
  $$\cP(Q)\to \cP'(Q)\;;\;P\mapsto \sig_P.$$

    \end{sbprop}
  
  \begin{pf}

  By \ref{Pcone2.5} and by the fact $\sig_P= \prod_w \sig_{P_w}$,  we can (and do) assume that we are in the pure situation of weight $w$. We denote $Q(w)$ by $Q$. 
  
 Take $p\in D_{\SL(2)}$ such that $\cW(p)=Q$, and take $\tau_p$ as a splitting $\alpha$ of $Q$. Let $n=\sharp(Q)$ be the rank of $p$. Let $N_1, \dots, N_n$ be the monodromy logarithms of $p$. We identify $Q$ with $\{1,\dots, n\}$.
    
  We prove that $\sig_P\in \cP'(Q)$ for $P\in \cP(Q)$. Let $R'=\{\chi\in L\;|\; \Lie(P)\cap \Lie(G_\R)_{\chi^{-1}}\neq 0\}$. Since $\alpha^{\star}(\prod_w \bG_m^{Q(w)})\subset P$ and since $P$ is parabolic, we have $R(Q)=R'\cup (R')^{-1}$. By the property (ii) of $P$ in \ref{Pcone2}, we have $\Lie(G_\R)_{\chi^{-1}}\subset \Lie(P)$ for $\chi\in R'$. Hence $\sig_P= \{h\in N_\R\;|\; h(\chi)\geq 0\;\text{for all} \; \chi\in R'\}$. 
  It remains to prove that $\sig_P \subset \R^Q_{\geq 0}$ in $N_\R=\R^Q$. Since $N_j\in \Lie(P)$ and 
   $\Ad(\tau^{\star}_p(t))(N_j)= t_j^{-2}N_j$ $(1\leq j\leq n)$, for any $\chi\in L^+$, $\chi^2$ is contained in the submonoid of $L$ generated by $R'$. This proves that $h(\chi)\geq 0$ for any $h\in \sig_P$ and $\chi\in L^+$. This implies $\sig_P\subset \R^Q_{\geq 0}$. 

  Thus we have a map $\cP(Q)\to \cP'(Q)$. 
  
  Next we define a map $\cP'(Q) \to \cP(Q)$.

  Let $\sig\in \cP'(Q)$ and let $\cS(\sig)\subset L$ be the corresponding fs submonoid of $L$. For $\chi\in L$, let $V[\chi]\subset H_{0,\R}$ be the sum of the $\chi'$-components $(H_{0,\R})_{\chi'}$ of $H_{0,\R}$ for all $\chi'\in L$ such that $\chi(\chi')^{-1}\in \cS(\sig)$. For $\chi, \chi'\in L$, we have $V[\chi]\supset V[\chi']$ if and only if $\chi(\chi')^{-1}\in \cS(\sig)$. 
  Let
    $P$ be the algebraic subgroup of $G_\R$ consisting of all elements which preserve $V[\chi]$ for all $\chi\in L$. We prove $P\in \cP(Q)$. 
  
  Since $L=\cS(\sig)\cup \cS(\sig)^{-1}$ (\ref{Pcone}), we have either $V[\chi]\supset V[\chi']$ or $V[\chi]\subset V[\chi']$. 
   As in \cite{KU0} 2.7,  this totally ordered property of the set $\{V[\chi]\;|\;\chi\in L\}$ shows that $P$ is a parabolic subgroup of $G_\R$. 
    We show that $P$ is defined over $\Q$.  
   For $\chi\in L$, let $U[\chi]= \sum_{\chi'} (H_{0,\R})_{\chi'}$ where $\chi'$ ranges over all elements of $L$ such that $\chi(\chi')^{-1}\in L^+$. Then $U[\chi]= \bigcap_{W'\in Q} W'_{m(W')}$ where 
    $m(W')\in \Z$ is the $W'$- component of $\chi\in L=\Z^Q$. 
   Since $W'$ are rational, $U[\chi]$ is rational. Since 
      $V[\chi]$ is the sum of $U[\chi']$ for all $\chi'$ such that $\chi(\chi')^{-1}\in \cS(\sig)$, $V[\chi]$ is also rational.
   Hence $P$ is rational.  The properties (i) and (ii) of $P$ in \ref{Pcone2} 
    are checked easily.

  As is easily seen, the maps $\cP(Q) \to \cP'(Q)$ and $\cP'(Q)\to \cP(Q)$ are the converses of each other.
    \end{pf}

 \begin{sbpara}\label{Pcone4} Let the notation be as in Proposition \ref{Pcone3}. Via the bijection in Proposition \ref{Pcone3}, we identify the fan $\cP'(Q)$ with the set $\cP(Q)$ of $\Q$-parabolic subgroups of $G_\R(\gr^W)$.

 Let $\Sig$ be the fan of all faces of the cone $\nu:=\prod_w \R^{Q(w)}_{\geq 0}\subset N_\R$. By the canonical homomorphism $\cS(\nu)=L^+=\prod_w \N^{Q(w)}\to M_{\frak E'}/\cO^\times_{\frak E'}$ where $\frak E'=D_{\SL(2)}(\gr^W)(Q)$ Proposition (\ref{logstalk}), we have a morphism $\Mor(-, D_{\SL(2)}(\gr^W)(Q) \to [\Sig]$. Consider the diagrams
  $$\Mor(-, D_{\SL(2)}(\gr^W)(Q))\to [\Sig]\leftarrow [\cP(Q)], \quad 
  D_{\SL(2)}(\gr^W)\to \Sig\leftarrow \cP(Q).$$

  \end{sbpara} 
  
  \begin{sblem}\label{Pcone5} Let $p\in D_{\SL(2)}(\gr^W)(Q)$. Then  $\cP(p)\subset \cP(Q)$. For $P\in \cP(Q)$, $P\in \cP(p)$ if and only if the image of $P$ in $\Sig$ 
  coincides with the image of $p$ in $\Sig$.

  \end{sblem}
  
  \begin{pf} It is clear that $\cP(p)\subset \cP(Q)$. To prove the rest, we may assume $Q(w)=\cW(p_w)$ $(w\in \Z)$. 
  It is sufficient to prove that in this case, for $P\in \cP(Q)$, $P_u\supset G_\R(\gr^W)_{Q,u}$ if and only if the image of $P$ under the map $\cP(Q)\to \Sig$ coincides with the face $\nu$ of $\nu$. Let $\sig\in \cP'(Q)$ be the cone in $N_\R$ corresponding to $P$ and let $\cS:=\cS(\sig)$ be the corresponding fs monoid in $L$. 
  Then the image of $\sig$ in $\Sig$ is $\nu$ if and only of $\cS^\times \cap L^+=\{1\}$. 
  By the proof of Proposition \ref{Pcone3}, we have $$\Lie(P_u)= \sum_{\chi\in \cS\smallsetminus \cS^\times} 
  \fg_\R(\gr^W)_{\chi^{-1}}, \quad \Lie(G_\R(\gr^W)_{Q,u}= \sum_{\chi\in L^+\smallsetminus 
  \{1\}} \fg_\R(\gr^W)_{\chi^{-1}}.$$
  Hence $P_u\supset G_\R(\gr^W)_{Q,u}$ if $L^+\smallsetminus \{1\}\subset \cS\smallsetminus \cS^{\times}$, 
  i.e., if $L^+\cap \cS^\times =\{1\}$. Let $w\in \Z$ and let $N_1, \dots, N_n\in \Lie(G_\R(\gr^W_w)_{Q(w),u})$ $(n=\sharp(Q(w)))$ be the monodromy logarithms of $p_w$. If $P_u\supset G_\R(\gr^W)_{Q,u}$, then $N_j \in \Lie(P_u)$. Since $\Ad(\tau^{\star}_p(t))N_j=t_j^{-2}N_j$ $(1\leq j\leq n)$,  this proves $L^+\smallsetminus \{1\}\subset \cS\smallsetminus \cS^\times$. 
  \end{pf}

  \begin{sbpara}\label{Pcone6} Let the notation be as in \ref{Pcone4}. We show that the object of $\cB'_\R(\log)$ which represents the fiber product of $\Mor(-,D^{\star,-}_{\SL(2)}(Q))\to [\Sig] \leftarrow [\cP(Q)]$ is identified, as a set, with the inverse image of $D_{\SL(2)}^{\star,\BS}(Q)$ of $\frak D':=D_{\SL(2)}^{\star,-}(Q)\subset \frak D:=D_{\SL(2)}^{\star,-}$ in $D^{\star,\BS}_{\SL(2)}$ (\ref{SB2}).  By Lemma \ref{gest8} and by Lemma \ref{Pcone5}, a point of this fiber product is identified with a triple $(x, P, Z)$ where  
  $x\in \frak D'$, $P\in \cP(p)$,  $Z\subset D$ satisfying the following condition (i). Let $\cS=\cS(\sig)$ be the 
  fs submonoid of $L$ corresponding to the cone $\sig\in \cP'(Q)$ which corresponds to $P$. Write $x=(p, Z')\in \frak D'$ and define a subgroup $T(x, P)$ 
  of $T(x)= \Hom((M^{\gp}_{\frak D}/\cO^\times_{\frak D})_x, \R^{\mult}_{>0})$ as follows. If $x$ is an $A$-orbit, let 
 $T(x,P)= \Hom(L/\cS^\times, \R^{\mult}_{>0})\subset \Hom(L, \R^{\mult}_{>0})=A_p= T(x)$. If $x$ is a $B$-orbit, let
 $T(x,P)= \R_{>0}\times \Hom(L/\cS^\times, \R^{\mult}_{>0})\subset \R_{>0}\times \Hom(L, \R^{\mult}_{>0})=B_p=T(x)$.
 
 \smallskip
  
  (i)  $Z$ is a $T(x, P)$-orbit in $Z'$.    
  
  \smallskip
  
  We prove this by showing the following claim.

  \smallskip
  
  {\bf Claim.} $T(x, P)= A_{p,P}$ if $x$ is an $A$-orbit and  $T(x,P)=B_{p,P}$ if $x$ is a $B$-orbit (\ref{2.6.3}).
  
  \smallskip
  
  Let $S_{p,P}= \Hom(L/\cS^\times, \bG_m)\subset \Hom(L, \bG_m)=S_p$. Then $S_{p,P}$ coincides with the
  part of $S_p$ consisting of all elements whose adjoint action on $\Lie(P/P_u)$ is trivial. That is, $S_{p,P}$ is the inverse image in $S_p$ of the center of $P/P_u$. Since $S_{p,P}$ is $\Q$-split, the image of $S_{p.P}$ in $P/P_u$ is contained in $S_P$. This proves that $A_{p,P}$ coincides with the connected component of $S_{p,P}(\R)$ containing the unit element. This proves the above Claim.
  
  Since $Z'=\tau_p^{\star}(A_p)Z$, a triple $(x,P, Z)$ as above corresponds to a point $(p, P, Z)$ of $D^{\star,\BS}_{\SL(2)}(Q)$ (\ref{2.6.3}) in one to one manner.

       \end{sbpara}

 \begin{sbpara} For $Q\in \prod_w \cW(\gr^W_w)$, we define the structure of $D_{\SL(2)}^{\star,\BS}(Q)$ as an object of $\cB'_\R(\log)$ by identifying it as a log modification of $D_{\SL(2)}^{\star,-}(Q)$ by \ref{Pcone6}. 
When $Q$ moves, these structures on $D^{\star,\BS}_{\SL(2)}(Q)$ glue globally to a structure of $D_{\SL(2)}^{\star,\BS}$ as an object of $\cB'_\R(\log)$. 
 
For a $\Q$-parabolic subgroup $P$ of $G_\R(\gr^W)$, let $$D^{\star,\BS}_{\SL(2)}(P)= \{(p,P', Z)\in D^{\star,\BS}_{\SL(2)}\;|\; P'\supset P\}.$$ Then $D^{\star,\BS}_{\SL(2)}(P)$ is an open set of $D^{\star,\BS}_{\SL(2)}$, and when $P$ moves, we have a covering of $D_{\SL(2)}^{\star,\BS}$ by these open sets. 

\end{sbpara}

\begin{sbprop}\label{S-B} The diagram
$$\begin{matrix}   D^{\star,\BS}_{\SL(2)} & \to & D^{\star,-}_{\SL(2)} \\ \downarrow && \downarrow \\
 D_{\SL(2)}(\gr^W)^{\BS} &\to &D_{\SL(2)}(\gr^W) \end{matrix}$$
 is cartesian in $\cB'_\R(\log)$ and also in the category of topological spaces.

 \end{sbprop}

 \begin{pf} This is because $D^{\star,\BS}_{\SL(2)}(Q)$ represents the fiber product of $\Mor(-, D_{\SL(2)}^{\star,-})(Q)\to [\Sig]\leftarrow [\cP(Q)]$ and $D_{\SL(2)}(\gr^W)^{\BS}(Q)$ represents the fiber product of $\Mor(-, D_{\SL(2)}(\gr^W))\to [\Sig]\leftarrow [\cP(Q)]$.  \end{pf}

  \begin{sbprop}

Let $F\in D(\gr^W)$, $\bar L= \bar \cL(F)$. 
Then $D^{\star,\BS}_{\SL(2)}$ is an $\bar L$-bundle over $D_{\SL(2)}(\gr^W)^{\BS}\times \spl(W)$.

 \end{sbprop}
  
  \begin{pf}
  This follows from Proposition \ref{S-B} and the corresponding result for $D^{\star,-}_{\SL(2)}$. 
  \end{pf}

. \begin{sbpara}
For $p\in D_{\SL(2)}(\gr^W)$ and $P\in \cP(p)$, let $S_{p,P}\subset S_p$ be the torus defined in \ref{Pcone6}, let $X(S_{p,P})$ be the character group of $S_{p,P}$, and let $X(S_{p,P})^+=\cS/\cS^\times$, where $\cS:=\cS(\sig_P)$ (\ref{Pcone}) with  $\sig_P$ the cone corresponding to $P$ (\ref{Pcone3}). Define a real toric variety $\bar A_{p,P}$ by
$$\bar A_{p,P}:= \Hom(X(S_{p,P})^+, \R^{\mult}_{\geq 0})\supset A_{p,P}=\Hom(X(S_{p,P}), \R^{\mult}_{>0}).$$
We have a canonical morphism
$$\bar A_{p,P}\to \bar A_p$$
induced from the homomorphism $X(S_p)^+\to X(S_{p,P})^+$ which is induced by the inclusion map $S_{p,P}\to S_p$. 

\end{sbpara}
     
\begin{sblem}\label{lemSB} (1) The homomorphism $X(S_P)\to X(S_{p,P})$ induced by $S_{p,P} \to S_P$ (\ref{Pcone6}) sends $X(S_P)^+$ to $X(S_{p,P})^+$. 

\medskip

(2) The map $A_{p, P}\to A_P$ extends uniquely to a morphism $\bar A_{p,P}\to \bar A_P$ in $\cB'_\R(\log)$.
  \end{sblem}
  
  \begin{pf} We prove (1). As a monoid, 
  $X(S_P)^+$ is  generated by $\Delta(P)$ (\ref{BS2}). For $\chi\in \Delta(P)$, $\chi^{-1}$ appears in $\Lie(P)$. Hence the image of $\chi^{-1}$ in $X(S_{p,P})$ appears in $\Lie(P)$. Hence the image of $\chi$ in $X(S_{p,P})$ belongs to $X(S_{p,P})^+$. 
  
  \smallskip
  
  (2) follows from (1). In fact, the homomorphism $X(S_P)^+\to X(S_{p,P})^+$ in (1) induces the morphism $\Hom(X(S_P)^+, \R^{\mult}_{\geq 0}) \to \Hom(X(S_{p,P})^+, \R^{\mult}_{\geq 0})$. 
  \end{pf}

\begin{sbpara}\label{lsSB1} In Theorem \ref{lsSB}, we will consider the local structure of $D_{\SL(2)}^{\star,\BS}$, comparing it with the local structure of $D_{\BS}$. Here we give preparations. We consider the following two situations (bd) and (d).

\smallskip

(bd)  $\frak D= D_{\SL(2)}^{\star,\BS}$ and $\frak E=D_{\SL(2)}(\gr^W)^{\BS}$.\smallskip

(d) $\frak D=D_{\BS}$ and $\frak E = D_{\BS}(\gr^W)$. 

\smallskip

Fix $p\in \frak E$ and $\br\in Z(p)$ (\ref{sim5}). In the situation (bd) (resp.\ (d)), fix $P\in \cP(p)$ (resp.\  fix a $\Q$-parabolic subgroup $P$ of $G_\R(\gr^W)$ such that $p\in \frak E(P)$).

Let $R$ be an $\R$-subspace of ${\frak g}_\R(\gr^W)$ satisfying the following conditions (C1) and (C2).

\smallskip

(C1) 
${\frak g}_\R(\gr^W)= \Lie(\tau^{\star}(A_{p,P})) \oplus R\oplus \Lie(K_{\br})$ (resp.\ ${\frak g}_\R(\gr^W)= \Lie((A_P)_{\br}) \oplus R\oplus \Lie(K_{\br})$, where $(A_P)_{\br}$ denotes the Borel-Serre lifting \ref{BS1} of $A_P$ at $\br$).

\smallskip

(C2)  $R\subset \Lie(P)$.

\smallskip
These conditions on $R$ are similar to those in \ref{thm+2}. 

 Like in \ref{thm+2}, let $S$ be an $\R$-subspace of $\Lie(K_{\br})$ such that $\Lie(K_{\br})=\Lie(K'_{\br})\oplus S$.

We define an object $Y$ of $\cB'_\R(\log)$ as follows. Let
$$Y= \bar A_P\times R \times S\;\;\text{in the situation (d)}.$$
In the situation (bd), we define $Y$ as follows. 
Let 
$$X= \bar A_{p,P}  \times R \times S.$$  
Let $Y$ be the subset of $X$ consisting of all elements $(t,f,k)$ satisfying the following conditions (i) and (ii).

\smallskip

(i)  If $\chi\in X(S_p)$ and $t(\chi_+)=0$, then $f_{\chi}=0$. In other words, if $m(w, j)$ denotes the $(w,j)$-component of $\chi\in X(S_p)=\prod_w \Z^{Q(w)}$,
$f_{\chi}=0$ unless $m(w,j)\leq 0$ for any $w\in \Z$ and $j\in J(w)$. Here $\chi_+$, $f_{\chi}$, and $J(w)$ are as in \ref{thm+2}. 
\smallskip

(ii) $k\in S_J$. Here $S_J$ is as in \ref{thm+2}. 

\smallskip

Regard $X$ as an object of $\cB'_\R(\log)$ in the natural way, and regard $Y\subset X$ as an object of $\cB'_\R(\log)$ by \ref{embstr}. 

Both in the situations (bd) and (d), let
$$Y_0=\{(t,f,k)\in Y\;|\; t\in A_{p.P}\}\subset Y.$$

\end{sbpara}

  \begin{sbthm}\label{lsSB} Let the notation be as in \ref{lsSB1}. Consider the situation (bd) $\frak D=D^{\star,\BS}_{\SL(2)}$ and $\frak E=D_{\SL(2)}(\gr^W)^{\BS}$ (resp.\ (d) $\frak D= D_{\BS}$ and $\frak E= D_{\BS}(\gr^W)$).

  \smallskip
  (1) For a sufficiently small open neighborhood $U$ of $(0,0,0)$ in $Y$, there exists a unique open immersion $U\to \frak E$ in $\cB'_\R(\log)$ which sends $(t, f,  k)\in U\cap Y_0$  to the element $$\exp(f){\tau}^{\star}_p(t) \exp(k) \br \quad (\text{resp.}\; t\circ \exp(f)\exp(k)\br)
 $$ of $D(\gr^W)\subset \frak E$.
  
    \smallskip
  (2) Let $\bar L=\bar \cL(\br)$ and $L=\cL(\br)$. Then for a sufficiently small open neighborhood $U$ of $(0,0,0)$ in $Y$, there exists a unique open immersion $U\times \spl(W) \times \bar L\to \frak D$ in $\cB'_\R(\log)$ having the following property.
   It sends
   $(t, f, k,s,\delta)\in Y\times \spl(W) \times L$, where $(t,f,k) \in U\cap Y_0$, $s\in \spl(W)$, and $\delta\in L$, to the element of $D$ 
    (resp.\ to the element $t\circ x$ where $x$ is the element of $D$) 
whose image in 
$D(\gr^W) \times \spl(W) \times \cL$ under the isomorphism \ref{grsd} is 
$$(\exp(f){\tau}^{\star}_p(t) \exp(k) \br, s, \Ad(\exp(f)\tau^{\star}_p(t)\exp(k))\delta)$$
$$(\text{resp.}\;
 (\exp(f)\exp(k) \br, s, \Ad(\exp(f)\exp(k))\delta)).$$ 

\smallskip

(3) For a sufficiently small open neighborhood $U$ of $(0,0,0)$ in $Y$, the diagram 
 $$\begin{matrix}  U\times \spl(W)\times \bar L &\to& \frak D\\
 \downarrow &&\downarrow \\
 U&\to& \frak E
 \end{matrix}$$
 is cartesian in $\cB'_\R(\log)$ and in the category of topological spaces.

\smallskip

(4) The image of the map in (1) is contained in $\frak E(Q)\cap \frak E(P)$ with $Q=(\cW(p_w))_w$ (resp.\ in $\frak E(P)$) and the image of the map in (2) is contained in $\frak D(Q)\cap \frak D(P)$ (resp.\ in $\frak D(P)$).
\smallskip

(5) The underlying maps of the morphisms in (1) and (2) are described as in \ref{lsSB2} below.

 \end{sbthm}
 
 \begin{sbpara}\label{lsSB2} The maps in (1) and (2) in Theorem \ref{lsSB} are induced from the maps
 $$Y\to \frak E, \quad Y \times \spl(W) \times \bar L\to \frak D,$$
  respectively, defined as follows. 
  
  We consider first the situation (bd) in \ref{lsSB1}.

Let $A'$ be the subset of $A_{p,P}=\Hom(X(S_{p,P}), \R^{\mult}_{>0})$ consisting of all elements whose restriction to $t^{-1}(\R_{>0})\subset X(S_{p,P})^+$ coincides with the restriction of $t: X(S_{p,P})^+\to \R_{\geq 0}$ where $t$ ranges over $\bar A_{p,P}$. 
Let $J=(J(w))_w$ for $t$ be as in \ref{lsSB1} and let $p_J\in D_{\SL(2)}(\gr^W)$ be as in \ref{thm+4} for $J$. Then the first map $Y\to \frak E$ sends $(t,f,k)$ to 
  $$p':=\exp(f) \tau^{\star}_p(t')\exp(k)p_J \quad \text{where} \;t'\in A'.$$

The second map $Y \times \spl(W) \times \bar L\to \frak D$ 
sends $(t,f,k,s,\delta)$ to the following element $(p', P', Z)$ of $\frak D=D^{\star,\BS}_{\SL(2)}$ (\ref{2.6.3}) where $P'$ and $Z$ are as follows.

Let $\bar A_{p,P}\to \bar A_P= \R_{\geq 0}^{\Delta(P)}$ be the morphism in Lemma \ref{lemSB}. Let $I=\{j\in \Delta(P)\;|\; t_j=0\}$ where $t_j$ denotes the $j$-component of the image of $t$ in $\R_{\geq 0}^{\Delta(P)}$. Then $P'$ is the $\Q$-parablolic subgroup of $G_\R(\gr^W)$ such that $P'\supset P$ which corresponds to the subset $I$ of $\Delta(P)$ (\ref{BS2}). 

If $\delta\in L$, $Z$ is the subset of $D$ whose image under the embedding $D\to D(\gr^W)\times\spl(W) \times \cL$ is the set 
$$\{(\exp(f) \tau_p^{\star}(t') \exp(k)\br, s, \Ad(\exp(f) \tau^{\star}_p(t') \exp(k))\delta)\;|\; t'\in A'\}.$$ 
If $\delta= 0\circ \delta^{(1)}\in \bar L\smallsetminus L$ $(\delta^{(1)}\in L\smallsetminus \{0\})$ (\ref{2.3ex} (4)), 
$Z$ is the subset of $D$ whose image under the embedding $D\to D(\gr^W)\times\spl(W) \times \cL$ is the set 
$$\{(\exp(f) \tau_p^{\star}(t') \exp(k)\br, s, \Ad(\exp(f) \tau^{\star}_p(t') \exp(k))(c\circ \delta^{(1)})\;|\; t'\in A', c\in \R_{>0}\}.$$ 

Next consider the situation (d) in \ref{lsSB1}. In this situation, the first map sends $(t,f,k)$ to $t\circ \exp(f)\exp(k)\br$. The second map sends $(t,f,k,s, \delta)$ with $\delta\in L$ to the element $t\circ x$ and sends $(t,f,k,s, 0\circ \delta)$ with $\delta\in L\smallsetminus \{0\}$ to the element $(0,t)\circ x$, where $x$ is the element of $D$ whose image in $D(\gr^W) \times \spl(W) \times \cL$ (\ref{grsd}) is $(\exp(f)\exp(k)\br, s, \Ad(\exp(f)\exp(k))\delta)$.
 Here we denote by $(t,x) \mapsto t\circ x$ the morphisms $\bar A_P\times D(\gr^W) \to D_{\BS}(\gr^W)$, $\bar A_P \times D\to D_{\BS}$, and $\bar B_P\times D_{\nspl}\to D_{\BS}$, which extends the morphisms $A_P\times D(\gr^W) \to D_{\BS}(\gr^W)$, $A_P\times D\to D_{\BS}$, and $B_P\times D\to D_{\BS}$, defined by $(t,x) \mapsto t \circ x$, respectively. 

\end{sbpara}

  \begin{sbpara}\label{2.6.21} We prove Theorem \ref{lsSB}.
  
  The theorem is clear in the situation (d) in \ref{lsSB1}.

  We consider the situation (bd) in \ref{lsSB1}. We reduce the theorem in this situation to Theorem \ref{ls1}. 
  
    It is easily seen that the validity of the theorem does not depend on the choices of $R$ and $S$. 
We take any $S$ satisfying the condition in \ref{thm+2}, and hence the condition in \ref{lsSB1}. 
We choose $R$ in the following way.
    
 Let $Q=(Q(w))_w$ where $Q(w)=\cW(\gr^W_w)$.  
 Take a splitting $\alpha$ of $Q$ and let $R(Q)$ be as in \ref{Pcone2.2}.
Let $\sig_P\in \cP'(Q)$ be the cone corresponding to $P\in \cP(Q)$ (\ref{Pcone3})
 and let $\cS:=\cS(\sig_P)$ be the corresponding fs submonoid of $X(S_p)$. Note that $R(Q)\subset \cS\cup \cS^{-1}$ (\ref{Pcone}). 
 
 Choose a subset $I_1$ of $R(Q)\cap \cS\cap \cS^{-1}$ such that $R(Q)\cap \cS\cap \cS^{-1}$ is the disjoint union of $\{1\}$, $I_1$, and $I_1^{-1}$. Let $I_2:=R(Q) \cap \cS^{-1}\smallsetminus R(Q)\cap \cS\cap \cS^{-1}$. Hence $R(Q)$ is the disjoint union of $\{1\}$, $I_1$, $I_1^{-1}$, $I_2$, and $I_2^{-1}$. 
 Choose an $\R$-subspace $C$ of $\fg_\R(\gr^W)$ such that the subspace $\fg_\R(\gr^W)_1= \{v\in \fg_\R(\gr^W)\;|\; \Ad(\tau^{\star}_p(t))v=v\;\text{for all}\;t\in A_p\}$ of $\fg_\R(\gr^W)$ coincides with the direct sum of $\Lie(\tau^{\star}_p(A_p))$, $C$, and $\frak g_\R(\gr^W)_1\cap \Lie(K_{\br})$. 
 Let 
 $$R'= C \oplus (\bigoplus_{\chi\in I_1\cup I_2} \fg_\R(\gr^W)_{\chi}).$$
 Then $R'\subset \Lie(P)$, and $R'$ satisfies the conditions (C1) and (C2) on $R$ of \ref{thm+2}. (Here we used the fact that the Cartan involution $\theta_{K_{\br}}$ associated to the maximal compact subgroup $K_{\br}$ of $G_\R(\gr^W)$ sends $\fg_\R(\gr^W)_{\chi}$ to $\fg_\R(\gr^W)_{\chi^{-1}}$ for any $\chi\in X(S_p)$, and $\Lie(K_{\br})$ coincides with  $\{v\in \fg_\R(\gr^W)\;|\; \theta_{K_{\br}}(v)=v\}$.)

     Take an $\R$-subspace $C'$ of $\fg_\R(\gr^W)$ such that $\Lie(\tau^{\star}_p(A_p))= \Lie(\tau^{\star}_p(A_{p,P})) \oplus C'$. We take $R:=C' \oplus R'$ as $R$ of \ref{lsSB1}.

    Define $X$ and $Y$ of
    \ref{lsSB1} by using these $R$ and $S$. Denote by $X'$ and $Y'$, respectively, the $X$ and $Y$ of \ref{thm+2} defined by taking $R'$ and $S$ as $R$ and $S$ of \ref{thm+2}.

 As in \ref{Pcone4}, let $\Sig$ be the fan of all faces of the cone $\Hom(X(S_p)^+, \R^{\add}_{\geq 0}) = \prod_w \R^{Q(w)}_{\geq 0}$. Let 
 $\Sig'$ be the fan of all faces of the cone $\sig_P$.

 Then $D^{\star,\BS}_{\SL(2)}(Q) \cap D^{\star,\BS}_{\SL(2)}(P)$ represents the fiber product of $\Mor(-, D^{\star}_{\SL(2)}(Q)) \to [\Sig] \leftarrow [\Sig']$. On the other hand, the fiber product of $\Mor(-, X') \to [\Sig]\leftrightarrow [\Sig']$ is represented by $X'':= \Hom(\cS, \R^{\mult}_{\geq 0}) \times \frak g_\R(\gr^W) \times \frak g_\R(\gr^W) \times \frak g_\R(\gr^W) \times S$ and the fiber product of $\Mor(-, Y') \to [\Sig]\leftarrow [\Sig']$ is represented by the inverse image $Y''$ 
  of $Y'$ in $X''$ under the canonical map $X''\to X'$, where $Y''$ is endowed with the 
  structure of an object of $\cB'_\R(\log)$ by using the embedding $Y''\to X''$ (\ref{embstr}). 
 
 We identify $X$ with $\Hom(\cS, \R^{\mult}_{\geq 0}) \times R'\times S$ via the isomorphism $\Hom(\cS, \R^{\mult}_{\geq 0})\cong \bar A_{p,P} \times C'$.

 To reduce Theorem \ref{lsSB} to Theorem \ref{ls1}, 
it is sufficient to prove the following (*).
 
 \smallskip
 
(*)   If $(t,f,g,h,k)\in Y''$, then $(t,f,k)\in Y$ in $X=\Hom(\cS, \R_{\geq 0}^{\mult})\times R' \times S$. We have an isomorphism 
$$Y''\overset{\cong}\to Y\;;\; (t,f,g,h,k) \mapsto (t, f,k)$$
 in $\cB'_\R(\log)$.
 
 \smallskip

 Before the proof of (*), we note the following (1) and (2).
 
 \smallskip
 
 (1) Let $(t,f,g,h,k)\in X''$ ($t\in \Hom(\cS, \R^{\mult}_{\geq 0})$, $f,g,h\in \fg_\R(\gr^W)$, $k\in S$). Then $(t,f,g,h,k)$ belongs to $Y''$ if and only if the conditions (i)--(iv) in \ref{thm+2}, among which (iii) and (iv) are modified as follows, are satisfied. We replace $R$ in (iii) in \ref{thm+2} by $R'$. In (iv) in \ref{thm+2}, we define $J= (J(w))_{w\in \Z}$ where $J(w)= \{j\in Q(w)\;|\; t_{w,j}=0\}$. Here $t_{w,j}\in \R_{\geq 0}$ denotes the $(w,j)$-component of the image of $t$ in $\bar A_p$. Then $k\in S_J$. 
 
 \smallskip
 
 (2) Let $(t,f,k)\in X$ ($t\in \Hom(\cS, \R^{\mult}_{\geq 0})$, $f\in R'$, $k\in S$). Then $(t,f,k)$ belongs to $Y$ if and only if the following conditions (2-i) and (2-ii) are satisfied.
 
 \smallskip
 (2-i) Let $\chi\in X(S_p)$. If $t(\chi_+)=0$, then $f_{\chi}=0$. 
 
 \smallskip
 
 (2-ii) The same as the form of (iv) in the above (1).

 \smallskip
Now we prove the assertion (*). Let $(t,f,g,h,k)\in Y''$. We first prove that $(t,f,k)\in Y$. To show this, it 
  is sufficient to prove $f\in R'$. 
  Let $\chi\in R(Q)$. If $t(\chi_-)\neq 0$, since $t(\chi_+)g_{\chi}=t(\chi_-)f_{\chi}$ 
  and $g_{\chi}\in R'$, we have $f_{\chi}=t(\chi_-)^{-1}t(\chi_+)g_{\chi}\in R'$. 
  Assume $t(\chi_-)=0$. If $\chi\in \cS$, then $t(\chi_+)=t(\chi_-)t(\chi)=0$. Hence $f_{\chi}=0$. If $\chi\notin \cS$, then $\chi\in \cS^{-1}$ and hence $f_{\chi}\in \fg_\R(\gr^W)_{\chi}\subset R'$. 
  
  We next prove that $Y'' \to Y$ is an isomorphism. 
  
  For this, we define a morphism 
$Y\to X''$ of the converse direction by
  $(t,f,k)\mapsto (t,f,g,h,k)$ with $g=\sum_{\chi\in \cS^{-1}} t(\chi^{-1})f_{\chi}$, $h= \sum_{\chi\in \cS^{-1}} t(\chi^{-1})^2f_{\chi}$. 
  
  We show that the image of this morphism is contained in $Y''$. Let $\chi\in R(Q)$.
   We prove $t(\chi_+)g_{\chi}= t(\chi_-)f_{\chi}$ and $t(\chi_+)h_{\chi}= t(\chi_-)g_{\chi}$. If $\chi\in \cS^{-1}$, 
   we have $t(\chi_+)g_{\chi}= t(\chi_+)t(\chi^{-1})f_{\chi}= t(\chi_-)f_{\chi}$ and 
   $t(\chi_+)h_{\chi}= t(\chi_+)t(\chi^{-1})^2f_{\chi}= t(\chi_-)g_{\chi}$. If $\chi\notin \cS^{-1}$, 
   we have $f_{\chi}=0$ by the definition of $R'$, and hence $g_{\chi}=h_{\chi}=0$ 
   by the definitions of $g$ and $h$. If $t(\chi_+)=0$, 
   then $f_{\chi}=0$ and hence $g_{\chi}=0$. We prove that if $t(\chi_-)=0$, then 
  $g_{\chi}=h_{\chi}=0$.
  In the case $t(\chi_+)=0$, then $f_{\chi}=0$ and hence $g_{\chi}=h_{\chi}=0$.  In the case $\chi\in \cS$, we have $t(\chi_+)=t(\chi_-)t(\chi)=0$. In the case $t(\chi_+)\neq 0$ and $\chi\notin \cS$, we have $\chi\in \cS^{-1}$ and $t(\chi_-)=t(\chi_+)t(\chi^{-1})$ and hence $t(\chi^{-1})=0$. Hence $g_{\chi}=t(\chi^{-1})^2f_{\chi}=0$ and $h_{\chi}=0$ similarly. We prove $g_{\chi}, h_{\chi}+f_{\chi^{-1}}\in R'$. 
   If $\chi\in \cS^{-1}$, $g_{\chi}=t(\chi)^{-1}f_{\chi}\in R'$ and $h_{\chi}= t(\chi^{-1})^2\in R'$ and hence $h_{\chi}+f_{\chi^{-1}}\in R'$. If $\chi\notin \cS^{-1}$, $g_{\chi}=h_{\chi}=0$ and hence $h_{\chi}+f_{\chi^{-1}}=f_{\chi^{-1}}\in R'$. 
   
   Thus we have a morphism $Y\to Y''$. It is clear that the composition $Y\to Y''\to Y$ is the identity morphism. We prove that the composition $Y''\to Y\to Y''$ is also the identify morphism. Let $(t,f,g,h,k)\in Y''$ and let $(t,f,g', h', k)$ be the image of $(t,f,k)\in Y$ under 
   $Y\to Y''$. We prove $g'_{\chi}=g_{\chi}$ and $h'_{\chi}=h_{\chi}$ for 
   any $\chi\in R(Q)$. Assume first $\chi\in \cS^{-1}$. If $t(\chi_+)\neq 0$, 
   then $g_{\chi}= t(\chi_+)^{-1}t(\chi_-)f_{\chi}= t(\chi^{-1})f_{\chi}=g'_{\chi}$, and we have similarly $h_{\chi}=h_{\chi}'$. 
   If $t(\chi_+)=0$, then $t(\chi_-)= t(\chi_+)t(\chi^{-1})=0$, and hence $f_{\chi}=g_{\chi}=h_{\chi}=0$, and we have $g'_{\chi}=0$ and $h'_{\chi}=0$ by $f_{\chi}=0$. Next assume $\chi\notin \cS^{-1}$.
    Then by the definition of $R'$, we have $a_{\chi}=0$ for any $a\in R'$. 
    Since $f_{\chi}, g_{\chi}, h_{\chi}+f_{\chi^{-1}}\in R'$, we have $f_{\chi}=g_{\chi}=h_{\chi}=0$, and we have $g'_{\chi}=h'_{\chi}=0$ by $f_{\chi}=0$. 
Theorem \ref{lsSB} is proved.

       \end{sbpara}

 \begin{sbthm}\label{thmBS}  (1) The identity map of $D$ extends uniquely to a morphism $D^{\star,\BS}_{\SL(2)}\to D_{\BS}$ in $\cB'_\R(\log)$. It sends $(p, P,, Z)\in 
 D^{\star,\BS}_{\SL(2)}$ to $(P, A_P\circ Z)\in D_{\BS}$.

 \smallskip
 
 (2) The diagram
 $$\begin{matrix} D^{\star,\BS}_{\SL(2)} & \to & D_{\BS}\\
 \downarrow &&\downarrow \\
 D_{\SL(2)}(\gr^W)^{\BS} &\to & D_{\BS}(\gr^W)\end{matrix}$$
 is cartesian in $\cB'_\R(\log)$ and also in the category of  topological spaces.
 
 \smallskip
 
(3) The inverse image of $D^{\mild}_{\BS}$ in $D^{\star,\BS}_{\SL(2)}$ coincides with $D^{\star,\BS,\mild}_{\SL(2)}$. 

   \end{sbthm}

\smallskip

 \begin{pf} Let $(p, P, Z)\in D_{\SL(2)}^{\star,\BS}$ and take $\br \in Z$. 
 We compare the situations (bd) and (d) in Theorem \ref{lsSB} by taking $p$, $\br$ 
 for both the situations (bd) and (d), and by taking $R$ and $S$ for these situations as follows. 
 Take $R$ and $S$ for the situation (d). Take this $S$ as $S$ for the situation (bd). 
 Let $C$ be an $\R$-subspace of $\Lie((A_P)_{\br})$ such that $\Lie((A_P)_{\br})= \Lie(\tau^{\star}(A_{p,P})) \oplus C$ and take $C\oplus R$ as the $R$ for the situation (bd). 
 Then Theorem \ref{thmBS} (1) and (2) follow from Theorem \ref{lsSB}, Lemma \ref{lemSB}, 
 and the fact
 $$(\tau^{\star}(t) \bmod P_u)\circ \exp(f)\exp(k)\br= \exp(f)\tau_p^{\star}(t) \exp(k)\br.$$
 Theorem \ref{thmBS} (3) is clear. 
    \end{pf}

\subsection{The category $\cB'_\R(\log)^+$}\label{ss:+}

 The aim of this Section \ref{ss:+} is to define a full subcategory $\cB'_\R(\log)^+$ 
of $\cB'_\R(\log)$, consisting of nice objects, and prove that the spaces of $\SL(2)$-orbits in this Section 2 belong to $\cB'_\R(\log)^+$ (Theorem \ref{slis+}). 

We discuss also full subcategories $\cB'_\R(\log)^{[+]}$ and $\cB'_\R(\log)^{[[+]]}$ of $\cB'_\R(\log)$ such that
 $$\cB'_\R(\log)\supset \cB'_\R(\log)^+\supset \cB'_\R(\log)^{[+]}\supset \cB'_\R(\log)^{[[+]]}.$$

\begin{sbpara}\label{b[[+]]}  

We first define a full subcategory  $\cB'_\R(\log)^{[[+]]}$ of $\cB'_\R(\log)$.

We define standard objects of $\cB'_\R(\log)^{[[+]]}$. 
Take $n\geq 0$, a real analytic manifold $A$, and a real analytic closed submanifold $A_J$ of $A$ for each subset $J$ of $\{1,\dots, n\}$ satisfying $A_{\emptyset}=A$, $A_J\subset A_{J'}$ if $J\supset J'$. Define
$$Y= \{(t, x)\in \R^n_{\geq 0}\times A\;|\; x\in A_{J(t)}\}$$
where $J(t)= \{j\;|\;1\leq j\leq n, t_j=0\}$. We regard $Y$ as an object of $\cB'_\R(\log)$  by taking
$\R^n_{\geq 0} \times A$ as $X$ in \ref{embstr} where the log structure of $X$ with sign is induced  from that of $\R^n_{\geq 0}$ (\ref{2.3ex} (1)).

Let $\cB'_\R(\log)^{[[+]]}$ be the full subcategory of $\cB'_\R(\log)$ consisting of all objects 
which are locally isomorphic to open subobjects of $Y$ as above.

Real analytic manifolds with corners belong to $\cB'_\R(\log)^{[[+]]}$. 

\end{sbpara}

\begin{sbpara}\label{b[+]}

We next define a full subcategory $\cB'_\R(\log)^{[+]}$ of $\cB'_\R(\log)$.

We define standard objects of $\cB'_\R(\log)^{[+]}$. Take an fs monoid $\cS$, a real analytic manifold $A$, and a real analytic closed submanifold $A_I$ of $A$ for each face $I$ of $\cS$ satisfying $A_{\cS}=A$, $A_I\subset A_{I'}$ if $I\subset I'$. Define
$$Y= \{(t, x)\in \Hom(\cS, \R^{\mult}_{\geq 0})\times A\;|\; x\in A_{I(t)}\}$$
where $I(t)$ is the face $ \{a\in \cS\;|\;t(a) \neq 0\}$ of $\cS$. We regard $Y$ as an object of $\cB'_\R(\log)$  by taking
$\Hom(\cS, \R^{\mult}_{\geq 0}) \times A$ as $X$ in \ref{embstr} where the log structure of $X$ with sign is induced  from that of $\Hom(\cS, \R^{\mult}_{\geq 0})$ (\ref{2.3ex} (3)).

Let $\cB'_\R(\log)^{[+]}$ be the full subcategory of $\cB'_\R(\log)$ consisting of all objects 
which are locally isomorphic to open subobjects of $Y$ as above. 

Since a standard object of $\cB'_\R(\log)^{[[+]]}$ is  the case $\cS=\N^n$ of a standard object of $\cB'_\R(\log)^{[+]}$, we have 
$\cB'_\R(\log)^{[+]}\supset \cB'_\R(\log)^{[[+]]}$.

\end{sbpara}

The following Lemmas \ref{lem++} and \ref{lem++2} are proved easily. 
\begin{sblem}\label{lem++} Let $S$ be an object of $\cB'_\R(\log)^{[+]}$. Then $S$ belongs to $\cB'_\R(\log)^{[[+]]}$ if and only if
for any $s\in S$, $(M_S/\cO^\times_S)_s$ is isomorphic to $\N^r$ for some $r\geq 0$ (which may depend on $s$). 

\end{sblem}

\begin{sblem}\label{lem++2}

  Let $S'\to S$ be a log modification in $\cB'_\R(\log)$. If $S$ belongs to $\cB'_\R(\log)^{[+]}$, then $S'$ also belongs to $\cB'_\R(\log)^{[+]}$.

\end{sblem}

\begin{sbpara}\label{b+}  We define a full subcategory $\cB'_\R(\log)^+$ of $\cB'_\R(\log)$.

Let $\cB'_\R(\log)^+$ be the full subcategory of $\cB'_\R(\log)$ consisting of all objects $S$ such that locally on $S$, there is a log modification $S'\to S$ such that $S'$ belongs to $\cB'_\R(\log)^{[[+]]}$.

We have clearly $\cB'_\R(\log)^{[[+]]}\subset \cB'_\R(\log)^+$. 
\end{sbpara}

 \begin{sblem}\label{leq1}
 Let $S$ be an object of $\cB'_\R(\log)^+$ and assume that $(M^{\gp}_S/\cO_S^\times)_s$ is of rank $\leq 1$ as an abelian group for any $s\in S$. Then $S$ belongs to $\cB'_\R(\log)^{[[+]]}$. 
 \end{sblem}
 
 \begin{pf}
 This is because any log modification $S' \to S$ is an isomorphism.
 \end{pf}

\begin{sbprop}\label{[+]+} $\cB'_\R(\log)^{[+]}\subset \cB'_\R(\log)^+$. 
\end{sbprop}

\begin{pf}  Let $S$ be an object of $\cB'_\R(\log)^{[+]}$. Locally on $S$, by the resolution of singularity in toric geometry (\cite{O} p.23), there exists a log modification $S'\to S$ such that for any $s\in S'$, $(M_{S'}/\cO^\times_{S'})_s\cong \N^r$ for some $r$. By Lemmas \ref{lem++} amd \ref{lem++2}, 
$S'$ belongs to $\cB'_\R(\log)^{[[+]]}$. 
\end{pf} 

\begin{sbprop}\label{lm+} Let $S'\to S$ be a log modification in $\cB'_\R(\log)$. Then, $S$ belongs to $\cB'_\R(\log)^+$ if and only if $S'$ belongs to $\cB'_\R(\log)^+$.

\end{sbprop}

\begin{pf} 
First assume that $S$ belongs to $\cB'_\R(\log)^+$. We prove that $S'$ belongs to $\cB'_\R(\log)^+$. We may assume that $S$ belongs to $\cB'_\R(\log)^{[[+]]}$. Locally on $S$, there is a log modification $S''\to S$  which is a composition $S''\to S'\to S$ where the first arrow is a log modification and the second arrow is the given morphism, such that for any $s\in S''$, $(M_{S''}/\cO^\times_{S''})_s \cong \N^r$ for some $r$. By Lemmas \ref{lem++} and \ref{lem++2}, $S''$ belongs to $\cB'_\R(\log)^{[[+]]}$. Hence $S'$ belongs to $\cB'_\R(\log)^+$. 

Next assume that $S'$ belongs to  $\cB'_\R(\log)^+$. 
We prove that $S$ belongs to $\cB'_\R(\log)^+$. By the
 assumption, there are an open covering $(U_{\lam})_{\lam}$ of $S'$ and a log modification $V_{\lam}\to U_{\lam}$ for each $\lam$ such that $V_{\lam}$ belongs to $\cB'_\R(\log)^{[[+]]}$. Since $S'\to S$ is proper, locally on $S$, we can take a finite covering $(U_{\lam})_{\lam}$. Hence locally on $S$, there is a log modification $S''\to S$ having the following properties (i)--(iii). (i) $S''\to S$ is a composition $S''\to S'\to S$ where the first arrow is a log modification and the second arrow is the given morphism. (ii) For each $\lam$, we have a morphism $U_{\lam} \times_{S'} S'' \to V_{\lam}$ over 
 $U_{\lam}$ which is a log modification. (iii) For any $s\in S''$, $(M_{S''}/\cO^\times_{S''})_s \cong \N^r$ for some $r\geq 0$. By Lemma \ref{lem++} and \ref{lem++2}, $U_{\lam} \times_{S'} S''$ belongs to $\cB'_\R(\log)^{[[+]]}$. Since $(U_{\lam} \times_{S'} S'')_{\lam}$ is an open covering of $S''$, $S''$ belongs to $\cB'(\log)^{[[+]]}$. Hence $S$ belongs to $\cB'_\R(\log)^+$. \end{pf}

\begin{sbprop}\label{cv++3} 
The category  $\cB'_\R(\log)^{+}$ (resp.\ $\cB'_\R(\log)^{[+]}$, resp.\ $\cB'_\R(\log)^{[[+]]}$) is stable in $\cB'_\R(\log)$ under taking finite products.

\end{sbprop}
\begin{pf} This is clear for  $\cB'_\R(\log)^{[+]}$ and $\cB'_\R(\log)^{[[+]]}$. The part for  $\cB'_\R(\log)^+$ follows from the part for $\cB'_\R(\log)^{[[+]]}$. 
\end{pf}

 \begin{sblem}\label{fiber+} Let  $Y\subset \Hom(\cS, \R^{\mult}_{\geq 0})\times A$  be a standard object of $\cB'_\R(\log)^{[+]}$ in \ref{b[+]}, let $S$ be an object 
 of $\cB'_\R(\log)^{[[+]]}$, and  let $S\to \Hom(\cS, \R^{\mult}_{\geq 0})$ be a morphism in $\cB'_\R(\log)$. 
 Then the fiber product of $S\to \Hom(\cS, \R^{\mult}_{\geq 0}) \leftarrow Y$ in $\cB'_\R(\log)$ belongs to $\cB'_\R(\log)^{[[+]]}$. 
 
  \end{sblem}
 
 \begin{pf} Working locally on $S$, we may assume that $S$ is an open set of the standard object $\R^n_{\geq 0} \times A'$ in \ref{b[[+]]} $(A'$ here plays the role of  $A$ in \ref{b[[+]]}), and that we have a commutative diagram of functors
 $$\begin{matrix} \Mor(-, S)&&\to&& \Mor(-, \Hom(\cS, \R^{\mult}_{\geq 0}))\\ \downarrow &&&& \downarrow\\
 \Mor(-, \R^n_{\geq 0}) & \to & [\Sig'] &\to &[\Sig]
 \end{matrix}$$ where $\Sig$ is the fan of all faces of the cone $\Hom(\cS, \R^{\add}_{\geq 0})$ and $\Sig'$ is the fan of all faces of the cone $\R^n_{\geq 0}\subset \R^n$. Then the fiber product in problem coincides with 
 the space $$\{(t, a, a') \in \R^n_{\geq 0} \times A\times A'\;|\; a\in A_{I(t)}, a'\in A'_{J(t)}\}$$
 where $J(t)=\{j\;|\; 1\leq j\leq n, \; t_j=0\}$ and $I(t)$ is the face of $\cS$ which corresponds to the image of $t$ under $\R^n_{\geq 0} \to \Sig' \to \Sig$. 
  \end{pf}
  
 \begin{sblem}\label{fiber+2} Let  $Y\subset \Hom(\cS, \R^{\mult}_{\geq 0})\times A$ be a standard object of $\cB'_\R(\log)^{[+]}$ in \ref{b[+]}, let $S$ be an object of $\cB'_\R(\log)^+$, and  let $S\to \Hom(\cS, \R^{\mult}_{\geq 0})$ be a strict morphism in $\cB'_\R(\log)$. Then the fiber product of $S\to \Hom(\cS, \R^{\mult}_{\geq 0}) \leftarrow Y$ in $\cB'_\R(\log)$ belongs to $\cB'_\R(\log)^+$.

 \end{sblem}

\begin{pf}
Since $S\to \Hom(\cS, \R^{\mult}_{\geq 0})$ is strict, working locally on $\Hom(\cS, \R^{\mult}_{\geq 0})$ and on $S$, we have a rational finite subdivision $\Sig'$ of the cone $\Hom(\cS, \R^{\add}_{\geq 0})$ such that 
the fiber product $S'$ of $S\to \Hom(\cS, \R^{\mult}_{\geq 0})\leftarrow |\toric|(\Sig')$ belongs to $\cB'_\R(\log)^{[[+]]}$ and such that $\cS(\sig')$ for all $\sig'\in \Sig'$ are isomorphic to $\N^r\times \Z^m$ for some $r, m$. 
Replacing $S$ by $S'$ and replacing $\cS$ by $\cS(\sig')$ $(\sig'\in \Sig')$, we are reduced to Lemma \ref{fiber+}. 
\end{pf}

\begin{sbprop}\label{VV+} Let $n\geq 0$, and let $V$ be a finite dimensional $\R$-vector space endowed with an action of $\bG_m^n$.
 Let $Y$ be the subset of $\R^n_{\geq 0}\times V \times V$ consisting of all elements $(t, u, v)$ satisfying the following conditions $(i)$ and $(ii)$  for any $\chi\in X(\bG_m^n)$. In the following, we write $\chi=\chi_+(\chi_-)^{-1}$ as in \ref{thm+2}. 
 
\smallskip
$(i)$
$t(\chi_+)v_{\chi}= t(\chi_-)u_{\chi}$.

\smallskip

$(ii)$ If $t(\chi_+)=0$, then $u_{\chi}=v_{\chi}=0$.

\smallskip
Endow $Y$ with the structure of an object of $\cB'_\R(\log)$ by the embedding $Y\to \R^n_{\geq 0} \times V\times V$ as in \ref{embstr}. 
 Let $S$ be an object of $\cB'_\R(\log)^+$ and assume that we are given a strict morphism $S\to \R^n_{\geq 0}$, and let $E$ be the fiber  product of $S\to \R^n_{\geq 0}\leftarrow Y$ in $\cB'_\R(\log)$. Then $E$ belongs to $\cB'_\R(\log)^+$.

\end{sbprop}

\begin{pf} In \ref{rvtoric}, we take $L=X(\bG_m^n)$. Let $L^+\subset L$ be the submonoid  corresponding to $\N^n$ in the identification $L=\Z^n$. Take a finite subset $R$ of $L$ such that $\{\chi\in L\;|\; V_{\chi}\neq 0\}\subset R$, $R=R^{-1}$, and $R^+:=R\cap L^+$ generates $L^+$ as a monoid.  Let $\Sig$ be fan of all faces of the cone $\R^n_{\geq 0}\subset \R^n=N_\R$, and let $\Sig'$ be the rational finite subdivision of $\Sig$ defined in \ref{Pcone} (3) with respect to $R$ and $L^+$. 

Let $Y'$, $S'$, $E'$ be  the fiber products of $Y\to \R^n_{\geq 0}\leftarrow |\toric|(\Sig')$, $S\to \R^n_{\geq 0}\leftarrow |\toric|(\Sig')$, $E\to \R^n_{\geq 0}\leftarrow |\toric|(\Sig')$, respectively
(we identify $\R^n_{\geq 0}$ with $|\toric|(\Sig)$).  For $\sig'\in \Sig'$, let $Y'(\sig')$, $S'(\sig')$, $E'(\sig')$ be the open sets of $Y'$, $S'$, $E'$, respectively, corresponding to $\sig'$. These are the fiber products of
$Y\to \R^n_{\geq 0}\leftarrow \Hom(\cS(\sig'), \R^{\mult}_{\geq 0})$, $S\to \R^n_{\geq 0}\leftarrow \Hom(\cS(\sig'), \R^{\mult}_{\geq 0})$, $E\to \R^n_{\geq 0}\leftarrow \Hom(\cS(\sig'), \R^{\mult}_{\geq 0})$, respectively. In particular, $Y'(\sig') \subset \Hom(\cS(\sig'), \R^{\mult}_{\geq 0}) \times V \times V$.

We prove that $Y'(\sig')$ is isomorphic to a standard object of the category $\cB'_\R(\log)^{[+]}$ (\ref{b[+]}). 

Since $R\subset \cS(\sig')\cup \cS(\sig')^{-1}$ (\ref{Pcone}), we can take  subsets $R_1$ and $R_2$ such that $R$ is the disjoint union of $R_1$ and $R_2$ and such that $R_1\subset \cS(\sig')$ and $R_2\subset \cS(\sig')^{-1}$. 
Consider the map 
$$Y'(\sig') \to \Hom(\cS(\sig'), \R^{\mult}_{\geq 0})\times V\;;\; (t, u, v) \mapsto (t, \sum_{\chi\in A_1} v_{\chi} + \sum_{\chi\in A_2} u_{\chi}).$$
This induces an isomorphism $$Y'(\sig') \overset{\cong}\to  \{(t,x)\in \Hom(\cS(\sig'), \R^{\mult}_{\geq 0}) \times V\;|\; x\in V_{I(t)}\}\leqno{(1)}$$
in $\cB'_\R(\log)$, where $I(t)$ denotes the face $\{\chi\in \cS(\sig')\;|\; t(\chi)\neq 0\}$ of $\cS(\sig')$, and for a face $I$ of $\cS(\sig')$, we define
$$V_I = \{x\in V\;|\; x_{\chi}=0\;\text{if $\chi\in L$ and $\chi_+\notin I$}\}.$$
The inverse map of (1) is given by $(t, x)\mapsto (t, u, v)$ where $u=\sum_{\chi\in R_1} t(\chi)x_{\chi} + \sum_{\chi\in R_2} x_{\chi}$ and $v= \sum_{\chi\in R_1} x_{\chi} + \sum_{\chi\in R_2} t(\chi^{-1})x_{\chi}$. We omit the more details of the proof of this isomorphism (1), for the argument is straightforwards and similar to the proof of $Y''\cong Y$ in the proof of Theorem \ref{lsSB} (\ref{2.6.21}). 
Note that the right-hand side of (1) is a standard object of $\cB'_\R(\log)^{[+]}$ (\ref{b[+]}). 

By Lemma \ref{fiber+2}, the fiber product $E'(\sig')$ of $S'(\sig') \to \Hom(\cS(\sig'), \R^{\mult}_{\geq 0}) \leftarrow Y'(\sig')$ belongs to $\cB'_\R(\log)^+$. Hence $E'$ belongs to $\cB'_\R(\log)^+$. By Proposition \ref{lm+}, this proves that $E$ belongs to $\cB'_\R(\log)^+$. 
\end{pf}

\begin{sbprop}\label{SB+}
$D^{\star,\BS}_{\SL(2)}$ belongs to $\cB'_\R(\log)^{[+]}$.

\end{sbprop}

This follows form Theorem \ref{lsSB} for the situation (bd) in \ref{lsSB1} and from \ref{VV+}. 

\begin{sbthm}\label{slis+}

The spaces $D^I_{\SL(2)}$, $D^{II}_{\SL(2)}$, $D^{\star}_{\SL(2)}$, $D^{\star,+}_{\SL(2)}$, $D^{\star,-}_{\SL(2)}$, $D^{\star,\BS}_{\SL(2)}$, $D_{\SL(2)}(\gr^W)^{\sim}$ and $D_{\SL(2)}(\gr^W)$  belong to $\cB'_\R(\log)^+$. 

\end{sbthm}

 \begin{sbrem} We think that this Theorem \ref{slis+} is a version for the spaces of $\SL(2)$-orbits, treated in this Section 2, of the following results (1), (2) on the spaces of  Borel-Serre orbits and of nilpotent orbits.
 
 \smallskip
 
(1)  The space $D_{\BS}$ of Borel-Serre orbits is a real analytic manifold with corners. (Part I.)
 
 \smallskip
 
 (2) For a weak rational fan $\Sig$ in $\fg_\R$ and for a neat subgroup $\Gamma$ of $G_{\Z}$ which is strongly compatible with $\Sig$, the space $\Gamma \bs D_{\Sig}$ is a log manifold (Part III, Theorem 2.5.2). 
 
 \smallskip
 
 These (1) and (2) tell that $D_{\BS}$ and $\Gamma \bs D_{\Sig}$ are beautiful spaces. Theorem \ref{slis+} also says that the spaces of $\SL(2)$-orbits are beautiful spaces. 
 
 \end{sbrem}

  \begin{sbpara} We prove Theorem \ref{slis+}.
  Theorem \ref{slis+}  for $D^{\star,\BS}_{\SL(2)}$ follows from \ref{SB+} and \ref{[+]+}. Theorem \ref{slis+}  for $D^{\star}_{\SL(2)}$, $D^{\star,+}_{\SL(2)}$, $D^{\star,-}_{\SL(2)}$ follows from that for $D_{\SL(2)}^{\star,\BS}$ by \ref{lm+}. In the pure situation, this implies that $D_{\SL(2)}(\gr^W_w)$ belongs to $\cB'_\R(\log)^+$ for any $w$, hence $D_{\SL(2)}(\gr^W)$ belong to $\cB'_\R(\log)^+$, and hence $D_{\SL(2)}(\gr^W)^{\sim}$ belongs to $\cB'_\R(\log)^+$ by \ref{lm+}. 
 Theorem \ref{slis+} for $D^{II}_{\SL(2)}$ follows from that for $D_{\SL(2)}(\gr^W)^{\sim}$ by \ref{Lbund} and \ref{cv++3}. 
  
  Finally we prove that $D^I_{\SL(2)}$ belongs to $\cB'_\R(\log)^+$. We apply Proposition \ref{VV+}. Let $x=(p, Z)\in D_{\SL(2)}$, fix $\br\in Z$, and let $\bar \br=\br(\gr^W)\in D(\gr^W)$. Let $n$ be  $\text{rank}(p)$ if $x$ is an $A$-orbit, and let $n=\text{rank}(p)+1$ if $x$ is a $B$-orbit. Let $V=\Lie(G_{\R,u})$ where $G_{\R,u}$ denotes the unipotent radical of $G_{\R}$. 
    We define the action of $\bG_m^n$ on $V$ as follows. In the case $x$ is an $A$-orbit (resp.\ a $B$-orbit), lift the homomorphism $\tau_p\; (\text{resp.}\;  \tilde \tau_p) : \bG_m^n \to G_\R(\gr^W)$ (\ref{sim5}) to $\tau_x: \bG_m^n\to G_\R$ by using the splitting $\spl_W(\br)$ of $W$. We consider the adjoint action of $\bG_m^n$ on $\Lie(G_{\R,u})$ via $\tau_x$. 
    Define $Y\subset \R^n_{\geq 0}\times V \times V$ as in \ref{VV+}. Then by Part II, Theorem 3.4.6, in the case where $x$ is an $A$-orbit (resp.\ a $B$-orbit), there are an open neighborhood $S$ of $y:=(p, \delta_W(\br))$ in $D_{\SL(2)}(\gr^W)^{\sim}\times \cL(\bar \br)$ (resp.\ $y:=(p, 0\circ \delta_W(\br))$ in $D_{\SL(2)}(\gr^W)^{\sim}\times (\bar \cL(\bar \br)\smallsetminus \{0\}))$ and a strict morphism $S\to \R^n_{\geq 0}$ which sends $y$ to $0=(0,\dots,0)$, an open neighborhood $U$ of $(y, 0,0,0)$ in the fiber product of $S\to \R^n_{\geq 0} \leftarrow Y$ (here $(0,0,0)\in \R^n_{\geq 0}\times V \times V$) and an open immersion $U\to D^I_{\SL(2)}$ which sends $(y,0,0,0)$ to $x$. By Proposition \ref{VV+}, $U$ is an object of $\cB'_\R(\log)^+$. This shows that $D^I_{\SL(2)}$ is an object of $\cB'_\R(\log)^+$. 
Theorem \ref{slis+} is proved.      
      
  \end{sbpara}

\begin{sblem} Let $n\geq 0$. Then the part of $D_{\SL(2)}(\gr^W)^{\sim}$ consisting of points of rank $\leq n$ is open in $D_{\SL(2)}(\gr^W)^{\sim}$.

\end{sblem}

\begin{pf} This part is the union of open sets $D_{\SL(2)}(\gr^W)^{\sim}(\Phi)$ where $\Phi$ ranges over all admissible sets of weight filtrations on $\gr^W$ associated to points of rank $\le n$.
\end{pf}

\begin{sbpara}
We denote the above part of $D_{\SL(2)}(\gr^W)^{\sim}$ by $(D_{\SL(2)}(\gr^W)^{\sim})_{\leq n}$. In the pure situation, this part is written as $D_{\SL(2),\leq n}$.
\end{sbpara}

\begin{sbprop} (1) Let $U$ be the inverse image of $(D_{\SL(2)}(\gr^W)^{\sim})_{\leq 1}$ in $D^{\star}_{\SL(2)}$ (resp.\ $D^{II}_{\SL(2)}$). Then $U$ is an object of $\cB'_\R(\log)^{[[+]]}$. 

\smallskip

(2) Let $U$ be the inverse image of $\prod_w D_{\SL(2)}(\gr^W_w)_{\leq 1}$ in $D^{\star,-}_{\SL(2)}$. Then $U$ is an object of $\cB'_\R(\log)^{[[+]]}$. 
\end{sbprop}
\begin{pf} We prove (1). By \ref{slis+}, \ref{logst1} (which describes the stalks of $M_S/\cO^\times_S$ for $S= D_{\SL(2)}(\gr^W)^{\sim}$)  and \ref{leq1}, $(D_{\SL(2)}(\gr^W)^{\sim})_{\leq 1}$ belongs to $\cB'_\R(\log)^{[[+]]}$. Hence $U$ belongs to $\cB'_\R(\log)^{[[+]]}$ by \ref{Lbund} and \ref{cv++3}. 

We prove (2). Similarly, $D_{SL(2)}(\gr^W_w)_{\leq 1}$ belongs to $\cB'_\R(\log)^{[[+]]}$ and hence $\prod_w D_{\SL(2)}(\gr^W_w)_{\leq 1}$ belongs to $\cB'_\R(\log)^{[[+]]}$ by \ref{cv++3}. Hence $U$ belongs to $\cB'_\R(\log)^{[[+]]}$ by \ref{Lbund} and \ref{cv++3}.
\end{pf}

\section{Valuative Borel--Serre orbits and valuative SL(2)-orbits}\label{s:val}
In this Section \ref{s:val}, we study the spaces $D_{\BS,\val}$, $D_{\SL(2),\val}$, and $D^{\star}_{\SL(2),\val}$, and their relations.

\subsection{The associated valuative spaces}\label{ss:valsp}
In this Section \ref{ss:valsp}, 

\smallskip

(1) 
for an object $S$ of $\cB'_\bR(\log)$, we define a locally ringed space $S_{\val}$ over $\bR$ with a \lq\lq valuative log structure with sign" , and

\smallskip

(2) more generally, for a field $K$ endowed with a non-trivial absolute value $|\;\;|\;:\;K\to \R$ and for a  locally ringed space $S$ over $K$ endowed with an fs log structure satisfying the conditions 
 in \ref{value}, we construct a topological space $S_{\val}$.  
 
\smallskip
In (2), $S_{\val}$ is merely a topological space and does not have more structures as in (1).  

(1) becomes important in the rest of this Section 3, and (2) will become important in Section 4. 

(1) is shortly explained  in 
Part II, 3.7. 

We call $S_{\val}$ {\it the valuative space associated to $S$}. 
\begin{sbpara}\label{V(S)}
  Let $L$ be an abelian group  whose group law is written multiplicatively. 
 A submonoid $V$ of $L$ is said to be {\it valuative} if $V\cup V^{-1}=L$. 

  An integral monoid $V$ is said to be {\it valuative} if it is a valuative submonoid of $V^{\gp}$.

  For an fs monoid $\cS$, let $V(\cS)$ be the set of all valuative submonoids $V$ of $\cS^{\gp}$ such that $V\supset \cS$ and $V^\times \cap \cS=\cS^\times$. 
\end{sbpara}

\begin{sbpara}\label{Sgen} Let $K$ be a field endowed with a non-trivial absolute value $|\;\;|: K\to \R$. 
Let $S$ be a locally ringed space over $K$ 
satisfying  the equivalent conditions in \ref{value} and endowed with an fs log structure.

Let $S_{\val}$  be the set of all triples $(s,V,h)$, where $s\in S$, $V\in V((M_S/\cO^\times_S)_s)$ (\ref{V(S)}), and writing by $\tilde V $ the inverse image of $V$ in $M^{\gp}_{S,s}$, $h$ is a homomorphism 
 $(\tilde V)^\times \to \R^{\mult}_{>0}$ extending 
$f\mapsto |f(s)|$ on $\cO_{S,s}^\times$. Here $\R^{\mult}_{>0}$ denotes the set $\R_{>0}$ regarded as a multiplicative group. 
\end{sbpara}

\begin{sbpara} There is a variant, which we denote by $S_{\val(K)}$,  of $S_{\val}$: Let $S_{\val(K)}$ be the set of all triples $(s, V, h)$ where $s$ and $V$ are as above but $h$ is a homomorphism $(\tilde V)^\times \to K^\times$ extending $f\mapsto f(s)$ on $\cO_{S,s}^\times$. In \cite{KU}, Section 3.6, in the case $K=\C$, this space $S_{\val,\C}$ was denoted by $S_{\val}$. But in this Part IV, we consider only $S_{\val}$ in the sense of \ref{Sgen} except in the proof of \ref{perth1}, and we hope no confusion occurs in Part IV.  In the case a confusion can happen in the future, we denote $S_{\val}$ in \ref{Sgen} by $S_{\val(|\;\;|)}$. We will call $S_{\val(K)}$ the valuative space of $K$-points associated to $S$, and 
$S_{\val(|\;\;|)}$ the valuative space of absolute values associated to $S$.

\end{sbpara}

\begin{sbpara}
In the case $K=\R$ and $M_S$ is a log structure with sign (as in the case $S\in \cB'_\R(\log)$) (\ref{lsign}), $S_{\val}$ is identified with the set of all triples $(s,V,h)$, where $s\in S$, $V$ is an element of $V((M_S/\cO^\times_S)_s)$, and writing by $\tilde V_{>0} $ the inverse image of $V$ in $M^{\gp}_{S,>0, s}$, $h$ is a homomorphism 
 $(\tilde V)_{>0}^\times \to \R^{\mult}_{>0}$ extending 
$f\mapsto f(s)$ on $\cO_{S,>0,s}^\times$.
\end{sbpara}

\begin{sbpara}\label{valtop}  Let $S$ be as in \ref{Sgen}.
The topology of $S_{\val}$ is defined as follows. 

Let $(s_0,V_0,h_0)\in S_{\val}$.
Assume that 
we are given a chart $\cS\to M_S$ near the point $s_0\in S$. 
We introduce a fundamental system of neighborhoods of the point $(s_0,V_0,h_0)\in S_{\val}$.

Let $U$ be a neighborhood of $s_0$ in $S$, $I$ a finite subset of $\cS^{\gp}$ such that, for any $f\in I$, the image ${\bar f}_{s_0}$ of $f$ in $(M_S^{\gp}/\cO^\times_S)_{s_0}$ is contained in $V_0$, and $\varepsilon>0$.
Let $B(U,I, \varepsilon)$ be the set of all points $(s,V, h)$ of $S_{\val}$ satisfying the following conditions (i)--(iii).

(i) $s\in U$. 

(ii) For any $f\in I$, the image ${\bar f}_s$ of $f$ in $(M_S^{\gp}/\cO^\times_S)_s$ belongs to $V$.

(iii) For any $f\in I$, $|h(f) - h_0(f)|<\varepsilon$. Here we define $h(f)$ (resp.\ $h_0(f)$) to be $0$ unless ${\bar f}_s \in V^\times$ (resp.\ ${\bar f}_{s_0}\in V_0^\times$). 

Define a topology of $S_{\val}$ so that the sets $B(U,I, \varepsilon)$, where $U$, $I$, and $\varepsilon$ vary, form a fundamental system of neighborhoods of the point $(s_0,V_0,h_0)$.
This topology is independent of the choice of a chart $\cS$, and hence is well defined globally. 

\end{sbpara}

We now consider the relation of $S_{\val}$ and the projective limit of toric varieties for sub-divisions of fans. This will be used to prove properties of $S_{\val}$,  and to endow $S_{\val}$ in the case $S\in \cB'_\R(\log)$ with a structure of a locally ringed space over $\R$ and a log structure with sign.

\begin{sbpara}\label{l:val0} Let the notation be as in \ref{rvtoric}.

Let $V$ be a valuative submonoid of $L$. 
For a submonoid $\cS$ of $L$, we say that $V$ {\it dominates} $\cS$ if $\cS\subset V$ and $\cS^\times = \cS\cap V^\times$. 
For a rational finitely generated sharp cone $\sig$ in $N_\bR$, we say that $V$ dominates $\sig$ if $V$ dominates $\cS(\sig):=\{l \in L\,|\, l(\sig) \geq 0\}$.

For a rational fan $\Sig$ in $N_\bR$, $V$ dominates some cone in $\Sig$ if and only if $\cS(\sig)\subset V$ for some $\sig\in \Sig$. If $V$ dominates a cone in $\Sig$, then such a cone is unique and is the smallest cone $\sig\in\Sig$ such that $\cS(\sig)\subset V$. 

If $\Sig'$ is a rational finite subdivision of $\Sig$, 
$V$ dominates some cone in $\Sig$ if and only if $V$ dominates some cone in $\Sig'$. In this case, if $V$ dominates $\sig'\in \Sig'$, $V$ dominates the smallest cone $\sig\in \Sig$ such that $\sig'\subset \sig$. 
\end{sbpara}

\begin{sblem}
\label{l:val}
Let $\Sig$ 
be a finite rational fan in $N_\bR$.
Then we have a bijection from the set of all valuative submonoids $V$ of $L$, which dominate some cone in $\Sig$, onto the projective limit $\varprojlim \;\Sig'$, where $\Sig'$ ranges over all finite rational subdivisions of $\Sig$. 
  This bijection sends $V$ to 
$(\sig_{\Sig'})_{\Sig'}$, where $\sig_{\Sig'}$ denotes the cone in $\Sig'$ dominated by $V$. 
The inverse map is given by $(\sig_{\Sig'})_{\Sig'}\mapsto \small\bigcup_{\Sig'}
\cS(\sig_{\Sig'})$. 
\end{sblem}

\begin{pf}
Straightforward.
\end{pf}

\begin{sblem}
\label{l:val2}
Let $\Sig$ be a finite rational fan in $N_\bR$. 
Then we have the following bijection from $\varprojlim_{\Sig'}  |\toric|(\Sig')$, where $\Sig'$ ranges over all finite rational subdivisions of $\Sig$, to the set of all pairs $(V, h)$ 
of a valuative submonoid $V$ of $L$ dominating some cone in $\Sig$ and a homomorphism $h\colon V^\times \to \R_{>0}$. If $(x_{\Sig'})_{\Sig'}$ is an element of $\varprojlim_{\Sig'} |\toric|(\Sig')$ and $(\sig_{\Sig'}, h_{\Sig'})$ $(\sig_{\Sig'}\in \Sig', h_{\Sig'}: \cS(\sig_{\Sig'})\to \R_{\geq 0})$ is the pair corresponding to $x_{\Sig'}$ (\ref{rvtoric}), then the pair $(V, h)$ corresponding 
to $(x_{\Sig'})_{\Sig'}$ is as follows. 
$V=\small\bigcup_{\Sig'} \cS(\sig_{\Sig'})$, and $h$ is the homomorphism $V^\times\to\R_{>0}$ whose restriction to $\cS(\sig_{\Sig'})^\times$ is $h_{\Sig'}$ for any $\Sig'$. 
\end{sblem}

\begin{pf}
  This can be shown by using \ref{l:val}. 
\end{pf}

\begin{sbprop}\label{4.1.9} Let $S$ be as in \ref{Sgen}, and assume that we are given  a chart $\cS\to M_S$ with $\cS$ an fs monoid, let $L=\cS^{\gp}$, let $N=\Hom(L, \Z)$, and let $\Sig$ be the fan in $N_{\R}$ of all faces of the cone $\Hom(\cS, \R^{\add}_{\geq 0})$. Here $\R^{\add}$ denotes $\R_{\geq 0}$ regarded as an additive monoid. Then 
we have a cartesian diagram of topological spaces
$$\begin{matrix}  S_{\val} &\to & \varprojlim_{\Sig'}  |\toric|(\Sig')\\
\downarrow && \downarrow \\
S &\to& |\toric|(\Sig)=\Hom(\cS, \R^{\mult}_{\geq 0})\end{matrix}$$
where $\Sig'$ ranges over all finite rational subdivisions of $\Sig$, and the lower row sends $s\in S$ to the homomorphism $f\mapsto |f(s)|$ $(f\in\cS)$.  

\end{sbprop}

\begin{pf} For $s\in S$, let $\cS(s)=\cS(\sig)$ where $\sig$ is the element of $\Sig$ such that the image of $s$ in $|\toric|(\Sig)$ corresponds to a pair $(\sig, h)$ for some $h:\cS(\sig)^\times \to \R^{\mult}_{>0}$ (\ref{rvtoric}). Then $\cS(s)^\times$ coincides with the inverse image of $\cO^\times_{S,s}$ under the canonical map $\cS^{\gp}\to M^{\gp}_{S,s}$, and $\cS(s)$ is generated by $\cS$ and $\cS(s)^\times$. We have $\cS(\sig)/\cS(\sig)^\times \overset{\cong}\to (M_S/\cO^\times_S)_s$. 

By 
\ref{l:val2}, the fiber product $S \times_{|\toric|(\Sig)} \varprojlim_{\Sig'} |\toric|(\Sig')$ is identified with the set of 
all triples $(s, V, h)$ where $s\in S$, $V$ is a valuative submonoid of $\cS^{\gp}$ such that $V\supset \cS$ and 
$V^\times \cap \cS=\cS(s)^\times$, and $h$ is a homomorphism $V^\times \to \R_{>0}^{\mult}$ whose restriction  to $\cS(s)^\times$ 
coincides with the composition $\cS(s)^\times \to \cO_{S,s}^\times \to \R^{\mult}_{>0}$ 
where the last map is $f\mapsto |f(s)|$. By the isomorphism $\cS(\sig)/\cS(\sig)^\times \overset{\cong}\to (M_S/\cO^\times_S)_s$, a valuative submonoid $V$ of $\cS^{\gp}$ such that $V\supset \cS$ and  $V^\times \cap \cS=\cS(s)^\times$ corresponds bijectively to a valuative submonoid $V'$ contains $(M_S/\cO^\times_S)_s$ and $(V')^\times \cap (M_S/\cO^\times_S)_s=\{1\}$. 
Furthermore, if $\tilde V'$ 
denotes the inverse image of $V'$ in $M^{\gp}_{S,s}$, 
$(\tilde V' )^\times$ is the pushout of $V^\times \leftarrow \cS(s)^\times \to \cO_{S,s}^\times$. Hence $h$ corresponds to a homomorphism $h': (\tilde V')^\times \to \R^{\mult}_{>0}$ whose restriction to $\cO_{S,s}^\times$ coincides with $f\mapsto |f(s)|$. Hence we have a bijection $(s, V, h) \mapsto (s,V',h')$ from the fiber product to $S_{\val}$. 

In the converse map $(s,V',h')\mapsto (s, V, h)$, $V$ is the inverse image of $V'$ under the canonical map $\cS^{\gp}\to (M_S^{\gp}/\cO^\times_S)_s$ and $h$ is the homomorphism $V^\times\to \R^{\mult}_{>0}$ induced by $h'$.

By using these explicit constructions of the bijection between $S_{\val}$ and the fiber product, it is easy to see that this bijection is a homeomorphism.
\end{pf}

\begin{sbcor}\label{vproper} For $S$ as in \ref{Sgen}, the map $S_{\val}\to S$ is proper.

\end{sbcor}

\begin{sblem}\label{valsrt} 

Let $S$ and $S'$ be as in \ref{Sgen} and assume that we are given a strict morphism $S'\to S$ of locally ringed spaces over $\R$ with log structures (for the word \lq\lq strict'', see \ref{gest6}). Then the canonical map $S'_{\val} \to S'\times_S S_{\val}$ is a homeomorphism. 
\end{sblem}

\begin{pf} For any $s'\in S'$ with image $s$ in $S$, the canonical map  $(M_S/\cO^\times_S)_s\to (M_{S'}/\cO^\times_{S'})_{s'}$ is an isomorphism from the assumption. From this, we see that the map $S'_{\val}\to S'\times_S S_{\val}$ is bijective. Since this map is continuous and since both $S'_{\val}$ and $S'\times_S S_{\val}$ are proper over $S'$ (\ref{vproper}), this map is a homeomorphism. 
\end{pf}
\begin{sblem}\label{abslog}  Let $S$ be as in \ref{Sgen} and let
 $|S|$ be the topological space $S$ with the sheaf of all $\R$-valued continuous functions. Endow $|S|$ with the log structure $M_{|S|}$ associated to the composition $M_S\to \cO_S\to \cO_{|S|}$, where the second arrow is $f\mapsto |f|$,  which we regard as a pre-log structure. 
Here $|f|$ denotes the function $s\mapsto|f(s)|$.

 Then $M_{|S|}$ is an fs log structure, and we have a canonical homeomoprhism $|S|_{\val}{\cong} S_{\val}$.

\end{sblem}
\begin{pf} If $\cS\to M_S$ is a chart with $\cS$ an fs monoid, then the composition $\cS\to M_S\to M_{|S|}$ is also a chart. Hence $M_{|S|}$ is an fs log structure. The canonical map $(M_S/\cO^\times_S)_s\to (M_{|S|}/\cO^\times_{|S|})_s$ is an isomorphism for any $s\in S$, and hence we have a canonical bijection $|S|_{\val}\to S_{\val}$. It is easy to see that this is a homeomorphism. 
\end{pf}

\begin{sbpara}\label{3.1.13}
Assume now
 $S$ is an object of  $\cB'_\R(\log)$ (\ref{cblog}). We endow $S_{\val}$ with a sheaf $\cO_{S_{\val}}$  of rings and a log structure $M_{S_{\val}}$ 
 with sign as follows. Locally on $S$, take a positive chart $\cS\to M_{S,>0}$ (\ref{pchar}), let $\Sig$ be the fan of faces of the cone $\Hom(\cS, \R^{\add}_{\geq 0})$,  and for a rational finite subdivision $\Sig'$ of $\Sig$, regard $S(\Sig') := S\times_{|\toric|(\Sig)} |\toric|(\Sig')$ as an object of $\cB'_\R(\log)$ by taking the fiber product in $\cB'_\R(\log)$. Here we use the fact that the underlying topological space of this fiber product is the same as the fiber product of the underlying topological spaces by \ref{gest6} (i) and by the fact that $S\to |\toric|(\Sig)$ is strict.
 
 We define $\cO_{S_{\val}}$ (resp.\ $M_{S_{\val}}$) as the inductive limit of $\cO_{S(\Sig')}$ (resp.\ $M_{S(\Sig')}$) by using Proposition \ref{4.1.9}. This sheaf of rings and the log structure with sign
 are independent of the choice of the chart and hence defined globally. In fact, if we have two charts $\cS\to M_{S,>0}$ and $\cS'\to M_{S,>0}$, there is a third chart $\cS''\to M_{S,>0}$ with homomorphisms $\cS\to \cS''$ and $\cS'\to \cS''$ of charts. It is easy to see that the sheaf of rings and the log structure with sign given by the chart $\cS$  (resp.\  $\cS'$) are isomorphic to the ones given by the chart $\cS''$, and that the composite isomorphisms between the ones given by the chart $\cS$ and the ones given by the chart $\cS'$ are independent of the choice of the third chart $\cS''$. 

We call $\cO_{S_{\val}}$ the {\it sheaf of real analytic functions}. 
\end{sbpara}

\begin{sbpara}\label{lmval}
A log modification  $S'\to S$ in $\cB'_\R(\log)$ induces an isomorphism $$(S')_{\val}\overset{\cong}\to S_{\val}$$ of locally ringed spaces over $\R$ with log structures with sign. 

\end{sbpara}

\begin{pf} This is clear. \end{pf}

\begin{sbpara}\label{Vval1} For $S\in \cB'_\R(\log)$ and for $x=(s, V, h)\in S_{\val}$, $V$ is identified with the inverse image of $(M_{S_{\val}}/\cO^\times_{S_{\val}})_x$ under the canonical map $(M^{\gp}_S/\cO_S^\times)_s
\to (M^{\gp}_{S_{\val}}/\cO^\times_{S_{\val}})_x$, and $h: \tilde V_{>0}^\times \to \R^{\mult}_{>0}$ coincides with the composition $\tilde V_{>0}^\times \to \cO_{S_{\val},>0,x}\to \R^{\mult}_{>0}$ where the first arrow is induced from $\tilde V_{>0}\subset M_{S,>0,s}\to M_{S_{\val},>0,x}$ and the second arrow is $f\mapsto f(x)$. 

\end{sbpara}

\begin{sbpara}\label{Vval2}
  Let $S$ be a locally ringed space.
  Then a log structure $M$ on $S$ is said to be {\it valuative} if it is integral and satisfies the following condition.
  
  \smallskip
  
For any local section $f$ of $M^{\gp}$, locally we have either $f\in M$ or $f^{-1}\in M$, 
that is, if every stalk of $M$ is valuative. 

\smallskip
  
 By \ref{Vval1}, for $S\in \cB'_\R(\val)$, the log structure of $S_{\val}$ is valuative.
   \end{sbpara}

    \begin{sbpara}\label{gest9} Let $\cS_1$, $T$, $\bar T$, $Z$, $\bar Z$ and $T(s)\subset T$, $Z(s)\subset Z$ (for $s\in \bar Z$)  be as in \ref{gest1}. We give a description (1) below of the valuative space $\bar Z_{\val}$ associated to $\bar Z$,  as a set. This will be used in \ref{orbittr}.

    For a valuative submonoid $V$ of $\cS_1^{\gp}$, let $T(V): =\Hom(\cS_1^{\gp}/V^{\times}, \R_{>0})\subset T=\Hom(\cS_1, \R^{\mult}_{>0})$. 
    Then we have:
    
    \smallskip
    
    (1) $\bar Z_{\val}$ is identified with the set of all triples $(s, V, Z')$ where $s\in \bar Z$, $V$ is a valuative submonoid of $\cS_1^{\gp}$ such that $V\supset \cS_1$ and such that $V^{\times}\cap \cS_1=\text{Ker}(\cS_1\to (M_{\bar Z}/\cO^\times_{\bar Z})_s)$, and $Z'$ is a $T(V)$-orbit in $Z(s)$. 
      (Note that $Z(s)$ is a $T(s)$-torsor and $T(V) \subset T(s)$.)   
 \smallskip
 
 This is proved as follows. Let $L=\cS_1^{\gp}$ and let $\Sig$ be the fan of all faces of the cone $\Hom(\cS_1, \R^{\add}_{\geq 0})\subset N_\R=\Hom(L, \R)$. Then by \ref{4.1.9}, $\bar Z_{\val}$ is the projective limit of the log modifications of $\bar Z$ corresponding to rational finite subdivisions $\Sig'$ of $\Sig$. By \ref{l:val}, the projective limit of the sets $\Sig'$ is identified with the set of valuative cones $V$ as above. Hence the above description (1) of $\bar Z_{\val}$ 
   follows from the descriptions of log modifications of $\bar Z$ in \ref{gest1} as sets 
    by taking the projective limit.

      \end{sbpara}

\subsection{The category $\cC_\R(\val)^+$}

We define categories $\cC_\R(\val)$ and $\cC_\R(\val)^+\subset \cC_\R(\val)$. In Section 3.3, we will see that the valuative spaces associated to the spaces of SL(2)-orbits and the space of Borel-Serre orbits belong to $\cC_\R(\val)^+$.

\begin{sbpara}\label{satval} 

Let $\cC_\R(\val)$ be the category of objects of $\cC_\R$ (\ref{2.5.1}) endowed with a valuative log structure (\ref{Vval2}) with sign.

We have $\cC_\R(\sat) \supset \cC_\R(\val)$ (\ref{cblog} (2)).

\end{sbpara}

\begin{sbprop}\label{B++C+} Let $S$ be an object of $\cB'_\R(\log)^+$. Then $S_{\val}$ belongs to $\cC_\R(\val)$.

\end{sbprop}

For the proof, we use the following lemma.

\begin{sblem}\label{prosys}  Let $(S_{\lam})_{\lam}$ be a directed projective system in $\cC_\R$,  let $S$ be  the projective limit  of the topological spaces $S_{\lam}$, and endow $S$ with the inductive limit of the inverse images of $\cO_{S_{\lam}}$.
Assume that there is an open set $S'$ of $S$  satisfying the following conditions {{\rm(i)}} and {{\rm (ii)}}.

\smallskip
{{\rm(i)}}
 $S'$ belongs to $\cC_\R$. 
 
 \smallskip
 
{{\rm (ii)}} For any open set $U$ of $S$, the map $\cO_S(U)\to \cO_S(U\cap S')$ is injective.

  \smallskip
  
Then $S\in \cC_\R$.
\end{sblem}

\begin{pf} Let $\cF$ be the sheaf on $S$ of morphisms to $\R^n$ of locally ringed spaces over $\R$, where $\R^n$ is endowed with the sheaf of all real analytic functions. We have a morphism $a:\cF\to \cO_S^n$ by $f\mapsto (f^*(t_j))_{1\leq j\leq n}$ where $t_j$ 
are the standard coordinate functions of $\R^n$. We have also a morphism $b: \cO_S^n \to \cF$, which comes from the fact that  
since $S_{\lambda}$ belong to $\cC_\R$,   $\cO_S^n$ is regarded as the inductive limit of the inverse images of sheaves on $S_{\lambda}$ of morphisms to $\R^n$. 
As is easily seen, the composition $ab: \cO_S^n\to \cO_S^n$ is the identity morphism.  We prove that $ba: \cF\to \cF$ is the identity morphism. Let $f\in \cF(U)$ with $U$ an open set of $S$. It is easy to see that $f$ and $ba(f)$  induce the same underlying continuous maps $U\to \R^n$ which we denote by $g$. It remains to prove that the homomorphisms $g^{-1}(\cO_{\R^n})\to \cO_U$ given by $f$ and $ba(f)$ coincide. Since $\cO_S(V) \to  \cO_S(V\cap S')$ is injective for any open set $V$ of $U$, it is sufficient to prove that the restrictions of $f$ and $ba(f)$ to $U\cap S'$ coincide. But $S'$ belongs to $\cC_\R$, and hence $ab$ gives the identity morphism of $\cF|_{S'}$. 
\end{pf}

\begin{sbpara} We prove Proposition \ref{B++C+}. By \ref{lmval}, it is sufficient to prove this for objects of $\cB'_\R(\log)^{[[+]]}$. 
As in \ref{b[[+]]}, let $Y\subset \R^n_{\geq 0}\times A$ be a standard object of $\cB'_\R(\log)^{[[+]]}$. 
It is sufficient to prove that $Y_{\val}$ belongs to $\cC_\R$. Let $L=\Z^n$, let 
$\Sig$ be the set of all faces of the cone $\R^n_{\geq 0}\subset \Hom(L, \R^{\add})= \R^n$, and for a rational finite subdivision $\Sig'$ of $\Sig$, let 
$Y(\Sig')= \{(t, x)\in |\toric|(\Sig') \times A\;|\; x\in A_{J(t)}\}$ 
where $J(t)= \{j\;|\; 1\leq j\leq n, t_j=0\}$ with $t_j$ the $j$-th component of the image of $t$ in $|\toric|(\Sig) =\R^n_{\geq 0}$. We apply Lemma \ref{prosys} by taking  the projective system $(Y(\Sig'))_{\Sig'}$ in $\cC_\R$ as $(S_{\lambda})_{\lambda}$ 
and by taking the open set $\R^n_{>0} \times A$ of $Y_{\val}$ as $S'$. Then the projective limit $S$ in \ref{prosys} is $Y_{\val}$. The injectivity of $\cO_S(U) \to \cO_S(U\cap S')$ for any open set $U$ of $Y_{\val}$ is seen easily. Hence $Y_{\val}$ belongs to $\cC_\R$ by Lemma \ref{prosys}.

\end{sbpara}

\begin{sbpara}\label{c+} We define a full subcategory $\cC_\R(\val)^+$ of $\cC_\R(\val)$. 

This is the category of all objects which are locally isomorphic to open subobjects of $S_{\val}$ with objects $S$  of $\cB'_\R(\log)^+$. We can replace $\cB'_\R(\log)^+$ by $\cB'_\R(\log)^{[[+]]}$ in this definition, to get the same category $\cC_\R(\val)^+$. Hence $\cC_\R(\val)^+$  is the category of objects which are locally an open subobject
of $$Y_{\val}= \{(t,x)\in (\R^n_{\geq 0})_{\val}\times A\;|\; x\in A_{J(t)}\}.$$
Here $n$, $A$, $(A_J)_J$ and $Y$ are as in \ref{b[[+]]}  and $J(t)=\{j\;|\; 1\leq j\leq n, t_j=0\}$ where $t_j$ denotes the $j$-th component of the image of $t$ in $\R^n_{\geq 0}$. 

\end{sbpara}

\begin{sbprop}\label{cv+2} For any object $S$ of $\cB'_\R(\log)$ and for any object $X$ of $\cC_\R(\val)$, the canonical map $\Mor(X, S_{\val})\to \Mor(X, S)$ is bijective.

 Consequently, if  $S$ is an object of $\cB'_\R(\log)^+$,  $S_{\val}$ represents the functor $X\mapsto \Mor(X, S)$ from $\cC_\R(\val)$ to (Sets). 
\end{sbprop}

\begin{pf} It is sufficient to prove that in the situation of \ref{4.1.9},  the canonical map 
$\Mor(X,  S\times_{|\toric|(\Sig)} |\toric|(\Sig')) \to \text{Mor}(X,S)$ is bijective. Here $S\times_{|\toric|(\Sig)} |\toric|(\Sig')$ denotes the fiber product in $\cB'_\R(\log)$. By \ref{fiberpr}, it is the fiber product in $\cC_\R(\sat)$. Hence it is sufficient to prove that the map $\Mor(X, |\toric|(\Sig')) \to \Mor(X, |\toric|(\Sig))$ is bijective. This last fact is reduced to Proposition \ref{cv+0}. 
\end{pf} 

\begin{sbprop}\label{cv+3} 

$(1)$
The category $\cC_\R(\val)^+$ has finite products. 

We will denote the product in $\cC_\R(\val)^+$ as $X\times_{\val} Y$.

\smallskip

$(2)$ A finite product in $\cC_\R(\val)^+$ is a finite product in $\cC_\R(\val)$. 

 \smallskip

$(3)$ The functor $\cB'_\R(\log)^+\to \cC_\R(\val)^+\;;\;S\mapsto S_{\val}$ preserves finite products. 

\smallskip

$(4)$ For objects $S_1,\dots,S_n$ of $\cC_\R(\val)^+$, the product of $S_1,\dots, S_n$ in the category $\cC_\R(\sat)$ exists. As a topological space, it is the product of the topological spaces $S_j$. 

We will denote this product in $\cC_\R(\sat)$ by $S_1\times_{\sat}\dots \times_{\sat} S_n$.
\end{sbprop}

\begin{pf}  If $Y$ and $Y'$ are objects of $\cB'_\R(\log)^+$, then by Proposition \ref{cv++3}, $(Y\times Y')_{\val}$ is an object of $\cC_\R(\val)^+$ and for any object $X$ of $\cC_\R(\val)^+$, we have
$$\Mor(X, (Y\times Y')_{\val})= \Mor(X, Y\times Y') = \Mor(X, Y) \times \Mor(X, Y') = \Mor(X, Y_{\val}) \times \Mor(X, (Y')_{\val})$$
where the first and the third equalities follow from Proposition \ref{cv+2}, and the second equality follows from Proposition \ref{fiberpr} (2). This proves (1), (2), (3).

We prove (4). Locally on each $S_j$, we have $S_j=(S'_j)_{\val}$ for an object $S'_j$ of $\cB'_\R(\log)^+$. Locally in each $S'_j$, take a chart $\cS_j\to M_{S'_j}$, let $\Sig_j$ be the fan of all faces of the cone $\Hom(\cS_j, \R^{\add}_{\geq 0})$, and consider $S= \varprojlim \prod_{j=1}^n  S'_j \times_{|\toric|(\Sig_j)} |\toric|(\Sig'_j)$
where $\Sig'_j$ ranges over all rational finite subdivisions of $\Sig_j$. Endow $S$ with the inductive limit of the inverse images of $\cO$ and the log structures with sign of $ \prod_{j=1}^n  S'_j \times_{|\toric|(\Sig_j)} |\toric|(\Sig'_j)$. Then $S$ belongs to $\cC_\R(\sat)$ by Lemma \ref{prosys},  and is the product of $S_j$ in $\cC_\R(\sat)$.  This locally constructed $S$ glues to a global $S$. 
 \end{pf}

The following Lemma will be used in Section 3.4.

\begin{sblem}\label{timesval} Let $n\geq 0$ and let $S_j\to S'_j$ $(1\leq j\leq n)$ be morphisms in $\cC_\R(\val)^+$ having the Kummer property of log structure in the sense (K) below.
Let $S$ (resp.\ $S'$) be the product $S_1\times_{\sat}\dots\times_{\sat} S_n$   (resp.\ $S_1'\times_{\sat} \dots\times_{\sat} S_n'$) in the category $\cC_\R(\sat)$, and let $S_{\val}=S_1\times_{\val} \dots \times_{\val} S_n$ (resp.\ $S'_{\val}= S'_1\times_{\val} \dots \times_{\val} S'_n$) be the product in the category $\cC_\R(\val)^+$ (\ref{cv+3}). 

\smallskip
(K) We say that a morphism $X\to Y$ of a locally ringed spaces with log structures has Kummer property of log structure if for any $x\in X$ and the image $y$ of $x$ in $Y$, 
 the homomorphism $(M_Y/\cO^\times_Y)_y \to (M_X/\cO_X^\times)_x$ is injective, and for any $a\in (M_X/\cO_X^\times)_x$, there is $m\geq 1$ such that $a^m$ belongs to the image of $(M_Y/\cO^\times_Y)_y$.

\smallskip

Then the diagram
$$\begin{matrix} S_{\val}&\to &S'_{\val}\\
\downarrow &&\downarrow\\
S &\to &S'\end{matrix}$$
is cartesian in the category of topological spaces.

\end{sblem}
\begin{pf} The set $S_{\val}$ is identified with the set of all triples $(s,V, h)$ where $s\in S$, $V$ is a valuative submonoid of $(M^{\gp}_S/\cO^\times_S)_s$ such that $V\supset (M_S/\cO_S^\times)_s$ and $V\cap (M_S/\cO_S^\times)_s=\{1\}$, and $h$ is a homomorphism $(\tilde V)^\times \to \R_{>0}$, where $\tilde V$ denotes the inverse image of $V$ in $M_{S,>0,s}$ such that the restriction of $h$ to $\cO_{S,>0,s}^\times$ coincides with $f\mapsto f(s)$. Furthermore, $(M_S/\cO_S^\times)_s \cong 
\prod_{j=1}^n (M^{\gp}_{S_j}/\cO^\times_{S_j})_{s_j}$ where $s_j$ denotes the image of $s$ in $S_j$. The similar things hold for $S'$. From these, we see that the diagram is cartesian in the category of sets. Since $S_{\val}$ and the fiber product $E$ of $S\to S'\leftarrow S'_{\val}$ in the category of topological spaces are proper over $S$, we see that the canonical map $S_{\val}\to E$ is a homeomorphism. 
\end{pf}

\subsection{$D_{\BS,\val}$, $D^I_{\SL(2),\val}$, $D^{II}_{\SL(2),\val}$, $D^{\star}_{\SL(2),\val}$}\label{4.3}

\begin{sbpara} We define 
$$D_{\BS,\val}, \quad D^I_{\SL(2),\val}, \quad D^{II}_{\SL(2),\val}, \quad D^{\star}_{\SL(2),\val}$$
as the valuative spaces associated to the objects 
$$D_{\BS}, \quad D^I_{\SL(2)}, \quad D^{II}_{\SL(2)}, \quad D^{\star}_{\SL(2)}$$
 of $\cB'_\R(\log)$ (\ref{3.1.13}), respectively. 
 
 By \ref{slis+}, \ref{B++C+} and \ref{c+}, they belong to $\cC_\R(\val)^+$. 
 
 We call $D_{\BS,\val}$ the space of valuative Borel--Serre orbits, and call the other spaces $D^I_{\SL(2),\val}$ etc., spaces of valuative $\SL(2)$-orbits. 

$D^I_{\SL(2),\val}$ and $D^{II}_{\SL(2),\val}$ are identified as sets because $D^I_{\SL(2)}$ and $D^{II}_{\SL(2)}$ are identified as sets and the morphism $D^I_{\SL(2)}\to D^{II}_{\SL(2)}$ is strict (\ref{valsrt}). They are denoted just by $D_{\SL(2),\val}$ when we regard them just as sets. 
\end{sbpara}

\begin{sbpara}
Since a log modification induces an isomorphism of associated valuative spaces (\ref{lmval}), we have
$$D^{\star,+}_{\SL(2),\val} \overset{\cong}\to D^{\star}_{\SL(2),\val} \overset{\cong}\to D^{\star,-}_{\SL(2),\val} \overset{\cong}{\leftarrow} D^{\star,\BS}_{\SL(2),\val}.$$
Hence the morphisms 
$$D^{\star,+}_{\SL(2)} \to D^{II}_{\SL(2)}, \quad D^{\star,\BS}_{\SL(2)}\to D_{\BS}$$ (Section 2.5, Section 2.6) induce morphisms  
$$D^{\star}_{\SL(2),\val} \to D^{II}_{\SL(2),\val},\quad D^{\star}_{\SL(2),\val}\to D_{\BS,\val}.$$

\end{sbpara}

\begin{sbpara}\label{orbittr} These valuative spaces are described as sets as follows. Let the situations (a)--(d)  and the notation be as in \ref{sit2} and \ref{objass}. 
By \ref{gest9}, we have: 

As a set, 
$\frak D_{\val}$ is identified with the set of all triples $(x, V, Z)$ where $x\in \frak D$, $V$ is a valuative submonoid of $(M^{\gp}_{\frak D}/\cO^\times_{\frak D})_x=X(S_x)$ such that 
$X(S_x)^+ \subset V$ and $X(S_x)^+ \cap V^\times =\{1\}$, and $Z$ is a $T(V)$-orbit in the $T(x)$-torsor $Z(x)$. 
Here 
$$T(V) := \Hom(X(S_x)/V^\times, \R^{\mult}_{>0})\subset T(x)=\Hom(X(S_x), \R^{\mult}_{>0}).$$

\end{sbpara}

\begin{sbpara}\label{orbitlog}  Let the notation be as in \ref{orbittr}. For a point $z=(x, V, Z)\in \frak D$, the stalk $(M_{\frak D}/\cO^\times_{\frak D})_z$ is described as follows. In the situations (a)--(c), $(M_{\frak D}/\cO^\times_{\frak D})_z= V/V^\times$. In the situation (d), $(M_{\frak D}/\cO^\times_{\frak D})_z= V'/(V')^\times$ where $V'=V\cap (X(S_x)^+)^{\gp}$.

\end{sbpara}

\begin{sbpara}\label{another} In \cite{KU0} 2.6 and \cite{KU} 5.1.6, which treated the pure case, we defined the set $D_{\BS,\val}$ in a different style. 
Following the style in \cite{KU0} 2.6 and \cite{KU} 5.1.6, we can define
$D_{\BS,\val}$ also as the set of all triples $(T, V, Z)$ where 
$T$ is an $\R$-split torus in $G_\R(\gr^W)$, $V$ is a valuative submonoid of the character group $X(T)$ of $T$, and $Z\subset D$, satisfying the following conditions (i)--(iv).

\smallskip

(i) Let $T_{>0}$ be the connected component of $T(\R)$ containing the unit element.
Then $Z$ is either a $T_{>0}$-orbit for the lifted action \ref{liftac} or an $\R_{>0}\times T$-orbit in $D_{\nspl}$ for the lifted action. Here $t\in \R_{>0}$ acts on $\gr^W$ by the multiplication by $t^w$ on $\gr^W_w$. 

\smallskip

(ii) Let $\br\in Z$, let $\bar \br:=\br(\gr^W)\in D(\gr^W)$, let $K_{\bar \br}$ be the maximal compact subgroup of $G_\R(\gr^W)$ associated to $\bar \br$, and let $\theta_{K_{\bar \br}}: G_\R(\gr^W)\to G_\R(\gr^W)$ be the Cartan involution associated to $K_{\bar\br}$. Then 
$\theta_{K_{\bar \br}}(t)=t^{-1}$ for any $t\in T$. 

\smallskip

(iii) $V^\times =\{1\}$.

\smallskip

(iv) Consider the direct sum  decomposition $\gr^W=\bigoplus_{\chi\in X(T)} (\gr^W)_{\chi}$ by the action of $T$. Then for any $\chi\in X(T)$, the subspace $\bigoplus_{\chi'\in V^{-1}\chi} (\gr^W)_{\chi'}$ is $\Q$-rational. 

\smallskip
The relation with the presentation \ref{orbittr} of $D_{\BS,\val}$ is as follows. $(P, V, Z)\in D_{\BS,\val}$ in the presentation in \ref{orbittr} corresponds to $(T, V', Z)$ in the above presentation, where $T\subset S_P$ is the annihilator of $V^{\times}$ in $S_P$ and $V'=V/V^\times \subset X(T)$.
The group $T(V)$ in \ref{orbittr} coincides with $T_{>0}$ in the above (i).

Conversely, for a triple $(T, V, Z)$ here, the corresponding triple in the presentation of $D_{\BS,\val}$ in \ref{orbittr} is $(P, V', Z)$ where $P$ is the $\Q$-parabolic subgroup of $G_\R(\gr^W)$ defined as the
connected component (as an algebraic group) of the algebraic subgroup of $G_\R(\gr^W)$ consisting of all elements which preserve the subspaces $\bigoplus_{\chi'\in V^{-1}\chi} (\gr^W)_{\chi'}$ of $\gr^W$, and $V'$ is the inverse image of $V$ under the homomorphism $X(S_P)\to X(T)$ induced by the canonical homomorphism $T\to S_P$.

\end{sbpara}

\begin{sbpara}\label{2slval} We describe the map   $D^{\star}_{\SL(2),\val}\to D^{II}_{\SL(2),\val}$.

This map is described as  $x=(p,V, Z)\mapsto (p, V', Z')$ using \ref{orbittr}  as follows.

\smallskip

(0) On $D$, this map is the identity map.

\smallskip

(1) For an $A$-orbit which does not belong to $D$, $V'=V$ and $Z'=Z_{\spl}$.

\smallskip

(2) Assume that $(p,V,Z)$ is a $B$-orbit, let $n$ be the rank of $p$, and identify $X(S_x)=\Z\times X(S_p)$ with $\Z\times \Z^n$. 
  Let $e= (1, -1, \dots, -1)\in \Z\times \Z^n$. 

\smallskip

(2.1) Assume $-e\notin V$ (hence $e\in V$). Then 
  $V'= \{a=(a_0, a_1,\dots, a_n)\in \Z^{n+1}\;|\; a - a_0e\in V\}$, and $Z'=Z$.

\smallskip

(2.2) Assume $e, -e\in V$. Then $V'=\{a\in \Z^n\;|\; (0, a)\in V\}$, and $Z'=Z$. 

\smallskip

(2.3) Assume $e\notin V$ (hence $-e\in V$). Then 
$V'=\{a\in \Z^n\;|\; (0, a)\in V\}$, and $Z'=Z_{\spl}$. 

\end{sbpara}

\subsection{The morphism $\eta^{\star}: D^{\star}_{\SL(2),\val}\to D_{\BS,\val}$}\label{ss:SB}

 \begin{sbpara}\label{eta1}
 
The map
$
\eta^{\star}: D^{\star}_{\SL(2),\val} \to D_{\BS, \val}
$
is described as follows. 
This description  is similar to the pure case in \cite{KU0} Theorem 3.11,   \cite{KU}  Theorem 5.2.11.

The map $\eta^{\star}$ sends $(p, V, Z)\in D_{\SL(2),\val}^{\star}$ in the presentation of $D_{\SL(2),\val}^{\star}$ in \ref{orbittr} to $(T, V', Z)\in D_{\BS,\val}$ in the presentation of $D_{\BS,\val}$ in \ref{another}, where $T$ and $V'$ are as follows. Let $T'\subset S_p$ be the annihilator of $V^\times\subset X(S_p)$. Then $T$ is the image of $T'\to G_\R(\gr^W)$ under $\tau^{\star}_p$. $V'$ is the inverse image of $V/V^\times\subset X(T')$ under the homomorphism $X(T)\to X(T')$ induced by the canonical homomorphism $T'\to T$. 
 
    \end{sbpara}

    \begin{sbpara}\label{eta2}

we have also the following description of $\eta^{\star}$ 
by regarding $D^{\star}_{\SL(2),\val}$ as $D^{\star,\BS}_{\SL(2),\val}$.

Let the notation be as in Section 2.6. By \ref{gest9}, an element of $D^{\star,\BS}_{\SL(2),\val}$ is written as $(p, P, V, Z)$ where $p\in D_{\SL(2)}(\gr^W)$, $P\in \cP(p)$, $V$ is a valuative submonoid of $X(S_{p,P})$ such that $X(S_{p,P})^+\subset V$ and $X(S_{p.P})^+\cap V^\times =\{1\}$, and $Z$ is either a $\tau^{\star}(\Hom(X(S_{p,P})/V^\times, \R_{>0}))$-orbit in $D$ or a $\tilde \tau^{\star}(\R_{>0}\times \Hom(X(S_{p,P})/V^\times, \R_{>0}))$-orbit in $D_{\nspl}$ for the lifted action, such that the image of $Z$ in $D(\gr^W)$ is contained in $Z(p)$. 

The map $\eta^{\star}$ sends $(p, P, V, Z)\in D^{\star,\BS}_{\SL(2),\val}$ to $(P, V', Z)\in D_{\BS,\val}$ in the presentation of $D_{\BS,\val}$ as a set in \ref{orbittr}, where $V'\subset X(S_P)$ is the inverse image of $V$ under the homomorphism $X(S_P)\to X(S_{p,P})$ induced by the canonical homomorphism $S_{p.P}\to S_P$. 
 
 \end{sbpara}

 \begin{sblem}\label{eta3} The morphism $\eta^{\star}: D^{\star}_{\SL(2),\val}\to D_{\BS,\val}$ has the Kummer property of log structure in the sense of \ref{timesval} (K).
 
 \end{sblem}
 
 \begin{pf} Let $x=(p,P,V,Z)\in D^{\star,\BS}_{\SL(2),\val}= D^{\star}_{\SL(2),\val}$ (\ref{eta2}) and let $y$ be the image of $x$ in $D_{\BS,\val}$. 
 By \ref{orbitlog}, in the case of $A$-orbit (resp.\ $B$-orbit), $V$ is a valuative submonoid of $X(S_{p,P})$ (resp.\ $\Z\times X(S_{p,P})$) and the stalk of $M/\cO^\times$ of $D^{\star}_{\SL(2),\val}$ at $x$ is identified with $V/V^\times$. On the other hand, the stalk of $M/\cO^\times$ of $D_{\BS,\val}$ at $y$ is identified with $V'/(V')^\times$ where in the case of $A$-orbit (resp.\ $B$-orbit), $V'$ is the inverse image of $V$ in $(X(S_P)^+)^{\gp}$ (resp.\ $\Z\times (X(S_P)^+)^{\gp}$) for the canonical map $(X(S_P)^+)^{\gp}\subset X(S_P) \to X(S_{p,P})$. Note that $(X(S_P)^+)^{\gp}$ is of finite index in $X(S_P)$. Furthermore, since the kernel of $S_{p,P}\to S_P$ is finite, the cokernel of $X(S_P)\to X(S_{p,P})$ is finite. Hence the map $V'/(V')^\times \to V/V^\times$ is injective, and for any element $a$ of $V/V^\times$, there is $m\geq 1$ such that $a^m$ belongs to the image of $V'/(V')^\times$.
  \end{pf}
 
 \begin{sbthm}
\label{SL2BS}

 The map $\eta^*:D^{\star}_{\SL(2),\val}\to D_{\BS,\val}$ in $\cC_\R(\val)^+$ has the following properties. 
 
 \smallskip
  
(1)  The map $\eta^{\star}:D^{\star}_{\SL(2),\val}\to D_{\BS,\val}$ is injective.
 
 \smallskip

(2)  Let $Q\in \prod_w \cW(\gr^W_w)$ and define the open set  $D^{\star}_{\SL(2),\val}(Q)$ of $D^{\star}_{\SL(2),\val}$ as the inverse image of the open set $D_{\SL(2)}(\gr^W)(Q)$ of $D_{\SL(2)}(\gr^W)$. Then the topology of $D^{\star}_{\SL(2),\val}(Q)$ coincides with the restriction of the topology of $D_{\BS,\val}$ through $\eta^{\star}$.

\smallskip
(3) The diagram
$$\begin{matrix} D^{\star}_{\SL(2),\val}  &\overset{\eta^{\star}}\to& D_{\BS,\val}\\
\downarrow &&\downarrow\\
\prod_w D_{\SL(2)}(\gr^W_w)_{\val}  &\overset{\eta}\to & \prod_w D_{\BS}(\gr^W_w)_{\val}
\end{matrix}$$
is cartesian in the category of topological spaces.

  \end{sbthm}

  \begin{pf} We prove (3) first. For each $x=(x_w)_w\in \prod_w D_{\SL(2)}(\gr^W_w)_{\val}$ and the image $y=(y_w)_w$ of $x$ in $\prod_w D_{\BS}(\gr^W_w)_{\val}$, and for each $w\in \Z$, there is an open neighborhood $U_w$ of $x_w$ and an open neighborhood $V_w$ of $y_w$ having the following properties (i) and (ii).

  \smallskip
  
  (i) The image of $U_w$ in $D_{\BS}(\gr^W_w)_{\val}$ is contained in $V_w$.
  
  \smallskip
  
  (ii) Let $U$ be the inverse image of $\prod_w U_w$ in $D^{\star}_{\SL(2),\val}$ and let $V$ be the inverse image of $\prod_w V_w$ in $D_{\BS,\val}$. Take any $F\in D(\gr^W)$ and let $\bar L= \bar \cL(F)$. Then we have a commutative diagram
  $$\begin{matrix}U &\cong &\prod_{\val,w} U_w\times_{\val} \spl(W) \times_{\val} \bar L\\
  \downarrow &&\downarrow\\
  V &\cong  &\prod_{\val,w} V_w \times_{\val} \spl(W) \times_{\val} \bar L
  \end{matrix}$$
  where $\prod_{\val,w}$ is the product in $\cC_\R(\val)^+$, the upper row is an isomorphism over $\prod_{\val,w} U_w$, and the lower row is an isomorphism over $\prod_{\val,w} V$. 
  
  \smallskip
  
  By \ref{timesval} and \ref{eta3}, the following diagram is cartesian in the category of topological spaces.
  $$\begin{matrix} \prod_{\val,w} U_w \times_{\val} \spl(W) \times_{\val} \bar L &\to & \prod_{\val,w} V_w\times_{\val} \spl(W) \times_{\val} \bar L\\
  \downarrow &&\downarrow\\
  \prod_w U_w \times \spl(W) \times \bar L &\to & \prod_w V_w \times \spl(W) \times \bar L.
  \end{matrix}$$
 (3) of Theorem \ref{SL2BS} follows from these two cartesian diagrams.
 
 Next we prove (1). The injectivity was proved in \cite{KU0} Theorem 3.11 in the pure case. Hence the map $\prod_w D_{\SL(2)}(\gr^W_w)_{\val} \to \prod_w D_{\BS}(\gr^W_w)_{\val}$ is injective. By (3), this proves the injectivity of $D_{\SL(2),\val}^{\star}\to D_{\BS,\val}$.

   We prove (2). 
  Assume first we are in the pure situation of weight $w$. Let ${\cal T}_1$ be the topology of $D_{\SL(2)}$ defined in Part II, and let ${\cal T}_{1,\val}$ be the topology of $D_{\SL(2),\val}$ defined in this Part IV. Let ${\cal T}_{2,\val}$ be the topology of $D_{\SL(2),\val}$ which is the weakest topology satisfying the following two conditions (i) and (ii). 
   
   \smallskip
   
   (i) For any open set $U$ of $D_{\BS,\val}$, the pull-back of $U$ in $D_{\SL(2),\val}$ is open.
   
   \smallskip
   
   (ii) For any $Q\in \prod_w \cW(\gr^W_w)$, $D_{\SL(2),\val}(Q)$ is open.
   
   \smallskip
   
   Let ${\cal T}_2$ be the topology of $D_{\SL(2)}$ as a quotient space of $D_{\SL(2),\val}$ which is endowed with the topology ${\cal T}_{2,\val}$. 
  Recall that in \cite{KU0} and \cite{KU} which treated the pure case, the topologies of $D_{\SL(2)}$ and $D_{\SL(2),\val}$ 
  were defined as ${\cal T}_2$ and ${\cal T}_{2,\val}$, respectively (not as in the present series of papers). 
  The study of ${\cal T}_2$ in \cite{KU} Section 10 and the study of ${\cal T}_1$ in Part II, Section 3.4  show that ${\cal T}_1={\cal T}_2$. Since the map $\eta^{\star}$ from $D_{\SL(2),\val}$ with ${\cal T}_{1,\val}$ to $D_{\BS,\val}$ is continuous as we have seen in Section 2.6, 
  we have that ${\cal T}_{1,\val} \geq {\cal T}_{2,\val}$. Since the map $D_{\SL(2),\val} \to D_{\SL(2)}$ is proper for ${\cal T}_{1,\val}$ (\ref{vproper}) and also for ${\cal T}_{2,\val}$ (\cite{KU} Theorem 3.14), we have ${\cal T}_{1,\val}={\cal T}_{2,\val}$.

  Thus we have proved (2) in the pure case. By (3), we have a cartesian diagram of topological spaces 
  $$\begin{matrix} D^{\star}_{\SL(2),\val} & \to & D_{\BS,\val}\\
  \downarrow && \downarrow\\
 (\prod_w D_{SL(2)}(\gr^W_w)_{\val}) \times \spl(W) &\to& (\prod_w D_{\BS}(\gr^W_w)_{\val})\times \spl(W).
    \end{matrix}$$
    The vertical arrows are proper by \ref{Lbund} and Part I, Cor. 8.5. Hence (3) is reduced to  the pure case. 
\end{pf}
 
 \begin{sbpara}
 As in (2) of Theorem \ref{SL2BS}, the topology of $D^{\star}_{\SL(2),\val}(Q)$ coincides with the induced topology from $D_{\BS,\val}$. We show an example in which the topology of $D^{\star}_{\SL(2),\val}$ is not the induced one from $D_{\BS,\val}$. This example is pure of weight $3$. So $D^{\star}_{\SL(2),\val}$ is written as $D_{\SL(2),\val}$ below. 
 
 Let $H_{0,\Z}=H'_{0,\Z} \otimes  \text{Sym}^2(H'_{0,\Z})$ where $H'_{0,\Z}$ is a free $\Z$-module of rank $2$ with basis $e_1,e_2$. Hence $H_{0,\Z}$ is of rank $6$. The intersection form $\langle -,-\rangle$ on $H_{0,\Z}$ is $b\otimes \text{Sym}^2(b)$ where $b$ is the
  anti-symmetric bilinear form on $H'_{0,\Z}$ characterized by $\langle e_1, e_2\rangle =-1$. We have the following 
  $\SL(2)$-orbit $(\rho, \varphi)$ in two variables: $$\rho(g_1,g_2)= g_1\otimes \text{Sym}^2(g_2), 
  \quad \varphi(z_1,z_2) =F(z_1)\otimes  \Sym^2F(z_2)$$ $(g_1,g_2\in \SL(2)$, $z_1,z_2\in \C)$ where $F(z)$ is the decreasing 
  filtration on $H'_{0,\C}$ defined by $$F(z)^2=0 \subset F(z)^1= \C\cdot (ze_1+e_2)\subset F(z)^0= H'_{0, \C}.$$
  The associated homomorphism $\tau^{\star}: \bG_m^2\to G_\R$ is as follows. 
  $\tau^{\star}(t_1, t_2)$ acts on $e_2\otimes e_2^2$ by $t_1t_2^3$, on $e_2\otimes e_1e_2$ by $t_1t_2$, 
  on $e_2\otimes e_1^2$ by $t_1t_2^{-1}$, on $e_1\otimes e_2^2$ by $t_1^{-1}t_2$, 
  on $e_1\otimes e_1e_2$
   by $t^{-1}_1t_2^{-1}$, and on $e_1\otimes  e_1^2$ by $t_1^{-1}t_2^{-3}$.
  The associated weight filtrations $W^{(1)}$ and $W^{(2)}$ are as follows.
  $$W^{(1)}_1=0 \subset W^{(1)}_2= e_1 \otimes \text{Sym}^2H'_{0,\R} =W^{(1)}_3\subset W^{(1)}_4=H_{0,\R}.$$
  $$W^{(2)}_{-1}=0\subset W^{(2)}_0=\R e_1 \otimes e_1^2=W^{(2)}_1\subset W^{(2)}_2 = W^{(2)}_1+ \R  e_1 \otimes e_1e_2 + \R  e_2 \otimes e_1^2= W^{(2)}_3$$ $$\subset W^{(2)}_4= W^{(2)}_3+ 
  \R e_1\otimes e_2^2 + \R e_2\otimes e_1e_2= W^{(2)}_5\subset W^{(2)}_6= H_{0,\R}.$$
 Let
 $$\Phi= \{W^{(1)}, W^{(2)}\}.$$
 
 We show that $D_{\SL(2),\val}(\Phi)$ is not open for the topology induced from the topology of $D_{\BS,\val}$. 
 
 Let $V$ be the valuative submonoid of $X(\bG_m^2)$ which is, under the identification $X(\bG_m^2)= \Z^2$, identified with the set of all $(a, b)\in \Z^2$ satisfying either $(a>0)$ or ($a=0$ and $b\geq 0$). 
 Consider the point $x:=(p, V, Z)\in D_{\SL(2),\val}$ where 
 $p$ is the class of this $\SL(2)$-orbit and $Z$ is the torus orbit $\{F(iy_1) \otimes \text{Sym}^2 F(iy_2) \;|\; y_1,y_2\in \R_{>0}\}$ of $p$. The map $D_{\SL(2),\val}\to D_{\BS,\val}$ sends $x$ to $y:=(T, V, Z)$ 
 in the presentation \ref{another} of $D_{\BS,\val}$ where $T$ is the image of  $\tau^{\star}=\tau_p^{\star}: \bG_m^2\to G_\R$ and we regard $V$ as 
 a submonoid of $X(T)$ via the canonical isomorphism $\bG_m^2\cong T$ given by $\tau_p^{\star}$. In the presentation of $D_{\BS,\val}$ in \ref{orbittr}, 
 this point  $y$ coincides with $(P, V', Z)$ where $P$ and $V'$ are as follows. $P$ is the $\Q$-parabolic subgroup of $G_\R$ consisting of all elements which preserve the following subspaces $W'_w$ $(w\in \Z)$. 
 $$W'_1=\R  e_1 \otimes e_1^2, \quad W'_2=W'_1+ \R  e_1 \otimes e_1e_2, \quad W'_3= W'_2+\R  e_1 \otimes e_2^2,$$
 $$W'_4= W'_3+\R e_2 \otimes e_1^2, \quad W'_5= W'_4+ \R  e_2\otimes e_1e_2.$$
 We have $S_P\cong \bG_m^3$. The inclusion map $T\to P$ induces a canonical homomorphism $T\to S_P$.  $V' \subset X(S_P)$ is the inverse image of $V$ under the canonical homomorphism $X(S_P)\to X(T) \cong X(\bG_m^2)$.

 Let $f$ be the element of $\Lie(P_u)$ which sends $e_2\otimes e_1^2$ to $e_1\otimes e_2^2$ and kills $e_1\otimes \text{Sym}^2H'_{0,\R}$ and $e_2\otimes e_1e_2$ and $e_2\otimes e_2^2$.  
  For $c\in \R$, we have the SL(2)-orbit in two variables $(\rho^{(c)}, \varphi^{(c)})$  defined by
 $$\rho^{(c)}(g_1, g_2)= \exp(cf) \rho(g_1, g_2) \exp(-cf), \quad \varphi^{(c)}(z_1,z_2)= \exp(cf)\varphi(z_1, z_2).$$
 The associated weight filtrations of $(\rho^{(c)}, \varphi^{(c)})$ are $\exp(cf)W^{(1)}, \exp(cf)W^{(2)}$. Since $f$ respects $W^{(1)}$ but not $W^{(2)}$, $\{\exp(cf)W^{(1)}, \exp(cf)W^{(2)}\}=\{W^{(1)}, \exp(cf)W^{(2)}\}$ is  
 not contained in $\Phi$ if $c\neq 0$. If $c\in \Q$, the filtration $\exp(cf)W^{(2)}$ is rational, and hence $(\rho^{(c)}, \varphi^{(c)})$ determines an element $p^{(c)}$ 
 of $D_{\SL(2)}$. Let $x^{(c)}:=(p^{(c)}, V, Z^{(c)})\in D_{\SL(2), \val}$ where $V$ is the same as above and $Z^{(c)}$ is the torus orbit of $p^{(c)}$.

 Now, when $c\in \Q\smallsetminus \{0\}$ converges to $0$ in $\R$, the image $y^{(c)}$ of $x^{(c)}$ in $D_{\BS,\val}$ converges to $y$. This is because $P$ 
 acts on $D_{\BS,\val}(P)$ continuously in the natural way, $y\in D_{\BS,\val}(P)$, and $y^{(c)}=\exp(cf)y$ for this action of $P$. Since $y^{(c)}\notin D_{\SL(2),\val}(\Phi)$, $D_{\SL(2),\val}(\Phi)$ is 
 not open for the topology induced by the topology of $D_{\BS,\val}$. 
 
 The set $\{x^{(c)}\;|\; c\in \Q\}$ is discrete in $D_{\SL(2),\val}$, though the image $\{y^{(c)}\;|\; c\in \Q\}$ in $D_{\BS,\val}$ has the topology of the subspace $\Q$ of $\R$ via the correspondence $y^{(c)}\leftrightarrow c$.  
 
 Thus the topology of $D_{\SL(2),\val}$ is not the induced topology from $D_{\BS,\val}$.

 \end{sbpara}

\subsection{The map $\eta: D_{\SL(2),\val}\to D_{\BS,\val}$}\label{4.6}
\label{ss:eta}
\begin{sbpara}\label{defeta} We define a canonical map 
$$
\eta: D_{\SL(2),\val} \to D_{\BS, \val}
$$
 following the method in the pure case \cite{KU}. But we will see that this map need not be continuous.

$\eta$ is the unique map such that for any $x\in D_{\SL(2)}$ and any $\tilde x\in D_{\SL(2),\val}$ lying over $x$, the restriction of $\eta$ to the subset $\bar Z(x)_{\val}$ of $D_{\SL(2),\val}$ (note $\tilde x\in \bar Z(x)_{\val}$)  is the unique morphism in $\cC_\R(\val)^+$ whose restriction to $Z(x)$ is the inclusion morphism $Z(x) \overset{\subset}\to D \subset D_{\BS,\val}$. 

The map $\eta$ coincides with the composition of the two maps $D_{\SL(2),\val}\to D^{\star}_{\SL(2),\val}\overset{\eta_{\val}^{\star}}\longrightarrow  D_{\BS,\val}$ where the first arrow is the following map $\lambda_{\val}$. The restriction of $\lambda_{\val}$ to $D_{\SL(2),\nspl,\val}$ is the morphism on the associated valuative spaces induced from the morphism $\lambda: D^{II}_{\SL(2),\nspl} \to D^{\star}_{\SL(2)}$
in  \ref{lam}. The restriction of $\lambda_{\val}$ to $D_{\SL(2),\spl,\val}$ is the morphism on the associated valuative spaces induced from the isomorphism $\eta:D^{II}_{\SL(2),\spl} \overset{\cong}\to D^{\star}_{\SL(2),\spl}$ 
 in \ref{lam}. 
 
 The composition $D_{\SL(2),\val} \overset{\lam_{\val}}\longrightarrow D^{\star}_{\SL(2),\val}\to D_{\SL(2),\val}$ is the identity map. By Theorem \ref{SL2BS} (1), the map $\eta: D_{\SL(2),\val}\to D_{\BS,\val}$ is injective. 

\end{sbpara}

\begin{sbprop}
$(1)$ The restriction of $\eta$ to the open set 
$D^{II}_{\SL(2),\nspl, \val}\cup D$ of $D^{II}_{\SL(2),\val}$ is a morphism in $\cC_\R(\val)$. 

$(2)$ For any $\Phi\in \overline{\cW}$, 
the topology of $D^{II}_{\SL(2),\nspl,\val}(\Phi)\cup D$ coincides with the topology induced from the topology of $D_{\BS,\val}$.
\end{sbprop}

\begin{pf} This restriction of $\eta$ to $D^{II}_{\SL(2),\nspl,\val}\cup D$ is the composition 
$D^{II}_{\SL(2), \nspl,\val}\cup D\to D^{\star}_{\SL(2),\val} \overset{\eta^{\star}_{\val}}\longrightarrow D_{\BS,\val}$
where the first arrow is the open immersion 
 induced from isomorphism $D^{II}_{\SL(2)}\cup D\cong D_{\SL(2)}^{\star,+}(\sig_2)$ in \ref{0thm} (2). This proves (1). By this, (2) follows from Theorem \ref{SL2BS} (2).
\end{pf}

\begin{sbprop}\label{twoSL23} The equivalent conditions {\rm (i)}--{\rm (vii)} of $\ref{twoSL2}$ 
 are equivalent to each of the following conditions.

\smallskip

{{\rm (viii)}} The identity map of $D$ extends to an isomorphism $D^{II}_{\SL(2),\val}\cong D^{\star}_{\SL(2),\val}$ in $\cC_\R(\val)$.

\smallskip

{{\rm (ix)}} The map $\lambda_{\val}: D^I_{\SL(2),\val} \to D^{\star}_{\SL(2),\val}$ (\ref{defeta}) is continuous.

\smallskip

{{\rm (x)}} The map $\eta: D^{II}_{\SL(2),\val}\to D_{\BS,\val}$ is continuous.

\smallskip

{{\rm (xi)}} The map $\eta: D^I_{\SL(2),\val}\to D_{\BS,\val}$ is continuous.

\smallskip

{{\rm (xii)}} The map $D^{\star,\mild}_{\SL(2),\val}\to D_{\SL(2),\val}$ is injective.

\end{sbprop}

\begin{pf} (ii) $\Rightarrow $ (viii). Take the associated valuative spaces.

The implications 
(viii) $\Rightarrow$ (ix) and 
(viii)  $\Rightarrow$ (xii) are clear. 

The implications (viii) $\Rightarrow$ (x) and 
(x) $\Rightarrow$ (xi) are easily seen.

(xi) $\Rightarrow$ (ix).   Use the fact that the topology of $D^{\star}_{\SL(2),\val}(\Phi)$ is the restriction of the topology of $D_{\BS,\val}$ (\ref{SL2BS} (2)).

(ix) $\Rightarrow $ (iv). This is because $D^I_{\SL(2),\val}\to D^I_{\SL(2)}$ is proper surjective (\ref{vproper}).

(xii)   $\Rightarrow $ (i). The proof of (vi) $\Rightarrow$ (i) of Proposition \ref{twoSL2} actually proves this. In that proof, assuming  (i) does not hold, 
we used $x= (p, Z)\in D^{\star, \mild}_{\SL(2)}$ with $p$ of rank $1$ 
such that $x\neq x_{\spl}$. Since $p$ is of rank one, these $x$ and $x_{\spl}$ 
are regarded canonically as elements of $D_{\SL(2),\val}^{\star,\mild}$ whose images in $D_{\SL(2),\val}$ coincide. 
\end{pf}

\section{New spaces $D^{\sharp}_{\Sig,[:]}$ and $D^{\sharp}_{\Sig,[\val]}$ of nilpotent orbits}\label{s:newval}

In this Section \ref{s:newval}, we define and consider the new spaces $D^{\sharp}_{\Sig, [:]}$ and  $D^{\sharp}_{\Sig,[\val]}$ of nilpotent orbits (nilpotent $i$-orbits, to be precise, in our terminology). 

In Sections \ref{ra1}--\ref{ra5}, for a topological space $S$ endowed with an fs log structure on the sheaf of  all $\R$-valued continuous functions, we define topological spaces $S_{[:]}$ (the space of ratios, Section \ref{ra2}) and $S_{[\val]}$ (Section \ref{ra5}), and we define proper surjective continuous maps $S_{[:]}\to S$, $S_{\val}\to S_{[:]}$, and $S_{[\val]}\to S_{[:]}$, where $S_{\val}$ is as in Section 3.1. As will be explained  in \ref{ra6}, in the case $S= D^{\sharp}_{\Sig}$,  we obtain the new spaces of nilpotent $i$-orbits $D^{\sharp}_{\Sig, [:]}$ as $S_{[:]}$ and $D^{\sharp}_{\Sig, [\val]}$ as $S_{[\val]}$, 
and $S_{\val}$ coincides with $D^{\sharp}_{\Sig, \val}$ which we have already defined in Part III. 
We construct CKS maps $D^{\sharp}_{\Sig, [:]}\to D_{\SL(2)}^I$ and $D^{\sharp}_{\Sig, [\val]}\to D_{\SL(2),\val}^I$ in \ref{ss:cks1}.
We have already constructed the CKS map $D^{\sharp}_{\Sig, \val}\to D_{\SL(2)}^I$ in Part III.

\subsection{The space of ratios in toric geometry}\label{ra1}
\begin{sbpara} The space of ratios which we consider appears in the following way. 

Consider $S=\Spec(k[T_1, T_2])$ with $k$ a field. Regarding $S$ as the toric variety associated  to the cone $\R^2_{\geq 0}\subset \R^2$, consider the toric varieties over $k$ associated to rational finite subdivisions of the cone $\R^2_{\geq 0}$ (\ref{rvtoric}), and let $X$ be the projective limit of these toric varieties regarded as topological spaces with Zariski topology.
  It is the projective limit obtained by blowing-up the origin $s=(0,0)\in S$ first and then continuing blowing-up  the intersections of irreducible components of the inverse image of $\Spec(k[T_1, T_2]/(T_1T_2))\subset S$ on the blowing-up. 

Let $X_0\subset X$ be the inverse image of $s$, and endow $X_0$ with the topology as a subspace of $X$. Then we have the following  continuous surjective map from $X_0$ to the interval $ [0, \infty] \supset \R_{>0}$ despite that the 
Zariski topology and the topology of real numbers are very much different in nature. If $x\in X_0$, the image of $x$ in $[0, \infty]$ is defined as 
$$  \sup \{a/b\;|\;(a, b)\in \N^2\smallsetminus \{(0,0)\}, T_1^b/T_2^a\in \cO_{X,x}\} $$ $$= \inf \{a/b \; |\; (a, b)\in \N^2\smallsetminus \{(0,0)\},T_2^a/T_1^b\in \cO_{X,x}\}.$$  
 Here $\N=\Z_{\ge0}$ and $\cO_X$ is the inductive limit of the inverse images on $X$ of the structural sheaves of the blowing-ups. The image of $x$ in $[0, \infty]$ is, roughly speaking,  something like the ratio  $\log(T_1)/\log(T_2)$ at $x$. 

In the definition below, this $[0, \infty]$ is {\it the space $R(\N^2)$ of ratios} of the fs monoid $\N^2= (M_S/\cO_S^\times)_s$ which is generated by the classes of $T_1$ and $T_2$. The above relation with the projective limit of blowing-ups is generalized  in \ref{zariski}. 
\end{sbpara}

\begin{sbpara}
In this Section \ref{ra1}, the notation $\cS$ is used for an fs monoid. We denote the semigroup law of $\cS$ multiplicatively unless we assume and state that $\cS=\N^n$. So the neutral element of $\cS$ is denoted by $1$.
\end{sbpara}

\begin{sbpara}\label{ratY}  
  For a sharp fs monoid $\cS$, let $R(\cS)$ be the set of all maps $r:(\cS\times \cS)\smallsetminus \{(1,1)\}\to [0,\infty]$ satisfying the following conditions (i)--(iii).  
  
(i) $r(g,f)=r(f,g)^{-1}$.

(ii) $r(f, g)r(g,h)=r(f,h)$ if $\{r(f,g), r(g,h)\}\neq \{0,\infty\}$. 

(iii) $r(fg, h)=r(f, h)+r(g,h)$.

We endow $R(\cS)$ with the topology of simple convergence. It is a closed subset of the product of copies of the compact set $[0, \infty]$ and hence is compact.

\end{sbpara}

\begin{sbrem} From the condition (i), we have $r(f,f)=1$. (Conversely, $r(f,f)=1$ and (ii) imply (i).) From this and from $r(1,f)+r(f,f)=r(f,f)$ which comes from (iii), we get
$$r(1,f)=0, \quad r(f,1)=\infty\quad \text{for any $f\in \cS\smallsetminus \{1\}$}. $$

\end{sbrem}
\begin{sbpara} For example, we have $R(\N^2)\cong [0,\infty]$, where $r\in R(\N^2)$ corresponds to $r(q_1,q_2)\in [0, \infty]$ with $q_j$ the standard $j$-th basis of $\N^2$. 

A description of $R(\N^n)$ for general $n$ is given in \ref{RNn}.
\end{sbpara} 

\begin{sbpara}\label{RandR'} We have a canonical bijection between $R(\cS)$ and the set $R'(\cS)$ of 
all equivalence classes of $((\cS^{(j)})_{0\leq j\leq n}, (N_j)_{1\leq j\leq n})$, where  $n\geq 0$, $\cS^{(j)}$ is a  face of $\cS$ such that
$$\cS=\cS^{(0)} \supsetneq \cS^{(1)}\supsetneq \dots \supsetneq \cS^{(n)}=\{1\},$$
and $N_j$  is a homomorphism $\cS^{(j-1)}\to \R^{\add}$ such that $N_j(\cS^{(j)})=0$ and such that $N_j(\cS^{(j-1)}\smallsetminus \cS^{(j)})\subset\R_{>0}$. 
The equivalence relation is given by multiplying each $N_j$ by an element of $\R_{>0}$ (which may depend on $j$). 

We define a map $R(\cS)\to R'(\cS)$ as follows. Let  $r\in R(\cS)$. We give the corresponding element of $R'(\cS)$.

For $f\in \cS\smallsetminus \{1\}$, let $\cS(r, f)=\{g\in \cS\;|\;r(g,f)\neq \infty\}$. Then the conditions (i)---(iii) on $r$ in \ref{ratY} show that $\cS(r,f)$ is a face of $\cS$. For $f, g\in \cS$, we have 
$\cS(r,f)\subset \cS(r,g)$ if and only if $r(f,g)\neq \infty$, and we have $\cS(r,f)\supset \cS(r,g)$ if and only if $r(f,g)\neq 0$. 
  Hence the faces of $\cS$ of the form $\cS(r,f)$ $(f\in \cS\smallsetminus\{1\})$ together with the face $\{1\}$ form a  totally ordered set for the inclusion relation. 
Let $\cS= \cS^{(0)}\supsetneq \cS^{(1)}\supsetneq \dots \supsetneq \cS^{(n)}=\{1\}$ be all the members of this set. 
 Take $q_j\in \cS^{(j-1)}\smallsetminus \cS^{(j)}$ $(1\leq j\leq n)$. We have a homomorphism $N_j: \cS^{(j-1)}\to \R$ defined by $N_j(f)=r(f, q_j)$. This $N_j$ kills $\cS^{(j)}$ and $N_j(\cS^{(j-1)}\smallsetminus \cS^{(j)})\subset \R_{>0}$. If we replace $q_j$ by another element $q'_j$, $N_j$ is multiplied by $r(q_j, q'_j)\in \R_{>0}$. 
  Thus we have the map  $R(\cS)\to R'(\cS)\;;\; r\mapsto  \text{class}((\cS^{(j)})_j,(N_j)_j)$. 
 
  Next we define the inverse map $R'(\cS)\to R(\cS)$. Let $\text{class}((\cS^{(j)})_{0\leq j\leq n}, (N_j)_{1\leq j\leq n})\in R'(\cS)$. 
 Let $(f,g) \in (\cS\times S) \smallsetminus \{(1,1)\}$. We define $r(f,g)$ as follows.
Let $j$ be the largest integer $\geq 0$ such that  $f$  belongs to $\cS^{(j)}$ and let $k$ be that of  $g$.

(1) If $j=k< n$, $r(f,g)= N_{j+1}(f)/N_{j+1}(g)$.

(2) If $j>k$, $r(f, g)= \infty$.

(3) If $j<k$, $r(f,g)=0$. 

 This gives the map $R'(\cS)\to R(\cS)$.

 It can be seen easily that the maps $R(\cS)\to R'(\cS)$ and $R'(\cS)\to R(\cS)$ are the inverses of each other.

\end{sbpara}

\begin{sbpara}\label{V(S)2}
As in \ref{V(S)}, for a sharp fs monoid $\cS$, let $V(\cS)$ 
 be the set of all valuative submonoids $V$ of $\cS^{\gp}$ such that $V\supset \cS$ and $V^\times \cap \cS=\{1\}$. We endow $V(\cS)$ with the following topology. For a finite set $I$ of $\cS^{\gp}$, let $U(I)=\{V\in V(\cS)\;|\;I\subset V\}$. Then these $U(I)$ form a basis of open sets of $V(\cS)$. 
\end{sbpara}

\begin{sbpara}\label{valcs}
We define a map $$V(\cS)\to R(\cS)\;;\;V\mapsto r_V.$$
For $V\in V(\cS)$, $r_V\in R(\cS)$ is the map $\cS\times\cS\smallsetminus\{(1,1)\}\to[0,\infty]$ defined by
$$r_V(f,g)=   \sup \{a/b\;|\;(a, b)\in \N^2\smallsetminus \{(0,0)\}, f^b/g^a\in V\} $$ $$= \inf \{a/b \; |\; (a, b)\in \N^2\smallsetminus \{(0,0)\}, g^a/f^b\in V\} $$  
$((f,g)\in (\cS\times \cS)\smallsetminus \{(1,1)\})$ (\ref{ratY}).

\end{sbpara}

\begin{sbprop}\label{valcs2}
The map $V(\cS)\to R(\cS)$ is continuous and surjective. 
\end{sbprop}

\begin{pf} We first prove the continuity of $V(\cS)\to R(\cS)$. 
 Let $f, g \in \cS\smallsetminus \{1\}$, and assume $r_V(f,g)> a/b$ where $a,b\in \N$ and $b>0$. We have $f^b/g^a\in V$. If $V'\in V(\cS)$ and $f^b/g^a\in V'$, we have $r_{V'}(f,g)\geq a/b$. This proves the continuity of $V(\cS)\to R(\cS)$ (\ref{V(S)2}, \ref{ratY}).

We next prove the surjectivity of $V(\cS)\to R(\cS)$. Let 
$\text{class}((\cS^{(j)})_{0\leq j\leq n}, (N_j)_{1\leq j\leq n})\in R'(\cS)$ (\ref{RandR'}). Then the corresponding element of $R(\cS)$  is the image in $R(\cS)$ of the following element $V\in V(\cS)$.
For $1\leq j\leq n$, define the $\Q$-vector subspace $Q^{(j)}$ of the $\Q$-vector space  $\cS_{\Q}=\Q\otimes \cS^{\gp}$ by
$Q^{(j)}:=\Ker(N_j: \cS^{(j-1)}_\Q\to \R)$. Then $Q^{(j)}\supset \cS^{(j)}_{\Q}$. Take an isomorphism 
of $\Q$-vector spaces $\lambda_j : Q^{(j)}/\cS^{(j)}_\Q\overset{\cong}\to \Q^{d(j)}$ where $d(j):=\dim(Q^{(j)}/ \cS^{(j)}_{\Q})$.  
Define $V$ by the following. Let $a\in \cS^{\gp}$. When there is $j$ such that $1\leq j\leq n$, $a\in \cS^{(j-1)}_{\Q}$ and $a\notin Q^{(j)}$, then $a\in V$ if and only if $N_j(a)>0$. When there is $j$ such that $a\in Q^{(j)}$ and $a\notin \cS^{(j)}_\Q$, then $a\in V$ if and only if the first non-zero entry of $\lambda_j(a)\in \Q^{d(j)}$ is $>0$. 
\end{pf}

\begin{sbpara} Let $k$ be a field, let $S$ be the toric variety $\Spec(k[\cS])$, and let 
$X$ be the projective limit as a topological space of the toric varieties over $k$ (with Zariski topology) which correspond to finite rational subdivisions of the cone $\Hom(\cS, \R^{\add}_{\geq 0})$ (\ref{rvtoric}). 
Let $\cO_X$ be the inductive limit of the inverse images on $X$ of the structural sheaves of these toric varieties. Let
$X_0\subset X$ be the inverse image of  $s\in S=\Spec(k[\cS])$ where $s$ is the $k$-rational point of $S$ at which all non-trivial elements of $\cS$ have value $0$. Endow $X_0$ with the topology as a subspace of $X$.

We have a continuous map $X_0\to V(\cS)$ which sends $x\in X_0$ to $\{f\in \cS^{\gp}\;|\; f\in \cO_{X,x}\}\in V(\cS)$. The induced map $X_0(k)\to V(\cS)$ is surjective. In fact, for each $V\in V(\cS)$, the inverse image of $V$ in $X_0$ under the map $X_0\to V(\cS)$ is identified with $\Spec(k[V^\times])$. It has a $k$-rational point which sends all elements of $V^\times$ to $1$. 

Composing with the map in \ref{valcs} as 
$$X_0(k)\subset X_0\to V(\cS)\to R(\cS),$$
and using Proposition \ref{valcs2}, we have

\end{sbpara}

\begin{sbprop}\label{zariski} $(1)$ The map $X_0\to R(\cS)$ is continuous.

$(2)$ The induced map $X_0(k)\to R(\cS)$ is surjective. 
\end{sbprop}

 \begin{sbcor}\label{zarcor} If we regard $R(\cS)$ as a quotient space of $V(\cS)$ or $X_0$, the topology of $R(\cS)$ coincides with the quotient topology.

 \end{sbcor}
 
 This is because $V(\cS)$ and $X_0$ are quasi-compact and $R(\cS)$ is Hausdorff. 
 
 Thus Zariski topology and the topology of real numbers are well connected here. 
 
\subsection{The space $S_{[:]}$ of ratios}\label{ra2}

\begin{sbpara}\label{4.2.1} For a locally ringed space $S$ endowed with an fs log structure, we define the set $S_{[:]}$ as the set of all pairs $(s, r)$ where $s\in S$ and $r \in R((M_S/\cO_S^\times)_s)$. 

We have the canonical surjection $S_{[:]}\to S\;;\;(s,r)\mapsto s$. 
\end{sbpara} 

\begin{sbpara}\label{absv} Let $K$ be a field 
endowed with a non-trivial absolute value $|\:\;|: K\to \R_{\geq 0}$. Let $S$ be a locally ringed space over $K$ satisfying the equivalent conditions in \ref{value}, and assume that we are given an fs log structure on $S$. 

We define a natural topology of $S_{[:]}$ for which the projection $S_{[:]}\to S$ is a proper continuous map and which induces on each fiber of this projection the topology of $R((M_S/\cO_S^\times)_s)$ defined in \ref{ratY}.

\end{sbpara}

\begin{sbpara}\label{ratS} Let $K$ and $S$ be as in \ref{absv}. To define the  topology on $S_{[:]}$, the method is, so to speak, to combine the topology of $S$ and the topologies of $R(\cS)$ (Section \ref{ra1}) for $\cS=(M_S/\cO_S^\times)_s$ $(s\in S)$ by using a chart of the log structure. 

Assume first that we are given a chart $\cS\to M_S$ of the log structure, where $\cS$ is an fs monoid. 
Fix $c\in \R_{>0}$. We have a map
$$S_{[:]}\to [0, \infty]^{\cS\times \cS}\;;\; (s,r)\mapsto r_c$$
 where $r_c: \cS\times \cS \to [0, \infty]$ is defined by the following (1) and (2).
 Let $f,g\in\cS$.
  
 \smallskip
 
 (1) If the images of $f$ and $g$ in $M_{S,s}$ belong to $\cO^\times_{S,s}$, then $$r_c(f,g)=\sup(c, - \log(|f(s)|))/\sup(c, -\log(|g(s)|)).$$

(2) Otherwise, $$r_c(f,g) = r({\bar f}_s, {\bar g}_s)$$
where ${\bar f}_s$ (resp.\ ${\bar g}_s$) denotes the image of $f$ (resp.\ $g$) in $(M_S/\cO_S^\times)_s$. 
\end{sbpara}

\begin{sblem}\label{5.2.4} 
$(1)$ The map $$S_{[:]}\to S\times [0, \infty]^{\cS\times \cS}\;;\; (s,r)\mapsto (s, r_c)$$ is injective.

$(2)$ The topology on $S_{[:]}$ induced by the embedding in  (1) is independent of the choices of the chart and of the constant $c>0$. 
\end{sblem}
\begin{pf} (1) follows from the fact that the map $\cS\to (M_S/\cO^\times_S)_s$ is surjective for any $s\in S$.

We prove (2). 
 If we have two charts $\cS\to M_S$ and $\cS'\to M_S$, we have locally on $S$ a third chart $\cS''\to M_S$ with homomorphisms of charts $\cS\to \cS''$ and $\cS'\to \cS''$. It is clear that if these third chart and  two homomorphisms of charts are given and if the constant $c>0$ is fixed, the topology given by the chart $\cS''\to M_S$  and $c$ is finer than the topology given by $\cS\to M_S$ or $\cS'\to M_S$ and $c$. Hence it is sufficient to prove that if we have a homomorphism $\cS'\to \cS$ from a chart $\cS'\to M_S$ to a chart $\cS\to M_S$, the topology given by the former chart and the constant $c'>0$ is finer than the topology given by the latter and $c>0$. 
It suffices to prove that for $f,g\in \cS$, the map $(s,r)\mapsto r_c(f,g)$ is continuous for the topology given by $\cS'\to M_S$ and $c'$. 
 
 \smallskip

{\bf Claim 1.} Let $f,g\in \cS$ and let $s\in S$, and assume that the images of $f$ and $g$ in $(M_S/\cO^\times_S)_s$ coincide. Let $c,c'>0$. Then for some neighborhood $U$ of $s$ in $S$, we have a continuous map $R_{c,c'}(f,g): U\to \R_{>0}$ whose value at $s'\in U$ is $\sup(c, -\log(|f(s')|))/\sup(c', -\log(|g(s')|))$ 
if the images of $f$ and $g$ in $M_{S,s'}$ belong to $\cO_{S,s'}^\times$, and is $1$ otherwise. 

\smallskip

This Claim 1 is proved easily.

We continue the proof of (2). Let $f,g\in \cS$. Then locally on $S$, we have $f',g'\in \cS'$ and sections $u,v$ of $\cO_S^\times$ such that $f=f'u$ and $g=g'v$ in $M_S$. We have
$$r_c(f,g)= r_{c'}(f', g')R_{c,c'}(f,f')(s)R_{c',c}(g', g)(s).$$ This proves the desired continuity of $r_c(f,g)$. \end{pf}
\begin{sbpara}\label{topcs}

By the independence (2) in \ref{5.2.4}, we have a canonical topology of $S_{[:]}$ (globally). 

\end{sbpara}

\begin{sbpara}\label{sharp:} Assume that $\cS$ is sharp and that for any $f\in \cS\smallsetminus \{1\}$ and any $s\in S$, we have $|f(s)|<1$. (Note that we have such a chart locally on $S$.)
Let $Y=(\cS\times \cS) \smallsetminus \{(1,1)\}$. 
Then we have a slightly different embedding 
$$S_{[:]}\to S\times [0, \infty]^Y;\;\;(s,r)\mapsto (s, r_*)$$
where $r_*: Y\to [0, \infty]$ is defined as follows.
Let $(f,g)\in Y$.

(1) If the images of $f$ and $g$ in $M_{S,s}$ belong to $\cO^\times_{S,s}$, then $$r_*(f,g)= \log(|f(s)|)/\log(|g(s)|).$$

(2) Otherwise, $$r_*(f,g)= r({\bar f}_s, {\bar g}_s)$$ where ${\bar f}_s$ (resp.\ ${\bar g}_s$) denotes the image of $f$ (resp.\ $g$) in $(M_S/\cO_S^\times)_s$.

\end{sbpara}

\begin{sblem}\label{ratL} 
  Let the assumption be as in \ref{sharp:}. 

$(1)$ The map $S_{[:]}\to S\times [0, \infty]^Y$ is injective.

$(2)$ The topology of $S_{[:]}$ induced by this embedding coincides with the topology defined in \ref{topcs}.

$(3)$ The image of the embedding (1) consists of all pairs $(s, r)\in S\times [0, \infty]^Y$ such that $r$ satisfies 
the conditions {\rm (i)}--{\rm (iii)} in \ref{ratY} and such that the following conditions {\rm (iv)} and {\rm (v)} are satisfied.
Let $(f,g)\in Y$.

\medskip

{\rm (iv)} If the images of $f$ and $g$ in $M_{S,s}$ belong to $\cO_{S,s}^\times$, $r(f,g)= \log(|f(s)|)/\log(|g(s)|)$. 

{\rm (v)} Otherwise, $r(f,g)$ depends only on the images of $f$ and $g$ in $(M_S/\cO^\times_S)_s$.

\medskip

$(4)$ The image of the embedding in $(1)$ is a closed set of $S\times [0, \infty]^Y$. 
\end{sblem}

\begin{pf}
(1) and (3) follow from the fact that the map $\cS\to (M_S/\cO^\times_S)_s$ is surjective for any $s\in S$.

(4) follows from (3). 

We prove (2). 
If $f\in \cS\smallsetminus \{1\}$, by the property $|f(s)|<1$ for any $s\in S$, we see that 
there is a continuous function 
$R_{c}(f): S\to \R_{>0}$ whose value at $s\in S$ is $-\sup(c, -\log(|f(s)|))/\log(|f(s)|)$ if the image of $f$ in $M_{S,s}$ belongs to $\cO_{S,s}^\times$  and is $1$ otherwise. For $f,g\in \cS\smallsetminus \{1\}$, we have  $$r_c(f,g)= r_*(f,g)R_{c}(f)(s) R_{c}(g)(s)^{-1}.$$ Furthermore,  for $f\in \cS$, $r_c(1,f)$  is the value of the continuous function $c/\sup(c, -\log(|f(s)|))$ at $s$, and $r_c(f,1)$ is the value of the continuous function $\sup(c, -\log(|f(s)|))/c$ at $s$, and for $f\in \cS\smallsetminus \{1\}$, we have $r_*(1,f)=0$, $r_*(1, f)=\infty$. 
\end{pf}

\begin{sbprop}
The canonical map $S_{[:]}\to S$ is proper.
\end{sbprop}

\begin{pf}
Since $[0, \infty]^Y$ is compact, this follows from (4) of Lemma \ref{ratL}. 
\end{pf}

\begin{sbpara} 
For each $s\in S$, the topology of $R((M_S/\cO^\times_S)_s)$ defined in Section \ref{ra1} coincides with the topology of the fiber $R((M_S/\cO^\times_S)_s)$ over $s$ of $S_{[:]}\to S$, as a subspace of $S_{[:]}$.
\end{sbpara}

\begin{sblem} Let $S$ and $S'$ be as in \ref{absv} and assume we are given a strict morphism $S'\to S$ of locally ringed spaces over $K$ with log structures. (For the word \lq\lq strict'', see \ref{gest6}.) Then the canonical map $S'_{[:]}\to S'\times_S S_{[:]}$ is a homeomorphism.

\end{sblem}

\begin{pf} This is proved in the same way as Lemma \ref{valsrt}. 
\end{pf} 
\begin{sbpara}\label{rat4} We consider $S_{[:]}$ more locally. 

Assume we are given a chart  $\cS\to M_S$.

Let $\Phi$ be a set of faces of $\cS$ which is totally ordered for the inclusion relation and which contains $\cS$. 
Let
$S_{[:]}(\Phi)$ be the subset of $S_{[:]}$ consisting of all $(s, r)$ such that the inverse images in $\cS$ of the faces of $(M_S/\cO_S^\times)_s$ associated to $r$ (\ref{RandR'}) belong to $\Phi$.

Then $S_{[:]}(\Phi)$ for all $\Phi$ forms an open covering of $S_{[:]}$.
\end{sbpara}

\begin{sbpara} Let the notation be as in \ref{rat4}. Assume further that for any $f\in \cS\smallsetminus \{1\}$, we have $|f(s)|<1$ for any $s\in S$. (Such a chart always exists locally on $S$.) In the following proposition \ref{Phistr}, we give a description of the topological space $S_{[:]}(\Phi)$.

Write $\Phi = \{\cS^{(j)} \;|\; 0\leq j \leq n\}$, $\cS= \cS^{(0)} \supsetneq \cS^{(1)}\supsetneq \dots \supsetneq \cS^{(n)}$. For each $1\leq j\leq n$, fix $q_j\in \cS^{(j-1)}\smallsetminus \cS^{(j)}$. Consider the topological subspace $$P\subset \R_{\geq 0}^n \times \prod_{j=0}^n \Hom(\cS^{(j)}, \R^{\add})$$ (here the Hom space is endowed with the topology of simple convergence) consisting of elements $(t,h)$ $(t=(t_j)_{1\leq j\leq n}, t_j\in \R_{\geq 0}, h=(h_j)_{0\leq j\leq n}, h_j: \cS^{(j)}\to \R)$ satisfying the following conditions (i) -- (iii) for $0\leq j<n$.

(i) $h_j(q_{j+1})=1$. 

(ii) $h_j(f) = t_{j+1}h_{j+1}(f)$ for any $f\in \cS^{(j+1)}$. 

(iii) $h_j(\cS^{(j)}\smallsetminus \cS^{(j+1)}) \subset \R_{>0}$. 

\end{sbpara}

\begin{sblem}\label{5lemA} We have a unique continuous map $P\to \Hom(\cS, \R^{\mult}_{\geq 0})$ which sends 
$(t,h)$ to the following $a\in \Hom(\cS, \R^{\mult}_{\geq 0})$. 
Let $j$ be the smallest integer such that $0\leq j\leq n$ and such that $t_k\neq 0$ if $j<k\leq n$. Then 
$$a(f) = \exp(- h_j(f) \prod_{k=j+1}^n t_k^{-1})\in \R_{>0}\quad \text{if} \; f \in \cS^{(j)},$$
$$a(f) = 0 \quad\text{if}\; f \in \cS\smallsetminus \cS^{(j)}.$$
\end{sblem}

\begin{pf}
The problem is the continuity of the map. This is shown as follows. Let $f\in \cS$. It is sufficient to prove that the map $P \to \R_{\geq 0}\;;\;(t,h)\mapsto a(f)$ $(f\in\cS)$ (with notation as above) is continuous. Let $j$ be the largest integer such that $0\leq j\leq n$ and such that $f\in \cS^{(j)}$. Then this map 
 is the composition of the continuous map $P\to \R_{\geq 0}$ which sends $((t_j)_j, (h_j)_j)\in P$ to $\prod_{k=j+1}^n t_k \cdot h_j(f)^{-1}$ (note $h_j(f)>0$) and the continuous map $\R_{\geq 0}\to \R_{\geq 0}$ which sends $t\in \R_{>0}$ to $\exp(-t^{-1})$ and $0$ to $0$.
\end{pf}

\begin{sbprop}\label{Phistr} 
Let the notation be as above. We have a cartesian diagram of topological spaces
$$\begin{matrix}
S_{[:]}(\Phi) &\to & P \\
\downarrow && \downarrow \\
S & \to& \Hom(\cS, \R^{\mult}_{\geq 0})
\end{matrix}$$
where the lower horizontal  arrow sends $s\in S$ to the map $f\mapsto |f(s)|$ $(f\in\cS)$, the right vertical arrow is as $a\mapsto a(f)$ $(f\in\cS)$ in \ref{5lemA}, the left vertical arrow is the canonical one, and the upper horizontal 
arrow sends $(s,r)\in S_{[:]}(\Phi)$ $(s\in S$, $r\in R((M_S/\cO_S^\times)_s))$ to $(s, ((t_j)_j, (h_j)_j))$ where 
 $t_j= \log(|q_{j+1}(s)|)/\log(|q_j(s)|)$ (resp.\ $t_j= r(q_{j+1}, q_j)$) if $1\leq j<n$ and if $q_jq_{j+1}$ is invertible (resp.\ not invertible) at $s$, $t_n= -1/\log(|q_n(s)|)$, $h_j(f)=r(f, q_{j+1})$ for $0\leq j<n$, and $h_n(f) = -\log(|f(s)|)$. (Note that if $(s,r)\in S_{[:]}(\Phi)$ and $f\in \cS^{(n)}$, the image of $f$ in $M_{S,s}$ belongs to $\cO_{S,s}^\times$ and hence $|f(s)|\in \R_{>0}$.)
\end{sbprop}

\begin{pf}
The converse map is given by $(s, (t, h))\mapsto (s,r)$ where 
$r$ is as follows. Let $(a, b)\in (M_S/\cO_S^\times)_s\times (M_S/\cO^\times_S)_s\smallsetminus \{(1,1)\}$
 and take $f, g\in \cS$ such that the image of $f$ (resp.\ $g$) in $(M_S/\cO_S^\times)_s$ is
  $a$ (resp.\ $b$). Take the largest $j$ such that $0\leq j\leq n-1$ and $f,g\in \cS^{(j)}$. 
  Then $r(f,g)= h_j(f)/h_j(g)\in [0,\infty]$. 
  (Note that at least one of $f, g$ is outside $\cS^{(j+1)}$ and hence at least one of  $h_j(f)$ and $h_j(g)$ is non-zero.) 

It is easy to see that this is the converse map and continuous. 
\end{pf}

\begin{sbrem}  In $(t,h)\in P$ $(t=(t_j)_{1\leq j\leq n}\in \R^n_{\geq 0})$, $t_j$ for $1\leq j\leq n-1$ is determined by $h$ as $t_j= h_{j-1}(q_{j+1})$. $t_n$ is determined by the image $a$ of $(t,h)$ in $\Hom(\cS, \R^{\mult}_{\geq 0})$ as $t_n=- 1/\log(a(q_n))$. 

These explain the fact that in the above proof of \ref{Phistr}, the converse map $(s, (t,h))\mapsto (s,r)$ is described without using $t$. 

\end{sbrem}

\begin{sbpara}\label{4.2.14} Let $\Hom(\cS, \R^{\mult}_{\geq 0})_{<1}$ be the open set of $\Hom(\cS, \R^{\mult}_{\geq 0})$ consisting of all elements $h$ such that $h(f)<1$ for any $f\in \cS\smallsetminus \{1\}$. Then the images of $S$ and $P$ in $\Hom(\cS, \R^{\mult}_{\geq 0})$, under the maps in \ref{Phistr}, are contained in $\Hom(\cS, \R^{\mult}_{\geq 0})_{<1}$. Hence, by \ref{Phistr}, we have 

\end{sbpara}

\begin{sbcor}\label{4.2.15} 
In the case $S=\Hom(\cS, \R^{\mult}_{\geq 0})_{<1}$ with the sheaf of all $\R$-valued continuous functions and with the natural log structure, $S_{[:]}(\Phi)$ is identified with $P$. 
\end{sbcor}

\begin{sbpara} We give a comment on this space $P$. 

For $1\leq j\leq n$, we fixed an element $q_j$ of $\cS^{(j-1)}\smallsetminus \cS^{(j)}$ (\ref{Phistr}). 
Let $m(j)= \dim_{\Q}(\cS^{(j-1)}_\Q/\cS^{(j)}_\Q) -1$ if $1\leq j\leq n$, and let $m(n+1)=\dim_{\Q}(\cS^{(n)}_\Q)$. 
For $1\leq j\leq n+1$, fix elements $q_{j,k}$ $(0\leq k\leq m(j))$ of $(\cS^{(j-1)})^{\gp}$ satisfying the following conditions (i)--(iii).

(i) $q_{j,0}=q_j$ if $1\leq j\leq n$. 

(ii)  For $1\leq j\leq n$, $(q_{j,k}\bmod \cS^{(j)}_\Q)_{0\leq k\leq m(j)}$ is  a $\Q$-basis of $\cS^{(j-1)}_\Q/\cS^{(j)}_\Q$.

(iii) $(q_{n+1,k})_{1\leq k\leq m(n+1)}$ is  a $\Q$-basis of $\cS^{(n)}_\Q$.

\end{sbpara}

\begin{sbprop}\label{qjk}

We have an injective open map  
$$P\overset{\subset}\to 
\R_{\geq 0}^n\times \prod_{j=1}^{n+1} \R^{m(j)}$$
which sends  $(t, h)\in P$ $(t\in \R^n_{\geq 0}$, $h\in \prod_{j=0}^n \Hom(\cS^{(j)}, \R^{\add}))$ to $(t,  a)$ where
$a=(a_j)_{1\leq j\leq n+1}$, $a_j=(a_{j,k})_{1\leq k\leq m(j)}$ with  
$$ a_{j,k}= h_{j-1}(q_{j,k})\in \R$$ for $1\leq j\leq n+1$.
Here we define $h_{j-1}(q_{j, k})$ by using the unique extension of $h_{j-1}:\cS^{(j-1)}\to \R$ to a homomorphism $(\cS^{(j-1)})^{\gp}\to \R^{\add}$.

\end{sbprop}
The proof is easy.

\begin{sbpara}\label{Nn:}

Consider the case $S=|\Delta|^n$ where $|\Delta|=\{t\in \R\;|\;0\leq t<1\}$ with the sheaf of all $\R$-valued continuous functions and with the fs log structure associated to $\N^n\to \cO_S\;;\; m\mapsto \prod_{j=1}^n q_j^{m(j)}$ where $q_j$ $(1\leq j\leq n)$ are the coordinate functions. Let $\cS$ be the multiplicative monoid generated by $q_j$ $(1\leq j\leq n)$ which is identified with $\N^n$. Then $|\Delta|^n$ is identified with $\Hom(\cS, \R^{\mult}_{\geq 0})_{<1}$ in
\ref{4.2.14}.

 Let
$\Phi=\{\cS^{(j)}\; |\; 0\leq j\leq n\}$ where $\cS^{(j)}$ is generated by $q_k$ $(j< k\leq n)$.   Then $S_{[:]}$ is covered by the open sets $S_{[:]}(g(\Phi))$ where $g$ ranges over elements of the permutation group $\frak S_n$ acting on $\cS$, and $g$ induces a homeomorphism $S_{[:]}(\Phi)\cong S_{[:]}(g(\Phi))$. 
We describe $S_{[:]}(\Phi)$. 

\end{sbpara}

\begin{sbprop}\label{Nnphi} Let the notation be as in \ref{Nn:}. Then 
we have a commutative diagram 
$$\begin{matrix}
S_{[:]}(\Phi)&\cong & \R^n_{\geq 0}\\
\downarrow &&\downarrow\\
S&=& |\Delta|^n
\end{matrix}$$
in which the upper horizontal isomorphism sends $(s,r)\in S_{[:]}(\Phi)$ to 
$(t_1, \dots, t_n)$ 
where $t_j=r(q_{j+1}, q_j)$ $(1\le j\le n-1)$ and $t_n=-1/\log(q_n(s))$, 
and the right vertical arrow is $(t_j)_{1\leq j\leq n}\mapsto (q_j)_{1\leq j\leq n}$ where $q_j=\exp(-\prod_{k=j}^n t_k^{-1})$.

\end{sbprop}

\begin{pf} This follows from \ref{4.2.15}.
\end{pf}

\begin{sbcor}\label{RNn} Let the notation be as in \ref{Nn:}. Regarding $R(\cS)$ as the fiber of $S_{[:]}\to S=|\Delta|^n$ over the point $(0,\dots,0)\in S$,
define $R(\cS)(\Phi)= R(\cS)\cap S_{[:]}(\Phi)$. Then 
we have a homeomorphism
$$R(\cS)(\Phi)\cong \R^{n-1}_{\geq 0}$$
which sends $r\in R(\cS)(\Phi)$ to $(t_1, \dots, t_{n-1})$ where $t_j= r(q_{j+1}, q_j)$. 

\end{sbcor}
\begin{pf}
This follows from  \ref{Nnphi}.
\end{pf}

\begin{sblem}\label{fsK} Let $S$ and $|S|$, $M_{|S|}$ be as in \ref{abslog}. Then we have a canonical homeomorphism $S_{[:]}\cong |S|_{[:]}$. 
\end{sblem}

\begin{pf} As in the proof of \ref{abslog}, we have a canonical isomorphism $(M_S/\cO^\times_S)_s\cong (M_{|S|}/\cO^\times_{|S|})_s$ for each $s\in S$. This gives a canonical bijection between $S_{[:]}$ and $|S|_{[:]}$. 
By 
 Proposition \ref{Phistr}, 
 they have the same topology. 
 \end{pf}

\subsection{$S_{[:]}$, $S_{[\val]}$, and $S_{\val}$}\label{ra5}

Let $K$ and $S$ be as in \ref{absv}. 

 We construct a topological space $S_{[\val]}$ and proper surjective continuous maps $$S_{\val}\to S_{[:]}, \quad S_{[\val]}\to S_{[:]}.$$ Here $S_{\val}$ is as in Section \ref{ss:valsp}.

\begin{sbpara}
Let $S_{\val}\to S_{[:]}$ be the map $(s, V, h)\mapsto (s, r_V)$ 
  where $V\mapsto r_V$ is the map 
$V(\cS)\to R(\cS)$ for $\cS=(M_S/\cO^\times_S)_s$ (\ref{Sgen}, \ref{valcs}).

\end{sbpara}

\begin{sbprop} The map $S_{\val}\to S_{[:]}$ is continuous, and proper and surjective. 
 
\end{sbprop}

\begin{pf}
The surjectivity follows from the surjectivity in \ref{valcs2}. Once we prove the continuity, properness follows from the properness of $S_{\val}\to S$ and of $S_{[:]}\to S$. We prove the continuity. Working locally on $S$, we may and do assume that we have a chart $\cS\to M_S$ with $\cS$ a sharp fs monoid such that for any $f\in \cS\smallsetminus \{1\}$ and $s\in S$, we have $|f(s)|<1$. 

Fix $(s_0,V_0,h_0)\in S_{\val}$ and let $(s_0, r_0)\in S_{[:]}$ be its image. 
 We show that when $(s, V, h)\in S_{\val}$ converges to $(s_0, V_0, h_0)$, its image $(s,r)\in S_{[:]}$ converges to $(s_0, r_0)$. 
Let $f,g\in \cS\smallsetminus \{1\}$. It is sufficient to prove that 
$r_*(f,g)\in [0,\infty]$ (\ref{sharp:}) converges to $(r_0)_*(f,g)\in [0,\infty]$. If at least one of $f$ and $g$ is invertible at $s_0$ (that is, if at least one of the images of $f$ and $g$ in $M_{S, s_0}$ belongs to $\cO_{S,s_0}^\times$), 
then the function $(s,r)\mapsto r_*(f,g)\in [0, \infty]$ on $S_{[:]}$ comes from the continuous function $s\mapsto \log(|f(s)|)/\log(|g(s)|)\in [0,\infty]$  on some neighborhood of $s_0$ in $S$. Hence we may assume that both $f$ and $g$ are not invertible at $s_0$. Assume $(r_0)_*(f,g) >a/b$, $a,b\in \N$, $b>0$. 
 It is sufficient to prove that $r_*(f,g) >a/b$ when $(s,V,h)$ is sufficiently near  $(s_0, V_0,h_0)$. Let $\varphi=f^b/g^a\in \cS^{\gp}$. Since the image $\bar \varphi_{s_0}$ of $\varphi$ in $(M_S/\cO_S^\times)_{s_0}$ belongs to $V_0$, there is a neighborhood $U$ of $(s_0,V_0, h_0)$ in $S_{\val}$ such that if $(s,V, h)\in U$, then $\bar \varphi_s\in V$. 
  If $(s,V,h)\in U$ and if at least one of $f$ and $g$ are not invertible at $s$, then $r_*(f,g)=r({\bar f}_s, {\bar g}_s)\geq a/b$ because $\bar \varphi_s\in V$. Consider points $(s,V,h)\in U$ such that  both $f$ and $g$ are invertible at $s$.  On $U$, the function $(s,V, h)\mapsto h(\varphi)$ is continuous. 
(Here $h(\varphi)$ is defined to be $0$ if $\bar \varphi_s\notin V^\times$.)
We have $$r_*(f,g)=b^{-1}r_*(f^b,g)= b^{-1}r_*(g^a\varphi, g)= (a/b) + b^{-1}\log(h(\varphi))/\log(|g(s)|).$$
When $(s, V, h)\in U$ converges to $(s_0, V_0, h_0)$, $h(\varphi)$ converges to $h_0(\varphi)\in\R$ and $g(s)$ converges to $0$. If $h_0(\varphi) =0$, then  when $(s,V,h)$ converges to $(s_0, V_0, h_0)$, we have $h(\varphi)<1$ and $|g(s)|<1$ and hence $\log(h(\varphi))/\log(|g(s)|)>0$. If $h_0(\varphi)>0$, then when $(s,V,h)$ converges to $(s_0,V_0,h_0)$, $\log(h(\varphi))/\log(|g(s)|)$ converges to $0$. 
 \end{pf}

\begin{sbpara}\label{rat6} We next discuss $S_{[\val]}$. To define it, we use the following {\it new log structure on $S_{[:]}$} which is endowed with the sheaf $\cO_{S_{[:]}}$ of all $\R$-valued continuous functions.  (We use the word \lq\lq new log structure'', to distinguish this log structure from  the \lq\lq old'' log structure on $S_{[:]}$ which is defined as the inverse image of the log structure of $S$ on the inverse image of $\cO_S$ on $S_{[:]}$.)

Assume that we are given a chart $\cS\to M_S$ with $\cS$ a sharp fs monoid such that $|f(s)| <1$ for any $f\in \cS\smallsetminus \{1\}$ and for any $s\in S$.  Let $\cS^{(j)}$ $(0\leq j\leq n)$ be faces of $\cS$ such that 
$\cS=\cS^{(0)}\supsetneq \cS^{(1)}\supsetneq \dots \supsetneq \cS^{(n)}$ and let $\Phi=\{\cS^{(j)}\;|\; 0\leq j\leq n\}$. Take $q_j\in \cS^{(j-1)}\smallsetminus \cS^{(j)}$ for $1\leq j\leq n$. Then we define the new log structure on $S_{[:]}(\Phi)$ as the fs log structure associated to 
$$\N^n \to \cO_{S_{[:]}}\;;\; m\mapsto 
(\prod_{j=1}^{n-1} r(q_{j+1}, q_j)^{m(j)/2})\cdot (-1/\log(|q_n|))^{m(n)/2}.$$
Then it is easy to see that this log structure glues to an fs log structure on $S_{[:]}$ which is independent of any choices. 

In the identification $S_{[:]}=|S|_{[:]}$ (\ref{fsK}), the new log structure of $S_{[:]}$ and that of $|S|_{[:]}$ coincide. 

\end{sbpara}

\begin{sbrem}
It may seem strange to take the 
square root $(-)^{m(j)/2}$ in the definition of this log structure. 
But this becomes important in Section 5 to have that the CKS map $D^{\sharp}_{\Sig,[:]}\to D_{\SL(2)}$ respects (and $D^{\sharp,\mild}_{\Sig,[:]}\to  D_{\SL(2)}^{\diamond}$ which appears later (\ref{diathm}) also respects) the log structures. 

\end{sbrem}

\begin{sbpara}
Let $S_{[\val]}$ be the valuative space $(S_{[:]})_{\val}$ (Section \ref{ss:valsp}) associated to $S_{[:]}$ endowed with this new log structure. 

By Section \ref{ss:valsp}, the map $S_{[\val]}\to S_{[:]}$ is proper and surjective.

\end{sbpara}

\begin{sblem} Let $S$ (resp.\ $S'$) be a topological space endowed with the sheaf of all $\R$-valued continuous functions and with an fs log structure, and let $S'\to S$ be a strict morphism (\ref{gest6}) of locally ringed spaces over $\R$ with log structures. Then the canonical map $S'_{[\val]}\to S'\times_S S_{[\val]}$ is a homeomorphism.

\end{sblem}

\begin{pf}
This is proved in the same way as \ref{valsrt}. 
\end{pf}

\begin{sbprop}\label{v[[v]]}

Assume $K=\R$. 

There is a unique homeomorphism  $$(|\Delta|^n)_{[\val]}\cong (\R^n_{\geq 0})_{\val}$$ in
which  $(q_j)_{1\leq j\leq n}\in (|\Delta|\smallsetminus \{0\})^n\subset (|\Delta|^n)_{[\val]}$ corresponds to $(-1/\log(q_j))_{1\leq j\leq n}\in \R^n_{>0}\subset (\R^n_{\geq 0})_{\val}$.
 
\end{sbprop}

\begin{pf}
This is deduced from Proposition \ref{Nnphi}. 
\end{pf}
\begin{sbpara}\label{rvalex}  {\bf Example.} We compare $S_{[:]}$, $S_{\val}$, and $S_{[\val]}$ in the case $K=\R$ and $S=\R^2_{\geq 0}$ with the standard log structure. The maps from these spaces to $S$ are homeomorphisms outside $(0,0)\in S$. We describe the fibers over $(0,0)$ explicitly. 
\medskip

(1)  The fiber of $S_{[:]}\to S$ over $(0,0)\in S$ is canonically homeomorphic to  the interval $[0, \infty]$.  It consists of points 
$r(a)$ with $a\in [0, \infty]$. $(q_1,q_2)\in \R^2_{>0}$ converges to $r(a)$ if and only if $q_j\to 0$ and $\log(q_2)/\log(q_1)\to a$.

\medskip
(2) A difference between the surjection $S_{\val}\to S_{[:]}$ and the surjection $S_{[\val]}\to S_{[:]}$ is that the fiber of the former surjection over $r(a)$ has cardinality $>1$ if and only if $a\in \Q_{>0}$ and the fiber of the latter surjection over $r(a)$ has cardinality $>1$ if and only if $a=0$ or $a=\infty$.

\medskip

(3) The fiber of $S_{\val}\to S$ over  $(0,0)\in S$ consists of  points $p(a)$ $(a\in [0, \infty]\smallsetminus \Q_{>0})$ and 
$p(a, c)$ $(a\in \Q_{>0}$, $c\in [0, \infty])$.

$(q_1,q_2)\in \R^2_{>0}$ converges to $p(a)$ if and only if $q_j\to 0$ and $\log(q_2)/\log(q_1) \to a$.

$(q_1,q_2)\in \R^2_{>0}$ converges to $p(a,c)$ if and only if $q_j\to 0$, $\log(q_2)/\log(q_1) \to a$, and $q_1^a/q_2\to c$.  

Under the map $S_{\val}\to S_{[:]}$, 
 $p(a)$ goes to $r(a)\in S_{[:]}$ and 
$p(a, c)$ goes to $r(a)\in S_{[:]}$. 

\medskip

(4) The fiber of $S_{[\val]}$ over $(0,0)\in S$ consists of points  $s(a)$ $(a\in [0,\infty]\smallsetminus \Q_{>0})$ and 
$s(a, c)$ $(a\in \Q_{>0}$, $c\in [0, \infty])$.

$(q_1,q_2)\in \R^2_{>0}$ converges to $s(a)$ if and only if $q_j\to 0$ and, for $t_j:= -1/\log(q_j)$ (so $t_j\to 0$), $\log(t_2)/\log(t_1)\to a$.

$(q_1,q_2)\in \R^2_{>0}$ converges to $s(a,c)$ if and only if $q_j\to 0$ and, for $t_j:= -1/\log(q_j)$ (so $t_j\to 0$), $\log(t_2)/\log(t_1) \to a$, and $t_1^a/t_2$ converges to $c$. 

Under the map $S_{[\val]}\to S_{[:]}$, 
$s(1,c)$ goes to $r(c)\in S_{[:]}$. 
$s(a)$ with $a<1$ and $s(a, c)$ with $a<1$ go to $r(0)$ in $S_{[:]}$, and 
$s(a)$ with $a>1$ and $s(a, c)$ with $a>1$ go to $r(\infty)$ in $S_{[:]}$.
\medskip

(5) Some examples of convergences. 

(5.1)  Fix $c\in \R_{>0}$. If $q\in \R_{>0}$ and $q\to 0$, $(cq, q)\in \R^2_{>0}$ converges to $r(1)$ in $S_{[:]}$, to $p(1, c)$ in $S_{\val}$, and to $s(1, 1)$ in $S_{[\val]}$. Thus the limit in $S_{[:]}$ and the limit in $S_{[\val]}$ are independent of $c$, but the limit in $S_{\val}$ depends on $c$.

(5.2)  Fix $a\in \R$ such that $0<a<1$. If $t\in \R_{>0}$ and $t\to 0$, $(\exp(-1/t), \exp(-1/t^a))\in \R_{>0}^2$ 
converges to $r(0)$ in $S_{[:]}$, to $p(0)$ in $S_{\val}$, and to $s(a)$ (resp.\ $s(a,1)$) in $S_{[\val]}$ if $a\notin \Q$ (resp.\ $a\in \Q$). Thus  the 
limit in $S_{[:]}$ and the limit in $S_{\val}$ are independent of $a$, but the limit in $S_{[\val]}$ depends on $a$.

\end{sbpara}

\subsection{The spaces $D^{\sharp}_{\Sig, [:]}$ and $D^{\sharp}_{\Sig,[\val]}$}\label{ra6}

\begin{sbpara} Let $\Sig$ be a weak fan in  $\fg_\Q$ (Part III, 2.2.3) and let $\Gamma$ be a neat subgroup of $G_{\Z}$ which is strongly compatible with $\Sig$. Then we have a space $\Gamma \bs D_{\Sig}$ which is endowed with a sheaf of holomorphic functions and an fs log structure. By taking $K=\C$ in Section 4.2, we have a
topological space  $(\Gamma \bs D_{\Sig})_{[:]}$ with a proper surjective map
$(\Gamma \bs D_{\Sig})_{[:]}\to \Gamma \bs D_{\Sig}$.

\end{sbpara}

\begin{sbpara} Let $\Sig$ be a weak fan in $\fg_\Q$  and let $D^{\sharp}_{\Sig}$ be the topological space defined in Part III, 2.2.5. We define topological spaces $D^{\sharp}_{\Sig,[:]}$ 
and $D^{\sharp}_{\Sig,[\val]}$, and proper surjective maps $D^{\sharp}_{\Sig,[:]}\to D^{\sharp}_{\Sig}$, $D^{\sharp}_{\Sig,\val}\to D^{\sharp}_{\Sig, [:]}$, and  $D^{\sharp}_{\Sig,[\val]}\to D^{\sharp}_{\Sig, [:]}$. Here $D^{\sharp}_{\Sig,\val}$ is the topological space defined in Part III, 3.2.

\end{sbpara}

\begin{sbpara}\label{4.4.3}
Let $\sig\in \Sig$ and consider the open set $D_{\sig}^{\sharp}$ of $D_{\Sig}^{\sharp}$. There is a neat subgroup $\Gamma$ of $G_\Z$ which is strongly compatible with the fan $\text{face}(\sig)$ of all faces of $\sig$. 

We define the topological space $D^{\sharp}_{\sig,[:]}$ as the fiber product of $D^{\sharp}_{\sig}\to \Gamma \bs D_{\sig}\leftarrow (\Gamma \bs D_{\sig})_{[:]}$. This is independent of the choice of $\Gamma$. 

Furthermore,  the inverse image of the new log structure of $(\Gamma \bs D_{\sig})_{[:]}$ on 
$D^{\sharp}_{\sig,[:]}$ (given on the sheaf of all $\R$-valued continuous functions), which we call the {\it new log structure of $D^{\sharp}_{\sig,[:]}$}, is independent of the choice of $\Gamma$. 

These $D^{\sharp}_{\sig,[:]}$ glue to a topological space $D^{\sharp}_{\Sig,[:]}$ over $D^{\sharp}_{\Sig}$, and the new log structures of 
$D^{\sharp}_{\sig,[:]}$ glue to an fs log structure on the sheaf of all $\R$-valued functions on $D^{\sharp}_{\Sig,[:]}$, which we call {\it the new log structure}. 

We define $D^{\sharp}_{\Sig,[\val]}$ as the valuative space associated to $D^{\sharp}_{\Sig,[:]}$ with the new log structure. 

We have a canonical proper surjective maps $D^{\sharp}_{\Sig,[:]}\to D^{\sharp}_{\Sig}$ and $D^{\sharp}_{\Sig,[\val]}\to D^{\sharp}_{\Sig,[:]}$. 
\end{sbpara}

\begin{sbpara}\label{trouble} Before we define the canonical map $D^{\sharp}_{\Sig, \val}\to D^{\sharp}_{\Sig, [:]}$, we remark that, though we have a canonical  new log structure on $D^{\sharp}_{\Sig,[:]}$, we do not have a canonical log structure on $D^{\sharp}_{\Sig}$. For $\sig\in \Sig$ and for a neat subgroup $\Gamma$ of $G_{\Z}$ which is strongly compatible with $\text{face}(\sig)$, the pull-back of the log structure of $\Gamma \bs D_{\sig}$ on $D^{\sharp}_{\sig}$ depends on the choice of $\Gamma$. Here we endow $D^{\sharp}_{\sig}$ the sheaf of all $\C$-valued continuous functions. 

For example, consider the classical case $H_{0,\Z}=\Z^2$ of pure of weight $1$ of Hodge type $(1,0)+(0,1)$, in which $D$ is the upper half plane. For the standard choice of $\sig$ and 
$\Gamma= \begin{pmatrix} 1& \Z \\
0&1\end{pmatrix}$,  $\Gamma \bs D_{\sig}$ is isomorphic to the unit disc and the log structure is generated by the coordinate function $q$.  $D^{\sharp}_{\sig}$ is identified with $\{x+iy\;|\; x\in \R, 0< y\leq \infty\}$ and the canonical projection $D^{\sharp}_{\sig}\to \Gamma\bs D_{\sig}$ is identified with $z\mapsto \exp(2\pi i z)$. We have $D^{\sharp}_{\sig,[:]}=D^{\sharp}_{\sig}$ and the new log structure on it is generated by $1/y^{1/2}$, or equivalently by $1/(\log|q|)^{1/2}$.

Take $n\geq 2$ and replace $\Gamma$ by $\Gamma':=\begin{pmatrix} 1&n\Z\\0&1\end{pmatrix}$. Then
 the log structure of $\Gamma' \bs D_{\sig}$ is generated by $q^{1/n}$. Hence the inverse image on $D^{\sharp}_{\sig}$ of the log structure of $\Gamma \bs D_{\sig}$ and that of $\Gamma'\bs D_{\sig}$ do not coincide.

This problem does not happen for the new log structure, for $1/(\log|q^{1/n}|)^{1/2}=n^{1/2}/(\log|q|)^{1/2}$ and $1/(\log|q|)^{1/2}$ generate the same log structure. 

\end{sbpara}

\begin{sbpara}\label{valD}  Endow $D^{\sharp}_{\Sig}$ with the sheaf of all $\C$-valued continuous functions. 

For $\sig\in \Sig$, take a neat sungroup $\Gamma$ of $G_{\Z}$ which is strongly compatible with $\sig$, and consider the inverse image on $D^{\sharp}_{\sig}$ of the log structure of $\Gamma \bs D_{\sig}$. 

We show that $D^{\sharp}_{\sig,\val}$ in Part III is identified with the valuative space $S_{\val}$ in Section 3.1 associated to $S:=D^{\sharp}_{\sig}$ with this log structure (with $K=\C$). 

In \ref{rvtoric}, let  $N= \{x\in \sig_\R\;|\; \exp(x)\in \Gamma\;\text{in}\;G_\R\}$, let $L= \Hom(N, \Z)$, and regard $\sig$ as a cone in $N_\R:=\R\otimes N$. Let $\Sig$ be the fan  of all faces of $\sig$, and denote $|\toric|(\Sig)$ by $|\toric|_{\sig}$. Then we have a commutative diagram
$$\begin{matrix}  D_{\sig,\val}^{\sharp}  &\leftarrow & E^{\sharp}_{\sig,\val} &\overset{\subset}\to & |\toric|_{\sig,\val}\times \Dc\\ \downarrow && \downarrow && \downarrow \\
D^{\sharp}_{\sig} &\leftarrow & E^{\sharp}_{\sig} & \overset{\subset}\to & |\toric|_{\sig}\times \Dc
\end{matrix}$$
where the squares are cartesian, 
 $E^{\sharp}_{\sig}$ is a $\sig_\R$-torsor over $D^{\sharp}_{\sig}$ and $E^{\sharp}_{\sig,\val}$ is a $\sig_\R$-torsor over $D^{\sharp}_{\sig,\val}$ for certain natural actions of $\sig_\R$ on $E^{\sharp}_{\sig}$ and $E^{\sharp}_{\sig,\val}$, and the pull-back of the log structure of $|D^{\sharp}_{\sig}|$ (\ref{abslog}) on $E^{\sharp}_{\sig}$ coincides with the pull-back of the canonical log structure of $|\toric|_{\sig}$. In the upper row, the space in the middle and the space on the right are the valuative  spaces associated to their lower spaces, respectively. Hence the valuative space associated to $D^{\sharp}_{\sig}$  coincides with the quotient $D^{\sharp}_{\sig,\val}$ of $E^{\sharp}_{\sig,\val}$ by $\sig_\R$, that is, $D^{\sharp}_{\sig,\val}$.

Here the problem of the dependence of the  log structure of $S=D^{\sharp}_{\sig}$ on $\Gamma$ (\ref{trouble}) does not affect for the following reason. For another choice $\Gamma'$ of $\Gamma$ such that $\Gamma'\subset \Gamma$, the identity map of $S$ is a morphism from $S$ with the log structure given by $\Gamma'$ to $S$ with the log structure given by $\Gamma$, and this morphism has the Kummer property \ref{timesval} of log structure. Hence the associated valuative space is independent of the choice of $\Gamma$. 

\end{sbpara}

\begin{sbpara}
By \ref{valD} and by Section \ref{ss:valsp} and Section \ref{ra2}, we have a proper surjective map $D^{\sharp}_{\sig,\val}\to D^{\sharp}_{\sig,[:]}$, and this glues to a proper surjective map $D^{\sharp}_{\Sig,\val}\to D^{\sharp}_{\Sig,[:]}$.

\end{sbpara}

\begin{sbpara}\label{rvalex2} {\bf Example.} We describe the differences of the topologies of $D^{\sharp}_{\Sig,[:]}$, $D^{\sharp}_{\Sig,\val}$ and $D^{\sharp}_{\Sig,[\val]}$.

Let $N_1, N_2\in \frak g_{\Q}$ and assume $N_1N_2=N_2N_1$ and that $N_1$ and $N_2$  are nilpotent and linearly independent over $\Q$. Let $F\in \Dc$ and assume that $(N_1, N_2, F)$ generates a nilpotent orbit in the sense of Part III, 2.2.2. Let $\Sig$ be the fan of all faces of the cone $\R_{\geq 0} N_1+\R_{\geq 0}N_2$. When $y_1, y_2\in\R$ tend to $ \infty$, $\exp(iy_1N_1+iy_2N_2)F$ converges in $D^{\sharp}_{\Sig}$.

(1)  Fix a constant $a\in \R$. When $y\to \infty$, 
 $\exp(iyN_1+i(y+a)N_2)F$ converges in $D^{\sharp}_{\Sig,\val}$, 
$D^{\sharp}_{\Sig,[\val]}$, $D^{\sharp}_{\Sig,[:]}$. 
The limit in $D^{\sharp}_{\Sig,[\val]}$ is independent of $a$ and hence the limit in $D^{\sharp}_{\Sig,[:]}$ is independent of $a$, but the limit in $D^{\sharp}_{\Sig,\val}$ depends on $a$. 

(2) Fix a constant $a\in \R$ such that $0<a<1$. Then when   $y\to \infty$, $\exp(iyN_1+iy^aN_2)F$ converges in $D^{\sharp}_{\Sig,\val}$, 
$D^{\sharp}_{\Sig,[\val]}$, $D^{\sharp}_{\Sig,[:]}$. 
The limit in $D^{\sharp}_{\Sig,\val}$ is independent of  $a$ and hence the limit in $D^{\sharp}_{\Sig,[:]}$ is independent of $a$, but the limit in 
$D^{\sharp}_{\Sig,[\val]}$ depends on $a$.

\end{sbpara}

\subsection{CKS maps to $D_{\SL(2)}$ and $D_{\SL(2),\val}$}\label{ss:cks1}

\begin{sbpara}\label{4.5.1} Recall that in Part III, Theorem 3.3.2, we proved that the identity map of $D$ extends uniquely to a continuous map
$$D^{\sharp}_{\Sig,\val}\to D_{\SL(2)}^I.$$
Part III, Section 3.3 is devoted to its proof.
The corresponding result in the pure case is \cite{KU} Theorem 5.4.4 whose full proof is given in [ibid] Chapter 6.

In this section \ref{ss:cks1}, we prove the following related theorems \ref{valper} and \ref{valper2}. 
\end{sbpara}

\begin{sbthm}\label{valper}
$(1)$ The identity map of $D$ extends uniquely to continuous maps
$$D^{\sharp}_{\Sig, [:]}\to D^I_{\SL(2)}, \quad D^{\sharp}_{\Sig,\lval}\to D^I_{\SL(2),\val}.$$ 
These maps respect the log structures on the sheaves of all $\R$-valued continuous functions.
Here we use the new log structures on $D^{\sharp}_{\Sig, [:]}$ in \ref{4.4.3} (cf.\ \ref{rat6}) and the log structure on $D^I_{\SL(2)}$ discussed in Theorem \ref{slis+}.

$(2)$ The CKS map $D_{\Sig,\val}^{\sharp}\to D^I_{\SL(2)}$ defined in Part III Theorem 3.3.2 coincides with the composition $D_{\Sig,\val}^{\sharp}\to D^{\sharp}_{\Sig,[:]}\to D^I_{\SL(2)}$.
\end{sbthm}

\begin{sbpara}\label{Dnilp}

 Let $D_{\nilp}$ be the set of $(N_1, \dots, N_n, F)$, where $n\geq 0$, $N_j\in \fg_\R$ and $F\in \Dc$, which 
 satisfies the following two conditions.

 (i) $(N_1,\dots, N_n, F)$ generates a nilpotent orbit in the sense of Part III, 2.2.2.

(ii)  For any $w\in \Z$ and for any $1\leq j\leq n$, let $W^{(j)}$ be 
the relative monodromy filtration 
of $y_1N_1+\dots +y_jN_j$ relative to $W$ $(W^{(j)}$ exists and does not depend on the choices of $y_j\in \R_{>0}$ by the condition (i)). Then the filtrations $W^{(j)}(\gr^W)$ on $\gr^W$ $(1\leq j\leq n)$ are  $\Q$-rational.

\end{sbpara}

\begin{sbpara}\label{assoc} We review the map $D_{\nilp} \to D_{\SL(2)}$ which sends $(N_1, \dots, N_n, F)\in D_{\nilp}$ to the class of the associated $\SL(2)$-orbit (Part II, 2.4). 

It is the map which sends $(N_1, \dots, N_n, F)$ to the limit of $\exp(\sum_{j=1}^n iy_jN_j)F$ where $y_j\in \R_{>0}$, $y_j/y_{j+1}\to \infty$ $(1\leq j\leq n$, $y_{n+1}$ denotes $1$) in $D_{\SL(2)}^I$. 

This map is also characterized as follows. Recall that an element $(p, Z)$ of  $D_{\SL(2)}$ is determined by the following (i) and (ii).

(i) Whether $(p, Z)$ is an $A$-orbit or a $B$-orbit. 

(ii) $(\Phi, \br)$ where $\Phi$ is the set of weight filtrations on $\gr^W$ associated to $p$ (Part II, 2.5.2 (ii)) and $\br$ is any element of $Z$. 

Let $(p,Z)\in D_{\SL(2)}$ be the image of $(N_1, \dots, N_n, F)$. Then $(p,Z)$ is a $B$-orbit if and only if there is $j$ such that $N_j\neq 0$, $N_k=0$ for $1\leq k< j$, and $\gr^W(N_j)=0$. $\Phi$ is the set of $W^{(j)}(\gr^W)$ for all $j$ such that $\gr^W(N_k) \neq 0$ for some $k\leq j$. $\br$ in the above (ii) is given as follows:

Since $(N_1,\dots, N_n, F)$ generates a nilpotent orbit, $(W^{(n)}, F)$ is an MHS. Let $(W^{(n)}, \hat F_{(n)})$ be $\R$-split MHS associated to it.  
Then $(N_1, \dots, N_{n-1}, \exp(iN_n)\hat F_{(n)})$ generates a nilpotent orbit and hence 
$(W^{(n-1)}, \exp(iN_n)\hat F_{(n)})$ is an MHS. Let $(W^{(n-1)}, \hat F_{(n-1)})$ be the $\R$-split MHS associated to it. Then $(N_1, \dots, N_{n-2}, \exp(iN_{n-1})\hat F_{(n-1})$ generates a nilpotent orbit and 
hence 
$(W^{(n-2)}, \exp(iN_{n-1})\hat F_{(n-1)})$ is an MHS. $\dots$ In this way, we have $\R$-split MHS $(W^{(j)}, \hat F_{(j)})$ for $1\leq j\leq n$ by a downward induction on $j$. 
 (See Part II 2.4.6). 
We obtain $\br\in D$ as  $\br=\exp(iN_k)\hat F_{(k)}$ if $k$ is the minimal $j$ such that $N_j\neq 0$, where in the case $N_j=0$ for all $j$, we define $\br = F$.

\end{sbpara}

\begin{sbpara}\label{Dsig:} Assume $(\Gamma, \Sig)$ is strongly compatible. By \ref{RandR'}, 
$D^{\sharp}_{\Sig, [:]}$ is identified with the set of $(\sig, Z, (\cS^{(j)})_{0\leq j\leq n}, (N_j)_{1\leq j\leq n})$ where $(\sig, Z)\in D^{\sharp}_{\Sig}$, and if $s$ denotes the image of $(\sig,Z)$ in $S:=\Gamma \bs D_{\Sig}$, $\cS^{(j)}$ are faces of $(M_S/\cO_S^\times)_s$ such that $(M_S/\cO_S^\times)_s=\cS^{(0)}\supsetneq \cS^{(1)}\supsetneq \dots \supsetneq \cS^{(n)}=\{1\}$ and $N_j$ is a homomorphism $\cS^{(j-1)}\to \R^{\add}$ such that $N_j(\cS^{(j)})=0$ and $N_j(\cS^{(j-1)}\smallsetminus \cS^{(j)})\subset \R_{>0}$. 

For $s= \text{class}(\sig, Z)\in S=\Gamma \bs D_{\Sig}$, $(M_S/\cO_S^\times)_s$ is canonically isomorphic to $\Hom(\Gamma(\sig), \N)$. 
Hence $\sig$ is identified with $\Hom((M_S/\cO_S^\times)_s, \R^{\add}_{\geq 0})$ and the face $\cS^{(j)}$ of $(M_S/\cO_S^\times)_s$ in the above corresponds to a face $\sig_j$ of $\sig$ consisting of all homomorphisms $(M_S/\cO_S^\times)_s\to \R^{\add}_{\geq 0}$ which kills $\cS^{(j)}$. 

Hence $D^{\sharp}_{\Sig,[:]}$ is identified with the set of $(\sig, Z, (\sig_j)_{0\leq j\leq n}, (N_j)_{1\leq j\leq n})$ where $(\sig,Z)\in D^{\sharp}_{\Sig}$, $\sig_j$ are faces of $\sig$ such that $0=\sig_0\subsetneq \sig_1\subsetneq \dots \subsetneq \sig_n=\sig$, and $N_j$ is an element of $\sig_{j,\R}/\sig_{j-1,\R}$ which belongs to the image of an element of the interior of $\sig_j$.

\end{sbpara}

\begin{sbpara}\label{Dsig:2} Let $(\sig, Z, (\sig_j)_{0\leq j\leq n}, (N_j)_{1\leq j\leq n})\in D^{\sharp}_{\Sig,[:]}$ (\ref{Dsig:}) and let $\tilde N_j$ be an element of the interior of $\sig_j$ whose image in $\sig_{j,\R}/\sig_{j-1,\R}$ coincides with $N_j$. Let $F\in Z$. Then $(\tilde N_1, \dots, \tilde N_n, F)$ generates a nilpotent orbit, as is easily seen. 

\end{sbpara}

\begin{sbthm}\label{valper2} The map $D^{\sharp}_{\Sig, [:]} \to D^I_{\SL(2)}$ (\ref{valper}) sends $(\sig, Z, (\sig_j)_j, (N_j)_j)$ to the image of 
$(\tilde N_1, \dots, \tilde N_n, F)\in D_{\nilp}$ in $D_{\SL(2)}$ (\ref{assoc}). Here $F\in Z$ and $\tilde N_j$ is any element of the interior of $\sig_j$ whose image in $\sig_{j,\R}/\sig_{j-1,\R}$ coincides with $N_j$.

\end{sbthm}

\begin{sblem}\label{dzs} Let $\sig\subset \fg_\R$ be a nilpotent cone, let $F\in \Dc$, and assume that $(\sig, F)$ generates a nilpotent orbit. Let $N\in \sig_{\R}$ and let $F'=\exp(iN)F$. Let $M(\sig, W)$ be the relative monodromy filtration of $\sig$ with respect to $W$. 

(1) $\delta(M(\sig, W), F')= \delta(M(\sig, W), F)+N$  where the last $N$ denotes the homomorphism $\gr^{M(\sig,W)}\to \gr^{M(\sig, W)}$ which is the sum of the maps  $\gr^{M(\sig,W)}_k\to \gr^{M(\sig, W)}_{k-2}$ $(k\in \Z)$ induced by $N$. 

(2) $\zeta(M(\sig, W), F')= \zeta(M(\sig, W), F)$. 

(3) $\spl_{M(\sig, W)}(F')= \spl_{M(\sig,W)}(F)$. 

\end{sblem}

\begin{pf} (1) follows from the definition of $\delta$. 

By (1), (2) follows from the facts that $\delta(M(\sig, W), F)$ and $N$ commute, $N$ is of Hodge type $(-1,-1)$ for $F(\gr^{M(\sig), W)})$, and $\zeta_{-1,-1}=0$ in general.

 (3) follows from (1) and (2). 
\end{pf}

\begin{sbpara}\label{Dsig:3} Let $(\sig, Z, (\sig_j)_{0\leq j\leq n}, (N_j)_{1\leq j\leq n})\in D^{\sharp}_{\Sig,[:]}$, $F\in Z$, $\tilde N_j$ be as in \ref{Dsig:2}. We show that the image of $(\tilde N_1, \dots, \tilde N_n, F)\in D_{\nilp}$ in $D_{\SL(2)}$ is independent of the choices of $\tilde N_j$ and the choice of $F\in Z$. 

We prove that the associated element of $D_{\SL(2)}$  does not depend on the choice of $F\in Z$. If
$F'\in Z$, $F'=\exp(iN)F$ for some $N\in \sig_\R$. Hence by \ref{dzs} (3) applied to $(\sig, F)$ which generates a nilpotent orbit,  $\hat F_{(n)}$ is independent  of the choice.

We prove that the associated element of $D_{\SL(2)}$ does not depend on the choice of a lifting $\tilde N_j$ of $N_j$. If $\tilde N'_j$ is another lifting of $N_j$, $\tilde N'_j= \tilde N_j + R_j$ for some $R_j\in \sig_{j-1,\R}$. By \ref{dzs} (3) applied to $(\sig_{j-1}, \exp(iN_j) \hat F_{(j)})$ which generates a nilpotent orbit,  $\hat F_{(j-1)}$ is independent of the choice by downward induction on $j$. 

\end{sbpara}

\begin{sbpara}\label{Dsig:4} By \ref{Dsig:3}, we have a map $D^{\sharp}_{\Sig,[:]}\to D_{\SL(2)}$ which sends $(\sig, Z, (\sig_j)_{0\leq j\leq n}, (N_j)_{1\leq j\leq n})\in D^{\sharp}_{\Sig,[:]}$ to the image of $(\tilde N_1, \dots, \tilde N_n, F)\in D_{\nilp}$ in $D_{\SL(2)}$. 

Comparing the definition of this map and the definition of the map $\psi: D^{\sharp}_{\Sig,\val}\to D_{\SL(2)}$  (\ref{4.5.1}) given in Part III, 3.3.1, we see that 
 the composition $D^{\sharp}_{\Sig,\val}\to D^{\sharp}_{\Sig,[:]}\to D_{\SL(2)}$ coincides with $\psi$. 

\end{sbpara}

\begin{sbpara} We complete the proofs of  \ref{valper}  (1) and \ref{valper2}.

This map $D^{\sharp}_{\Sig,[:]}\to D_{\SL(2)}$ in \ref{Dsig:4} is continuous because $D^{\sharp}_{\Sigma,\val}\to D^{I}_{\SL(2)}$ is continuous and $D^{\sharp}_{\Sigma,\val}\to D^{\sharp}_{\Sigma,[:]}$ is proper and surjective.

 \end{sbpara}

\begin{sbpara}\label{4.5log}

Endow $D^{\sharp}_{\Sigma,[:]}$ with the the new log structure in \ref{rat6} on
the sheaf of all $\R$-valued continuous functions. 
Consider the log structure on $D^I_{\SL(2)}$ in \ref{slis+}.
 We show that the continuous map $D^{\sharp}_{\Sigma,[:]}\to D^I_{\SL(2)}$ respects these log structures. 
  We check this on $E^{\sharp}_{\sig,[:]}$. On the toric component of $E^{\sharp}_{\sig,[:]}$, 
the log structure is generated by $t_j:=(y_{j+1}/y_j)^{1/2}$ ($1\leq j\leq n$, $y_{n+1}$ denotes $1$). Let $\beta$ be a distance to the boundary for $\Phi$, where $\Phi$ is the set of $W(\sigma_j,W)$ and let $\tau: \bG_{m,\R}^n\to \prod_w \Aut_\R(\gr^W_w)$ be the homomorphism whose $w$-component is the 
$\tau$ (\ref{pure2}) of the $\SL(2)$-orbit in $n$ variables associated to $(N_1(\gr^W_w), \dots, N_n(\gr^W_w), F(\gr^W_w))$.   Then $\beta(\exp(\sum_{j=1}^n iy_jN_j)F) = t u$ where $u:=\beta(\tau(t)^{-1}\exp(\sum_{j=1}^n iy_jN_j)F)$ is invertible in the ring of real analytic functions. 
  Hence $D^{\sharp}_{\Sigma,[:]}\to D^I_{\SL(2)}$ respects the log structures. 

\end{sbpara}

\begin{sbpara} 

By \ref{4.5log},  the map  $D^{\sharp}_{\Sigma,[:]} \to D^I_{\SL(2)}$ induces the continuous map  $D^{\sharp}_{\Sig,\lval}\to D^I_{\SL(2),\val}$ of associated valuative spaces.
This proves \ref{valper} (2).

\end{sbpara}

\begin{sbpara} Consequently, in the pure case, we have an amplified fundamental diagram
$$\begin{matrix}
&&&&D^{\sharp}_{\Sigma,\lval}&
\overset{\psi}\to & D_{\SL(2),\val}&\overset{\eta}\to  &D_{\BS,\val}\\
&&&&&&\downarrow &&\downarrow \\
&&&&\downarrow && D^{\BS}_{\SL(2)} &\to &D_{\BS}\\
&&&&&&\downarrow&&\\
 \Gamma \bs D_{\Sig,\val}& \leftarrow & D_{\Sig,\val}^{\sharp}&\to &D^{\sharp}_{\Sig,[:]} & \overset {\psi} \to &D_{\SL(2)}&&\\
\downarrow &&&&\downarrow&&&&\\
 \Gamma \bs D_{\Sig}&&\leftarrow &&D_{\Sig}^{\sharp}&&&&
\end{matrix}$$
  which is commutative and in which the maps respect the structures of the spaces.

\end{sbpara}

\section{Mild nilpotent orbits and the space $D^{\diamond}_{\SL(2)}$ of $\SL(2)$-orbits}\label{s:dia}

In this Section 5, we consider the spaces of mild nilpotent orbits, and the space $D^{\diamond}_{\SL(2)}$ which is closely related to mild nilpotent orbits.

In Section 5.1, we give main definitions and main results of Section 5. In the rest of Section 5, we give the proofs of the results in Section 5.1. These results in Section 5.1 were obtained in our joint efforts with Spencer Bloch.

\subsection{Mild nilpotent orbits and the space $D^{\diamond}_{\SL(2)}$}

Let $\cL= W_{-2}\End_\R(\gr^W)$ be as in \ref{cL(F)}.
 
 \begin{sbpara}\label{cD} Let $D_{\nilp}^{\mild}$ be the subset of $D_{\nilp}$ (\ref{Dnilp}) consisting of all elements $(N_1,\dots, N_n, F)$ satisfying the following condition.

 For any $y_j\geq 0$ $(1\leq j\leq n)$, there is a splitting (which may depend on $(y_j)_j$) of $W$ which is compatible with $\sum_{j=1}^n y_jN_j$. 
  \end{sbpara}
    
    We have the following \lq\lq\SL(2)-orbit theorem for mild degeneration''.
  
\begin{sbthm}\label{mnilp1}
 Let $(N_1, \dots, N_n, F)\in D_{\nilp}^{\mild}$. 
 
 \medskip
 
(1) If $y_j/y_{j+1}\to \infty$ $(1\leq j\leq n$, $y_{n+1}$ denotes $1$),  $\delta_W(\exp(\sum_{j=1}^n iy_jN_j)F)$ converges in $\cL$. Moreover, there are $a_m \in \cL$ for $m\in \N^n$ and $\varepsilon\in \R_{>0}$ satisfying the following (i) and (ii).

\medskip

(i)  $\sum_{m \in \N^n} \; (\prod_{j=1}^n \;x_j^{m(j)})a_m$ absolutely converges for $x_j\in \R$, $|t_j|<\varepsilon$ ($1\leq j \leq n$).

(ii)  For $y_j\in \R_{>0}$ ($1\leq j\leq n$) such that $t_j:=(y_{j+1}/y_j)^{1/2}<\varepsilon$ ($1\leq j\leq n$, $y_{n+1}$ denotes $1$), we have $\exp(\sum_{j=1}^n iy_jN_j)F\in D$ and   
$$\delta_W 
(\exp
 (\sum_{j=1}^n iy_jN_j)F) =
 \sum_{m\in \N^n}\;
 (\prod_{j=1}^n 
 t_j^{m(j)}) a_m.$$

\medskip 
(2) Let $\tau^{\star}: \bG^n_{m,\R}\to G(\gr^W)$ be the homomorphism whose $G_\R(\gr^W_w)$ component is the $\tau^{\star}$ (\ref{pure2}) of the $\SL(2)$-orbit in $n$ variables on $\gr^W_w$ associated to $(N_1(\gr^W_w), \dots, N_n(\gr^W_w), F(\gr^W_w))$. Then 
 there are $a_m \in \cL$ for $m\in \N^n$ and $\varepsilon\in \R_{>0}$ satisfying the above condition  (i) and the modification of the above condition (ii) by replacing $\delta_W 
(\exp
 (\sum_{j=1}^n iy_jN_j)F)$ with $\Ad(\tau^{\star}(t))^{-1}\delta_W 
(\exp
 (\sum_{j=1}^n iy_jN_j)F)$ where $t=(t_1,\dots, t_n), \; t_j=(y_{j+1}/y_j)^{1/2}$.

 \medskip

(3) If $y_j/y_{j+1}\to \infty$ $(1\leq j\leq n$, $y_{n+1}$ denotes $1$),  $\exp(\sum_{j=1}^n iy_jN_j)F$ converges in $D_{\SL(2)}^{\star,\mild}$.

\end{sbthm}

In fact, (3) follows from (2) by the definition of the structure of $D^{\star}_{\SL(2)}$ as an object of $\cB'_\R(\log)$ given in \ref{2.3.11}.

\begin{sbpara}\label{5.1.3}
By \ref{mnilp1}, we have maps 
$$D_{\nilp}^{\mild} \to \cL, \quad D_{\nilp}^{\mild} \to D^{\star,\mild}_{\SL(2)}$$
by taking the limit of the convergence in \ref{mnilp1}. 
\end{sbpara}
 
\begin{sbpara}\label{mild}
 We define the {\it mild part $D_{\Sigma}^{\mild}$ of the set of nilpotent orbits $D_{\Sigma}$} as the part of points $(\sigma,Z)$ which satisfy the following condition:
 
(C)  For each $N$ in the cone $\sigma$, there is a splitting of $W$ (which can depend on $N$) which is compatible with $N$. 

For the other spaces of nilpotent orbits $D_{\Sigma}^{\sharp}$, $D_{\Sigma,[:]}^{\sharp}$, $D_{\Sigma,[\val]}^{\sharp}$ etc.\ we define their mild parts $D_{\Sigma}^{\sharp,\mild}$, $D_{\Sigma,[:]}^{\sharp,\mild}$, $D_{\Sigma,[\val]}^{\sharp,\mild}$ etc.\ as the inverse images of $D_{\Sigma}^{\mild}$.

\end{sbpara}
 
\begin{sbpara} In the above definition \ref{mild} of the mildness, the following stronger condition (C${}'$) need not be satisfied.
 
 $({\rm C}')$ There is a splitting of $W$ which is compatible with any element $N$ of the cone $\sigma$. 
\end{sbpara}
 
\begin{sbpara}\label{ExII}
 For example, in the case of Example II in Part I and Part II (the case of $0\to H^1(E)(1)\to * \to \Z\to 0$, where $E$ varies over elliptic curves), we had a nilpotent orbit of rank $2$, and that is a mild degeneration in the sense of \ref{mild} (that is, it satisfies (C)) but it does not satisfy $({\rm C}')$.
 \end{sbpara}

\begin{sbthm}\label{Lthm} 
 (1) There is a unique continuous map $D^{\sharp,\mild}_{\Sig, [:]}\to \cL$ which extends the map $D\to \cL\;;\;x\mapsto \delta_W(x)$.

 (2) There is a unique continuous map $D^{\sharp,\mild}_{\Sig, [:]}\to D^{\star,\mild}_{\SL(2)}$ which extends the identity map of $D$.

(3) The map in (1) (resp.\ (2)) sends $(\sig, Z, (\sig_j)_j, (N_j)_j) \in D^{\sharp,\mild}_{\Sig,[:]}$ (\ref{Dsig:}) to the image of $(\tilde N_1, \dots, \tilde N_n, F)\in D^{\mild}_{\nilp}$ in $\cL$ (resp.\ $D^{\star,\mild}_{\SL(2)}$) under the map in \ref{5.1.3}. Here $\tilde N_j$ is as in \ref{valper2} and $F$ is any element of $Z$. 
\end{sbthm}

(1) of \ref{Lthm} shows 
the convergence of Beilinson regulators in a family with mild degeneration. See 
 Section 7.2.

\begin{sbpara}\label{diamond}

We define a space $D^{\diamond}_{\SL(2)}$.

Let $D^{\diamond}_{\SL(2)}$  be the subset of $D^{\star,\mild}_{\SL(2)}\times \cL$ consisting of all elements $(p, Z, \delta)$ $((p,Z)\in D^{\star,\mild}_{\SL(2)}$ with $p\in D_{\SL(2)}(\gr^W)^{\sim}$ and $Z\sub D$
(\ref{slmap}), $\delta\in \cL)$ satisfying the following conditions (i) and (ii). 

\medskip

(i) Let $n$ be the rank of $p$ and let ${\bold 0}:=(0,\dots,0) \in \Z^n$. Then $\delta$ is of $\Ad(\tau^{\star}_p)$-weight $\leq \bold 0$.

(ii) For any $F\in Z$, $\delta_W(F)$ coincides with the component of $\delta$ of $\Ad(\tau^{\star}_p)$-weight $\bold 0$. 

\medskip

We define the structure of 
  $D^{\diamond}_{\SL(2)}$ as an object of $\cB'_{\R}(\log)$ by regarding $D^{\diamond}_{\SL(2)}$ (resp.\ $D^{\star,\mild}_{\SL(2)}\times \cL$) as $Y$ (resp.\ $X$)  in \ref{embstr}.
  
   We have the evident morphism $$D^{\diamond}_{\SL(2)}\to D^{\star,\mild}_{\SL(2)}\;; \;
(p,Z,\delta) \mapsto (p,Z)$$ of $\cB'_\R(\log)$.

\end{sbpara}

\begin{sbpara}

Via the map  $$D\to D^{\star,\mild}_{\SL(2)}\times \cL\;;\;F\mapsto (F, \delta_W(F)),$$
we regard $D$ as a subset of $D^{\diamond}_{\SL(2)}$.

\end{sbpara}

\begin{sbthm}\label{diathm} $(1)$  Let $D_{\nilp}^{\mild}\to 
D^{\star,\mild}_{\SL(2)}\times\cL$ be the map which sends $(N_1, \dots, N_n, F)\in D_{\nilp}^{\mild}$ to the limit of $(F_y, \delta_W(F_y))$ where $y=(y_j)_{1\leq j\leq n} \in \R_{>0}^n$, $F_y:=\exp(\sum_{j=1}^n iy_jN_j)F$, $y_j/y_{j+1}\to \infty$ $(1\leq j\leq n$, $y_{n+1}$ denotes $1$). Then, the image of this map 
 is contained in $D^{\diamond}_{\SL(2)}$.

$(2)$ There is a unique continuous map $D^{\sharp,\mild}_{\Sig, [:]}\to D_{\SL(2)}^{\diamond}$ which extends the identity map of $D$.

$(3)$ There is a unique continuous map $D^{\sharp,\mild}_{\Sig, [\val]}\to D_{\SL(2),\val}^{\diamond}$ which extends the identity map of $D$.
\end{sbthm}

 \begin{sbprop}\label{diatostar} (1) The map $D^{\diamond}_{\SL(2)}\to D_{\SL(2)}(\gr^W)^{\sim} \times \spl(W) \times \cL$ induced by $D^{\star,\mild}_{\SL(2)}\to D_{\SL(2)}(\gr^W)^{\sim} \times \spl(W)$ is injective, and the image of this map consists of all elements $(p,s,\delta)$ satisfying the following conditions (i) and (ii).
 
 (i) $\delta$ is of $\Ad(\tau^{\star}_p)$-weight $\leq \bold0$.
 
 (ii) Let $(\rho_w,\varphi_w)_w$ be an $\SL(2)$-orbit on $\gr^W$ which represents $p$. Then the component of $\delta$ of $\Ad(\tau^{\star}_p)$-weight $\bold 0$ is of Hodge type $(\leq -1,\leq -1)$ with respect to $(\varphi_w(i,\dots,i))_w$. 

(2) If $(p, Z, \delta)\in D^{\diamond}_{\SL(2)}$ and if $(p,s, \delta)$ is its image in $D_{\SL(2)}(\gr^W)^{\sim} \times \spl(W)\times \cL$, $Z$ is recovered from $(p,s,\delta)$ as follows.
   Under the embedding
 $D\to D(\gr^W) \times \spl(W) \times \cL$ in \ref{grsd}, the image of  $Z\subset D$ coincides with $(Z(p), s, \delta_{\bold 0})\subset D(\gr^W)\times \spl(W)\times \cL$. Here $\delta_{\bold 0}$ denotes the component of $\delta$ of $\Ad(\tau_p^{\star})$-weight $\bold 0$.

 \end{sbprop}

  \begin{sbpara}\label{weakstar} 
  By the {\it weak topology of $D^{\diamond}_{\SL(2)}$}, we mean the topology of $D^{\diamond}_{\SL(2)}$ as a subspace of $D_{\SL(2)}(\gr^W)^{\sim}\times \spl(W) \times \cL$. We denote  the topological space $D^{\diamond}_{\SL(2)}$ endowed with the weak topology by 
  $D^{\diamond,\text{weak}}_{\SL(2)}$.
This weak topology need not coincide with the  topology defined in \ref{diamond}.  See \ref{noIIstar}.
\end{sbpara}

\begin{sbrem} (1) Unlike other spaces of $\SL(2)$-orbits $(D^{\star}_{\SL(2)}$, $D^I_{\SL(2)}$, $D^{II}_{\SL(2)}$, ...), $D$ is not necessarily dense in $D^{\diamond}_{\SL(2)}$ (even for the weak topology). 

(2) The authors believe that $D^{\diamond}_{\SL(2)}$ belongs to $\cB'_\R(\log)^+$ and that this can be proved by using the methods in Section 2.7, but they have not yet proved it. 
\end{sbrem}

   \begin{sbpara} The above results show that we have commutative diagrams
     $$\begin{matrix} D^{\sharp,\mild}_{\Sig, [:]} &\to & D_{\SL(2)}^{\diamond} &\to & D_{\SL(2)}^{\star,\mild}
     \\
  \cap &&&& \downarrow\\
  D^{\sharp}_{\Sig,[:]} &&\to && D^{II}_{\SL(2)}
   \end{matrix}\quad \quad\quad 
  \begin{matrix}D^{\sharp,\mild}_{\Sig, [\val]} &\to & D_{\SL(2),\val}^{\diamond} &\to & D_{\SL(2),\val}^{\star,\mild}\\
  \cap &&&& \downarrow\\
   D^{\sharp}_{\Sig,[\val]} &&\to && D^{II}_{\SL(2),\val}. \end{matrix}$$

The rest of Section 5 is devoted to the proofs of the above results.

\end{sbpara}

\subsection{Preparations on pure $\SL(2)$-orbits}

We review pure $\SL(2)$-orbits in one variable more.

\begin{sbpara}

Assume that we are in the pure case of weight $w$, and assume that we are given  an $\SL(2)$-orbit $(\rho, \varphi)$ in one variable.

    Let $$N, N^+\in {\frak g}_\R$$
    be as follows. Let $\rho_* : {\frak s}{\frak l}(2, \R) \to {\frak g}_\R$ be the Lie algebra homomorphism induced by $\rho$. Then $N$ (resp.\ $N^+$) is the image of 
    $\begin{pmatrix} 0&1\\ 0&0\end{pmatrix}$ (resp.\ $\begin{pmatrix} 0&0\\ 1&0\end{pmatrix}$) in  ${\frak s}{\frak l}(2,\R)$. 
    
    \end{sbpara}

    \begin{sbpara}\label{pdec}
    
     We have a direct sum decomposition $$H_{0,\R}=\bigoplus_{k,r\geq 0}\;  H_{0,\R, (k,r)}$$
     defined as follows. Let $Z=  \Ker(N: H_{0,\R}\to H_{0,\R})$. Then $Z=\bigoplus_{k\geq 0} Z_{(-k)}$ where $Z_{(-k)}$ is the part of $Z$ of $\tau^{\star}$-weight $-k$. 
     Let       $$H_{0,\R, (k,r)}:=(N^+)^r Z_{(-k)}.$$
In particular, $Z_{(-k)}=H_{0,\R,(k,0)}$. 
\end{sbpara}

\begin{sbpara}\label{pdec2}
We have 

\smallskip

(1) Elements of $H_{0,\R, (k,r)}$ have $\tau^{\star}$-weight $2r-k$.

\smallskip

(2) For each $k\geq 0$,  $H_{0,\R, (k,\bullet)}:=\bigoplus_r H_{0,\R, (k,r)}$ is stable under the action of $\SL(2,\R)$ by $\rho$. As a representation of $\SL(2,\R)$, we have a unique isomorphism
$$H_{0,\R, (k, \bullet)}\cong \text{Sym}^k(A) \otimes Z_{(-k)}$$
where $A=\R^2=\R e_1\oplus \R e_2$,  on which $\SL(2,\R)$ acts via the natural action on $A$ and the trivial action on $Z_{(-k)}$, which sends $v\in Z_{(-k)}$ on the left-hand side to $e_1^k\otimes v\in \text{Sym}^k(A)\otimes Z_{(-k)}$ on the right-hand side.

\smallskip

(3) For $e\geq 0$, the kernel of $N^e: H_{0,\R}\to H_{0,\R}$ coincides with the direct sum of $H_{0,\R,(k,r)}$ for $k,r\geq 0$  such that $r < e$.

\smallskip

(4) The filtration $\varphi(0)$ is the direct sum of its restrictions $\varphi(0)_{(k,r)}$ to $H_{0, \C, (k,r)}$ for all $(k,r)$. $H_{0,\R,(k,r)}$ with Hodge filtration $\varphi(0)_{(k,r)}$ on $H_{0,\C, (k,r)}$ is an $\R$-Hodge structure of weight $w+2r-k$.

\smallskip

(5) For any $z\in {\mathbb P}^1(\C)$, the filtration $\varphi(z)$ on $H_{0,\C}$ is the direct sum of its restrictions $\varphi(z)_{(k,\bullet)}$ to $H_{0,\C,(k,\bullet)}$ for $k\geq 0$.

The filtration $\varphi(z)_{(k,\bullet)}$ is described as follows. In the isomorphism in (2), it is given by $\varphi(0)_{(k,0)}$ on $Z_{(-k),\C}$ and the filtration on $A_{\C}$ whose $F^0$ is $A_{\C}$, whose $F^2$ is $0$, and whose $F^1$ is $\C \cdot (ze_1+e_2)$ if $z\in \C$, and is $\C e_1$ if $z=\infty$.

  \end{sbpara}

   \subsection{More preparations on $\SL(2)$-orbits}

   \begin{sbpara}\label{presl2} 
Assume that we are given an $\SL(2)$-orbit $(\rho_w,\varphi_w)$ on $\gr^W_w$ in one variable for each $w\in \Z$. Let
   $$E=W_0\End_\R(\gr^W)= \bigoplus_{w\leq 0} E_w, \quad E_w=\bigoplus_{a\in \Z}\Hom_\R(\gr^W_a, \gr^W_{a+w}).$$ We apply our preparations in Section 5.2 to the $\SL(2)$-orbit in one variable of pure weight $w$ induced on each $E_w$ 
 by $(\rho_a, \varphi_a)$ and $(\rho_{a+w}, \varphi_{a+w})$ $(a\in \Z)$. By \ref{pdec}, we have a direct sum decomposition
$$E_w= \bigoplus_{k,r\geq 0} \; E_{w,(k,r)}.$$ 
\end{sbpara}
    
    \begin{sblem}  $E_{w, (k, r)}E_{w',(k',r')}\subset \bigoplus_{k'',r''} E_{w+w', (k'', r'')}$, 
where $(k'',r'')$ ranges over all elements of $\N\times \N$ such that $r''\leq r+r'$ and $k''-2r''=(k+k')-2(r+r')$. 
    
    \end{sblem}
    
    \begin{pf} This follows from (1) and (3) of \ref{pdec2}. 
    \end{pf}

  \begin{sbpara}\label{frakA} 
  Let $\R\{\{t\}\}$ be the ring of power series in $t$ over $\R$ which absolutely converge when $|t|$ is small.
  We define subrings   $\frak A_0$, $\frak A$, $\frak B_0$, $\frak B$ of $\R\{\{t\}\}\otimes_{\R} E$ as follows. 
$$\frak A_0 =E_{\bullet,(\bullet,0)}
  \subset \frak A=\sum_{r\geq 0}  t^{2r}\R\{\{t^2\}\} \otimes_{\R} E_{\bullet,(\bullet,r)},$$
$$ \frak B_0=\sum_{k\geq 0} t^k E_{\bullet,(k,0)}
\subset \frak B=\sum_{k\geq 0}  t^k \R\{\{t^2\}\}\otimes_{\R} E_{\bullet,(k,\bullet)}.$$
  For $w\leq 0$, define the two-sided ideals of these rings as
 $$W_w\frak A_0 =W_wE_{\bullet,(\bullet,0)}
  \subset W_w\frak A=\sum_{r\geq 0}  t^{2r}\R\{\{t^2\}\} \otimes_{\R} W_wE_{\bullet,(\bullet,r)},$$
$$ W_w\frak B_0=\sum_{k\geq 0} t^k W_wE_{\bullet,(k,0)}
\subset W_w\frak B=\sum_{k\geq 0} t^k \R\{\{t^2\}\}\otimes_{\R} W_wE_{\bullet,(k,\bullet)}.$$
  
  \end{sbpara}   
  
  \begin{sblem}\label{flem}
We have, for $w\leq 0$,
  $$\Ad(\tau^{\star}(t)) W_w\frak B= W_w\frak A, \quad \Ad(\tau^{\star}(t))W_w\frak B_0= W_w\frak A_0.$$
 
  \end{sblem} 
  
  These are direct consequences from the definitions in \ref{frakA}.

   We will apply the following \ref{Qalg2} in \ref{pfnm1} (resp. in the proof of \ref{r1eq}) by taking  $A= \C\otimes_\R \frak B$ (resp. $A=W_0\End_{\C}(H_{\C})$).

  \begin{sbpara}\label{Qalg1}
  Let $A$ be a $\Q$-algebra. 
  
  For a nilpotent ideal $I$ of $A$, we have
  bijections
  $$\exp: I \to 1+I, \quad \log : 1 +I \to I, \quad \exp(x)= \sum_{n=0}^{\infty} \frac{x^n}{n!}, \;\;\log(1-x) = \sum_{n=1}^{\infty} \frac{x^n}{n}$$
  (these are finite sums, for $x\in I$ are nilpotent) which are the inverse of each other. 

Let
  $I^{(r)}$ $(r\geq 1)$ be two-sided ideals of $A$ such that $I^{(1)}\supset I^{(2)}\supset I^{(3)}\supset \dots $, $I^{(r)}I^{(s)}\subset I^{(r+s)}$ for any $r,s\geq 1$,  and  $I^{(r)}=0$ for $r\gg 1$. Let $I=I^{(1)}$. Then $I$ is a nilpotent two-sided ideal. 
    \end{sbpara}
  
  \begin{sblem}\label{Qalg2} Let the notation be as in \ref{Qalg1}. Let
  $M_j$ $(1\leq j\leq m)$ be $\Q$-submodules of $I$ such that 
  $$I^{(r)}= \bigoplus_{j=1}^m\;  (M_j\cap I^{(r)})$$ for any $r\geq 1$. 
  Then if $x\in I$, there is a unique family $(x_j)_{1\leq j\leq m}$ of elements $x_j$ of $M_j$ such that 
  $$\exp(x)= \exp(x_1)\dots \exp(x_m).$$
  \end{sblem}

 \begin{pf}
 Easy induction on $r$ such that $I^{(r)}=0$. 
 \end{pf}

\subsection{Proof of Theorem \ref{mnilp1}}\label{ss:cks2}

\begin{sbpara}\label{Sch}
We first prove Theorem \ref{mnilp1} in the case $n=1$. 
We use the following part of the SL(2)-orbit theorem in one variable of Schmid (\cite{Scw}). 

Assume that we are in the pure case. 
 Then for $y\gg 0$, we have
$$\exp(iyN)F= \exp(g(y)){\tau}(y^{-1/2})\br\quad \text{with}\;\; \br=\exp(iN)\hat F$$
for some convergent power series $g(y)=\sum_{m=0}^{\infty} y^{-m}a_m$ in $y^{-1}$  with $a_m\in \End_{\R}(H_{0,\R})$ such that $a_0=0$ and 
such that $N^{m+1}a_m=0$ for any $m$. 

\end{sbpara}

\begin{sbprop}\label{5.4.1} Let $(N,F)\in D^{\mild}_{\nilp}$ (\ref{cD}) with one $N$. Let $(W^{(1)},\hat F)$ be the $\R$-split MHS associated to the MHS $(W^{(1)},F)$ and let $\br := \exp(iN)\hat F$. Then $(W,\br)$ is an $\R$-split MHS and the splitting $\spl_W(\br)$ of $W$ is compatible with $N$.

\end{sbprop}

\begin{pf}
This follows from \cite{BP} Lemma 2.2. 
\end{pf}

\begin{sbpara}\label{3.6.2} Let $(N, F)\in D^{\mild}_{\nilp}$ with one $N$. Let $\br=\exp(iN)\hat F$ as in \ref{5.4.1} and let $s=\spl_W(\br):\gr^W \overset{\cong}\to H_{0,\R}$, $s^{(1)}=\spl_{W^{(1)}}(F)=\spl_{W^{(1)}}(\hat F): \gr^{W^{(1)}}\overset{\cong}\to H_{0,\R}$.  By \ref{5.4.1}, $N$ is of weight $0$ for $s$. 
Let $\tau^{\star}: \bG_m \to \Aut(\gr^W)$ be the homomorphism associated to the $\SL(2)$-orbit on $\gr^W$ in one variable associated to $(N,F)(\gr^W)$. (In the case $N(\gr^W)=0$, $\tau^{\star}$ is defined to be the trivial homomorphism.)  Let $y\in \R_{>0}$ and let $t= y^{-1/2}$. 

For the proof of the case $n=1$ of \ref{mnilp1} (1) and (2), it is sufficient to prove that
$\delta_W(\exp(iyN)F)$ and $\Ad(\tau^{\star}(t))^{-1}\delta_W(\exp(iyN)F)$ converge in $\cL$ when $y\to \infty$.
We prove it.

Note  that the actions of $\tau(t)$ and $\tau^{\star}(t)$ on $D(\gr^W)$ are the same. 
\end{sbpara}

\begin{sbpara}\label{pfnm1} Let the notation be as in \ref{3.6.2}.

For $y\gg 0$, let  $g_w(y)$ for each $w\in \Z$ be as in the above result of Schmid in \ref{Sch} for $(N, F)(\gr^W_w)$, and let $g(y)=\bigoplus_w g_w(y) \in E= W_0\End_\R(\gr^W)$. 
By the above result of Schmid, we have 

\medskip

(1) $\exp(iyN(\gr^W))F(\gr^W) = \exp(g(y)){\tau}^{\star}(t) \br (\gr^W),  \quad g(y), \exp(g(y))\in \frak A$

\medskip
\noindent
 where $N(\gr^W)$ is the map $\gr^W\to \gr^W$ induced by $N$ and $\frak A$ is as in  \ref{frakA}. Let $h(y) = \Ad(\tau^{\star}(t))^{-1}g(y)$. Then 
 
 \medskip
 
 (2) $h(y)\in \frak B$

 \medskip
 \noindent
 by \ref{flem}.

Let $\delta^{(1)}=\delta_{W^{(1)}}(F)$ and let $\zeta^{(1)}$ be the corresponding $\zeta$ (\ref{II,1.2.3}), so that 

\medskip

(3) $F=s^{(1)} \exp(-\zeta^{(1)})\exp(i\delta^{(1)})(s^{(1)})^{-1}\hat F.$

\medskip

  Write $s^{(1)} \exp(-\zeta^{(1)})\exp(i\delta^{(1)})(s^{(1)})^{-1}= \exp(\alpha)\exp(\beta)$ where 
$\alpha, \beta\in W_0\End_\C(H_{0,\C})\cap W_{-2}^{(1)}\End_\C(H_{0,\C})$, $\alpha$ is of 
$s$-weight $\leq -1$ and $\beta$ is of $s$-weight $0$. By (3), we have

\medskip
(4)  $F(\gr^W) = \exp(\beta(\gr^W))\hat F(\gr^W)$

\medskip
\noindent
where $\beta(\gr^W)$ is the map $\gr^W\to \gr^W$ induced by $\beta$. We have 

\medskip
(5) $\exp(iyN)\exp(\beta)\hat F= s\exp(iyN(\gr^W))\exp(\beta(\gr^W)) \hat F(\gr^W)= s\exp(iyN(\gr^W))F(\gr^W)= s\exp(g(y))\tau^{\star}(t)\br (\gr^W)$

\medskip
\noindent
where the first $=$ follows from $N=sN(\gr^W)s^{-1}$ (\ref{5.4.1}), $\beta=s\beta(\gr^W)s^{-1}$ and $\hat F=s\hat F(\gr^W)$, the second $=$ follows from (4), and the last $=$ follows from (1). 

Since $s^{(1)}\delta^{(1)}(s^{(1)})^{-1}$ and 
$s^{(1)}\zeta^{(1)} (s^{(1)})^{-1}$ commute with $N$, $\alpha$ and $\beta$ commute with $N$. Hence we have $$s^{-1}\alpha s \in W_{-1}\frak A_0.$$ We have $$\exp(iyN)F= \exp(iyN)\exp(\alpha)\exp(\beta)\hat F = \exp(\alpha) \exp(iyN)\exp(\beta)\hat F$$
$$=\exp(\alpha) s \exp(g(y))\tau^{\star}(t) \br (\gr^W)  = s \exp(g(y))\exp(\alpha(y))\tau^{\star}(t)\br(\gr^W)$$
where $\alpha(y):= \Ad(\exp(g(y)))^{-1}(s^{-1}\alpha s)$. Here the third $=$ follows from (5).

Since $s^{-1}\alpha s\in W_{-1}\frak A_0$ and $\exp(g(y))\in \frak A$, we have $\alpha(y) \in W_{-1}\frak A$. Hence 
$$\Ad(\tau^{\star}(t))^{-1} \alpha(y)  \in W_{-1}\frak B.$$

To apply \ref{Qalg2}, we use the direct sum decomposition
$$\C\otimes_\R W_{-1}\End_\R(\gr^W)= M'_1\oplus M'_2 \oplus M'_3$$
where $M'_1= W_{-1}\End_\R(\gr^W)$, $M'_2$ is the $-1$-eigen space of the complex conjugation acting  on the $(\leq -1, \leq -1)$-Hodge component of $\C\otimes_\R W_{-1}\End_\R(\gr^W)$ with respect to $\br(\gr^W)$, and $M'_3= F^0(\C\otimes_\R W_{-1}\End_\R(\gr^W))$ for the Hodge filtration $\br(\gr^W)$.
 In \ref{Qalg2}, consider the case
$$A= \frak B, \quad I^{(r)}=W_{-r}A, \quad I=I^{(1)},$$
$$M_j = \bigoplus_{k\geq 0} t^k\R\{\{t\}\}\otimes_\R M'_{j,(k, \bullet)} \quad  (j=1,2,3).$$
Then the assumption of \ref{Qalg2} is satisfied. By \ref{Qalg2}, we have
$$\exp(\Ad(\tau^{\star}(t))^{-1}\alpha(y)) = \exp(a(y)) \exp(ib(y)) \exp(c(y))$$
where $a(y) \in M_1$, $ib(y)\in M_2$, $c(y) \in M_3$. Then 
$$\exp(iyN)F = s \exp(g(y)) \tau^{\star}(t) \exp(a(y)) \exp(i b(y)) \br(\gr^W).$$
Hence 
$$\delta_W(\exp(iyN)F)=\Ad(\exp(g(y)) \Ad(\tau^{\star}(t)) b(y) \in \frak A,$$
$$\Ad(\tau^{\star}(t))^{-1}\delta_W(\exp(iyN)F)= \Ad(\exp(h(y)) b(y)  \in \frak B.$$
Hence $\delta_W(\exp(iyN)F)$ and $\Ad(\tau^{\star}(t))^{-1}\delta_W(\exp(iyN)F)$ converge when $y\to \infty$.

\end{sbpara}

\begin{sbpara}\label{5.4.5}

We prove \ref{mnilp1} in general. 

Let $(N_1, \dots, N_n, F)\in D_{\nilp}$.
Let $\tau$ (resp.\ $\tau^{\star}): \bG_{m,\R}^n\to \prod_w \Aut_\R(\gr^W_w)$ be the homomorphism whose $w$-component is the $\tau$ (resp.\ $\tau^{\star}$)  (\ref{pure2}) of the $\SL(2)$-orbit in $n$-variables on $\gr^W_w$ associated to $(N_1(\gr^W_w), \dots, N_n(\gr^W_w), F(\gr^W_w))$. 
Note the action of $\bG_{m,\R}^n$ on $D(\gr^W)$ via $\tau$ and that via $\tau^{\star}$ are the same.

By SL(2)-orbit theorem in $n$ variables (\cite{KNU1}), 
 $$\Ad(\tau(t))^{-1}\delta_W(\exp(\sum_{j=1}^n iy_jN_j)F)\quad  (t= (t_1,\dots, t_n), \;t_j=(y_{j+1}/y_j)^{1/2}, \;y_{n+1}=1)$$ is a convergent series in $t_1, \dots, t_n$. Hence 
   $\delta_W(\exp(\sum_{j=1}^n iy_jN_j)F)$ and \newline 
   $\Ad(\tau^{\star}(t))^{-1}\delta_W(\exp(\sum_{j=1}^n iy_jN_j)F)$ have the shapes of Laurent series
$$\delta_W(\exp(\sum_{j=1}^n iy_jN_j)F)=(\prod_{j=1}^n t_j)^{-r} \cdot \sum_{m\in \N^n} (\prod_{1\leq j\leq n} t_j^{m(j)})a_m,$$
$$\Ad(\tau^{\star}(t))^{-1}\delta_W(\exp(\sum_{j=1}^n iy_jN_j)F)=(\prod_{j=1}^n t_j)^{-s} \cdot \sum_{m\in \N^n} (\prod_{1\leq j\leq n} t_j^{m(j)})b_m$$
for some $r,s\in \N$ and $a_m, b_m\in \cL$  
where the sums $\sum_{m\in \N^n}$ are convergent series.

Now assume $(N_1, \dots, N_n, F)\in D_{\nilp}^{\mild}$. We prove that we can take $r=s=0$ (that is, these series are actually taylor series). It is sufficient to prove that when we fix $j$ and  fix a sufficiently small $t_k>0$ for $k\neq j$, then these series become Taylor series in one variable in $t_j$. 

But in this situation, the first Laurent series becomes $\delta_W(\exp( iy'N')F')$ with $(N', F')\in D_{\nilp}$, where 
$$ y'= t_j^{-2}, \quad N'=t_j^2\sum_{k=1}^j  y_kN_k= \sum_{k=1}^j \;(\prod_{k\leq \ell\leq n, \ell\neq j} t_{\ell}^{-2})N_k,$$ $$F'=\exp(\sum_{k=j+1}^n  iy_k N_k)F=\exp(i\sum_{k=j+1}^n   \;(\prod_{\ell=k}^n t_{\ell}^{-2}) N_k)F.$$

We consider the second Laurent series. 
Let $\tau^{\star}_j: \bG_{m,\R}\to G_\R(\gr^W)$ be the restriction of $\tau^{\star}$ to the $j$-th $\bG_{m,\R}$.
It is sufficient to prove that when $t_k$ for $k\neq j$ are fixed, $\delta(t):=\Ad(\tau^{\star}_j(t_j))^{-1}\delta_W(\exp(\sum_{j=1}^n iy_jN_j)F)$ is a Taylor series in $t_j$. 
Let $\tau^{\star,\prime}_j: \bG_{m,\R}\to G_\R(\gr^W)$ be the $\tau^{\star}$ of the $\SL(2)$-orbit in one variable associated to $(N', F')$ where $N'$ and $F'$ are as above. By the case $n=1$ applied to $(N', F')$,  $\delta'(t):=\Ad(\tau_j^{\star,\prime}(t_j))^{-1}\delta_W(\exp(\sum_{j=1}^n iy_jN_j)F)$ is a Taylor series in $t_j$.
Let $W^{(j)}$ be the relative monodromy filtration $M(N_1+\dots+N_j, W)$. 
By \cite{KNU1} Proposition 4.2, there is a convergent Taylor series $u=\sum_m (\prod_{k=j+1}^n t_k^{m(k)})u_m$ in $t_{j+1}, \dots, t_n$ with $u_m\in W^{(j)}_{-1}\fg_\R$ such that $u_0=0$ and such that
$$\tau_j^{\star,\prime}(t_j)= \exp(u)\tau_j^{\star}(t_j)\exp(-u).$$
 We have
$$\delta(t)= \Ad(\exp(v)\exp(-u))^{-1}\delta'(t)\quad{where}\quad 
v= \Ad(\tau^{\star}_j(t_j))^{-1}u.$$
Since $u_m\in W^{(j)}_{-1}\fg_\R$, $v$ is a Taylor series in $t_j$. Hence $\delta(t)$ is a 
Taylor series in $t_j$.

\end{sbpara}

\begin{sbpara}\label{Coef}

In the mild $\SL(2)$-orbit theorem Theorem \ref{mnilp1} (1) and (2), 
the power series depend real analytically on  $(N_1, \dots, N_n, F)$ in the following sense. Let $A$ be a real analytic manifold and let $A\to \fg_\R\;;\;\alpha \mapsto N_{j, \alpha}$ ($1\leq j \leq n$) and $A\to \Dc\;;\;\alpha\mapsto F_{\alpha}$ be real analytic functions.
Assume that $N_{j,\alpha}$ are nilpotent and commute with each other, and assume that $(N_{1,\alpha}, \dots, N_{n,\alpha},F_{\alpha})$ generates a nilpotent orbit for any $\alpha$. Assume further that for each $1\leq j\leq n$, the relative monodromy filtration $M(N_1+\dots+N_j, W)$ is independent of $\alpha$. Then the $\varepsilon$ in \ref{mnilp1} (1) and (2) can be taken constant locally on $A$, and the coefficients of the power series in \ref{mnilp1} (1) and (2) are real analytic functions on $A$. 

This follows from the corresponding result \cite{KNU1} Proposition 10.8 (\cite{CKS} Remark 4.65 (ii) in the pure case)  for the original $\SL(2)$-orbit theorem and from the above proof in \ref{5.4.5} to reduce the mild $\SL(2)$-orbit theorem to the original one.

\end{sbpara}

\subsection{Proof of Theorem \ref{Lthm}}\label{ss:Lthm}
We prove Theorem \ref{Lthm}.
We first prove

\begin{sbprop}\label{5.5.1} Let $(\sig, Z, (\sig_j)_{0\leq j\leq n}, (N_j)_{1\leq j\leq n})\in D^{\sharp,\mild}_{\Sig,[:]}$ (\ref{Dsig:}). Then for  $\tilde N_j$  as in \ref{Dsig:2} and for $F\in Z$,  the image of $(\tilde N_1, \dots, \tilde N_n, F)\in D^{\mild}_{\nilp}$  in $D^{\star,\mild}_{\SL(2)}\times \cL$ (\ref{mnilp1}) is independent of the choices of $\tilde N_j$ and $F$. 

\end{sbprop}

\begin{pf}
For another choice $(\tilde N'_1,\dots, \tilde N_n', F')$ of $(\tilde N_1,\dots, \tilde N_n,F)$, we have
$\tilde N_1'=\tilde N_1$, $\tilde N'_j=\tilde N_j+R_{j-1}$ for $2\leq j\leq n$ and $F'= \exp(iR_n)F$ for some $R_j\in \sig_{j,\R}$. We have
$$\exp(\sum_{j=1}^n iy_j\tilde N'_j)F'= \exp(\sum_{j=1}^n iy_j(\tilde N_j + (y_{j+1}/y_j)R_j))F$$
($y_{n+1}$ denotes $1$). The limit of this for $y_j/y_{j+1}\to \infty$ coincides with the limit of that for  $R_j=0$ by \ref{Coef}. 
\end{pf}

\begin{sbpara}\label{5.4.2}

 By \ref{5.5.1}, we have a map $D^{\sharp,\mild}_{\Sig,[:]}\to D^{\star,\mild}_{\SL(2)}\times \cL$.

Let $D^{\sharp,\mild}_{\Sig,\val}\to D^{\star,\mild}_{\SL(2)}\times \cL$ be the composition with $D^{\sharp,\mild}_{\Sig,\val}\to D^{\sharp,\mild}_{\Sig,[:]}$ (4.4.6). Since the last map is proper surjective (4.4.6),  \ref{Lthm} is reduced to 

\end{sbpara}

\begin{sbprop}\label{5.4.3}
The map $D^{\sharp, \mild}_{\Sig, \val}\to D^{\star,\mild}_{\SL(2)}\times \cL$ is continuous. 

\end{sbprop}

Just as Part III, Theorem 3.3.2 was reduced to the case $y^*_{\lam,t}= y_{\lam,t}$ of Part III, Proposition 3.3.4 (see \ref{8.3.b}),  Proposition \ref{5.4.3}  is reduced to $(A_0)$ of  the following Proposition \ref{forLthm}. 

The proof of Proposition \ref{forLthm} given below is similar to the proof of Part III, Proposition 3.3.4.

\begin{sbprop}\label{forLthm}
 Let the situation and the assumption be as in Part III, 3.3.3 with $y^*_{\lam,t}=y_{\lam,t}$ there. 
  
 Assume that there is $\varepsilon\in \R_{>0}$ such that for any $(y_s)_{s\in S}\in \R^S$ satisfying the following condition (C), 
 there is a splitting of $W$ (which may depend on $(y_s)_s$) which is compatible with $\sum_{s\in S} y_sN_s$. 
 
 \medskip
 
 (C) If $1\leq j\leq n$, $s\in S_j$, and $y_s\neq 0$, then $y_ty_s^{-1}<\varepsilon$ for any $t \in S_{\geq j+1}$ and $|y_ty_s^{-1}-a_ta_s^{-1}|<\varepsilon$ for any $t\in S_j$. 
 
 \medskip
 Note that $(N_1, \dots, N_n, F)\in D_{\nilp}^{\mild}$ by this assumption. Let $\tau, \tau^{\star}: \bG_{m,\R}^n \to \Aut_\R(H_{0,\R})$ be the homomorphisms given by the $\SL(2)$-orbit in $n$ variables associated to $(N_1, \dots, N_n, F)$. Let 
 $$\delta= \lim 
 \delta_W(\exp(\sum_{j=1}^n iy_jN_j)F), \quad \delta'= \lim \Ad(\tau^{\star}(t))^{-1}\delta_W(\exp(\sum_{j=1}^n iy_jN_j)F)$$ where $y_j/y_{j+1}\to \infty$ ($1\leq j\leq n$, $y_{n+1}=1$) and where $t= (t_1,\dots, t_n)$, $t_j= (y_{j+1}/y_j)^{1/2}$.

 For $1\leq j\leq n+1$, let $e_{\lam,\geq j}:=\exp(\sum_{s\in S_{\geq j}} iy_{\lam,s}N_s) \in G_\C$.

 Then we have the following $(A_j)$ for $0\leq j\leq n$.
 
\medskip

$(A_j)$ for $1\leq j\leq n$:  Let $e\geq 1$. Then if $\lam$ is sufficiently large, there are $F^{(j)}_{\lam}\in \Dc$ satisfying the following (1)--(4). 

\medskip

(1) $y_{\lam,c_j}^ed(F_{\lam}, F^{(j)}_{\lam}) \to 0$.

(2) $((N_s)_{s\in S_{\leq j}}, e_{\lam,\geq j+1}F^{(j)}_{\lam})$ generates a nilpotent orbit. 

(3) 
$\delta_W(\exp(\sum_{s\in S} iy_{\lam,s} N_s)F^{(j)}_{\lam})$ converges to $\delta$.

(4) 
$\Ad(\tau^{\star}(t))^{-1} \delta_W(\exp(\sum_{s\in S} iy_{\lam,s} N_s)F^{(j)}_{\lam})$ with 
$t=(t_1, \dots, t_n)$, $t_j= (y_{\lam, c_{j+1}}/y_{\lam,c_j})^{1/2}$ converges to $\delta'$.

 \medskip
 
 $(A_0)$:  Let $e\geq 1$. Then if $\lam$ is sufficiently large, we have the following (3) and (4).
 
 \medskip

(3) $\delta_W(\exp(\sum_{s\in S} iy_{\lam,s} N_s)F_{\lam})$ converges to $\delta$.

(4) 
$\Ad(\tau^{\star}(t))^{-1} \delta_W(\exp(\sum_{s\in S} iy_{\lam,s} N_s)F_{\lam})$ with 
$t$ as in above $(A_j)$ (4).

\end{sbprop}

\begin{sbpara}
We prove Proposition 
\ref{forLthm} by downward induction on $j$. 

For $1\leq j\leq n$, let $\tau_j$ be the restriction of $\tau$ to the $j$-th factor of $\bG_{m,\R}$ and let 
$\tau_{\geq j} = \prod_{k=j}^n \tau_k((y_{\lam, c_{k+1}}/y_{\lam, c_k})^{1/2})\in G_\R$.

 $(A_n)$ follows from the condition (5) in Part III 3.3.3 for $j=n$ with $y^*_{\lam,t}=y_{\lam,t}$ and from \ref{Coef} 
 
Assume $0\leq j<n$.  We prove $(A_j)$ assuming $(A_{j+1})$. Take sufficiently large integers $e, e', e''\geq 0$. Take $F^{(j+1)}_{\lam}$ as in $(A_{j+1})$ with $e$ replaced by $e+e'+e''$. 
In the case $1\leq j<n$ (resp. $j=0$), let $F^{(j)}_{\lam}$ be $F_{\lam}^*$ in Part III, 3.3.3 (5) with $e$ there replaced by $e+e'+e''$ (resp. let 
$F^{(0)}_{\lam}= F_{\lam}$). 

We have 

\medskip

(5) $y_{\lam, c_{j+1}}^{e+e'+e''}d(F^{(j)}_{\lam}, F^{(j+1)}_{\lam})\to 0$.

\medskip
By Part III, Lemma 3.3.6, 
 $\tau_{\geq j+1}^{-1} e_{\lam,\geq j+1}F^{(j+1)}_{\lam}$ converges. Hence by (5), 
 $\tau_{\geq j+1}^{-1} e_{\lam,\geq j+1}F^{(j)}_{\lam}$ converges and
 we have

\medskip

(6) $y_{\lam,c_{j+1}}^{e+e'} d(\tau_{\geq j+1}^{-1}e_{\lam,\geq j+1}F^{(j)}_{\lam}, \tau_{\geq j+1}^{-1} e_{\lam,\geq j+1}F^{(j+1)}_{\lam})\to 0.$

\medskip

By the mild SL(2)-orbit theorem \ref{mnilp1} for 

\medskip

$((N_s)_{s\in S_{\leq j}}, \tau_{\geq j+1}^{-1}e_{\lam,\geq j+1}F^{(j)}_{\lam})$  and 
$((N_s)_{s\in S_{\leq j}}, \tau_{\geq j+1}^{-1}e_{\lam,\geq j+1}F^{(j+1)}_{\lam})$, 

\medskip
\noindent 
and by \ref{Coef}, 
we have:

\medskip

(7) The four sequences 

 $a_{\lam}:=\delta_W(\tau_{\geq j+1}^{-1} \exp(\sum_{s\in S} iy_{\lam,s}N_s)F^{(j)}_{\lam})$,

$b_{\lam}:=\delta_W(\tau_{\geq j+1}^{-1} \exp(\sum_{s\in S} iy_{\lam,s}N_s)F^{(j+1)}_{\lam})$,

$a'_{\lam}:=\Ad(\tau^{\star}(t))^{-1}\delta_W(\tau_{\geq j+1}^{-1} \exp(\sum_{s\in S} iy_{\lam,s}N_s)F^{(j)}_{\lam})$, 

$b'_{\lam}:=\Ad(\tau^{\star}(t))^{-1}\delta_W(\tau_{\geq j+1}^{-1} \exp(\sum_{s\in S} iy_{\lam,s}N_s)F^{(j+1)}_{\lam})$

\noindent
 converge in $\cL$ and we have  
$$y_{\lam, c_{j+1}}^e(a_{\lam}-b_{\lam})\to 0, \quad 
y_{\lam, c_{j+1}}^e(a'_{\lam}-b'_{\lam})\to 0.$$

By the induction assumption on $j$,  $\exp(\sum_{s\in S} iy_{\lam,s}N_s)F^{(j+1)}_{\lam}$ converges to $\delta$ and
\newline $\Ad(\tau^{\star}(t))^{-1}\exp(\sum_{s\in S} iy_{\lam,s}N_s)F^{(j+1)}_{\lam}$ converges to $\delta'$. Hence by (7), 
$\exp(\sum_{s\in S} iy_{\lam,s}N_s)F^{(j)}_{\lam}$ converges to $\delta$ and  
$\Ad(\tau^{\star}(t))^{-1}\exp(\sum_{s\in S} iy_{\lam,s}N_s)F^{(j)}_{\lam}$ converges to $\delta'$. 

\end{sbpara}

  \subsection{Proofs of other results in Section 5.1}\label{ss:cks0}

  \begin{sblem}\label{5.6.1} Let $x=(N_1, \dots, N_n, F)\in D^{\mild}_{\nilp}$ and  let $p\in D_{\SL(2)}(\gr^W)^{\sim}$ be the image of $x$. Let ${\bold 0}=(0,\dots,0)\in \Z^n$. 
  
  (1) Let $\delta$ be the image of $x$ in $\cL$. Then $\delta$ is of $\Ad(\tau^{\star}_p)$-weight $\leq {\bold 0}$.
  
  (2) Let $\delta'\in \cL$ be the limit of $\Ad(\tau^{\star}_p(t))^{-1}\delta_W(\exp(\sum_{j=1}^n iy_jN_j)F)$ $(t=(t_1,\dots, t_n), t_j=(y_{j+1}/y_j)^{1/2}$, $y_{n+1}$ denotes $1$, $t_j\to 0$). Then $\delta'$ coincides with the component of $\delta$ of $\Ad(\tau^{\star}_p)$-weight $\bold 0$.

  \end{sblem}
  
  \begin{pf} Let $\delta(y):=\delta_W(\exp(\sum_{j=1}^n iy_jN_j)F)$ and let $\delta'(y)= \Ad(\tau^{\star}_p(t))^{-1}\delta_W(\exp(\sum_{j=1}^n iy_jN_j)F))$ where $t_j$ is as above. Then $\delta(y)$ and $\delta'(y)$ are convergent series in $t_1,\dots, t_n$. For $a\in \Z^n$, let $\delta_a$ (resp.\ $\delta'_a$, resp.\ $\delta(y)_a$, resp.\ $\delta'(y)_a$) be the component of $\delta$ (resp.\ $\delta'$, resp.\ $\delta(y)$, resp.\ $\delta'(y)$) of $\Ad(\tau^{\star}_p)$-weight $a$. Then 
  $\delta(y)_a = (\prod_{j=1}^n t_j^{a(j)}) \delta'(y)_a$. Hence $\delta(y)_a$ is divisible by $\prod_{j=1}^n t_j^{\max(a(j), 0)}$. Hence the constant term $\delta_a$ of $\delta(y)_a$ is $0$ unless $a\leq {\bold 0}$. On the other hand, by the reduction to the case of one $N$, we have
  $$\delta'(y)\in \sum_{k\in \N^n} (\prod_{j=1}^n t_j^{k(j)}) \R\{\{t_1,\dots, t_n\}\}\cdot (\bigcap_{j=1}^n W^{(j)}_{k(j)}\cL).$$
  Hence the constant term of $\delta'(y)$ belongs to $\bigcap_{j=1}^n W^{(j)}_0\cL$. That is, $\delta'$ is of $\Ad(\tau^{\star}_p)$-weight $\leq {\bf 0}$. For $a\in \Z^n$ such that $a\leq {\bold 0}$, the constant term $\delta'_a$ of $\delta'(y)_a= \prod_{j=1}^n t_j^{-a(j)}\delta(y)_a$ is $0$ unless $a={\bold 0}$. This argument shows also $\delta_{\bold 0}=\delta'_{\bold 0}$. 
   \end{pf}
   
   \begin{sbpara}\label{5.6.2} We prove \ref{diathm} (1). By \ref{5.6.1}, it is sufficient to prove that the element $\delta'\in \cL$ in \ref{5.6.1} (2) 
 belongs to $\cL(\br)$ where $\br=(\varphi_w(i,\dots,i))_w$. 
   
   Since $\br(y):=\tau_p^{\star}(t)^{-1}\exp(\sum_{j=1}^n iy_jN_j)F(\gr^W)$ converges to $\br$, $\cL(\br(y))$ converges to $\cL(\br)$. Since $\delta'(y)\in \cL(\br(y))$, its limit $\delta'$ belongs to $\cL(\br)$.

   \end{sbpara}
   
   \begin{sbpara}\label{5.6.3} We prove \ref{diathm} (2).
   Let $s\in \spl(W)$ be the image of $x$ (\ref{5.6.1}) in $\spl(W)$. Consider  the element $(p,Z)$ of $D^{\star,\mild}_{\SL(2)}$ (\ref{slmap}) where $Z$ is the subset of $D$ whose image 
   in $D(\gr^W)\times \spl(W)\times \cL$ is $(Z(p), s, \delta_{\bold 0})$. Such element exists uniquely by \ref{5.6.2}. We show that
  $F_y:=\exp(\sum_{j=1}^n iy_jN_j)F$ converges to $(p, Z)$ in $D^{\star,\mild}_{\SL(2)}$. 
  
  Let $\Phi$ be the set $\{W^{(1)}(\gr^W), \dots, W^{(n)}(\gr^W)\}$ of weight
   filtrations on $\gr^W$ associated to $p$. Fix a distance $\beta: D(\gr^W) \to \R_{\geq 0}^n$ to $\Phi$-boundary. Let $D_{\SL(2)}^{\star,\mild}(\Phi)\to D_{\SL(2)}(\gr^W)^{\sim} \times \spl(W) \times \cL$ be 
   the map $\nu_{\alpha,\beta}$ in \ref{staran1} where $\alpha=\tau_p$. Then $\nu_{\alpha,\beta}(p,Z)=(p,s, \Ad(\tau_p(\beta(\br)))^{-1}\delta_{\bold 0})$. Hence it is sufficient to prove 
   that
   $\nu_{\alpha,\beta}(F_y)\in D(\gr^W) \times \spl(W) \times \cL$ converges to $(p, s, \Ad(\tau^{\star}_p(\beta(\br)))^{-1}\delta_{\bold 0})$. It is sufficient to prove that $\Ad(\tau^{\star}_p(\beta(F_y(\gr^W))))^{-1} \delta(y)$ 
   converges to $\Ad(\tau_p^{\star}(\beta(\br)))^{-1}\delta_{\bold 0}$. But this is deduced from the fact that $\beta(F_y(\gr^W))t^{-1}$ converges to $\beta(\br)$.
    \end{sbpara}
  
\begin{sbpara}\label{5.6.4} We prove \ref{diathm} (3). By (2), it is sufficient to prove  
the compatibility of the map $D^{\sharp,\mild}_{\Sig,[:]} \to D^{\diamond}_{\SL(2)}$ with log structures. 

This is reduced to the pure case treated in \ref{4.5log} because the log structure of $D^{\diamond}_{\SL(2)}$ is the inverse image of that of $D_{\SL(2)}(\gr^W)^{\sim}$.
\end{sbpara}

\begin{sbpara} Theorem \ref{diatostar} follows from Lemma \ref{5.6.1} and Theorem \ref{diathm}(1).

\end{sbpara} 

These complete the proofs of the results in Section 5.1.

\section{Complements}\label{s:NSB}

In Section \ref{ss:+3}, we give properties of the extended period domains. 

In Section \ref{ss:cks4}, we show that for nilpotent orbits in one variable, we have stronger results (\ref{rk1}) and (\ref{r1eq}) which connect the world of nilpotent orbits with the world of $\SL(2)$-orbits and Borel--Serre orbits. 

In Section \ref{ss:per},   we consider extended period maps.  

\subsection{Global properties of the extended period domains}\label{ss:+3}
\begin{sbthm}\label{1.4.3} 
  Let $X$ be one of $D^{\star}_{\SL(2)}$, $D^{\star,+}_{\SL(2)}$, $D^{\star,-}_{\SL(2)}$, $D^{\star,\BS}_{\SL(2)}$, 
$D^{\diamond}_{\SL(2)}$, $D_{\BS,\val}$,  $D_{\SL(2),\val}^I$, $D_{\SL(2),\val}^{II}$, $D_{\SL(2),\val}^{\star}$, $D^{\diamond}_{\SL(2),\val}$,  $D^{\sharp}_{\Sig,[:]}$, and 
$D^{\sharp}_{\Sigma,[\val]}$.
Let $\G$ be a subgroup of $G_\bZ$. 

{\rm(1)} The action of $\G$ on $X$ is proper, and the quotient
space $\G\bs X$ is Hausdorff.  
 
{\rm (2)} Assume that $\G$ is neat. 
Let $\g\in \G$,  $p\in X$, and assume $\g p=p$. Then $\g=1$. 

{\rm (3)} 
Assume that $\G$ is neat. 
Then the
projection $X\to\G\bs X$ is a local homeomorphism. 
  Further, for $X=
D^{\star}_{\SL(2)}$, $D^{\star,+}_{\SL(2)}$, $D^{\star,-}_{\SL(2)}$, $D^{\star,\BS}_{\SL(2)}$, 
there is a structure on the quotient such that the projection is a local isomorphism in $\cB'_{\bR}(\log)$. 
\end{sbthm}

Note that the corresponding results for $D_{\BS}$, $D_{\SL(2)}^I$ and $D_{\SL(2)}^{II}$, and $D^{\sharp}_{\Sig}$ and $D^{\sharp}_{\Sig,\val}$ were already proved in Part I Theorem 9.1, Part II Theorem 3.5.17, and Part III, Theorem 4.3.6, respectively.

   \begin{pf} (3) for $X$ follows from (1) and (2) for $X$. Hence it is sufficient to prove (1) and (2). Since we have continuous maps $D_{\SL(2),\val}^{\star}\to D_{\BS,\val} \to D_{\BS}$ and $D_{\Sig, [\val]}\to D_{\Sig, [:]}\to \Gamma\bs D_{\Sig}$ which are compatible with the actions of $\Gamma$, the results for $D_{\SL(2),\val}^{\star}$, $ D_{\BS,\val}$, $D_{\Sig, [\val]}$, $D_{\Sig, [:]}$ follow from the results for $D_{\BS}$ and $\Gamma\bs D_{\Sig}$. Since $D^{\star}_{\SL(2),\val}\to D^{\star}_{\SL(2)}$ is proper and surjective, 
the properness of the action of $\Gamma$ on $D^{\star}_{\SL(2)}$ follows from that for $D^{\star}_{\SL(2),\val}$. (2) for $D^{\star}_{\SL(2)}$ follows from the $\bar L$-property (i.e., Theorem \ref{ls1} for the situation (a) in \ref{sit})
and the result for the pure case. Since there are continuous maps $D^{\diamond}_{\SL(2),\val}\to D^{\diamond}_{\SL(2)}\to D^{\star}_{\SL(2)}$ which are compatible with the actions of $\Gamma$, the results for $D^{\diamond}_{\SL(2)}$ and $D^{\diamond}_{\SL(2),\val}$ follow from the result for $D^{\star}_{\SL(2)}$.
   \end{pf}
   \begin{sbcor} The spaces in the above theorem are Hausdorff.

\end{sbcor}

This is obtained from the above theorem by taking $\Gamma=\{1\}$.
 
 \begin{sbcor}  Let $X= D^I_{\SL(2)}$, $D^{II}_{\SL(2)}$, $D^{\star}_{\SL(2)}$ or $D^{\diamond}_{\SL(2)}$. Let $\Gamma$ be a neat subgroup of $G_{\Z}$. Then there is a unique structure on $\Gamma\bs X$ as an object of $\cB'_\R(\log)^+$ (\ref{b[+]}) such that the projection $X\to \Gamma\bs X$ is a morphism in $\cB'_\R(\log)^+$ which is locally an isomorphism.  
   \end{sbcor}

\begin{pf} 
This follows from (3) of Theorem \ref{1.4.3} and the corresponding results for $D^I_{\SL(2)}$ and $D^{II}_{\SL(2)}$ in Part II Theorem 3.5.17. 
\end{pf}

\subsection{Results on nilpotent orbits in one variable}\label{ss:cks4}

We prove results Theorem \ref{rk1} and Theorem \ref{r1eq} on nilpotent orbits in one variable. In \ref{rk21}--\ref{rk22}, give a counter-example for the extension of Theorem \ref{rk1} to nilpotent orbits in many variables. 
\begin{sbpara}
 Let $(D^{\sharp}_{\Sig,[:]})'\subset D^{\sharp}_{\Sig,[:]}$ be the union of the two open sets $D_{\Sig,[:]}^{\sharp,\mild}$ (\ref{mild}) and 
the inverse image of $D_{\SL(2),\nspl}$ by $D^{\sharp}_{\Sig,[:]}\to D_{\SL(2)}^I$ in \ref{valper} (1).

Then $(D^{\sharp}_{\Sig,[:]})'$ is the union of  $D_{\Sig,[:]}^{\sharp,\mild}$ and the set of the 
points $p$ of $D^{\sharp}_{\Sig,[:]}$ such that if 
 $N_1, \dots, N_n$ (ordered) is the monodromy logarithms associated to $p$, then 
 $(W, N_1)$ does not split. 
 
 The morphisms  $D_{\Sig, [:]}^{\sharp,\mild}\to  D^{\star}_{\SL(2)}$ (\ref{diathm}, \ref{afd2}) and $D_{\SL(2),\nspl}\to D^{\star}_{\SL(2)}$ (\ref{lam}) induce a morphism $(D_{\Sig,[:]}^{\sharp})' \to D_{\SL(2)}^{\star}$. 
 
 Let $(D^{\sharp}_{\Sigma,\lval})'$ be the inverse image of $ (D^{\sharp}_{\Sig,[:]})' $ in $D^{\sharp}_{\Sigma,\lval}$. Then we obtain the induced morphism 
$(D^{\sharp}_{\Sigma,\lval})'\to D^{\star}_{\SL(2),\val}$
and a commutative diagram 
$$\begin{matrix}
  (D^{\sharp}_{\Sigma,\lval})'&
\overset{\psi}\to & D^{\star}_{\SL(2),\val}\\
\downarrow &&\downarrow\\
  (D^{\sharp}_{\Sig,[:]})' & \overset {\psi} \to &D^{\star}_{\SL(2)}.
\end{matrix}$$
\end{sbpara}

Let $\Xi$ be as in \ref{afd3}.
Since $(D_{\Xi}^{\sharp})'=D^{\sharp}_{\Xi}$, we have

\begin{sbthm}\label{rk1} The identity map of $D$ extends uniquely to a continuous map
$$D_{\Xi}^{\sharp}\to D^{\star}_{\SL(2),\val}$$ and hence extends uniquely to a continuous map $D^{\sharp}_{\Xi} \to D_{\BS,\val}$.

\end{sbthm}

\begin{sbrem}\label{r1rem} 

(1) The image of $D^{\sharp}_{\Xi}$ in $D_{\SL(2)}$ is contained in $D_{\SL(2),\leq 1}$ for both structures $I$, $II$ of $D_{\SL(2)}$. (We denote by $\leq 1$ the part where the log structure is of rank $\leq 1$.) 

(2) However, the image of  $D^{\sharp}_{\Xi}$ in $\LbD$ is not necessarily contained in $D^{\star}_{\SL(2),\leq 1}$. (This is seen in \ref{c.ex} below.) 
Hence the morphism in \ref{rk1} cannot be obtained as the composition $D^{\sharp}_{\Xi}\to D^{\star}_{\SL(2),\leq 1}\cong D^{\star}_{\SL(2),\leq 1,\val}$ (the first arrow here need not exist). For $p\in D_{\Xi}^{\sharp}$, it can happen that the image of $p$ in $D^{\star}_{\SL(2),\val}$ has some information about $p$ which the image of $p$ in $D^{\star}_{\SL(2)}$ does not have (see \ref{info1} and \ref{ex1} below). 

(3) The image of $D^{\sharp,\mild}_{\Xi}\to D^{\star,\mild}_{\SL(2)}$ (2.1.4) is contained in  $D^{\star,\mild}_{\SL(2),\leq 1}$.  
\end{sbrem}

\begin{sbthm}\label{r1eq} Let 
 $p=(\R_{\geq 0}N, \exp(i\R N)F)\in D^{\sharp}_{\Xi}$ with $N\ne0$. 
Let $W'=W^{(1)}$ be the relative monodromy filtration of $N$ with respect to $W$. 
Let $(W',\hat F)$ be the $\R$-split mixed Hodge structure associated to the mixed Hodge structure $(W',F)$, i.e., $\spl_{W'}(F)(F(\gr^{W'}))$ $(1.2)$.
Then the following conditions are equivalent.

{\rm (i)} $p$ belongs to $D^{\sharp,\mild}_{\Xi}$. 

{\rm (ii)} $\exp(iyN)F$ converges in $D^{\diamond}_{\SL(2)}$ when $y\to \infty$.

{\rm (iii)} $\delta_W(\exp(iyN)F)$ converges in $\cL$ when $y\to \infty$.

{\rm (iv)} The image of $p$ in $D^{\star}_{\SL(2)}$ (\ref{rk1}) belongs to $D^{\star,\mild}_{\SL(2)}$.

{\rm (v)} The image of $p$ in $D_{\BS}$ (\ref{rk1}) belongs to $D^{\mild}_{\BS}$. 

{\rm (vi)} The image of $p$ in $D_{\SL(2)}$ belongs to $D_{\SL(2),\spl}$ (\ref{lam}).

{\rm (vii)} $\delta_W(\exp(iN)\hat F) =0$.

{\rm (viii)} The splitting $\spl_W(\exp(iN)\hat F)$ of $W$ is compatible with $N$. 

\end{sbthm}

\begin{pf}
We have proved (i) $\Rightarrow$ (ii). (ii) $\Rightarrow$ (iii) is clear. We know (ii) $ \Rightarrow$ (iv) 
$\Leftrightarrow$ (v), (v) $\Rightarrow$ (vi) $\Leftrightarrow$ (vii). (viii) $\Rightarrow$ (i) is clear.

It is sufficient to prove the implications (iii) $\Rightarrow$ (vii) and (vii) $\Rightarrow$ (viii). 

 Let $s= \spl_W(\exp(iN)\hat F)$, $\bar N=\gr^W(N)\in \bigoplus_w \Hom(\gr^W_w, \gr^W_w)$, $N_0=s\bar N s^{-1}$.

 We prove (vii) $\Rightarrow$ (viii). Assume (vii). Then $\exp(iN)\hat F= s(\exp(i\bar N){\hat F}(\gr^W)) = \exp(iN_0) \hat F$. 
For the mixed Hodge structure $(W', \exp(iN)\hat F)=(W', \exp(iN_0)\hat F)$, we have $N=\delta_{W'}(\exp(iN)\hat F) = \delta_{W'}(\exp(iN_0)\hat F) = N_0$ and (viii) holds.

For the proof of  (iii) $\Rightarrow$ (vii), we  first prove the following claim.

{\bf Claim.} $\delta_W(\exp(iN)\hat F)$ is of $W'$-weight $\leq -1$. 

Proof of Claim. Let $A= W_0 \End_{\C}(H_{0,\C})$. For $r\geq 1$, let $I^{(r)}$ be the 
two-sided ideal $W_{-1}A\cap W'_{-r}A$ of $A$, and let $I=I^{(1)}$. Let $M_1= I \cap \End_\R(H_\R)$. Let $M_2$ be the part of $iM_1\subset I$ consisting of all elements which belong to the $(\leq -1, \leq -1)$-Hodge component of $A$ with respect to $\exp(iN_0)\hat F$. Let $M_3$ be the part of $I$ consisting of all elements which belong to $F^0A$  with respect to $\exp(iN)\hat F$. 
Then we have $I^{(r)}= (I^{(r)}\cap M_1)\oplus (I^{(r)}\cap M_2)\oplus (I^{(r)}\cap M_3)$ for any $r\geq 1$. We have $\exp(iN)\hat F=\exp(x) \exp(iN_0)F$ for some $x\in I$. Hence by \ref{Qalg2}, there are $x_j\in M_j$ $(j=1,2,3)$ such that $\exp(iN)\hat F= \exp(x_1)\exp(x_2)\exp(x_3)\exp(iN_0)\hat F= \exp(x_1)\exp(x_2)\hat F= s' \exp(i\delta) \exp(i\bar N){\hat F}(\gr^W)$ where 
$s' = \exp(x_1)s\in \spl(W)$ and $i\delta =s^{-1}x_2s$. We have $\delta_W(\exp(iN)\hat F)=\delta\in W'_{-1}\End_\R(\gr^W_\R)$.  

Now we prove (iii) $\Rightarrow$ (vii). Assume that $\delta_W(\exp(iN)\hat F)\neq 0$. By Claim, there is $w\leq -1$ such that the component of $\delta_W(\exp(iN)\hat F)$ of  $\tau$-weight $w$ 
is non-zero. Since $\Ad(\tau(\sqrt{y})) \delta_W(\exp(iyN)F)$ converges to $\delta_W(\exp(iN)\hat F)$ when $y\to \infty$, $\delta_W(\exp(iyN)F)$ is $\Ad(\tau(\sqrt{y})^{-1})B(y)$ where $B(y)$ converges to an element whose $w$-part is non-zero. Hence the $w$-part of $\tau$-weight of $\delta_W(\exp(iyN)F)$ is $ y^{-w/2}C(y)$ where $C(y)$ converges to a non-zero element, and hence diverges. 
\end{pf}

\begin{sbpara}\label{rk21}
We have constructed CKS maps $D^{\sharp,\mild}_{\Sig, [:]}\to D^{\star,\mild}_{\SL(2)}$ (\ref{diathm}) and $D^{\sharp}_{\Xi}\to D^{\star}_{\SL(2)}$ (\ref{rk1}). In the rest of this subsection, 
we show an example of $\sigma$ of rank $2$ such that there is no continuous map $D^{\sharp}_{\sig,\val}\to D^{\star}_{\SL(2)}$ which extends the identity map of $D$. 

\end{sbpara}
\begin{sbpara} Take an integer $m\geq 1$. (The case  $m\geq 3$ will be a crucial example.)

Let $H_0$ be of rank $2m+1$ with base $e_j'$ $(1\leq j\leq m)$, $e_j$ $(1\leq j\leq m)$, and $e$.

The weight filtration is as follows. $W_{-m-1}=0$. $W_{-m}$ is generated by $e_j'$ and $e_j$ $(1\leq j\leq m)$. $W_{-1}=W_{-m}$. $W_0$ is the total space.

We have $N_1, N_2$ defined as follows.

$N_1e=0$, $N_1e_j=e'_j$. $N_1e'_j=0$.

$N_2e=e'_m$, $N_2e_j=e_{j-1}$ and $N_2e'_j=e'_{j-1}$ for $2\leq j\leq m$, $N_2e_1=N_2e'_1=0$. 

Let $\sig$ be the cone generated by $N_1$ and $N_2$.

Note that $(W, N_1)$ splits, but $(W, N_2)$ does not split.
\end{sbpara}

\begin{sbpara}

For $j=1,2$, let $W^{(j)}$ be the $W$-relative $N_j$-filtration.

We give a splitting of $W^{(1)}$, which is compatible with $N_1$,  as follows.  $e$ is of weight $0$. $e_j$ is of weight $-m+1$, and $e'_j$ is of weight $-m-1$ $(1\le j\le m)$.

We give a splitting of $W^{(2)}$, which is compatible with $N_2$, as follows. $e$ is of weight $0$. $e_j$ is of weight $-2(m-j)$, and $e'_j$ is of weight $-2(m-j+1)$ $(1\le j\le m)$.

\end{sbpara}

\begin{sbpara}
Define $\alpha_1,\alpha_2: {\mathbb G}_{m,\R}\to \Aut_\R(H_{0,\R}) $ by using the above splittings of $W^{(1)}$ and $W^{(2)}$, respectively. Then $\alpha_1$ and $\alpha_2$ commute. Define $\alpha^{\star}_1, \alpha^{\star}_2: {\mathbb G}_{m,\R}\to \Aut_\R(H_{0,\R}) $ by
$\alpha^{\star}_j(t) e=e$ and $\alpha_j^{\star}(t) x= t^m \alpha_j(t) x$ for $x\in W_{-m}$.

Let $$t(y)=
\alpha_1((y_2/y_1)^{1/2})\alpha_2((1/y_2)^{1/2}), 
\quad t^{\star}(y)=\alpha^{\star}_1((y_2/y_1)^{1/2})\alpha^{\star}_2((1/y_2)^{1/2}).$$ 
\end{sbpara}

\begin{sbpara} We have
$$\Ad(t(y))^{-1}(y_1N_1+y_2N_2)=N_1+N_{2,y},$$ 
where $N_{2,y}$ coincides with $N_2$ on $e_j$ and $e'_j$ $(1\le j\le m)$, but $N_{2,y}e=(y_2/y_1)^{(m+1)/2}e_m'$ and the last element converges to $0$ when $y_1/y_2\to \infty$.

$$\Ad(t^{\star}(y))^{-1}(y_1N_1+y_2N_2)=N_1+ N'_{2,y},$$ where
$N'_{2,y}$ coincides with $N_2$ on $e_j$ and $e'_j$ $(1\le j\le m)$, but $N'_{2,y}e=u_ye_m'$, where $$u_y=(y_2/y_1)^{1/2}y_2^{m/2}= y_1^{-1/2}y_2^{(m+1)/2}.$$

\end{sbpara}

\begin{sbpara}

Note that 
$u_y$ need not converge when $y_2\to \infty$ and $y_1/y_2\to \infty$. 

\end{sbpara}

\begin{sbpara} Let $F$ be as follows. $F^1=0$. $F^0$ is generated by $e$ and $e_m$. $F^{-j}$, for $1\leq j\leq m-1$, is generated by $F^{-j+1}$ and $e_{m-j}$, $e'_{m-j+1}$. $F^{-m}$ is the total space. 

Then $(N_1, N_2, F)$ generates a nilpotent orbit.

Hence 
$\exp(iy_1N_1+iy_2N_2)F$, as $y_2\to \infty$ and $y_1/y_2\to \infty$, converges in $D^{\sharp}_{\sig,\val}$.

\end{sbpara}

\begin{sblem}\label{noCKS}  Let the notation be as above. If $m\geq 3$, $\exp(iy_1N_1+iy_2N_2)F$, as $y_2\to \infty$ and $y_1/y_2\to \infty$, need not converge in $D^{\star}_{\SL(2)}$. 

\end{sblem}

This follows from

\begin{sblem}\label{noCKS2} Let the notation and the assumption be as above. 
Let $F_y:=t^{\star}(y)^{-1}\exp(iy_1N_1+iy_2N_2)F$.
Then, 
$\delta_W(F_y)$ does not converge in $\bar\cL$ when  $y_2\to \infty$ and $y_1/y_2\to \infty$.
\end{sblem}

\begin{sbpara}\label{rk22}

We prove Lemma \ref{noCKS2}.

Since $F_y= \exp(\Ad(t^{\star}(y))^{-1}(iy_1N_1+iyN_2)) F$, $F_y$ is described as follows.
$F_y^1=0$, $F_y^0$ is generated by $e+\sum_{k=1}^m i^k\cdot k!^{-1}\cdot u_y \cdot e_{m-j+1}'$ and $\exp(iN_1+iN_2)e_m$,  $F_y^{-j}$ $(1\leq j\leq m-1)$ is generated by $F_y^{-j+1}$ and
 $\exp(iN_1+iN_2)e_{m-j}$ and  $\exp(iN_1+iN_2)e_{m-j+1}'$, 
and $F_y^{-m}$ is the total space. 

The Hodge type of $\gr^W_{-m}$ of $F_y$ is that the $(j, -m-j)$-Hodge component is of dimension one if $j=0, -m$, two dimensional if $-1\geq j\geq 1-m$, and zero otherwise.

$\delta_W(F_y)$ sends $e$ to the sum of the $(j, -m-j)$ components of $v_y:=\sum_{1\leq k\leq m,k:\text{odd}} (-1)^{(k-1)/2}\cdot k!^{-1}\cdot t \cdot e_{m-k+1}$ for $-1\geq j\geq 1-m$. 

\medskip

{\bf Claim.} $v_y$ does not belong to the $((0,-m)+(-m,0))$-Hodge component of $\gr^W_{-m}$.

\medskip

By the claim, 
$v_y$ is $u_y$ times a {\it non-zero} element which is independent of $y_1, y_2$. Hence, when $y_1/y_2, y_2\to \infty$, $v_y$  need not converge in $\bar \cL$. 

We prove the claim. Assume that $v_y$ belongs to the 
 $((0,-m)+(-m,0))$-Hodge component of $\gr^W_{-m}$. Then we should have
$$\sum_{1\leq k\leq m, k:\text{odd}} (-1)^{(k-1)/2}\cdot k!^{-1}\cdot e'_{m-k+1} =  a\exp(iN_1+iN_2)e_m + b\exp(-iN_1-iN_2)e_m$$
for some $a,b\in \C$. If $V$ denotes the $\C$-vector space generated by $e_j$ $(1\leq j\leq m)$ and $e'_j$ $(1\leq j\leq m-3)$, we should have
$$e'_m - (1/6)e'_{m-2} \equiv a(ie'_m - e'_{m-1} -(i/2)e'_{m-2}) + b(-ie'_m - e'_{m-1}+(i/2) e'_{m-2})\; \bmod V.$$
(To get this, use $(N_1+N_2)^k=kN_1N_2^{k-1}+ N_2^k$ and hence $\exp(iN_1+iN_2) = 1+ \sum_{k=1}^\infty (i^k \cdot (k-1)!^{-1} \cdot N_1N_2^{k-1}+ i^k \cdot k!^{-1}\cdot N_2^k)$.) By comparing the coefficients of $e'_{m-1}$, we have $a+b=0$. Hence
$$e'_m-(1/6)e'_{m-2}\equiv a\cdot 2i\cdot (e'_m - (1/2)e'_{m-2})\;\bmod V.$$
This is impossible. 
\end{sbpara}

\begin{sbrem}
We do not know whether the identity map of $D$ always extends to a continuous map $D_{\Sig, [\val]}^{\sharp}\to D^{\star}_{\SL(2),\val}$ or not.
\end{sbrem}

\subsection{Extended period maps}\label{ss:per}
The following is a modified version of Part III, 7.5.1 (1).

\begin{sbthm}\label{perth1} 
Let $S$ be a connected, log smooth, fs log analytic space, and let $U$ be the open subspace of $S$ consisting of all points of $S$ at which the log structure of $S$ is trivial.
Let $H$ be a variation of mixed Hodge structure on $U$ with polarized graded quotients for the weight filtration, and with unipotent local monodromy along $S\smallsetminus U$. 
Assume that $H$ extends to a log mixed Hodge structure (Part III, \S1.3) on $S$ $($that is, $H$ is admissible along $S\smallsetminus U$ as a variation of mixed Hodge structure$)$. 
Fix a base point $u\in U$ and let $\Lambda=(H_0,W, (\langle\;,\;\rangle_w)_w, (h^{p,q})_{p,q})$ be $(H_{\Z,u}, W, (\langle\;,\;\rangle_{w,u})_w, (\text{the Hodge numbers of $H$}))$.
Let $\Gamma$ be a subgroup of $G_\Z$ which contains the global monodromy group $\text{Image}(\pi_1(U,u)\to G_\Z)$ and assume that $\G$ is neat. 
Let $\varphi: U\to \Gamma\bs D$ be the associated period map. Let $S^{\log}_{[:]}= S^{\log} \times_S S_{[:]}$ and let $S^{\log}_{[\val]}=S^{\log}\times_S S_{[\val]}$, and regard $U$ as open sets of these spaces. 

Then:

(1) 
The map $\varphi:U\to\G\bs D$ extends uniquely to a continuous maps
$$
S_{[:]}^{\log}\to\G\bs D_{\SL(2)}^I, \quad S_{[\val]}^{\log} \to \G \bs D^I_{\SL(2),\val}
$$

(2) Assume that the complement $S\smallsetminus U$ of $U$ is a smooth divisor on $S$. Then 
the map $\varphi:U\to\G\bs D$ extends uniquely to a continuous map
$$
S^{\log} \to \G \bs D^{\star}_{\SL(2),\val}
$$
and hence extends uniquely to a continuous map $S^{\log} \to \G \bs D_{\BS,\val}$. 

\end{sbthm}

\begin{pf} (1) is a modified version of Part III, 7.5.1 (1) which treated the extended period map $S^{\log}_{\val} \to D^I_{\SL(2)}$ where $S^{\log}_{\val}$ is the topological space defined in \cite{KU} 3.6.26. This map factors through the quotient space $S_{[:]}^{\log}$ of $S^{\log}_{\val}$ as is seen by the arguments in \ref{Dsig:3}.
Since  $S^{\log}_{\val}\to S_{[:]}$ is a proper surjective continuous map, the map $S_{[:]}\to \G \bs D_{\SL(2)}^I$ is continuous. 

The last map is compatible with log structures as is seen by the arguments in \ref{5.6.4} and hence induces a continuous map $S^{\log}_{[\val]}\to \G\bs D^I_{\SL(2),\val}$. 

(2) is proved similarly to (1) by using Theorem \ref{rk1}. 
\end{pf}

In the rest of this Section 6.3, we consider mild log mixed Hodge structures.

\begin{sbprop}\label{WNQ}
Let $\sigma$ be a rational nilpotent cone (it is an $\R_{\ge0}$-cone generated by rational elements) in $\fg_\R$. Assume that 
 there is $F\in \Dc$ such that $(\sigma, F)$ generates a nilpotent orbit. If $(W, N)$ splits for any rational element $N$ of $\sigma$, $(W, N)$ splits for any element $N$ of the cone $\sigma$. 
\end{sbprop}

\begin{pf} 
We may assume that $N$ is in the interior $\sigma_{>0}$ of $\sigma$. This is because if we denote by $\sig'$ the face of $\sig$ such that $N$ belongs to the interior of $\sig'$, $(\sig', \exp(iN')F)$ generates a nilpotent orbit for some $N'\in \sig_{>0}$ and hence we can replace $\sig$ by $\sig'$. 

Assume $N\in \sig_{>0}$. Let $\hat F$ be the $\R$-split MHS associated to the MHS $(M(W, \sig), F)$. Then $\exp(iN')\hat F\in D$ for any $N'\in \sig_{>0}$. Hence as the composition of the continuous map $D\to \cL\;;\;x\mapsto \delta_W(x)$ and the continuous map $\sig_{>0}\to D\;;\; N' \mapsto \exp(iN')\hat F$, the map $\sig_{>0}\to \cL\;;\;N'\mapsto \delta_W(\exp(iN')\hat F)$ is continuous. 
By the part (i) $\Rightarrow$ (vii) of \ref{r1eq}, 
the last map sends all rational elements of $\sigma_{>0}$ to $0$. 
 Hence  it sends $N$ also to $0$. By the part (vii) $\Rightarrow$ (i) of \ref{r1eq}, this shows that $(W, N)$ splits. 
\end{pf}

\begin{sbpara} Let $\cB(\log)$ be the category of locally ringed spaces over $\C$ endowed with fs log structures satisfying a certain condition,  defined in \cite{KU} (see \cite{KNU2} Part III, \S1.1 for the review).

Let $S$ be an object of $\cB(\log)$ and let $H$ be an LMH on $S$ with polarized graded quotients for the weight filtration. By \ref{WNQ}, for $s\in S$ and for $t\in S^{\log}$ lying over $s$, the following two conditions (i) and (ii) below are equivalent. Let $$\pi_1^+(s^{\log}):= \Hom((M_S/\cO_S^\times)_s, \N) \subset \pi_1(s^{\log})=\Hom((M_S/\cO^\times_S)_s, \Z),$$
$$\pi_1(s^{\log},\R_{\geq 0}) := \Hom((M_S/\cO^\times_S)_s, \R_{\geq 0}^{\add})\subset \R\otimes \pi_1(s^{\log})= \Hom((M_S/\cO_S^\times)_s, \R^{\add}).$$  Consider the action $\rho$ of $\pi_1(s^{\log})$ on $H_{\Z,t}$, and consider the homomorphism 
$$\log(\rho): \R\otimes  \pi_1(s^{\log})\to \End_{\R}(H_{\R, t})\;;\; a\otimes \gamma \mapsto a\log(\rho(\gamma)).$$ Let $W$ be the weight filtration on $H_{\R,t}$. 

\medskip

(i) For any $\gamma\in \pi_1(s^{\log}, \R_{\geq 0})$, $(W, \log(\rho)(\gamma))$ splits.

(ii) For any $\gamma\in \pi_1^+(s^{\log})$, $(W, \log(\rho(\gamma)))$ splits.

\medskip

 We say that $H$ is {\it mild} (we say also $H$ is {\it of mild degeneration})  if the equivalent conditions (i) and (ii) are satisfied for any $s$ and $t$.

\end{sbpara}

\begin{sblem}\label{pullback} Let $S$ and $H$ be as above and assume $H$ is mild. Let $S'\to S$ be a morphism in $\cB(\log)$. Then the pull back of $H$ to $S'$ is mild. 

\end{sblem}

This is clear.

\begin{sbprop}\label{Ccut} Let $S$ be a log smooth fs log analytic space over $\C$ and let $H$ be a log mixed Hodge structure  on $S$ with polarized graded quotients for the weight filtration $W$. Then the following two conditions (i) and (ii) are equivalent. 

\smallskip
 (i) $H$ is mild.
 
 (ii) For any smooth analytic curve $C$ over $\C$ and any analytic map $f: C\to S$ such that the subset $f^{-1}(S\smallsetminus U)$ of $S$ is finite, the pull back $f^*H$ on $C$ is mild. Here we endow $C$ with the log structure associated to the finite subset $f^{-1}(S\smallsetminus U)$. 
 
 \smallskip
 
If $S$ is an algebraic variety over $\C$, these conditions are equivalent to the modified version  of the condition (ii) in which we take only smooth algebraic curves $C$ in it. 
\end{sbprop}

\begin{pf}  By \ref{pullback}, we have (i)$ \Rightarrow$ (ii). We prove (ii) $\Rightarrow$ (i).  Assume (ii). Let $s\in S\smallsetminus U$ and let $t$ be a point of $S^{\log}$ lying over $s$. Let $\gamma\in \pi_1^+(s^{\log})$. We prove that $(W, \log(\rho)(\gamma))$ splits. 

Let $\sig$ be the face of $\pi_1^+(s^{\log})$, regarded as a monoid, such that $\gamma$ belongs to the interior of $\sig$. Then there are $s'\in S$ and $t'\in S^{\log}$ lying over $s'$ and isomomorphisms $\pi^+((s')^{\log})\cong \sig$ and $H_{\R, t'}\cong H_{\R,t}$ such that the action of $\pi_1^+(s^{\log})$ on $H_{\R,t}$ and the action of $\pi_1^+((s')^{\log})$ on $H_{\R,t'}$ are compatible via these isomorphisms. 
By this we are reduced to the case where $\gamma$ belongs to the interior of $\pi_1^+(s^{\log})$. 

Assume $\gamma$ belongs to the interior of $\pi_1^+(s^{\log})$. Then there are a smooth analytic curve over $\C$, a morphism $f:C \to S$ and  $s' \in C$ satisfying the following conditions (1)--(3). (1) $f(s')=s$. (2) $f^{-1}(S\smallsetminus U)$ is finite. (3) The image of $\pi_1^+((s')^{\log})\to \pi_1^+(s^{\log})$ contains $\gamma$. 
By the condition (ii), this proves that $(W, \log(\rho)(\gamma))$ splits. 

In the case where $S$ is an algebraic variety, the same arguments show that the modified version of (ii) implies (i). 
\end{pf}

\begin{sbthm}\label{perth2} Let the assumptions be as in \ref{perth1}. Assume furthermore that $H$ is mild. 

{\rm(1)} The period map $\varphi: S \to \G \bs D$ extends uniquely to continuous maps $$S^{\log}_{[:]}\to \Gamma\bs D^{\diamond}_{\SL(2)},\quad S^{\log}_{[:]}\to \Gamma\bs D^{\star,\mild}_{\SL(2)},$$ 
$$ S^{\log}_{[\val]}\to \Gamma \bs D^{\diamond}_{\SL(2),\val}, \quad S^{\log}_{\lval}\to \Gamma\bs D^{\star,\mild}_{\SL(2),\val}, \quad S^{\log}_{\lval}\to \Gamma\bs D^{\mild}_{\BS,\val}.$$ 

\medskip

{\rm(2)} For any point $s\in S$, there exist an open neighborhood
$V$ of $s$, a log modification $V'$ of $V$ $(\text{\cite{KU}}\ 3.6.12)$, a
commutative subgroup $\Gamma'$ of $\G$, and a fan $\Sig$ in $\fg_\bQ$
which is strongly compatible with $\Gamma'$ such that the period map
$\varphi|_{U\cap V}$ lifts to a morphism $U\cap V \to \Gamma'\bs D$ which
extends uniquely to a morphism $V'\to\Gamma'\bs D^{\mild}_\Sig$ of log manifolds.
$$
\begin{matrix}
U & \supset & U\cap V & \sub & V' \\
{{}^\varphi}\downarrow\;\;\; & & \downarrow &  & \downarrow\\
\Gamma \bs D & \leftarrow  & \Gamma' \bs D &  \sub  & \Gamma' \bs D^{\mild}_\Sig.
\end{matrix}
$$
Furthermore, we have
\medskip

{\rm(2.1)} Assume $S\smallsetminus U$ is a smooth divisor. 
Then we can take $V=V'=S$ and $\G'=\G$. 
That is, we have a commutative diagram 
$$
\begin{matrix} U&\sub & S \\ 
{{}^\varphi}\downarrow\;\;\;&&\downarrow\\
\Gamma\bs D&\sub&\Gamma\bs D^{\mild}_{\Sig}.\end{matrix}
$$

{\rm(2.2)} Assume $\Gamma$ is commutative.
Then we can take $\Gamma'=\Gamma$.

\medskip

{\rm (2.3)} 
Assume that $\Gamma$ is commutative and that the following condition {\rm(i)} is satisfied.
\smallskip

{\rm(i)} There is a finite family $(S_j)_{1\leq j\leq n}$ of connected
locally closed analytic subspaces of $S$ such that $S=\bigcup_{j=1}^n S_j$
as a set and such that, for each $j$, the inverse image of the sheaf
$M_S/\cO^\times_S$ on $S_j$ is locally constant.
\smallskip

Then we can take $\Gamma'=\Gamma$ and $V=S$. 
\end{sbthm}

(1) and (2) are modified version of Part III, Theorem 7.5.1 (1) and (2), respectively. (2) is proved in the same way as Part III, Theorem 7.5.1 (2). We can deduce (1) from (2) by using $D^{\sharp,\mild}_{\Sig}\to D^{\diamond}_{\SL(2)}$ (\ref{diathm})  by the arguments in the above proof of Theorem \ref{perth1} (1).

\section{Relations with asymptotic behaviors of regulators and local height pairings}\label{s:Ex}
In this Section \ref{s:Ex}, we show examples to describe the relations of this work to the work \cite{BK} on the asymptotic behaviors of regulators and local height pairings.

\subsection{Example III}\label{III}

This is Example III of Parts I, II. It
appeared in Part III as the case $b=2$ of 7.1.3. As in Section \ref{ss:reg} below, this example is related to the regulator of $K_2$ of a degenerating elliptic curve. 

In this Example III, and also in Example IV in Section 7.3 below, we compare $D_{\BS}$, $D_{\SL(2)}^I$, $D_{\SL(2)}^{II}$,  $D^{\star}_{\SL(2)}$, $D_{\SL(2)}^{\diamond}$,  and their associated valuative spaces, by regarding them as topological spaces, that is, 
we forget the real analytic structures. 

\begin{sbpara}
Let $H_0=\Z^3$ with basis $e_1,e_2,e_3$. The weight filtration is given by $$W_{-4}=0\subset W_{-3}=\R e_1+\R e_2=W_{-1}\subset W_0=H_{0,\R}.$$
The intersection form on $\gr^W_{-3}$ is the anti-symmetric form characterized by $\langle e_2, e_1\rangle=1$.

\end{sbpara}

\begin{sbpara} $D(\gr^W)\cong \frak h$, the upper half plane. 

$D_{\SL(2)}(\gr^W)=D_{\BS}(\gr^W)={\frak h}_{\BS}$.

\end{sbpara}

\begin{sbpara}
We have a homeomorphism
$$D^{\star}_{\SL(2),\val}\overset{\cong}\to D_{\BS,\val}$$
 and this induces a homeomorphism $$D^{\star}_{\SL(2)}\cong D_{\BS}$$ of quotient spaces. 

Let $W'$ be the increasing filtration on $\gr^W$ given by 
$$
W'_{-5}=0\subset W'_{-4}=\bR e_1=W'_{-3}
\subset W'_{-2}=\gr^W_{-3}
\subset W'_0=\gr^W,
$$
and let $\Phi=\{W'\}$.
Let $P$ be the parabolic subgroup of $G_\bR$ consisting of elements which preserve  $W'$.
Then, 
$D_{\BS}(P)=D^{\star}_{\SL(2)}(\Phi)$ and it is the inverse image of the open set $\{x+iy\;|\;x\in \R,\; y\in (0, \infty]\}$ of ${\frak h}_{\BS}$ under the projection $D_{\BS}=D^{\star}_{\SL(2)}\to D_{\BS}(\gr^W)=D_{\SL(2)}(\gr^W)={\frak h}_{\BS}$.

We have $$D^I_{\SL(2)}= D^{II}_{\SL(2)}$$ (Part II, 3.6.2).
So we denote $D_{\SL(2)}^I$ and $D_{\SL(2)}^{II}$ simply by $D_{\SL(2)}$.  
\end{sbpara}

\begin{sbpara}\label{homeos1} 

Let $$V:=\R e_1+\R e_2.$$ We have
$$\spl(W)\cong V, \quad \cL\cong V,$$  where $v\in V$ corresponds in the first isomorphism  to the splitting of $W$ given by $e_3+v$, i.e., $s\in\spl(W)$ such that $s(e_3(\gr^W_0))=e_3+v$, 
and  $v\in V$ corresponds in the second isomorphism to $\delta\in \cL$ such that $\delta(e_3(\gr^W_0))=v$. We have $\cL(F)=\cL$ for any $F\in D(\gr^W)$.

\end{sbpara}

\begin{sbpara}\label{IIIdia} 
We have homeomorphisms
$$D\cong  {\frak h}\times \cL \times \spl(W)  \cong \R_{>0} \times V  \times \R\times V,$$
where the first isomorphism  is  $F\mapsto (F(\gr^W), \delta_W(F), \spl_W(F))$, and the second isomorphism 
sends $(x+iy, \delta, s)$ to $(t, \delta, x, s)$, where $x, y\in \R$, $y>0$, $t:=1/\sqrt{y}$,  and we identify both $\cL$ and $\spl(W)$ with $V$ via the  isomorphisms in \ref{homeos1}. We call the composition $D\cong \R_{>0}\times V \times \R \times V$ the {\it standard isomorphism for $D$}.
Let $\bar V=\bar\cL$ be as in \ref{2.3ex} (4).

We have a commutative diagram of homeomorphisms
$$\begin{matrix} 
D^{\diamond,\text{weak}}_{\SL(2)}(\Phi) && \cong && (\R_{\geq 0} \times V  \times \R\times V)'\\
\uparrow &&&& \uparrow (1)\\
D^{\diamond}_{\SL(2)}(\Phi)&&  \cong && (\R_{\geq 0} \times V  \times \R\times V)'\\
\downarrow &&&& \downarrow (2)\\
D^{\star}_{\SL(2)}(\Phi) && \cong && \R_{\geq 0}\times {\bar V}\times \R \times V\\
\uparrow &&&& \uparrow \;\;\;\;\;\\
D^{\star}_{\SL(2),\val}(\Phi) && \cong && (\R_{\geq 0}\times {\bar V})_{\val} \times \R\times V\\
\downarrow &&&& \downarrow (3)\\
D_{\SL(2),\val}(\Phi) && \cong && (\R_{\geq 0}\times {\bar V})_{\val} \times \R\times V\\
\downarrow &&&& \downarrow\;\;\;\;\;\\
D_{\SL(2)}(\Phi) && \cong && \R_{\geq 0} \times {\bar V} \times \R\times V,
\end{matrix}$$
where 
$$(\R_{\geq 0}\times V  \times \R\times V)':=\{(t, \delta, x, s)\in \R_{\geq 0} \times V\times \R \times V\;|\;\delta\in \R e_1\;\text{if}\;t=0\}, $$ 
and where  $(\R_{\geq 0}\times {\bar V})_{\val}$ is the valuative space of $\R_{\geq 0}\times {\bar V}$ associated to the canonical log structure (see a description below and also \cite{KU}, 0.5.21).

The homeomorphism for $D^{\diamond,\text{weak}}_{\SL(2)}(\Phi)$ is compatible with the standard isomorphism for $D$, but other homeomorphisms are {\it not} compatible with the standard isomorphism. 

The homeomorphism for $D^{\diamond}_{\SL(2)}(\Phi)$ (resp.\ $D^{\star}_{\SL(2)}(\Phi)$, resp.\ $D_{\SL(2)}(\Phi))$ sends 
 a point of $D$ corresponding to $(t,c_1e_1+c_2e_2, x, u)\in \R_{>0}\times V \times \R \times V$ under the standard isomorphism to 
 $(t, c_1e_1+t^{-1}c_2e_2, x, u)$ (resp.\
$(t, tc_1e_1+t^{-1}c_2e_2, x, u)$, resp.\  $ (t, t^4c_1e_1+t^2c_2e_2, x, u))$. 

The homeomorphism for $D^{\star}_{\SL(2),\val}(\Phi)$ (resp.\ $D_{\SL(2),\val}(\Phi))$ is compatible with the homeomorphism for $D^{\star}_{\SL(2)}(\Phi)$ (resp.\ $D_{\SL(2)}(\Phi))$. 

Concerning the vertical arrows on the right-hand side, they are described as follows.
The arrows without labels are the canonical projections (\ref{ss:valsp}). The map (1) sends $(t, c_1e_1+c_2e_2,x,u)$ to $(t, c_1e_1+tc_2e_2,x,u)$. The map (2) sends $(t,c_1e_1+c_2e_2, x,u)$ to $(t,tc_1e_1+c_2e_2, x,u)$. The map (3) is explained below.

The valuative space  $(\R_{\geq 0}\times {\bar V})_{\val}$ is described as follows. Over $U= (\R_{>0}\times \bar V)\cup (\R_{\geq 0}\times V)\subset \R_{\geq 0}\times \bar V$,  it is $U$. 
The inverse image of  $\{0\}\times (\bar V\smallsetminus V)$ in $(\R_{\ge0}\times {\bar V})_{\val}$ consists of points 

(a) $p(0, \lambda)$ $(\lambda\in \bar V\smallsetminus V)$, 

(b) $p(c, \lambda)$ $(c\in \R_{>0}\smallsetminus \Q_{>0}$, $\lambda\in \bar V\smallsetminus V)$, 

(c) $p(c+, \lambda)$ $(c\in \Q_{\geq 0}$, $\lambda \in {\bar V}\smallsetminus V)$,

(d) $p(c-, \lambda)$ $(c\in \Q_{>0}$, $\lambda\in \bar V\smallsetminus V)$,

(e) 
$p(c, \mu)$ $(c\in \Q_{>0}$, $\mu \in V\smallsetminus \{0\})$.

Write $\lambda = 0\circ\mu$ with $\mu \in V \smallsetminus \{0\}$ (\ref{2.3ex} (4)). Then the above point is the limit of $t^{c'}\mu$,  where $t>0$ and $t\to 0$ and, in the cases of (b) and (e) (resp.\ case (a), resp.\ case (c), resp.\ case (d)), $c'=c$ (resp.\ $c'\to \infty$, resp.\ $c'>c$ and $c' \to c$, resp.\ $c'<c$ and $c'\to c$).

The map (3) sends 

$(t, \delta, x, u)$ $((t,\delta)\in (\R_{>0}\times {\bar V})\cup (\R_{\geq 0}\times V))$ to $(t, t^3\delta,x,u)$, 

$(p(0, \lambda), x,u)$ to $(p(0, \lambda), x,u)$,

$(p(c, \alpha),x,u)$ ($\alpha \in \bar V$) to $(p(c-3,\alpha),x,u)$ if $c>3$, to $(0, \alpha,x,u)$ if $c=3$, and to $(0,0,x,u)$ if $0<c<3$, 

$(p(c+, \lambda),x,u)$ to $(p((c-3)+, \lambda),x,u)$ if $c\geq 3$, and to $(0,0,x,u)$ if $0\leq c<3$, 

$(p(c-, \lambda),x,u)$ to $(p((c-3)-, \lambda), x,u)$ if $c>3$, and to $(0,0,x,u)$ if $0<c\leq 3$.

\end{sbpara}

\begin{sbpara} We describe for Example III

(1) that there is no continuous map $D_{\SL(2),\val}(\Phi)\to D_{\BS}(P)=D_{\SL(2)}^{\star}(\Phi)$ which extends the identity map of $D$, and 

(2) how $\eta: D_{\SL(2),\val}(\Phi)\to D_{\BS,\val}(P)=D^{\star}_{\SL(2),\val}(\Phi)$ is not continuous.

Fixing  $c_1, c_2\in \R$, for $t>0$, let $p(t)$ be the point of $D$ corresponding to $(t, c_1e_1+c_2e_2, 0,0)$ via the homeomorphism for $D^{\star}_{\SL(2),\val}(\Phi)$ in 7.1.5. Then, via the homeomorphism for $D_{\SL(2),\val}(\Phi)$, $p(t)$ corresponds to $(t, t^3c_1e_1+t^3c_2e_2, 0,0)$. 

Hence, when $t\to 0$, $p(t)$ converges in $D_{\SL(2),\val}(\Phi)$ to the point $p$ corresponding to $(0,0,0,0)$, but it converges in $D_{\SL(2)}^{\star}(\Phi)$ to the point 
which corresponds to $(0, c_1e_1+c_2e_2,0,0)$ which depends on $(c_1,c_2)$. This explains (1). Concerning (2), the image of $p$ under $\eta$ is the point $p'$ of $D^{\star}_{\SL(2)}(\Phi)$ corresponding to $(0,0,0,0)$. If $(c_1, c_2)\neq (0,0)$, $p(t)$ does not converge to $p'$ in $D^{\star}_{\SL(2)}(\Phi)$.

\end{sbpara}

\begin{sbpara}\label{noIIstar} As is mentioned in \ref{weakstar}, the topology of $D^{\diamond,\text{weak}}_{\SL(2)}$ does not coincide with the one of $D^{\diamond}_{\SL(2)}$.

Fixing $c\in \R$, for $t>0$, let $p(t)$ be the point of $D$  corresponding to $(t, tce_2, 0,0)$ via the homeomorphism for $D^{\diamond,\text{weak}}_{\SL(2)}(\Phi)$ in 7.1.5. Then,  when $t\to 0$, $p(t)$ converges in $D^{\diamond,\text{weak}}_{\SL(2)}(\Phi)$ to the point corresponding to $(0,0,0,0)$. On the other hand, $p(t)$ corresponds to $(t, ce_2, 0,0)$ under the homeomorphism for $D^{\star}_{\SL(2)}(\Phi)$ and hence converges in $D^{\star}_{\SL(2)}(\Phi)$ but the limit depends on the choice of $c$.

\end{sbpara}

\begin{sbpara}

The open set $D^{\star,\mild}_{\SL(2)}(\Phi)$ of $D^{\star}_{\SL(2)}(\Phi)$ is the part consisting of elements corresponding to $(t,\delta,x,u)$ such that $\delta\in V\subset \bar V$.

The map $D^{\star,\mild}_{\SL(2)}(\Phi)\to D_{\SL(2)}(\Phi)$ corresponds to $(t, \delta, x, u)\mapsto (t, t^3\delta, x, u)$. It does not extend to a continuous map $D^{\star}_{\SL(2)}\to D_{\SL(2)}$. In fact,  fixing $v\in V\smallsetminus \{0\}$, let $p(t)$ for $t>0$ be the point of $D$ corresponding to $(t, t^{-3}v, 0,0)$ via the homeomorphism for $D^{\star}_{\SL(2)}(\Phi)$. Then when $t\to 0$,  $p(t)$ converges to the point of $D^{\star}_{\SL(2)}(\Phi)$ corresponding to 
$(0, 0\circ v, 0,0)$ (\ref{2.3ex} (4)). But $p(t)$ converges to the point of $D_{\SL(2)}(\Phi)$ corresponding to $(0,v, 0,0)$ which depends on the choice of $v$.

\end{sbpara}

\begin{sbpara}\label{ex3}
 Let $a\in \Q_{>0}$ and define  $N_a\in \fg_\Q$ by
$N_a(e_3)=ae_2$, $N_a(e_2)=e_1$, and $N_a(e_1)=0$. 

For $b\in \R$, 
let $F_b\in \Dc$ be the decreasing filtration defined as follows: $F_b^1=0$, $F_b^0$ is generated by $e_3+ibe_1$, $F_b^{-1}$ is generated by $F_b^0$ and $e_2$, and $F_b^{-2}$ is the total space. 

Then $(N_a, F_b)$ generates a nilpotent orbit. 
Let $\sigma_a=\R_{\geq 0}N_a$.
Then $(\sig_a, \exp(i\sig_{a,\R})F_b)\in D^{\sharp}_{\sig_a}$ is the limit of $\exp(iyN_a)F_b$ for $y\to \infty$. This 
$(\sig_a,\exp(i\sig_{a,\R})F_b)$ belongs to $D^{\sharp,\mild}_{\sig_a}$ (\ref{mild}) if and only if $a=0$. 

We consider the image of $\exp(iyN_a)F_b\in D$ in $\R_{>0} \times V \times \R \times V$ under the isomorphism in \ref{IIIdia}. Let $t= 1/\sqrt{y}$. 

In the standard isomorphism for $D$, the image is $(t, at^{-2}e_2+be_1, 0, -(b/2)t^2e_2)$. (The last component is computed by using the relation of $\delta$ and $\zeta$ (\ref{II,1.2.3}).)

In the homeomorphsm for $D^{\diamond}_{\SL(2)}$, the image is $(t, at^{-3}e_2+be_1, 0, -(b/2)t^2e_2)$.

In the homeomorphsm for $D^{\star}_{\SL(2)}$, the image is $(t at^{-3}e_2+bte_1, 0, -(b/2)t^2e_2)$.

In the homeomorphism for $D_{\SL(2)}$, the image is $(t, ae_2+ bt^4e_1, 0, -(b/2)t^2e_2)$.

\end{sbpara}
By taking the limit for $t\to0$, we have:

\begin{sblem}\label{IIIconv}  
$(1)$ If $a\neq 0$, the image of $(\sig_a, \exp(i\sig_{a,\R})F_b)\in D^{\sharp}_{\sigma_a}$  in 
$D^{\star}_{\SL(2)}(\Phi)$ (resp.\ $D_{\SL(2)}(\Phi)$, resp.\ $D^{\star}_{\SL(2),\val}(\Phi)$, resp.\ $D_{\SL(2),\val}(\Phi))$ 
has the  coordinate 
$(0,\infty e_2, 0, 0)$ (resp.\ $(0, ae_2, 0, 0)$, resp.\ $(p(3, ae_2), 0,0)$, resp.\ $(0, ae_2, 0,0)$). 

$(2)$ If $a=0$, the image of $(\sig_0, \exp(i\sig_{0,\R})F_b)\in D^{\sharp,\mild}_{\sigma_0}$ in $D^{\diamond}_{\SL(2)}(\Phi)$ (resp.\ $D^{\star}_{\SL(2)}(\Phi)$, resp.\ $D_{\SL(2)}(\Phi)$, resp.\ $D^{\star}_{\SL(2)}(\Phi)$, resp.\ $D_{\SL(2)}(\Phi)$) has the  coordinate  $(0, be_1, 0, 0)$ (resp.\  $(0, 0, 0, 0)$, resp.\ $(0, 0, 0, 0)$, resp.\ $(0,0, 0,0)$ if $b\neq 0$ and $(0,0,0,0)$ if $b=0$, resp.\ $(0,0,0,0)$). 
\end{sblem}

\begin{sbpara}\label{info1} By Lemma \ref{IIIconv}, we have the following.
Consider the image $p$ of $(\sig_a,\exp(i\sig_{a,\R})F_b)\in D_{\sigma_a}^{\sharp}$ in one of $D_{\SL(2)}$,  $D_{\SL(2),\val}$, $D^{\star}_{\SL(2)}$,  $D^{\star}_{\SL(2),\val}$.  In the case $a=0$, consider also the image in $D^{\diamond}_{\SL(2)}$.

(1) $p$ remembers $a$ in the cases of  $D_{\SL(2)}$, $D_{\SL(2),\val}$, $D^{\star}_{\SL(2),\val}$, but $p$ does not remember $a$ in the other cases.
 In the case $a\neq 0$, $p$ does not remember $b$ in any of these cases.

(2) Assume $a=0$. Then $p$ remembers $b$ in the case of $D^{\diamond}_{\SL(2)}$, but $p$ does not remember $b$ in all other cases. 
\end{sbpara}

$$\begin{matrix} D^{\sharp,\mild}_{\Xi} &\to & D^{\diamond}_{\SL(2)}\\
\downarrow &&\downarrow\\
D^{\sharp}_{\Xi} &\to& D_{\SL(2)}
\end{matrix}$$

\subsection{Degeneration and regulator maps}\label{ss:reg} 

\begin{sbpara}\label{regdelta}

Let $X$ be a proper smooth variety over $\C$. Let $n\geq 1, r\geq 0$. Then we have the ($r$-th) regulator map (\cite{Be1})
$$\text{reg}_X: K_n(X) \to \bigoplus_{p,q} \; (H^m(X)(r)_{\C, p,q})^{-},$$ 
where $m=2r-n-1$ and $(p,q)$ ranges over all elements of $\Z^2$ such that $p+q=m-2r$ and $p<0,q<0$, $H^m(X)(r)_{\C,p,q}$ is the $(p,q)$-Hodge component of $H^m(X, \C)$ with respect to the Hodge structure $H^m(X)(r)$, and $(-)^-$ denotes the minus part for the complex conjugation which fixes the image of $H^m(X, \Z(r))= H^m(X, \Z)\otimes (2\pi i)^r\Z$.  

This regulator map is understood as $\delta$ (Section 1.2) of a mixed Hodge structure as follows. 
  An element $Z\in K_n(X)$ determines a mixed Hodge structure $H_Z$ with an exact sequence $0\to H^m(X)(r) \to H_Z \to \Z\to 0$. We have 
$$\text{reg}_X(Z)= \delta_W(H_Z),$$
where $W$ is the weight filtration of $H_Z$.
\end{sbpara}

\begin{sbpara} Let $X\to S$, $0\in S$, $n, r, m$  be as in 0.3, and let $Z\in K_n(X\smallsetminus X_0)$.

For $t\in S\smallsetminus \{0\}$, let $Z(t)\in K_n(X_t)$ be the pull back of $Z$. 

Then the regulator $\text{reg}(Z(t))$  is understood as $\delta_W(H_Z(t))$ where  $H_Z$ denotes the variation of mixed Hodge structure on $S\smallsetminus \{0\}$ defined by $Z$ which has an exact sequence 
$0\to H^m(X/S)(r) \to H_Z \to \Z \to 0$ with $H^m(X/S)$ as in 0.3 and whose fiber $H_Z(t)$ at $t$ is the mixed Hodge structure in \ref{regdelta} associated to $Z(t)$.  This $H_Z$  is admissible along $X_0$ and extends uniquely to a log mixed Hodge structure on $S$ which we denote by the same letter $H_Z$. 

Hence the  
behavior of $t\mapsto \text{reg}(Z(t))$ in the degeneration is explained by the theory of degeneration of mixed Hodge structure as in this paper. 

For the details of what follows, see \cite{BK}. 
 
\end{sbpara}

\begin{sbprop}\label{KXmild} Assume $Z$ comes from $K_n(X)$. Then the log mixed Hodge structure $H_Z$ on $S$ is mild. 

\end{sbprop}

\begin{pf} The
Clemens-Schmid sequence  
$H^m(X_0,\Q)\to H^m(X/S)_{\Q,t}\overset{N}\to H^m(X/S)_{\Q,t} \to H_{2d-m}(X_0,\Q)$ $(t\in S \smallsetminus\{0\}$ is near to $0$) induces an injection  
$H^m(X/S)_{\Q,t}/NH^m(X/S)_{\Q,t} \to H_{2d-m}(X_0,\Q)$. Here $N$ is the monodromy logarithm of $H_{Z,\Q}$ at $0\in S$. 
We have a commutative diagram 
$$\begin{matrix} K_n(X\smallsetminus X_0) &\overset{\partial}\to & K_{n-1}'(X_0) \\
\downarrow &&\downarrow \\
H^m(X/S)_{\Q,t}/NH^m(X/S)_{\Q,t} &\overset{\subset}\to& H_{2d-m}(X_0,\Q).\end{matrix}$$
Here the left vertical arrow sends $Z\in K_n(X\smallsetminus X_0)$ to $Ne$ where $e$ is the lifting of $1\in \Q$ to $H_{Z,\Q,t}$ under the exact sequence $0\to H^m(X/S)_{\Q,t}\to H_{Z,\Q,t}\to \Q\to 0$. $K'_{n-1}$ denotes the $K$-group of coherent sheaves. 
The right vertical arrow is the topological Chern class map. 
By the localization theory of $K$-theory, we have an exact sequence 
$K_n(X) \to K_n(X\smallsetminus X_0) \overset{\partial}\to K'_{n-1}(X_0)$.

Assume $Z\in K_n(X\smallsetminus X_0)$ comes from $K_n(X)$. Then $\partial(Z)=0$ and hence the above diagram shows that the image of $Z$ in $H^m(X/S)_{\Q,t}/NH^m(X/S)_{\Q,t}$ is zero. This proves that $(W, N)$ splits. 
\end{pf}

By \ref{KXmild} and by \ref{diathm}, we have

\begin{sbthm}\label{thm2} If $Z\in K_n(X\smallsetminus X_0)$ comes from $K_n(X)$, the regulator $\text{reg}_{X_t}Z(t)$ $(t\in S\smallsetminus \{0\})$ converges when $t\to 0$.

\end{sbthm}

\begin{sbrem} In \cite{BK}, this result \ref{thm2} will be generalized to the situation $S$ need not be of dimension $\leq 1$. This generalization will be reduced to \ref{thm2} by using \ref{Ccut}. 

\end{sbrem}

\begin{sbpara}\label{7.2.2}
Let $X\to S$ and $0\in S$ be as in Section 0.3. Take $H_0= \Z \oplus H^m(X/S)(r)_{\Z,t}$. The extended period maps in \S6.3 give  a commutative diagram
$$\begin{matrix} S^{\log} \times K_n(X) &\to& \Gamma\bs D^{\diamond}_{\SL(2)}\\
\downarrow &&\downarrow\\
S^{\log}\times K_n(X\smallsetminus X_0) &\to& \Gamma\bs D_{\SL(2)}
\end{matrix}$$
Here $\Gamma$ is the group of all elements $\gamma$ of $\Aut_\Z(H_0)$ satisfying the following conditions. (i) $\gamma$ preserves $H^m(X/S)(r)_{\Z,t}$. (ii) $\gamma e \equiv e \bmod H^m(X/S)(r)_{\Z,t}$ where $e$ denotes $(1,0)\in H_0$. (iii) The action of $\gamma$ on $H^m(X/S)_{\Z,t}$ is contained in the monodromy group of $H^m(X/S)_{\Z}$. 
\end{sbpara}

\begin{sbpara}\label{ell1} We give an explicit example. 
Assume that $X\to S$ is a family of elliptic curves which degenerates at $0\in S$, and assume $n=r=2$. 

Then the period domain and the extended period domains which appear here are those of Example III (Section \ref{III}). Let the notation be as in \ref{III}. 

We discuss an explicit example of $Z\in K_2(X\smallsetminus X_0)$. 

Let $\alpha_j$ and $\beta_k$ be a finite number of torsion sections of $X\smallsetminus X_0$ over $S\smallsetminus \{0\}$, let $m_j, n_k\in \Z$ such that $\sum_j m_j=\sum_k n_k=0$, and consider the divisors $\alpha=\sum_j m_j(\alpha_j)$, $\beta=\sum_k n_k (\beta_k)$ on $X\smallsetminus X_0$ of degree $0$. 
Then we have an element $Z_{\alpha,\beta}\in K_2(X\smallsetminus 
X_0)$ (see \cite{Bl1}, \cite{Sct}).  It is essentially the Steinberg symbol $\{f_{\alpha}, f_{\beta}\}$, where $f_{\alpha}$ (resp.\ $f_{\beta}$) is an element of $\Q\otimes \C(X)^\times$ whose divisor is $\alpha$ (resp.\ $\beta$). 
When $t$ tends to $0$ in $S$, there are $a,b\in W_{-3}H_{0,\R}$ such that we have
$$\text{reg}_{X_t}(Z_{\alpha,\beta}(t))= ay+b +O(y^{-1}),$$
where $y$ is defined by $q(t)=e^{2\pi i(x+iy)}$ $(x,y\in \R)$ with $q(t)$ the $q$-invariant of the elliptic curve $X_t$.

We have
$$a\equiv \sum_{j,k} m_jn_k B_3(\{r(\alpha_j)-r(\beta_k)\})e_2 \bmod \R e_1,$$
where:

$B_3$ is the Bernoulli polynomial of degree $3$, 

$r(\mu)$ for a torsion section $\mu$ is the element of $\Q/\Z$ such that as a section of the Tate curve ${\mathbb G_m}/q^{\Z}$, $\mu$ is expressed as $s q^{r(\mu)}\bmod q^{\Z}$ with $s$ a root of $1$, 

$\{-\}: \Q/\Z\to [0, 1)\subset \Q$ is the lifting.

Assume now that $r(\alpha_j)=r(\beta_k)=0$ for any $j,k$, that is, these torsion sections $\alpha_j$ and $\beta_k$ are roots of $1$ in the Tate curve ${\mathbb G}_m/q^{\Z}$. Then $Z_{\alpha, \beta}$ comes from $K_2(X)$, $a=0$, and the degeneration is mild. In this case, 
$$b= \sum_{j,k} m_jn_k D(\alpha_j/\beta_k)e_1,$$
where we regard $\alpha_j$ and $\beta_k$ as roots of $1$ and $D$ is the real analytic modified dilogarithm function (\cite{Bl1} of Bloch-Wigner. 

These things will be explained in \cite{BK} by using results in \cite{Sct},  \cite{Bl1}, \cite{GL} and using the results of this paper. 
\end{sbpara}

\subsection{Example IV}\label{IV}
This is Example IV of Part II. As is explained in Part II, 4.4 and also in Section \ref{ss:reBK} below, this example is related to the local height pairing of points of a degenerating elliptic curve. 

\begin{sbpara}
Let $H_0=\Z^4$ with basis $e_1,e_2,e_3, e_4$. The weight filtration is given by $$W_{-3}=0\subset W_{-2}=\R e_1\subset W_{-1}=
\bigoplus_{j=1}^3 \R e_j\subset W_0=H_{0,\R}.$$ 
The intersection form on $\gr^W_{-1}$ is the anti-symmetric form characterized by $\langle e_3, e_2\rangle=1$. 
\end{sbpara}

\begin{sbpara} We have 
$D(\gr^W)\cong \frak h$, the upper half plane, and 
$D_{\SL(2)}(\gr^W)=D_{\BS}(\gr^W)={\frak h}_{\BS}$. 

\end{sbpara}

\begin{sbpara}
We have a homeomorphism
$$D^{\star}_{\SL(2),\val}\overset{\cong}\to D_{\BS,\val}$$
and this induces a homeomorphism $$D^{\star}_{\SL(2)}\cong D_{\BS}$$ of quotient spaces. 

Let $W'$ be the increasing filtration on $\gr^W$ given by 
$$
W'_{-3}=0\subset W'_{-2}=\gr^{W}_{-2}+\bR (e_2\bmod W_{-2})  =W'_{-1}
\subset W'_0=\gr^W,
$$
and let $\Phi=\{W'\}$.
Let $P$ be the parabolic subgroup of $G_\bR$ consisting of elements which preserve $W'$.
Then, 
$D_{\BS}(P)=D^{\star}_{\SL(2)}(\Phi)$ and it is the inverse image of the open set $\{x+iy\;|\;x\in \R, y\in (0, \infty]\}$ of ${\frak h}_{\BS}$ under the projection $D_{\BS}=D^{\star}_{\SL(2)}\to D_{\BS}(\gr^W)=D_{\SL(2)}(\gr^W)={\frak h}_{\BS}$.

We have $$D^I_{\SL(2)}= D^{II}_{\SL(2)}$$ 
(Part II, 4.4).  We will denote both of them by $D_{\SL(2)}$.

We have a canonical homeomorphism 
$$D^{\diamond}_{\SL(2)}\overset{\cong}\to D^{\star,\mild}_{\SL(2)}.$$

\end{sbpara}

\begin{sbpara} We have
$$\spl(W)\cong \R^5, \quad \cL\cong \R,$$
where in the first isomorphism $(s_{3,4}, s_{2,4}, s_{1,4}, s_{1,3}, s_{1,2})\in \R^5$ corresponds to the splitting of $W$ which is given by $e_4+\sum_{j=1}^3 s_{j,4}e_j$ and  $e_k+s_{1,k}e_1$  $(k=2,3)$, and the second isomorphism is given by  $\delta\mapsto r$ $(\delta\in \cL$, $r\in \R)$, $\delta e_4=re_1$.
We have $\cL(F)=\cL$ for any $F\in D(\gr^W)$.

\end{sbpara}

\begin{sbpara}\label{IVdia}
We have homeomorphisms
$$D\cong  {\frak h}\times \cL \times \spl(W) \cong  \R_{>0}\times \R \times \R^6,$$
where 
the left isomorphism  is $F\mapsto (F(\gr^W), \delta_W(F), \spl_W(F))$, and the second isomorphism
 sends $(x+iy, \delta, s)$ to $(1/\sqrt{y}, \delta, x, s)$, where $x, y\in \R$, $y>0$. 
We call the composite homeomorphsm $D\cong \R_{>0}\times \R \times \R^6$ the {\it standard isomorphism for $D$}.

We have a commutative diagrams of homeomorphisms
$$\begin{matrix} D^{\star}_{\SL(2)}(\Phi) && \cong && \R_{\geq 0} \times [-\infty,\infty] \times \R^6\\
\uparrow &&&& \uparrow\;\;\;\;\; \\
D^{\star}_{\SL(2),\val}(\Phi) && \cong && (\R_{\geq 0}\times [-\infty, \infty])_{\val} \times \R^6\\
\downarrow &&&& \downarrow (1) \\
D_{\SL(2),\val}(\Phi) && \cong && (\R_{\geq 0} \times [-\infty, \infty])_{\val} \times \R^6\\
\downarrow &&&& \downarrow\;\;\;\;\; \\
D_{\SL(2)}(\Phi) && \cong && \R_{\geq 0}\times [-\infty,\infty] \times \R^6,
\end{matrix}$$
where the upper two homeomorphisms are compatible with the standard isomorphism for $D$, but via the lower two
homeomorphisms, 
a point of $D$ corresponding to $(t, \delta, u)\in \R_{>0}\times \R \times \R^6$ under the standard isomorphism for $D$ is sent to 
$ (t, t^2\delta, u)$. 

Concerning the vertical arrows on the right-hand side, the arrows except (1) are the canonical projections, and the arrow (1) is as follows:

The map (1)
sends  

$(t, \delta, u)$ $((t,\delta)\in (\R_{\geq 0}\times \R)\cup (\R_{>0}\times [-\infty,\infty]))$ to $(t, t^2\delta, u)$, 

$(p(c,\pm \infty), u)$ $(c\in \R_{>0}\smallsetminus \Q_{>0})$ to $ (p(c-2,\pm \infty), u)$ if $c>2$ and to $(0, 0, u)$ if $c<2$, 

$(p(c+, \pm \infty),u)$ $(c\in \Q_{\geq 0})$ to $(p((c-2)+,\pm \infty), u)$ if $c\geq 2$ and to $(0,0,u)$ if $c<2$,

 $(p(c-, \pm\infty), u)$ $(c\in \Q_{>0})$ to $(p((c-2)-, \pm \infty), u)$ if $c>2$ and to $(0,0, u)$ if $c\leq 2$,
 
 $(p(c, \delta),u)$ $(c\in \Q_{>0}$, $\delta \in \R\smallsetminus \{0\})$ to $(p(c-2, \delta),u)$ if $c>2$, to $(0, \delta, u)$ if $c=2$, and to $(0,0,u)$ if $c<2$. 
 
 Here the notation $p(c, \delta)$ etc.\ are understood as in \ref{IIIdia} by replacing $\bar V\supset V$ by $[-\infty, \infty]\supset \R$. 
  \end{sbpara}

\begin{sbpara} 
Let $a\in \Q_{\ge0}$ and define $N_a\in \fg_\Q$ by
$N_a(e_4)=ae_1$, $N_a(e_3)=e_2$, and $N_a(e_1)=N_a(e_2)=0$. Let $\sig_a$ be the cone generated by $N_a$.

For $b\in \R$, let  $F_b\in \Dc$ be the decreasing filtration defined as follows: $F_b^1=0$, $F_b^0$ is generated by $e_3$ and $e_4+ibe_1$, and $F_b^{-1}$ is the total space.

In $D^{\sharp}_{\sig_a}$, we have the limit $(\sig_a, \exp(i\sig_{a,\R})F_b)\in D^{\sharp}_{\sigma_a}$ of $\exp(iyN_a)F_b$ for $y\to \infty$. 
This $(\sig_a, \exp(i\sig_{a,\R})F_b)$ belongs to $D^{\sharp,\mild}_{\sig_a}$ if and only if $a=0$. 

For $y\in \R_{>0}$, via the first homeomorphism in the diagram in \ref{IVdia}, 
$\exp(iyN_a)F_b\in D$ is sent to $(1/\sqrt{y}, ay+b, \bold0)$ and hence
the limit $(\sig_a, \exp(i\sig_{a,\R})F_b)\in D_{\sig_a}^{\sharp}$ is sent to $(0, \infty, \bold0)$ in the case $a\neq 0$, and to $(0, b, \bold0)$ in the case $a=0$. Here $\bold0$ denotes $(0,\dots,0)\in \R^6$. 

On the other hand,  for $y\in \R_{>0}$, via the last homeomorphism in the diagram in \ref{IVdia}, $\exp(iyN_a)F_b\in D$ is sent to $(1/\sqrt{y}, a+y^{-1}b, \bold0)$ and hence 
the limit $(\sig_a, \exp(i\sig_{a,\R})F_b)\in D_{\sig_a}^{\sharp}$ is sent to $(0, a, \bold0)$. 

By taking the limit for $y\to \infty$, we have:
\end{sbpara}

\begin{sblem}\label{IVconv}
The limit of $\exp(iyN_a)F_b$ for $y\to \infty$ exists also in $D^{\star}_{\SL(2),\val}(\Phi)$,  and in $D_{\SL(2),\val}(\Phi)$. In the case $D^{\star}_{\SL(2),\val}(\Phi)$, the $(\R_{\geq 0}\times [-\infty, \infty])_{\val}$-component of the limit is $p(2, a)$ if $a\neq 0$ 
and is $(0, b)$ if $a=0$. In the case $D_{\SL(2),\val}(\Phi)$, the $(\R_{\geq 0}\times [-\infty, \infty])_{\val}$-component of the limit is $p(0, a)$.

\end{sblem}

\begin{sbpara}\label{ex1} By Lemma \ref{IVconv}, we have 
  the following.

Let $p$ be the image of $(\sig_a, \exp(i\sig_{a,\R})F_b)\in D_{\sigma_a}^{\sharp}$ in one of $D_{\SL(2)}$,  $D_{\SL(2),\val}$, $D^{\star}_{\SL(2)}$,  $D^{\star}_{\SL(2),\val}$. 
We have:

(1) $p$ remembers $a$ in the cases of $D_{\SL(2)}$, $D_{\SL(2),\val}$ and $D^{\star}_{\SL(2),\val}$, whereas $p$ does not remember $a$ in the case of $D^{\star}_{\SL(2)}$.

In the case $a\neq 0$,  $p$ does not remember $b$ in any of these cases.

(2) Assume $a=0$. Then $p$ remembers $b$ in the cases of $D^{\star}_{\SL(2)}$ and $D^{\star}_{\SL(2),\val}$, but $p$ does not remember $b$ in the cases of $D_{\SL(2)}$ and $D_{\SL(2),\val}$.

\end{sbpara}

\begin{sbpara}\label{c.ex} 
The following is mentioned in \ref{r1rem}. 

Though $\sig_a$ is of rank one, 
the image $p$ of $D^{\sharp}_{\sig_a}$ in $D^{\star}_{\SL(2)}$ is not contained in $D^{\star}_{\SL(2),\leq 1}$ 
($=$ the part of $D^{\star}_{\SL(2)}$ at which  the log structure $M$ satisfies $\text{rank}(M^{\gp}/\cO^\times)_p \leq 1)$  if $a\neq 0$. Indeed, in the case $a\neq 0$, the image of class\,$(N_a, F_b)\in D^{\sharp}_{\sig_a}$ in $D^{\star}_{\SL(2)}$ has the coordinate $(0, \infty, \bold0)$, which shows that $(M/\cO^\times)_p \cong \bN^2$.
\end{sbpara}

\subsection{Degeneration and height pairings}\label{ss:reBK}
\begin{sbpara}\label{hpdelta}

Let $X$ be a proper smooth algebraic variety over $\C$ of dimension $d$, and let $Y$ and $Z$ be algebraic cycles on $X$ of codimenson $r$ and $s$, respectively. We assume that $r+s=d+1$, that their supports are disjoint $|Y|\cap |Z|=\emptyset$, and that both $Y$ and $Z$ are homologically equivalent to $0$. Then we have a height pairing $\langle Y,Z \rangle_X\in \R$
(the local version of the height pairing for number field, at an Archimedean place). See \cite{Be2}, \cite{Bl2}. 

This height pairing is understood as $\delta$ (Section 1.2) of a mixed Hodge structure.  
We have $$\langle Y, Z\rangle_X = \delta_W(H_{Y,Z}),$$
where $H_{Y,Z}$ is the mixed Hodge structure whose weight filtration $W$ has the following properties:
$W_0H_{Y,Z}=H_{Y,Z}$, $W_{-3}=0$, $\gr^W_0=\Z$, $\gr^W_{-2}=\Z(1)$, $\gr^W_{-1}= H^{2r-1}(X)(r)$, constructed in \cite{Be2}, \cite{Bl2}. 
  The exact sequence $0\to H^{2r-1}(X)(r) \to W_0/W_{-2}\to \Z\to 0$ is given by the class of $Y$, and the exact sequence $0\to \Z(1) \to W_{-1}\to H^{2r-1}(X)(r)\to 0$ is given by the class of $Z$. 
\end{sbpara}

\begin{sbpara}
Let $X\to S$ and $0\in S$ be as in Section \ref{relBK}. 
  Let $Y$ and $Z$ be algebraic cycles on $X$ of codimension $r$ and $s$, respectively, such that $r+s=d+1$, where $d$ is the relative dimension of $X\to S$, such that $|Y|\cap |Z|=\emptyset$, and such that both $Y(t)$ and $Z(t)$ are homologically equivalent to $0$ for any $t\in S\smallsetminus \{0\}$. 
\end{sbpara}

\begin{sbpara} Since the height paring $\langle Y(t), Z(t)\rangle$ is understood as $\delta$ of mixed Hodge structure(\ref{hpdelta}), its behavior in the degeneration is explained by the theory of degeneration of mixed Hodge structure as in this paper. 

When $t\to 0$ with $x$ fixed, there are $a,b\in \R$ such that we have
$$\langle Y(t), Z(t)\rangle_{X_t} = ay+b+O(y^{-1})$$
where taking a local coordinate $q$ on $S$ at $0$ such that $q(0)=0$, we define $y$ by $q=e^{2\pi i (x+iy)}$ $(x,y\in \R,\, y>0)$. This follows by a general theory of degeneration of mixed Hodge structure as studied in Section \ref{IV}. 
  Here, $a=0$ if and only if the degeneration of $H_{Y,Z}$ at $0\in S$ is mild.

In \cite{BK}, it is shown that $a$ is the local geometric  intersection number of $Y$ and $Z$ over $0\in S$. 
\end{sbpara}

\begin{sbpara}\label{ell1b}  We give an explicit example. 
Assume  that $X\to S$ is a family of degenerating elliptic curves, and assume $r=s=1$. Let $Y=\sum_j m_j(\alpha_j)$, $Z=\sum_k n_k(\beta_k)$ where  $\alpha_j$ and $\beta_k$ are closures in $X$ of torsion sections of $X\smallsetminus X_0 \to S\smallsetminus \{0\}$, $m_j, n_k\in \Z$, $\sum_j m_j=\sum_k n_k=0$. We assume that the divisors $\alpha_j$ and $\beta_k$ of $X$ do not intersect for any pair $(j,k)$. 

This is an example discussed at the end of Part II, 4.4.
The extended period domains which appear here are those of Example IV (Section \ref{IV}).

  We have
$$a= \sum_{j,k} m_jn_k B_2(\{r(\alpha_j)-r(\beta_k)\}),$$
where $B_2$ is the Bernoulli polynomial of degree $2$. (The notation is as in \ref{ell1}.)
This was explained in Part II, Proposition 4.4.8.

If $r(\alpha_j)=r(\beta_k)=0$ for any $j,k$, then $a=0$ and the degeneration is mild. In this case, 
$$b= \sum_{j,k} m_jn_k l(\alpha_j/\beta_k),$$
where we regard $\alpha_j$ and $\beta_k$ as roots of $1$ and $l(t)=\log(|1-t|)$. 

These things are surprisingly similar to \ref{ell1}.

These things will be explained in \cite{BK}  more using the results of this paper. 

\end{sbpara}
\setcounter{section}{0}
\def\thesection{\Alph{section}}

\section{Corrections to \cite{KU}, supplements to Part III}\label{s:co}

Corrections of errors in the book \cite{KU} have been put in the home page of Princeton University Press. 
In this Section A, we update them. 
In Section A.1 and Section A.2, we describe important corrections to \cite{KU}.
In Section A.3, we give supplements to Part III. 
Other errors described in the above home page are minor ones.

Section \ref{ss:8.1} and Section A.3 are important for \ref{ss:Lthm} in the text. 

\subsection{Change on \cite{KU} \S6.4}\label{ss:8.1}

\begin{sbpara} The following errors (1) and (2) are in \cite{KU} Sections 6.4 and 7.1, respectively.

(1) 
Proposition 3.1.6 is used in 6.4.12 (line 2 from the end), but Proposition 3.1.6 is not strong enough for the arguments in 6.4.12.

(2) We can not have the second convergence in  7.1.2 (3).

\end{sbpara}

\begin{sbpara}

We make the following changes 1, 2, 3 on the book  \cite{KU}. The change 1 solves the above problem (1). The changes 2 and 3 solve the above problem (2).

Change 1. We replace \cite{KU} 7.1.2 by \ref{8.1.a}--\ref{8.1.b} below.

Change 2. We move \cite{KU} Section 7.1, revised as in the above Change 1, to the place just before \cite{KU} Section 6.4. That is, we exchange the order of the Sections 6.4 and 7.1.

Change 3.  We make the change on \cite{KU} 6.4.12 explained in \ref{8.1.c}.

 \end{sbpara}

 \begin{sbpara}\label{8.1.a}

We will prove Theorem A (i), that
 $E_\sig$ is open in $\tilde E_{\sig}$ for the strong topology. 

Since $\tilde E_{\sig,\val}\to\tilde E_\sig$ is proper surjective and $E_{\sig,\val}\sub\tilde E_{\sig,\val}$ is the inverse image of $E_\sig\sub\tilde E_\sig$, 
it is sufficient to prove that $E_{\sig,\val}$ is open in $\tilde E_{\sig,\val}$.

Assume $x_{\lam}=(q_\lam, F'_\lam)\in \tilde E_{\sig,\val}$ converges in $\tilde E_{\sig,\val}$ to $x=(q, F')\in E_{\sig,\val}$. We prove that $x_{\lam}\in  E_{\sig,\val}$ for any sufficiently large $\lam$. 

\end{sbpara}

\begin{sbpara}

We fix notation.

Let $|\;\;|:\toric_{\sig,\val}\to\abtoric_{\sig,\val}$
be the canonical projection induced by $\bC\to\bR$,
$z\mapsto|z|$.

Let $(A,V,Z)\in D^\sharp_{\sig,\val}$ be the image
of $(|q|, F')\in E^{\sharp}_{\sig,\val}$ under $E^{\sharp}_{\sig,\val}\to D^{\sharp}_{\sig,\val}$ (5.3.7),  and take an excellent basis $(N_s)_{s\in S}$ for
$(A,V,Z)$ such that $N_s\in\sig(q)$ for any $s$ (6.3.9).
Let $S_j$ $(1\le j\le n)$ be as in 6.3.3.
Take an $\bR$-subspace $B$ of $\sig_\bR$ such that
$\sig_\bR=A_\bR\op B$.

We have a unique injective open continuous map
$$
(\bR_{\ge0}^S)_\val\x B\to\abtoric_{\sig,\val}
$$
which sends $((e^{-2\pi y_s})_{s\in S},b)$ $(y_s\in\bR,
b\in B)$ to $\be((\ts_{s\in S}iy_sN_s)+ib)$ (cf. 3.3.5).
Let $U$ be the image of this map.
Define the maps $t_s:U\to\bR_{\ge0}$ $(s\in S)$ and
$b:U\to B$ by
$$
(t,b)=((t_s)_{s\in S},b):
U\simeq(\bR_{\ge0}^S)_\val\x B\to\bR_{\ge0}^S\x B
$$
(an abuse of notation $b$).
We have  $|q|\in U$ and $t(|q|):=(t_s(|q|))_{s\in S}=0$.
Since $|q_\lam|\to |q|$, we may assume $|q_\lam|\in U$.
Let
 $$F_{\lam}=\exp(i b(q_{\lam}))F'_{\lam}, \quad F=\exp(ib(q))F'.$$
Then $((N_s)_{s\in S}, F)$ generates a nilpotent orbit.
\end{sbpara}

\begin{sbpara}

We may assume that, for some $m$ $(1\leq m\leq n+1)$, $t_s(q_\lam)=0$ for any $\lam$ and $s\in S_{\leq m-1}$ and $t_s(q_\lam)\neq 0$ for any $\lam$ and $s\in S_{\geq m}$ (6.3.11), ($S_{\leq 0}$ and $S_{\geq n+1}$ are defined as the empty set). 

Take $c_j\in S_j$ for each $j$. 
For $s\in S_{\geq m}$, define $y_{\lam,s}\in \bR$ by $t_s(q_\lam)=e^{-2\pi y_{\lam,s}}$.
For each $j\in \bZ$ such that $m\leq j \leq n$, let $N_j= \ts_{s\in S_j} a_sN_s$ 
where $a_s\in \bR$ is the limit of $y_{\lam,s}/y_{\lam,c_j}$. 

Then $(N_1, \dots, N_n,F)$ generates a nilpotent orbit. 
Let $\tilde \rho: \bG_{m,\bR}^n\to G_\bR$ be the homomorphism of the $\SL(2)$-orbit $(5.2.2)$ associated to $(N_1, \dots, N_n,  F)$. 
For $m \leq j\leq n$, let
$$
e_{\lam,\geq j}=\exp(\ts_{s\in S_{\geq j}} iy_{\lam,s}N_s)\in G_\bC,
$$ 
$$
\tau_{\lam,j}= 
\tilde \rho_j\left(\sqrt{y_{\lam,c_{j+1}}/y_{\lam,c_j}}\right)\in G_\bR, 
\quad  \tau_{\lam, \geq j}
 = \tp_{k= j}^n \tau_{\lam,k}\in G_\bR
$$
$(y_{\lam,c_{n+1}}$ denotes $1$). 
Here $\tilde\rho_j$ is the restriction of $\tilde\rho$ to the $j$-th factor of $\bG_{m,\R}^n$.
Let $\hat F_{(j)}$ ($1\leq j\leq n$) be as in 6.1.3 associated to $(N_1,\dots, N_n, F)$.

By 3.1.6 applied to $S=E_{\sig}\subset X=\Ec_{\sig}$, we have

\end{sbpara}

\begin{sblem}\label{lem712} Let the situation and the notation be as above. 
Let $m\leq j \leq n$ and let $e\geq 0$. Then for any sufficiently large $\lam$, there exist $F^*_{\lam}\in \Dc$ satisfying the following (i) and (ii).

\medskip

$(i)$ $y_{\lam,s}^ed(F_{\lam}, F^*_{\lam})\to 0 \;(\forall s\in S_j)$. 

\smallskip

$(ii)$ $(N_s, F^*_{\lam})$ {\it satisfies Griffiths transversality for any}\; $s\in S_{\leq j}.$

\medskip
\noindent
Furthermore, in the case $j=n$, there is $F^*_{\lam}$ as above satisfying the following condition $(ii)^*$ which is stronger than the above  condition (ii).

\medskip

$(ii)^*$ $((N_s)_{s\in S}, F^*_{\lam})$ generates a nilpotent orbit.

\end{sblem}

\begin{sbprop}\label{prop712}
Let the situation and the assumption be as above.
 Then the following assertions $(A_j)$ $(m-1\leq j\leq n)$, $(B_j)$ $(m\leq j\leq n)$, $(C_j)$ $(m\leq j\leq n)$ are true.

$(A_j)$ (resp.\ $(B_j)$, resp.\ $(C_j)$) for $m\leq j\leq n$:
Let $e\geq 1$. Then for any sufficiently large $\lam$, there are $F^{(j)}_{\lam}\in \Dc$ satisfying the following (1)--(3).

$(1)$ $y^e_{\lam,j}d(F_{\lam}, F^{(j)}_{\lam}) \to 0$.

$(2)$ $((N_s)_{s\in S_{\leq j}}, e_{\lam,\geq j+1}F^{(j)}_{\lam})$ generates a nilpotent orbit.  

$(3)$  $\tau_{\lam,\geq j+1}^{-1}e_{\lam,\geq j+1} F^{(j)}_{\lam}\to \exp(iN_{j+1})\hat F_{(j+1)}$. 

(resp.\ $(3)$ $\tau_{\lam,\geq j}^{-1}e_{\lam,\geq j+1} F^{(j)}_{\lam}\to \hat F_{(j)}$. 

resp.\ $(3)$ $\tau_{\lam,\geq j}^{-1}e_{\lam,\geq j}F^{(j)}_{\lam}\to \exp(iN_j)\hat F_{(j)}$.)

\noindent
Here $(A_n)$ is formulated by understanding $N_{n+1}=0$ and $\hat F_{(n+1)}=F$. 

\medskip

$(A_{m-1})$: For any sufficiently large $\lam$, we have the following (2) and (3). 

$(2)$ $((N_s)_{s\in S_{\leq m-1}}, e_{\lam,\geq m}F_{\lam})$ generates a nilpotent orbit.  

$(3)$ $\tau_{\lam,\geq m}^{-1}e_{\lam,\geq m} F_{\lam}\to \exp(iN_m)\hat F_{(m)}$. 

\end{sbprop}

\begin{sbpara}
We prove Proposition \ref{prop712} by using the downward induction of the form

$(A_j)$ $\Rightarrow$ $(B_j)$ $\Rightarrow$ $(C_j)$ $\Rightarrow$ $(A_{j-1})$. (Here $m\leq j\leq n$.)

$(B_j)$ $\Rightarrow$ $(C_j)$ is clear. 

$(A_j)$ $\Rightarrow$ $(B_j)$ is easy.

$(A_n)$ follows from \ref{lem712}. 

We prove $(C_{j+1})$ $\Rightarrow$ $(A_j)$.

By \ref{lem712}, if $m\leq j\leq n$ (resp. $j=m-1$), there are $F^{(j)}_{\lam}\in \Dc$ satisfying (1) and (resp. $F^{(j)}_{\lam}:=F_{\lam}$) satisfies the condition 

\medskip

$(2')$ $(N_s, F^{(j)}_{\lam})$ satisfies Griffiths transversality for any $s\in S_{\le j}$.

\medskip

By $(C_{j+1})$, there are $F^{(j+1)}_{\lam}\in \Dc$ satisfying

\medskip

$(1'')$ $y^e_{\lam,j+1}d(F_{\lam}, F^{(j+1)}_{\lam}) \to 0$.

$(2'')$ $((N_s)_{s\in S_{\leq j+1}}, e_{\lam,\geq j+2}F^{(j+1)}_{\lam})$ generates a nilpotent orbit.  

$(3'')$  $\tau_{\lam,\geq j+1}^{-1}e_{\lam,\geq j+1} F^{(j+1)}_{\lam}\to \exp(iN_{j+1})\hat F_{(j+1)}$.

\medskip

By $(1'')$ and $(3'')$, we have

\medskip

(4) $\tau_{\lam,\geq j+1}^{-1}e_{\lam,\geq j+1} F_{\lam}\to \exp(iN_{j+1})\hat F_{(j+1)}$.

\medskip

By (4) and by (1), we have 

\medskip

(5) $\tau_{\lam,\geq j+1}^{-1}e_{\lam,\geq j+1} F^{(j)}_{\lam}\to \exp(iN_{j+1})\hat F_{(j+1)}$.

\medskip

Concerning to the left-hand side of (5), by $(2')$, $(N_s, \tau_{\lam,\geq j+1}^{-1}e_{\lam,\geq j+1} F^{(j)}_{\lam})$ satisfies Griffiths transversality for any $s\in S_{\leq j}$. On the other hand, concerning the right-hand side of (5), 
$((N_s)_{s\in S_{\leq j}}, \exp(iN_{j+1})\hat F_{(j+1)})$ generates a nilpotent orbit. Hence (5) and 7.1.1 show that
$((N_s)_{s\in S_{\leq j}}, \tau_{\lam,\geq j+1}^{-1}e_{\lam,\geq j+1} F^{(j)}_{\lam})$ generates a nilpotent orbit. This proves that $((N_s)_{s\in S_{\leq j}}, e_{\lam,\geq j+1} F^{(j)}_{\lam})$ generates a nilpotent orbit for any sufficiently large $\lam$. Hence for any sufficiently large $\lam$, $(W^{(j)}, e_{\lam,\geq j+1} F^{(j)}_{\lam})$ is a mixed Hodge structure where $W^{(j)}$ denotes the relative monodromy filtration of $N_1+\dots+N_j$ with respect to $W$. By this and by (5), we have 
$(A_j)$. 

\end{sbpara}

\begin{sbpara}\label{8.1.b} By
$(A_{m-1})$ (2) of Proposition \ref{prop712}, $x_{\lam}$ belongs to $E_{\sig,\val}$ if $\lam$ is sufficiently large. This proves that $E_{\sig,\val}$ is open in $\tilde E_{\sig,\val}$, and hence proves that $E_{\sig}$ is open in $\tilde E_{\sig}$. 
\end{sbpara}

\begin{sbpara}\label{8.1.c}
In \cite{KU} 6.4.12, in the line 2 from the end, we replace the part 

\medskip

\lq\lq by Proposition 3.1.6"

\medskip
\noindent
with 

\medskip

\lq\lq by the case $m=0$ and $x_{\lam}\in E^{\sharp}_{\sig,\val}$ of Proposition \ref{prop712}"

\medskip
\noindent
of the present paper.

\end{sbpara}

\begin{sbrem}\label{8.1.d}
By these changes, we have the following simplification in \cite{KU}.

We can assume $y^*_{\lam,t}=y_{\lam,t}$ in \cite{KU} Proposition 6.4.1. It is claimed in the original \cite{KU} \S6.4 that 
the proofs of \cite{KU} Theorem 5.4.3 (ii) and Theorem 5.4.4 are reduced to \cite{KU} Proposition 6.4.1, but actually they are reduced to the case $y^*_{\lam,t}=y_{\lam,t}$ of \cite{KU} Proposition 6.4.1.
\end{sbrem}

\subsection{Change on \cite{KU} \S7.2}\label{ss:8.3}

\begin{sbpara} Professor J.-P. Serre kindly pointed out that our book should not use \cite{BS} \S10, for errors are there (cf.\ 10.10. Remark in the version of \cite{BS} contained in 
\lq\lq Armand Borel oevres collected papers, Vol.\ III, Springer-Verlag, 1983").
We used a result \cite{BS} 10.4 in the proof of Lemma 7.2.12 of our book. 
In order to correct our argument, we change as follows.
\medskip

We put the following assumption in Theorem 7.2.2 (i):

\lq\lq Assume that $\sig$ is a nilpotent cone associated to a nilpotent orbit."
\medskip

We replace Lemma 7.2.12 and its proof in our book by the following proposition and its proof which does not use \cite{BS} \S10.

\end{sbpara}

\begin{sbprop}
 Let $\sig$ be a nilpotent cone associated to a nilpotent orbit and let $W(\sig)$ be the associated weight filtration.
Then, by assigning the Borel-Serre splitting, we have a continuous map $E^\sharp_{\sig,\val}\to\spl(W(\sig))$.

\end{sbprop}

\begin{pf}

The composite map $E^\sharp_{\sig,\val}\to D^\sharp_{\sig,\val}\overset{\psi}\to D_{\SL(2)}$ is continuous by the definition of the first map and by 6.4.1 for the CKS map $\psi$.
Let $N_1,\dots,N_n$ be a generator of the cone $\sig$.
Let $s$ be a bijection $\{1, \dots, n\}\to\{1, \dots, n\}$.
Then the image of the map $E^\sharp_{\sig,\val}\to D_{\SL(2)}$ is contained in the union $U$ of $D_{\SL(2)}(\{W(N_{s(1)}+\dots+N_{s(j)})\,|\,j=1,\dots,n\})$ where $s$ runs over all bijections $\{1, \dots, n\}\to\{1, \dots, n\}$.
Since $N_{s(1)}+\dots+N_{s(n)}= N_1+\dots+N_n$, the filtration $W(N_1+\dots+N_n)=W(\sig)$ appears for any $s$.
By \cite{KNU2} Part II, Proposition 3.2.12, the Borel-Serre splitting gives a continuous map $U\to\spl(W(\sig))$.
Thus we get our assertion.
\end{pf}

\begin{sbpara} We replace the third paragraph in 7.2.13 by the following:

\lq\lq Since the action of $\sig_\bR$ on $\spl(W(\sig))$ is proper and $E^\sharp_{\sig,\val}$ is Housdorff, the action of $\sig_\bR$ on $E^\sharp_{\sig,\val}$ is proper by applying Lemma 7.2.6 (ii) to the continuous map $E^\sharp_{\sig,\val}\to\spl(W(\sig))$ in Proposition 7.2.12.
Hence $\text{Re}(h_\lam)$ converges in $\sig_\bR$ by Lemma 7.2.7."

\end{sbpara}

\begin{sbpara}
Add the following sentence at the top of the fourth paragraph in 7.2.13:

\lq\lq Let $|\;\;| : E_{\sig,\val} \to E^\sharp_{\sig,\val}$ be the continuous map $(q,F) \to (|q|,F)$ in 7.1.3."

\end{sbpara}
\subsection{Supplements to Part III}

We add explanations to Part III, Section 3.3.

\begin{sbpara}
We put the following explanation \ref{8.3.a}
just after the statement of Part III, Theorem 3.3.1.

\end{sbpara}

\begin{sbpara}\label{8.3.a}

This Part III, Theorem 3.3.1 is the mixed Hodge version of \cite{KU} Theorem 5.4.3 of the pure case, and is proved in the same way. 

\end{sbpara}

\begin{sbpara} We replace the two lines 

\medskip

\lq\lq As in [KU09] 6.4, a key step ..... We only
prove this proposition."

\medskip
\noindent
just before Part III, 3.3.3 by the following \ref{8.3.b}.

\end{sbpara}

\begin{sbpara}\label{8.3.b}

This Part III, Theorem 3.3.2 is proved in the following way.

We can prove the evident mixed Hodge version of Proposition A.1.6 by
using the same arguments in the proof of Proposition A.1.6.

Just as  \cite{KU} Theorem 5.4.4 was reduced to the case $y^*_{\lam,t}=
y_{\lam,t}$ of \cite{KU} Proposition 6.4.1 by using the case $m=0$ and $x_{\lam}\in E^{\sharp}_{\sig,\val}$ of A.1.6  (see
\ref{8.1.c}, \ref{8.1.d}), Part III, Theorem 3.3.2 is  reduced to the case $y^*_{\lam,t}=
y_{\lam,t}$ of Part III, Proposition 3.3.4 by using the case $m=0$ of this mixed Hodge version of
Proposition A.1.6.

\end{sbpara}

\begin{sbpara}\label{8.3.c} We replace the part

\medskip
\lq\lq proposition implies"
 
 \medskip
 \noindent
 in line 2 of Remark after Part III, Proposition 3.3.4 with 
 
 \medskip

\lq\lq proposition and \cite{KU} 6.4.1 of the pure case imply".

\end{sbpara}

\noindent {\rm Kazuya KATO
\\
Department of mathematics
\\
University of Chicago
\\
Chicago, Illinois, 60637, USA}
\\
{\tt kkato@math.uchicago.edu}

\bigskip

\noindent {\rm Chikara NAKAYAMA
\\
Department of Economics 
\\
Hitotsubashi University 
\\
2-1 Naka, Kunitachi, Tokyo 186-8601, Japan
\\
{\tt c.nakayama@r.hit-u.ac.jp}

\bigskip

\noindent
{\rm Sampei USUI
\\
Graduate School of Science
\\
Osaka University
\\
Toyonaka, Osaka, 560-0043, Japan}
\\
{\tt usui@math.sci.osaka-u.ac.jp}
\end{document}